\input amstex
\input amsppt.sty
\NoBlackBoxes
\nologo
\topmatter
\title Notes \linebreak
on the inductive algorithm of \linebreak
resolution of singularities \linebreak
by S. Encinas and O. Villamayor
\endtitle
\author Kenji Matsuki
\endauthor
\endtopmatter
\leftheadtext{}
\rightheadtext{}

\document

$$\bold{CONTENTS}$$

\vskip.2in

Chapter 0. Introduction

\vskip.1in

Chapter 1. Basic objects and invariants

\vskip.1in
 
Chapter 2. Resolution of singularities of monomial basic objects

\vskip.1in
 
Chapter 3. Key inductive lemma

\vskip.1in
 
Chapter 4. General basic objects and invariants

\vskip.1in

Chapter 5. Inductive algorithm for resolution of singularities of general basic 

\hskip.7in objects

\vskip.1in
 
Chapter 6. A more down-to-earth approach to the inductive algorithm

\vskip.1in
 
Chapter 7. Embedded resolution of singularities

\vskip.1in
 
Chapter 8. Equivariance and resolution of singularities over base fields

\hskip.7in (of characteristic zero) which are possibly not algebraically closed

\vskip.1in

Chapter 9. Invariants revisited

\vskip.1in
 
Chapter 10. Non-embedded resolution of singularities

\vskip.1in

Chapter 11. Examples

\vskip.1in

References

\newpage

$$\bold{CHAPTER\ 0.\ INTRODUCTION}$$

\vskip.1in

The purpose of these notes is simply to record the regurgitation of the beautiful and
elegant ideas of Encinas and Villamayor on the problem of resolution of singularties in the
papers: 

``A course on constructive desingularization and equivariance"

``A new theorem of desingularization over fields of characteristic zero"

``On properties of constructive desingularization" (by Encinas).

The notes are the results of seminars held at Purdue University, organized by A. Gabrielov and
the author, in the Fall semester of 2000 and continued in the Spring semester of 2001.  

After the
first draft of these notes was written, we had the fortune of Villamayor himself visiting Purdue
University to give a series of lectures titled ``Constructive Desingularization". 
Consequently we added Chapter 6, which should explain the origin of the ingeneous
$t$-invariant, to the revised version based upon one of his lectures.  Some of the examples
presented in Chapter 11 are also taken from his lectures.  We thank Prof. Villamayor for
his generous permission to include these in this revised version.

\vskip.1in
 
The following are the main themes of these notes.

\proclaim{Main Theme 0-1 (Resolution of singularities)} Understand the solution by Encinas and
Villamayor (extending of course some of the original ideas of Hironaka) to the problem of
resolution of singularities: 

Let $X$ be a variety over a field $k$ of characteristic zero.  Establish an algorithm
to construct a sequence of blowups
$$X = X_0 \overset{\pi_1}\to{\leftarrow} X_1 \overset{\pi_2}\to{\leftarrow} \cdot\cdot\cdot
\overset{\pi_{l-1}}\to{\leftarrow} X_{l-1} \overset{\pi_l}\to{\leftarrow} X_l$$
so that 

(i) the centers $Y_{i-1} \subset X_{i-1}$ of the blowups $\pi_i \hskip.1in (i = 1,
... , l)$ are over
\linebreak
$\roman{Sing}(X) = X \setminus \roman{Reg}(X)$, 

(ii) the centers $Y_{i-1} \subset X_{i-1}$ are closed subschemes, which may be reducible and may
NOT be smooth or reduced in general\footnote"${}^1$"{We want to emphasize that we
do NOT require that the centers $Y_{i-1}$ be smooth or even reduced, or that they be contained
in the singular loci of the varieties $X_{i-1}$, i.e., $Y_{i-1}
\subset
\roman{Sing}(X_{i-1})$, as Hironaka or Bierstone-Milman does in their presentation of resolution
of singularities.  It seems that this slight weakening of the statement is a price we have to pay
for dealing only with the order function and weak transforms, and not with the
Hilbert-Samuel function, which is better suited for detecting the strict transforms.  The
author would like to thank Prof. Bierstone, who brought this fact to the attention of the
author and pointed out the mistakes in the earlier manuscript.},

(iii) $X_l$ is a variety smooth over $k$ and the induced morphism $X = X_0
\overset{\pi}\to{\leftarrow} X_l$, where $\pi = \pi_1 \circ \pi_2 \circ \cdot\cdot\cdot
\circ \pi_{l-1}
\circ
\pi_l$, is a projective birational morphism isomorphic over
$\roman{Reg}(X)$.
\endproclaim

The main body of these notes, Chapter 1 through Chapter 7, will be devoted to the solution to the
following problem of ``embedded" resolution of singularities, from which the solution to the
original problem of resolution of singularities immediately follows.  (See Chapter 10 for detail.)

\proclaim{Main Theme 0-2 (Embedded resolution of singularities)} Understand the solution to the
problem of ``embedded" resolution of singularities:

Let $X \subset W$ be a variety, embedded as a closed subscheme of another variety $W$ smooth over a
field $k$ of characteristic zero.  Establish an algorithm to construct a sequence of blowups
$$X = X_0 \subset W = W_0 \overset{\pi_1}\to{\leftarrow} X_1 \subset
W_1 \overset{\pi_2}\to{\leftarrow} \cdot\cdot\cdot \overset{\pi_{l-1}}\to{\leftarrow} X_{l-1}
\subset W_{l-1} \overset{\pi_l}\to{\leftarrow} X_l
\subset W_l$$
so that

(i) the centers $Y_{i-1} \subset W_{i-1}$ of the blowups $\pi_i \hskip.1in (i = 1,
... , l)$ are over $\roman{Sing}(X) = X \setminus \roman{Reg}(X)$, 

(ii) the centers $Y_{i-1} \subset W_{i-1}$\footnote"${}^2$"{We want to emphasize that we do NOT
require that $Y_{i-1} \subset X_{i-1}$, i.e., the centers $Y_{i-1}$ be contained in the strict
transforms
$X_{i-1}$ of $X = X_0$, or that they be contained in the singular loci of the strict transforms,
i.e., $Y_{i-1} \subset \roman{Sing}(X_{i-1})$, as Hironaka or Bierstone-Milman does in their
presentation of embedded resolution of singularities.  Therefore, though the centers $Y_{i-1}$ are
smooth in the ambient varieties $W_{i-1}$, their restrictions $Y_{i-1} \cap X_{i-1}$ to the strict
transforms may not be smooth or reduced in general.} are permissible with respect to the
exceptional divisors
$E_{i-1} \subset W_{i-1}$ for the morphisms $\psi_{i-1} = \pi_1 \circ \pi_2 \circ \cdot\cdot\cdot
\circ \pi_{i-2} \circ \pi_{i-1}$ (which are simple normal crossing divisors),

(iii) the strict transform $X_l$ (of $X_0$) $\subset W_l$ is a variety smooth over $k$, permissible
with respect to $E_l$, and the induced morphism $X = X_0 \overset{\pi}\to{\leftarrow} X_l$,
where \linebreak
$\pi = \psi_l = \pi_1
\circ
\pi_2 \circ \cdot\cdot\cdot \circ \pi_{l-1} \circ \pi_l$, is a projective
birational morphism isomorphic over $\roman{Reg}(X)$.
\endproclaim

Note that we say the center $Y_{i-1} \subset W_{i-1}$ is permissible with respect to $E_{i-1}$ if
$Y_{i-1}$ are smooth, at each closed point
$p
\in W_{i-1}$ there exists an open neighborhood $U_p$ with a system of regular parameters $(x_1,
... , x_d)$ such that \linebreak
$Y_{i-1} \cap U_p = \cap_{i \in M}\{x_i = 0\}$ and $E_{i-1} \cap
U_p =
\{\prod_{i \in N}x_i = 0\}$ for some subsets $M, N \subset \{1, ..., d = \dim W_{i-1}\}$, and that we
say $E_{i-1}$ is a simple normal crossing divisor where the irreducible
components of
$E_{i-1}$ are required to be smooth without self-intersection, in contrast to the condition of being a
normal crossing divisor where only the local requirement $E_{i-1} \cap U_p =
\{\prod_{i \in N}x_i = 0\}$ is posed with the system of regular parameters $(x_1, ..., x_d)$
chosen analytically.

\vskip.1in

We remark that the solution to the problem of embedded resolution of singularities is
derived from looking at our specific \footnote"${}^3$"{Without condition (i) imposed on
our formulation of embedded resolution, which requires the centers to be taken over
$\roman{Sing}(X)$, a solution to the problem of embedded resolution of singularities
follows immediately as a corollary to the solution to the problem of principalization,
if we apply the latter to the defining ideal ${\Cal I}_X \subset {\Cal O}_W$ and look
at the stage where the strict transform becomes the center of blowup.  However, in
order to satisfy condition (i), we need more requirements on the algorithm of
principalization.  This is why we have to look at the ``specific" algorithm as we
discuss in these notes.} algorithm to solve the problem of ``principalization" of
ideals.

\proclaim{Main Theme 0-3 (Principalization of ideals)} Understand the solution to the problem of
``principalization" of ideals: Let $W$ be a variety smooth over a field $k$ of characteristic zero
and ${\Cal I} \subset {\Cal O}_W$ be a coherent sheaf of ideals.  Establish an algorithm to
construct a sequence of blowups
$$W = W_0 \overset{\pi_1}\to{\leftarrow} W_1 \overset{\pi_2}\to{\leftarrow} \cdot\cdot\cdot
\overset{\pi_{l-1}}\to{\leftarrow} W_{l-1} \overset{\pi_l}\to{\leftarrow} W_l$$
so that

(i) the centers $Y_{i-1} \subset W_{i-1}$ of the blowups $\pi_i \hskip.1in (i = 1, ... , l)$ are
over the support of ${\Cal O}_W/{\Cal I}$,

(ii) the centers $Y_{i-1} \subset W_{i-1}$ are permissible with respect to the exceptional divisors
$E_{i-1} \subset W_{i-1}$ for the morphisms $\psi_{i-1} = \pi_1 \circ \pi_2 \circ \cdot\cdot\cdot
\circ \pi_{i-2} \circ \pi_{i-1}$,

(iii) the total transform ${\Cal I}_l = {\Cal I}{\Cal O}_{W_l}$ of the ideal
${\Cal I}$ (We write ${\Cal I}{\Cal O}_{W_l}$ for the ideal of ${\Cal O}_{W_l}$ generated by
$\psi_l^{-1}({\Cal I})$ by abuse of notation.) is a product of the principal ideals defining
divisors
$H_j$
$${\Cal I}_l = I(H_1)^{a_1} \cdot\cdot\cdot I(H_m)^{a_m}$$
where the divisors $H_j$ and the exceptional divisor $E_l$ for $\psi_l$ form a divisor
with only simple normal crossings.
\endproclaim

In our formulation of the probelm of principalization, it should be emphasized and warned against the
common usage of the word ``principal", we require not only the (total transform of the) ideal to be
locally generated by one element but also to be a product of the defining ideals of the irreducible
components of a simple normal crossing divisor.

\vskip.1in
 
We note that, for the sake of simplicity of presentation, we assume that the base field $k$ is
algebraically closed in Chapter 1 through Chapter 7, aside from the basic assumption that $k$ is
of characteristic zero.  (The general case where
$k$ may not be algebraically closed is discussed in Chapter 8 and it can be settled rather easily
after the discussion of equivariance under the action of the Galois group
$\roman{Gal}(\overline{k}/k)$ on the process prescribed over $\overline{k}$.)
\vskip.1in

The key strategy of Encinas and Villamayor is to reduce the problem of (embedded) resolution of
singularities, which is reformulated as the problem of principalization, to that of resolution of
singularities of ``basic objects", the notion we introduce in Chapter 1.  The basic objects
are designed to extract the inductive nature of the problem.  The elegance of their ideas is
condensed in the definition of the
$t$-invariant attached to a (sequence of) basic object(s).

Chapter 2 discusses a solution to the problem of resolution of singularities of monomial basic
objects, where the given ideals (of the basic objects) are already
products of the principal ideals defining the irreducible components of the boundary divisors. 
This turns out to be the easiest case where the solution can be given in a concise
combinatorial manner.  When (the maximum of) the invariant ``$w\text{-}\roman{ord}$" of a basic
object is equal to 0, resolution of singularities of the basic object is reduced to that of the
monomial ones.  Therefore, in the later chapters, we consider an algorithm for resolution of
singularities of (general) basic objects to be complete as soon as (the maximum of) the
invariant
$w\text{-}\roman{ord}$ is 0.

Chapter 3 reveals the key inductive lemma, which reduces the problem of resolution of
singularities of a basic object of dimension $d$ to that of resolution of
singularities of charts consisting of basic objects of dimension
$d-1$, and hence realizing our inductive strategy via the notion of basic objects.  But there
is a catch.  We have to assume that the basic object to start with to be ``simple" and
also have to assume the existence of smooth hypersurfaces (inside of the open subsets which give
rise to the charts) which cover the singularities of the simple basic object and which cross
transversally with the specified boundary divisor of the original simple basic object of
dimension
$d$.  The lemma forms the basis of
our inductive argument.

The key inductive lemma leads us naturally to the notion of ``general basic objects",
generalizing the notion of basic objects, so that we can carry out the inductive argument,
suggested by the lemma, in a more natural framework.  This is done in Chapter 4, clarifying some
minor obscure points in the original papers.  As a general basic object consists of
(local) charts of basic objects, there arises a problem of patching up the processes of
resolution of singularities of various (local) charts of basic objects to form a unique process
of resolution of singularities of the (global) general basic object.  This problem will be
solved by showing via Hironaka's trick that the invariants defined on individual (local) charts
patch up to provide well-defined (global) invariants on the general basic object, which in turn
determine the global centers of blowups in the process of resolution of singularities.  This is
another subject of Chapter 4. 

The key inductive lemma, however, falls short of completing the inductive
process of resolution of singularities of (general) basic objects for the following two reasons
(difficulties) (which was described as a ``catch" in the previous paragraph):

1. It requires the original basic object of dimension $d$ to be simple, though the resulting
(general) basic objects of dimension $d-1$ may not be (and in most cases actually are not) simple.

2. It requires the existence of smooth hypersurfaces (inside of the open subsets which give rise
to the charts) satisfying the conditions
(including the transversality) mentioned above.

The elegant and brilliant theorem of Encinas and Villamayor, discussed in Chapter 5, overcomes these
two difficulties in one stroke with the use of the ingeneous $t$-invariant, and hence provides a
complete inductive algorithm for resolution of singularities of general basic objects.

However ingeneous it may be, nonetheless, the use of the $t$-invariant to complete the inductive
algorithm in Chapter 5 may look ``slick" in the untrained eyes and seems as though it came ``out
of blue".  Chapter 6 presents a more down-to-earth approach to the inductive algorithm, which, by
decomposing the inductive algorithm into a few reduction steps, tries to explain where the
$t$-invariant comes from and how natural it is.

The inductive algorithm for resolution of singularities of general basic objects
provides a solution to the problem of principalization of ideals, achieving Main Theme
0-3.  Now a solution to the problem of embedded resolution of singularities follows as
an easy corollary, if we apply this specific algorithm for principalization to the
defining ideal ${\Cal I}_X \subset {\Cal O}_W$ of an embedding $X \subset W$.  The
argument is presented in Chapter 7, achieving Main Theme 0-2.

We observe in Chapter 8 that the inductive algorithm is equivariant under any group
action.  This implies, in particular, that the process of the algorithm prescribed over
$\overline{k}$, where $\overline{k}$ is the algebraic closure of the base field $k$, is
equivariant under the action of the Galois group
$\roman{Gal}(\overline{k}/k)$, and hence that the process is actually defined over
$k$.  This observation provides an inductive algorithm for embedded resolution of
singularities over any field of characteristic zero.

In Chapter 9, we construct an invariant, based upon the $w\text{-}\roman{ord}$,
$\Gamma$- and $t$-invariants, of general basic objects, so that the centers of blowups
in our inductive algorithm for resolution of singularities are exactly the loci where
the values of this invariant attain maxima.

A variety $X$ is covered by a finite number of open subsets $U$ which can be embedded into smooth
varieties $W_U$.  By choosing a number $d$ sufficiently large and replacing $W_U$ with $W_U \times
{\Bbb A}^{d - \dim W_U}$ if necessary, we may assume that all the ambient smooth varieties $W_U$
are of the same dimension $d$.  We observe then that the processes of embedded resolution of
singularities of $U \subset W_U$ prescribed by our inductive algorithm patch up and give rise
to a sequence representing non-embedded resolution of singularities of $X$ as stated in Main
Theme 0-1.  This is done in Chapter 10 via the analysis of the invariants constructed in
Chapter 9.

In Chapter 11, we give examples demonstrating the mechanism and some subtleties of our
inductive algorithm for embedded and non-embedded resolution of singularities.
 
\vskip.2in

The elementary nature of the inductive algorithm by Encinas and Villamayor, which does not even
make an explicit use of the Hilbert-Samuel function and builds its key invariants upon the order
(multiplicity) function, allowed us to try to make these notes self-contained.  We provide complete
proofs for embedded and non-embedded resolution of singularities over any field of characteristic
zero (as formulated in Main Themes 0-1 and 0-2, which are slightly weaker than the formulation by
Hironaka or Bierstone-Milman), with little reference to the other literature, for an easy
understanding on the side of the reader.  We even try to avoid referring to the original papers by
Encinas and Villamayor, though almost all the proofs are taken verbatim from them.

\vskip.1in

We are very much aware of the other important developments on the subject of canonical
and constructive resolution of singularities, especially the monumental paper by E.
Bierstone and P. Milman:

``Canonical desingularization in characteristic zero by blowing up the maximum strata of a local
invariant", Inventiones Mathematicae 128, 207-302 (1997).

The restricted attention in these notes to the inductive algorithm by Encinas and Villamayor,
with no discussion on the above-mentioned developments or on the more classical papers including
Hironaka's, is merely a result of the lack of resource and time to run the seminars but mainly
caused by the incompetence of the author, who is responsible for any mistakes in these notes.

\vskip.2in

$\bold{Note\ to\ the \ reader}$: In the process of revision, the size of these notes became
much bigger than what would not scare off a reader wishing for a concise and minimal
understanding of the subject.  For such a reader, we would like to recommend reading only of
Chapter 1 through Chapter 7 (Chapter 6 is not necessary for the logic of the development of the
argument but is of great help in order to understand the core ideas behind all the technical
details.), where, when he finishes, a self-contained proof for an algorithm of embedded
resolution of singularities is obtained.

\newpage

$$\bold{CHAPTER\ 1.\ BASIC\ OBJECTS\ AND\ INVARIANTS}$$

\vskip.1in

In Chapter 1 through Chapter 7, the base field $k$ is assumed to be algebraically
closed and of characteristic zero.

\vskip.1in

Let $W$ be a variety smooth over $k$ of dimension $d$ and $J \subset {\Cal O}_W$ a coherent sheaf
of ideals (which we simply call an ideal by abuse of language).

\proclaim{Definition 1-1 (Order of an ideal)} Let $p \in W$ be a point.  The order
$\nu_p(J)$ of an ideal $J \subset {\Cal O}_W$ at $p$ is defined to be
$$\nu_p(J) := \nu(J_p) = \max\{n \in {\Bbb Z}_{\geq 0};J_p \subset m_p^n\}$$
where $m_p$ is the maximal ideal of the local ring ${\Cal O}_{W,p}$ and where $J_p \subset {\Cal
O}_{W,p}$ is the stalk of
$J$ at $p$.
\endproclaim

\vskip.1in

\proclaim{Remark 1-2 (Some properties of order)}\endproclaim

(i) Let $\widehat{{\Cal O}_{W,p}}$ be the ($m_p$-adic) completion of ${\Cal O}_{W,p}$,
$\widehat{J_p} = J_p \otimes_{{\Cal O}_{W,p}} \widehat{{\Cal O}_{W,p}}$ and $\widehat{m_p} = m_p
\otimes_{{\Cal O}_{W,p}}
\widehat{{\Cal O}_{W,p}}$ the completions of the ideals $J_p$ and $m_p$, respectively, in
$\widehat{{\Cal O}_{W,p}} = {\Cal O}_{W,p} \otimes_{{\Cal O}_{W,p}}\widehat{{\Cal O}_{W,p}}$.  The
order
$\nu(\widehat{J_p})$ of
$\widehat{J_p}$ coincides with the order $\nu_p(J) = \nu(J_p)$ of $J_p$, i.e., 
$$\nu(J_p) = \nu(\widehat{J_p}) = \max\{n \in {\Bbb Z}_{\geq 0};\widehat{J_p} \subset
\widehat{m_p}^n\},$$  
since $\widehat{{\Cal O}_{W,p}}$ is faithfully flat over ${\Cal O}_{W,p}$.

Observe that $\widehat{{\Cal O}_{W,p}}$ is isomorphic to a power series ring over $k$ at a closed
point
$p
\in W$, i.e., once we fix a system of regular parameters $(x_1, ... , x_d)$ we have a
$k$-algebra isomorphism 
$$\widehat{{\Cal O}_{W,p}} \cong k[[x_1, ... , x_d]],$$  
sending $x_1, ... , x_d$ of $\widehat{{\Cal O}_{W,p}}$ to the corresponding variables in
$k[[x_1, ... , x_d]]$.

Therefore,
$$\nu(\widehat{J_p}) = \min\{\nu(f);f \in \widehat{J_p}\}$$
where $\nu(f)$ is the lowest degree of the Taylor expansion of $f$ considered as an element of
the power series ring.  In particular, if $\{f_i\}$ is a set of generators for $J_p$ over ${\Cal
O}_{W,p}$ and hence for
$\widehat{J_p}$ over $\widehat{{\Cal O}_{W,p}}$, then
$$\nu(J_p) = \nu(\widehat{J_p}) = \min\{\nu(f_i)\}.$$

(ii) Let $I, J \subset {\Cal O}_W$ be ideals.  Then
$$\align
\nu_p(I + J) &= \min\{\nu_p(I), \nu_p(J)\} \\
\nu_p(I \cdot J) &= \nu_p(I) \cdot \nu_p(J).\\
\endalign$$

\vskip.1in

The order of an ideal can be analyzed using ``derivatives".  The analysis naturally leads to
the following notion of the ``extension" of an ideal.

\proclaim{Definition-Proposition 1-3 (Extension of an ideal)} Let $J \subset {\Cal O}_W$ be an ideal.

(i) The extension $\Delta(\widehat{J_p})$ of $\widehat{J_p} \subset \widehat{{\Cal
O}_{W,p}}$, where $p \in W$ is a closed point, is defined to be the ideal generated by the elements
$f$ of
$\widehat{J_p}$ and their (partial) derivatives
$\frac{\partial f}{\partial x_i}$ via $k$-algebra isomorphism  $\widehat{{\Cal O}_{W,p}} \cong
k[[x_1,
... , x_d]]$ (cf. Remark 1-2 (i)), i.e.,
$$\Delta(\widehat{J_p}) = \langle \widehat{J_p}, \frac{\partial f}{\partial x_1},
... , \frac{\partial f}{\partial x_d}; f \in \widehat{J_p} \rangle.$$   
The extension
$\Delta(\widehat{J_p})$ is determined independently of the choice of the isomorphism.

(ii) There uniquely exists an ideal $\Delta(J) \subset {\Cal O}_W$, called the extension of $J$, such
that
$$\Delta(J)_p \otimes_{{\Cal O}_{W,p}} \widehat{{\Cal O}_{W,p}} = \Delta(\widehat{J_p})$$
for all closed points $p \in W$. 
\endproclaim 

\demo{Proof}\enddemo (i) It follows from the chain rule that the extension is independent of the
choice of the isomorphism.

(ii) Take an affine open covering $\{U\}$ of $W$ together with a system of regular parameters
$(x_1, ... , x_d)$ over $U$ so that $(dx_1, ... , dx_d)$ provide generators of the
locally free sheaf $\Omega_W^1$ of rank $d$ over $U$.  Take the dual generators
$(\frac{\partial}{\partial x_1}, ... , \frac{\partial}{\partial x_d})$ of the tangent
sheaf $T_W$ over $U$ so that $(\frac{\partial}{\partial x_i}, dx_j) =
\delta_{ij}.$  

We only have to take the ideal generated by the elements $f \in J(U)$ and their (partial)
derivatives
$\frac{\partial f}{\partial x_i} = (\frac{\partial}{\partial x_i}, df)$ over $U$ in
order to define and obtain the extension $\Delta(J)|_U$.

The characterization as described in (ii) can be easily checked and implies the uniqueness at any
closed point and hence of the sheaf.  This also implies that the collection $\{\Delta(J)|_U\}$ patch
up to provide the extension $\Delta(J)$ over $W$. 

\vskip.1in

The relation between the order of an ideal and its extension(s) is described by the following lemma.

\proclaim{Lemma 1-4 (Characterization of order in terms of extensions)} Let $V(I)$ denote the zero
locus of an ideal
$I
\subset {\Cal O}_W$.

(i) Let $b \in {\Bbb N}$ be a positive integer.  Then
$$p \in V(\Delta^{b-1}(J)) \Longleftrightarrow \nu_p(J) \geq b$$
and hence
$$p \in V(\Delta^{b-1}(J)) \setminus V(\Delta^b(J)) \Longleftrightarrow \nu_p(J) = b,$$
where $\Delta^b$ represents the $b$-iterations of the operation $\Delta$ of taking the
extension of an ideal.

In particular, the function $\nu_J:W \rightarrow {\Bbb Z}_{\geq 0}$ defined by $\nu_J(p) =
\nu_p(J)$ is upper semi-continuous.

(ii) Let $p \in W$ be a point.  Then
$$\nu_p(\Delta^{b-i}(J)) = i \Longleftrightarrow \nu_p(J) = b$$
for $i = 1, \cdot\cdot\cdot, b$.

In particular,
$$\nu_p(\Delta^{b-1}(J)) = 1 \Longleftrightarrow \nu_p(J) = b$$
\endproclaim 

\demo{Proof}\enddemo Assertions (i) and (ii) are immediate consequences of Remark 1-2 (i) for
closed points.  For arbitrary points, one has only to argue taking a general closed point in its
closure.  

\vskip.1in

We emphasize:
$$\boxed{\roman{The\ assumption\ of\ characteristic\ being\ equal\ to\ 0\ is\ essential\ for\
this\ lemma.}}$$ 

\vskip.1in

\proclaim{Remark 1-5 (Primitive but fundamental idea toward the inductional argument)}\endproclaim

As a consequence of Lemma 1-4, we come to the following primitive but fundamental observation, which
forms the core of our idea toward the inductive argument:

Let $J \subset {\Cal O}_W$ be an ideal and $b_{\max} = \max\{\nu_p(J);p \in W\}$.  Then the locus
$$S = \{p \in W; \nu_p(J) = b_{\max}\} = V(\Delta^{b_{\max}-1}(J))$$
is closed.  For any point $p \in S$, we can find a
neighborhood $U_p$ of $p$ in $W$ and a smooth hypersurface $H_p \subset U_p$ such that  
$$p \in S \cap U_p \subset H_p \subset U_p,$$
since $\nu_p(\Delta^{b_{\max}-1}(J)) = 1$.  In fact, we have only to take an element \linebreak
$f_p
\in
\Delta^{b_{\max}-1}(J)_p$ with $\nu(f_p) = 1$, and set $H_p = \{f_p = 0\} \subset U_p$ with $U_p$
an open neighborhood of $p$ where $f_p$ is regular and where the order of $f_p$ remains 1.

\vskip.1in

This observation suggests the possibility that the analysis of the ``worst" locus
$S$ of the ideal $J$ on a $d$-dimensional smooth variety $W$ may be reduced, at least locally, to the
one on a
$(d-1)$-dimensional smooth variety $H_p$, which is sometimes called a $\bold{hypersurface\
of\ maximal\ contact}$ (at $p$).

\vskip.1in

Now we introduce the notion of a basic object.

\proclaim{Definition 1-6 (Basic object)} A basic object is a triplet $(W, (J,b), E)$ where $W$ is a
variety smooth over $k$, $(J,b)$ is a pair consisting of an ideal $J \subset {\Cal O}_W$ and a
positive integer $b
\in {\Bbb N}$, and where $E = \{H_1, ... , H_r\}$ is a divisor with
simple normal crossings.  (We sometimes call $E$ the boundary divisor of the basic object.)  

We define the singular locus of the basic object to be
$$\roman{Sing}(J,b) = \{p \in W;\nu_p(J) \geq b\}.$$ 
We call a couple $(W,E)$, consisting of $W$ and $E$ as above, a pair.
\endproclaim

\vskip.1in

\proclaim{Note 1-7}\endproclaim

We apply the same slightly abusive notation $E = \{H_1, ... , H_r\}$ in the above as was
used in the original papers of Encinas and Villamayor: $H_i$ actually consists of smooth
irreducible components $H_{i,1}, ... , H_{i,l_i}$ disjoint from each other.  If $i \neq
j$, then $H_i$ and $H_j$ share no common irreducible components.  We require that $\cup_{i = 1,
... , r, l = 1, ... , l_i}H_{i,l}$ is a divisor with simple normal
crossings.  

(Thus, strictly speaking, if we want to avoid the abuse, we should write 
$$E = \{H_{1,1}, H_{1,2},
... , H_{1,l_1}, H_{2,1}, ... , H_{2,l_2}, ... , H_{r,1},
... , H_{r,l_r}\}.)$$

\vskip.1in

\proclaim{Definition 1-8 (Sequence of transformations and smooth morphisms of basic
objects)}  We define a sequence of transformations and smooth morphisms of basic objects
$$\align
(W, (J,b),E) = (W_0, (J_0,b), E_0) &\overset{\pi_1}\to{\leftarrow} (W_1, (J_1,b), E_1)
\overset{\pi_2}\to{\leftarrow} \cdot\cdot\cdot \\
(W_{i-1}, (J_{i-1},b), E_{i-1}) &\overset{\pi_i}\to{\leftarrow}
(W_i, (J_i,b), E_i) \\
\cdot\cdot\cdot \overset{\pi_{k-1}}\to{\leftarrow} (W_{k-1}, (J_{k-1},b),
E_{k-1}) &\overset{\pi_k}\to{\leftarrow} (W_k, (J_k,b), E_k) \\
\endalign$$
to satisfy the following conditions:

\vskip.1in

Case T: $(W_{i-1}, (J_{i-1},b), E_{i-1}) \overset{\pi_i}\to{\leftarrow} (W_i, (J_i,b), E_i)$ is a
transformation.

\vskip.1in

(i) $W_{i-1} \overset{\pi_i}\to{\leftarrow} W_i$ is the blowup
with a center
$Y_{i-1}
\subset W_{i-1}$ which is permissible with respect to $E_{i-1}$, i.e., $Y_{i-1}$ is smooth (maybe
reducible) and for any closed point $p \in W_{i-1}$ there exists an open neighborhood $U_p$ with a
system of regular parameters $(x_1, ... , x_d)$ such that
$$\align
Y_{i-1} \cap U_p &= \cap_{m \in M}\{x_m = 0\} \\
E_{i-1} \cap U_p &= \{\prod_{m \in N}x_m = 0\} \\
\endalign$$
for some subsets $M, N \subset \{1, ... , d = \dim W_{i-1}\}$. 

(ii) The center $Y_{i-1} \subset W_{i-1}$ is contained in the singular locus of the basic object
$(W_{i-1}, (J_{i-1},b), E_{i-1})$, i.e.,
$$Y_{i-1} \subset \roman{Sing}(J_{i-1},b).$$

(iii) $H_{r+i} = \pi_i^{-1}(Y_{i-1})$ and $E_i = \{H_1, ... , H_r, H_{r+1},
... , H_{r+i-1}, H_{r+i}\}$ where
$H_j \hskip.1in (j = 1, ... , r, r+1, ... , r+i-1)$ in $E_i$ denotes by abuse
of notation the strict transform of
$H_j$ in $E_{i-1} = \{H_1, ... , H_r, H_{r+1},
... , H_{r+i-1}\}$.  

(We also use the convention that if $H_{j,l}$ is a smooth irreducible component belonging
to $H_j$ in $E_{i-1}$ and if $H_{j,l} \subset Y_{i-1}$, then we exclude $H_{j,l}$ from $H_j$ in
$E_i$ and consider it as an element belonging to $H_{r+i}$ in $E_i$.)

(iv) $J_i \subset {\Cal O}_{W_i}$ is the unique ideal such that
$$J_{i-1}{\Cal O}_{W_i} = I(H_{r+i})^b \cdot J_i.$$

(We note that the existence of such an ideal $J_i$, i.e., the fact that $J_{i-1}{\Cal
O}_{W_i}$ is divisible by $I(H_{r+i})^b$, is guaranteed by the condition $Y_{i-1} \subset
\roman{Sing}(J_{i-1},b)$ and can be checked, e.g., by Lemma 1-13.)

\vskip.1in

Case S: $(W_{i-1}, (J_{i-1},b), E_{i-1}) \overset{\pi_i}\to{\leftarrow} (W_i, (J_i,b), E_i)$ is a
smooth morphism.

\vskip.1in

(i) $W_{i-1} \overset{\pi_i}\to{\leftarrow} W_i$ is a smooth morphism.

(ii) $E_i = \{\pi_i^{-1}(H_1), ... , \pi_i^{-1}(H_r), \pi_i^{-1}(H_{r+1}),
... , \pi_i^{-1}(H_{r+i-1})\}$ where \linebreak
$E_{i-1} = \{H_1, ... , H_r,
H_{r+1}, ... , H_{r+i-1}\}$.  By abuse of notation and for consistency in notaion with
Case T, we write
$E_i = \{H_1,
... , H_r, H_{r+1}, ... , H_{r+i-1}, H_{r+i}\}$ with the understanding that $H_1,
... , H_r, H_{r+1}, ..., H_{r+i-1}$ in $E_i$ denote the corresponding
pull-backs $\pi_i^{-1}(H_1), ... ,
\pi_i^{-1}(H_r),
\pi_i^{-1}(H_{r+1}), ... ,
\pi_i^{-1}(H_{r+i-1})$ and that
$H_{r+i} =
\emptyset$.

(iii) $J_i = J_{i-1}{\Cal O}_{W_i}$.

\vskip.1in

We define a sequence of transformations and smooth morphisms of pairs
$$(W,E) = (W_0,E_0) \overset{\pi_1}\to{\leftarrow} (W_1,E_1) \overset{\pi_2}\to{\leftarrow}
\cdot\cdot\cdot \overset{\pi_{k-1}}\to{\leftarrow} (W_{k-1},E_{k-1})
\overset{\pi_k}\to{\leftarrow} (W_k,E_k)$$
to satisfy conditions (i), (iii) in Case T and (i), (ii) in Case S.
\endproclaim

\vskip.1in

\proclaim{Definition 1-9 (Resolution of singularities of a basic object)} We call a sequence of
transformations only (i.e., all the $\pi_i$ are in Case T) of basic objects 
$$\align
(W, (J,b), E) = (W_0, (J_0,b), E_0) &\overset{\pi_1}\to{\leftarrow} (W_1, (J_1,b), E_1)
\overset{\pi_2}\to{\leftarrow} \cdot\cdot\cdot \\
(W_{i-1}, (J_{i-1},b), E_{i-1}) &\overset{\pi_i}\to{\leftarrow}
(W_i, (J_i,b), E_i) \\
\cdot\cdot\cdot \overset{\pi_{k-1}}\to{\leftarrow} (W_{k-1}, (J_{k-1},b),
E_{k-1}) &\overset{\pi_k}\to{\leftarrow} (W_k, (J_k,b), E_k) \\
\endalign$$ 
resolution of singularities of a basic object $(W, (J,b), E)$ if 
$$\roman{Sing}(J_k,b) = \emptyset.$$
\endproclaim 

As will be seen in Chapter 7 through Chapter 9, the problem of (embedded) resolution of
singularities, as well as the problem of principalization, can be readily reduced to the
problem of  resolution of singularities of (general) basic objects.   

\vskip.1in

One of the keys to solve the problem of resolution of singularities of (general)
basic objects is to define the following invariants
$\roman{ord}_k, w\text{-}\roman{ord}_k$ and
$t_k$ on a basic object
$(W_k,(J_k,b), E_k)$ appearing in a sequence of transformations and smooth morphisms as above (with
one extra condition on the sequence in order to define the invariant $t_k$).

\proclaim{Definition 1-10 (Key invariants of basic objects)} Let 
$$\align
(W, (J,b),E) = (W_0, (J_0,b), E_0) &\overset{\pi_1}\to{\leftarrow} (W_1, (J_1,b), E_1)
\overset{\pi_2}\to{\leftarrow} \cdot\cdot\cdot \\
(W_{i-1}, (J_{i-1},b), E_{i-1}) &\overset{\pi_i}\to{\leftarrow}
(W_i, (J_i,b), E_i) \\
\cdot\cdot\cdot \overset{\pi_{k-1}}\to{\leftarrow} (W_{k-1}, (J_{k-1},b),
E_{k-1}) &\overset{\pi_k}\to{\leftarrow} (W_k, (J_k,b), E_k) \\
\endalign$$
be a sequence of transformations and smooth morphisms of basic objects as defined in
Definition 1-8.  

\vskip.1in

(i) The invariant $\roman{ord}_k:\roman{Sing}(J_k,b) \rightarrow \frac{1}{b}{\Bbb Z}_{\geq 0}$ is a
function defined over $\roman{Sing}(J_k,b)$ such that
$$\roman{ord}_k(p) = \frac{\nu_p(J_k)}{b} \text{\ for\ }p \in \roman{Sing}(J_k,b).$$

(ii) The invariant $w\text{-}\roman{ord}_k:\roman{Sing}(J_k,b) \rightarrow \frac{1}{b}{\Bbb Z}_{\geq
0}$ is a function defined over $\roman{Sing}(J_k,b)$ such that
$$w\text{-}\roman{ord}_k(p) = \frac{\nu_p(\overline{J_k})}{b} \text{\ for\ }p \in
\roman{Sing}(J_k,b)$$ where $\overline{J_k} \subset {\Cal O}_{W_k}$ is the unique ideal
characterized by
$$J_k = I(H_{r+1})^{a_{r+1}} \cdot\cdot\cdot I(H_{r+k})^{a_{r+k}} \cdot
\overline{J_k}.$$
(We note that
$I(H_j)^{a_j}$ is a multi-index notation and denotes
$$I(H_j)^{a_j} = \prod_l I(H_{j,l})^{a_{j,l}}$$
where the $H_{j,l}$ are the smooth irreducible components in $H_j$ and where $a_{j,l} =
\nu_{\eta_{j,l}}(J_k)$ is the order of $J_k$ at the generic point $\eta_{j,l}$ of $H_{j,l}$.) 

(iii) First note that in order to define the invariant $t_k$ we require the following extra
condition $(\heartsuit)$ on the sequence of transformations and smooth morphisms of basic objects:
$$(\heartsuit) \hskip.1in \left\{\aligned &Y_{i-1} \subset
\underline{\roman{Max}}\ w\text{-}\roman{ord}_{i-1} (\subset
\roman{Sing}(J_{i-1},b)) \\
&\text{whenever\ } \pi_i \text{\ is\ a\ transformation\ with\ center\
}Y_{i-1}\\
\endaligned\right\} \text{\ for\ }i = 1, ... , k$$  
where 
$$\align
\underline{\roman{Max}}\ w\text{-}\roman{ord}_{i-1} &= \{p \in
\roman{Sing}(J_{i-1},b);w\text{-}\roman{ord}_{i-1}(p) =
\max w\text{-}\roman{ord}_{i-1}\}\\
\max w\text{-}\roman{ord}_{i-1} &= \max\{w\text{-}\roman{ord}_{i-1}(p);p \in
\roman{Sing}(J_{i-1},b)\}.\\
\endalign$$
Under condition $(\heartsuit)$ it follows that we have inequalities (See Proposition 1-12.)
$$\align
\max w\text{-}\roman{ord}_0 &\geq \max w\text{-}\roman{ord}_1 \geq \cdot\cdot\cdot \\
\max w\text{-}\roman{ord}_{i-1}
&\geq \max w\text{-}\roman{ord}_i \\
\cdot\cdot\cdot \geq \max w\text{-}\roman{ord}_{k-1} &\geq \max
w\text{-}\roman{ord}_k.\\
\endalign$$
Let $k_o$ be the index so that
$$\max w\text{-}\roman{ord}_{k_o-1} > \max w\text{-}\roman{ord}_{k_o} = \cdot\cdot\cdot = \max
w\text{-}\roman{ord}_k.$$ (We let $k_o = 0$ if $\max w\text{-}\roman{ord}_0 = \cdot\cdot\cdot = \max
w\text{-}\roman{ord}_k$.)  Set $E_k = E_k^{-} \cup E_k^{+}$ where $E_k^{-} = \{H_1, ... , H_r, ... ,
H_{r + k_o}\}$ as a subset of $E_k = \{H_1, ... , H_r, ... , H_{r + k_o},
... , H_{r + k}\}$ and where
$E_k^+$ is the complement of
$E_k^-$ in
$E_k$.  (Look at the convention explained in Definition 1-8 (iii).)  

\vskip.1in

The invariant $t_k:\roman{Sing}(J_k,b) \rightarrow \frac{1}{b}{\Bbb Z}_{\geq 0} \times {\Bbb Z}_{\geq
0}$ is a function defined over $\roman{Sing}(J_k,b)$ such that
$$t_k(p) = (w\text{-}\roman{ord}_k(p),n_k(p)) \text{\ for\ }p \in \roman{Sing}(J_k,b)$$
where 
$$n_k(p) = \left\{\aligned
\#\{H_i \in E_k;p \in H_i\} &\text{\ if\ }w\text{-}\roman{ord}_k(p) < \max\ w\text{-}\roman{ord}_k \\
\#\{H_i \in E_k^{-}; p \in H_i\} &\text{\ if\ }w\text{-}\roman{ord}_k(p) = \max\
w\text{-}\roman{ord}_k.\\
\endaligned\right.$$
\endproclaim

\vskip.1in

\proclaim{Remark 1-11}\endproclaim

(i) Both the invariants $w\text{-}\roman{ord}_k$ and $t_k$ depend not only on the basic object $(W_k,
(J_k,b), E_k)$ but also on the sequence, while the invariant $\roman{ord}_k$ is solely determined by
the basic object $(W_k, (J_k,b), E_k)$.

(ii) ($\bold{Why\ exclude\ }I(H_1), ... , I(H_r)\bold{\ from\ the\ definition\ of\
}\overline{J_k}\ ?$) It may look more natural to define
$\overline{J_k}$ to be the unique ideal so that
$$J_k = I(H_1)^{a_1} \cdot\cdot\cdot I(H_r)^{a_r}I(H_{r+1})^{a_{r+1}} \cdot\cdot\cdot
I(H_{r+k})^{a_{r+k}} \cdot \overline{J_k}$$
including $I(H_1), ... , I(H_r)$ in the right hand side of the equation, especially when
one realizes that this definition would make the invariant $w\text{-}\roman{ord}_k$ independent of
the sequence and solely determined by the basic object $(W_k, (J_k,b), E_k)$.  However, it will be
clear in Chapter 4, where we introduce the notion of a general basic object, that the natural domain
of definition for the invariant
$\roman{ord}_k$ is the singular locus of a general basic object, which restricts to the singular loci
of the basic objects forming the charts.  Therefore, if we adopt the above definition allowing
$I(H_1), ... , I(H_r)$ to affect $\overline{J_k}$, some of which may not lie over the singular loci,
then we will not be able to conclude that $w\text{-}\roman{ord}_k$ is a well-defined invariant on the
singular locus of the general basic object.  More concretely, the proof in Definition-Proposition
4-5 (ii) showing that the invariant $w\text{-}\roman{ord}_k$ is determined only in terms of
$\roman{ord}_i \hskip.1in (i = 0, ... , k)$, which are verified to be well-defined invariants on
the singular loci of the general basic objects $({\Cal F}_i,(W_i,E_i))$ and hence that so is the
invariant
$w\text{-}\roman{ord}_k$, will NOT work. 

On a more historical account, $\overline{J_k}$ is called the weak transform (of the ideal $J =
J_0$), in contrast to the strict transform or total transform.  It seems that the letter ``$w$" of
the invariant ``$w\text{-}\roman{ord}$" comes from the word ``weak" transform, indicating
it is the order of the weak transform. 

(iii) ($\bold{Why\ choose\ the\ domain\ to\ be\ }$$\roman{Sing}(J_k,b)$$\bold{\ and\ not\ the\
entire\ }$$W_k$$\bold{\ ?}$) \linebreak
As $\roman{ord}_k$ and
$w\text{-}\roman{ord}_k$ are determined by the orders of the ideals
$J_k$ and
$\overline{J_k}$, respectively, which are defined over the whole variety $W_k$, it may look
artificial at this point to restrict their domains of definition to the singular locus
$\roman{Sing}(J_k,b)$ of the basic object.  It may also look unnecessary to devide the orders by the
positive integer
$b$ to obtain these invariants.  However, when we introduce the notion of a ``general" basic object
(generalizing that of a basic object) where it consists of ``charts" given by many basic objects, it
becomes clear that these invariants are naturally defined only over the singular locus of the
general basic object, which restricts to the singular loci of the basic objects forming the charts,
and that, since the positive integers $b$ may vary from chart to chart, it is necessary to devide
the orders of the ideals by these integers in order for the invariants given on the individual
charts to patch.  We will discuss these issues more in detail in Chapter 4.  

(iv) We are NOT going to use the invariant $\roman{ord}_k$ explicitly in our algorithm for resolution
of singularities, though the order (multiplicity) function is the basis of almost all of our
analysis, and though $\roman{ord}_k$ and $w\text{-}\roman{ord}_k$ coincide for a sequence of
transformations and smooth morphisms of simple basic objects.  (See Remark 3-2 (ii) for the
definition of a simple basic object.)  It is used for the purpose of verifying that the
invariants
$w\text{-}\roman{ord}_k$ given on the individual charts, consisting of basic objects, for a general
basic object patch together in Chapter 4.  It is also used for the purpose of verifying that the
$\Gamma$-invariant can be defined purely in terms of the collection ${\goth C}$ of sequences of
transformations and smooth morphisms represented by a general basic object, free of its presentation
using charts.

(v) ($\bold{Local\ description\ of\ the\ \text{$t$}\text{-}invariant}$) The invariant $t_k$ is indeed
a local one, though superficially it depends on such global information as
$\max w\text{-}\roman{ord}_i$ \linebreak
$(i = 0, ... , k)$.  In fact, under condition
$(\heartsuit)$ it has the following description:  Let $k_{op}$, depending on $p \in W_k$, be the
index so that
$$w\text{-}\roman{ord}_{k_{op}-1}(p_{k_{op}-1}) > w\text{-}\roman{ord}_{k_{op}}(p_{k_{op}}) =
\cdot\cdot\cdot = w\text{-}\roman{ord}_k(p_k) = w\text{-}\roman{ord}_k(p)$$ 
where $p_i$ is the image on $W_i$ of $p = p_k \in W_k$.  (We let $k_{op} = 0$ if \linebreak
$w\text{-}\roman{ord}_0(p_0)
=
\cdot\cdot\cdot = w\text{-}\roman{ord}_k(p_k) = w\text{-}\roman{ord}_k(p)$.)  Then
$$t_k(p) = (w\text{-}\roman{ord}_k(p), n_k(p))$$
where $n_k(p)$ has the description
$$n_k(p) = \# \{H_i \in E_{k,p}^-\}$$
where $E_{k,p}^{-} = \{H_1, ... , H_r, H_{r+1}, ... , H_{r + k_{op}}\}$ as a subset of
\hfill\hfill\linebreak
$E_k =
\{H_1, ... , H_r, H_{r+1}, ... , H_{r + k_{op}}, H_{r + k_{op} + 1}, ... , H_{r + k}\}$, and
hence $t_k(p)$ is locally determined.

In order to see that the global definition given as in Definition 1-10 (iii) and the local
definition given as above coincide, one has only to observe under condition $(\heartsuit)$
that 

if $w\text{-}\roman{ord}_k(p) < \max\ w\text{-}\roman{ord}_k$, then $\pi_i \hskip.1in (i = k_{op} + 1, ... ,
k)$ is either a smooth morphism or a transformation whose center $Y_{i-1}$ is disjoint from
$p_{i-1}$ and hence $E_k = E_{k,p}^-$ in a neighborhood of $p$,

if $w\text{-}\roman{ord}_k(p) = \max\ w\text{-}\roman{ord}_k$, then $k_{op}
\leq k_o$ and in case $k_{op} < k_o$ the morphism $\pi_i \hskip.1in (i = k_{op} + 1, ... , k_o)$
is either a smooth morphism or a transformation whose center $Y_{i-1}$ is disjoint from $p_{i-1}$
and hence $E_k^- = E_{k,p}^-$ in a neighborhood of $p$.

(vi) The definition of $E_k^-$ here is slightly different from the one in the original papers by
Encinas and Villamayor,  where they say ``$E_k^-$ consists of the strict transforms of the
hypersurfaces (irreducible components) in
$E_{k_o}$".  After the index $k_o$, if we blow up along a divisor which is the strict
transform of an irreducible component in $E_{k_o}$, the strict transform of this divisor belongs to
$E_k^-$ according to their definition, while it does not belong to $E_k^-$ according to our
definition and convention in Definition 1.8 (iii).   

However, the difference between the two definitions occurs only
when $w\text{-}\roman{ord}_k = 0$.  In our algorithm for resolution of singularities for (general)
basic objects where we have condition $(\heartsuit)$, as soon as (the maximum of) the invariant
$w\text{-}\roman{ord}_k$ is 0, we apply the method of resolution of singularities specifically
prescribed for monomial (general) basic objects, where the invariant
$t$ plays no role.  Thus this difference has no effect as long as our algorithm for resolution
of singularities is concerned.   

Our choice of the definition is only justified for the virtue of consistency with the
convention in Definition 1.8 (iii) and consistency of notation for a sequence consisting
both of transformations and smooth morphisms.

(vii) ($\bold{The\ origin\ of\ the\ invariant\ }$$t_k$) The definition of $w\text{-}\roman{ord}_k$ is
natural from a view point of achieving principalization by ``extracting" the exceptional
divisors from the total transforms of the original ideals, though the definition of $t_k$ may baffle
us at this point.  The true ingenuity and power of the invariant $t_k$, especially in regard to
choosing appropriate permissible centers, will be revealed in Chapter 5.  A more
down-to-earth explanation of the ``origin" of the invariant
$t_k$ may be found in Chapter 6.

\vskip.1in

\proclaim{Proposition 1-12 (Properties of key invariants)} Let 
$$\align
(W, (J,b), E) = (W_0, (J_0,b), E_0) &\overset{\pi_1}\to{\leftarrow} (W_1, (J_1,b), E_1)
\overset{\pi_2}\to{\leftarrow} \cdot\cdot\cdot \\
(W_{i-1}, (J_{i-1},b), E_{i-1}) &\overset{\pi_i}\to{\leftarrow}
(W_i, (J_i,b), E_i) \\
\cdot\cdot\cdot \overset{\pi_{k-1}}\to{\leftarrow} (W_{k-1}, (J_{k-1},b),
E_{k-1}) &\overset{\pi_k}\to{\leftarrow} (W_k, (J_k,b), E_k) \\
\endalign$$
be a sequence of transformations and smooth morphisms of basic objects.

(i) The invariants $\roman{ord}_k$ and $w\text{-}\roman{ord}_k$ are upper semi-continuous functions.

(ii) Note first that (by condition (ii) in Case T in Definition 1-8)
$$\roman{Sing}(J_{i-1},b) \supset \pi_i(\roman{Sing}(J_i,b)) \text{\ for\ }i = 1, ... ,
k.$$

Suppose that the sequence satisfies condition $(\heartsuit)$.  (See Definition 1-10 (iii).) 

Then for $i = 1, ... , k$ we have inequalities
$$w\text{-}\roman{ord}_{i-1}(\xi_{i-1})
\geq w\text{-}\roman{ord}_i(\xi_i)$$
where $\xi_i \in \roman{Sing}(J_i,b)$ and $\xi_{i-1} =
\pi_i(\xi_i) \in \roman{Sing}(J_{i-1},b)$, which imply
$$\max\ w\text{-}\roman{ord}_{i-1}
\geq \max\ w\text{-}\roman{ord}_i.$$
That is to say, we have
$$\align
\max\ w\text{-}\roman{ord}_0 &\geq \max\ w\text{-}\roman{ord}_1 \geq \cdot\cdot\cdot \\
\max\ w\text{-}\roman{ord}_{i-1} &\geq \max\ w\text{-}\roman{ord}_i \\
\cdot\cdot\cdot \geq \max\ w\text{-}\roman{ord}_{k-1} &\geq \max\ w\text{-}\roman{ord}_k. \\
\endalign$$ 

The invariant $t_k$ is an upper semi-continuous function and for $i = 1, ... , k$ we have
inequalities
$$t_{i-1}(\xi_{i-1}) \geq t_i(\xi_i)$$
where $\xi_i \in \roman{Sing}(J_i,b)$ and $\xi_{i-1} =
\pi_i(\xi_i) \in \roman{Sing}(J_{i-1},b)$, which imply
$$\max\ t_{i-1}
\geq \max\ t_i.$$
That is to say, we have
$$\align
\max\ t_0 &\geq \max\ t_1 \geq \cdot\cdot\cdot \\
\max\ t_{i-1} &\geq \max\ t_i \\
\cdot\cdot\cdot \geq \max\ t_{k-1} &\geq \max\ t_k. \\
\endalign$$ 

\endproclaim

\demo{Proof}\enddemo (i) This is obvious, since the functions $\nu_{J_k}$ and $\nu_{\overline{J_k}}$
are upper semi-continuous functions (cf. Lemma 1-4 (i)).

(ii) First we note the following Lemma.

\proclaim{Lemma 1-13 (Behavior of extensions under transformation)} Let
$$(W_0, (J_0,b), E_0) \overset{\pi_1}\to{\leftarrow} (W_1, (J_1,b), E_1)$$
be a transformation of basic objects, which is the blowup of a center \linebreak
$Y_0 \subset
\roman{Sing}(J_0,b)
\subset W_0$ permissible with respect to $E_0 = \{H_1, ... , H_r\}$.  Let \linebreak
$H_{r+1} =
\pi_1^{-1}(Y_0)$ be the (exceptional) divisor defined by the ideal $I(H_{r+1}) = I(Y_0){\Cal
O}_{W_1}$.

Then we have for $i = 0, 1, ... , b$

$$\align
\Delta^{b-i}(J_0){\Cal O}_{W_1} &\subset I(H_{r+1})^i \\
\frac{1}{I(H_{r+1})^i}\Delta^{b-i}(J_0){\Cal O}_{W_1} &\subset \Delta^{b-i}(J_1).\\
\endalign$$
\endproclaim

\demo{Proof}\enddemo We prove the statements by decreasing induction on $i$.

Suppose $i = b$.  Let $\eta_{r+1,l}$ be the generic point of an irreducible component $H_{r+1,l}$
in $H_{r+1}$, which maps onto the generic point $\theta_{0,l}$ of an irreducible component
$Y_{0,l}$ of the center $Y_0$.

Then
$$\Delta^{b-b}(J_0){\Cal O}_{W_1,\eta_{r+1,l}} = J_0 {\Cal O}_{W_1,\eta_{r+1,l}} \subset
m_{\theta_{0,l}}^b {\Cal O}_{W_1,\eta_{r+1,l}} = I(H_{r+1})^b {\Cal O}_{W_1,\eta_{r+1,l}}$$
since $Y_{0,l} \subset Y_0 \subset \roman{Sing}(J_0,b)$, which implies (since $W_1$ is nonsingular
and hence factorial)
$$\Delta^{b-b}(J_0){\Cal O}_{W_1} = J_0{\Cal O}_{W_1} \subset I(H_{r+1})^b.$$   
Moreover, by definition we have
$$\frac{1}{I(H_{r+1})^b}\Delta^{b-b}(J_0){\Cal O}_{W_1} = \frac{1}{I(H_{r+1})^b}J_0{\Cal O}_{W_1} =
J_1 =
\Delta^{b-b}(J_1).$$

Now suppose that the statements hold for $i \geq j$.  

For any point $\xi_1 \not\in H_{r+1}$, the statements
for $i = j - 1$ clearly hold for the stalks at the point, since $\pi_1$ is an isomorphism at
$\xi_1$ and since $I(H_{r+1})_{\xi_1} = {\Cal O}_{W_1,\xi_1}$.  
 
So we have only to consider the statements for $i = j - 1$ for the stalks at a closed point
$\xi_1 \in H_{r+1}$.  

Set
$\xi_0 =
\pi_1(\xi_1)$.

Consider the completions $\widehat{{\Cal O}_{W,\xi_0}}$ and $\widehat{{\Cal O}_{W_1,\xi_1}}$
with systems of regular parameters $(y = x_1, ... , x_s, x_{s+1} , ... , x_d)$ and
$(y, \frac{x_2}{y}, ... , \frac{x_s}{y}, x_{s+1}, ... , x_d)$, where
$x_1, ... , x_s$ in the former ring define the center $Y_0$ and $y$ in the latter ring
defines the (exceptional) divisor $H_{r+1}$.

It suffices to show that for some generators $\{f\}$ of $\Delta^{b-(j-1)}(J_0)$ in $\widehat{{\Cal
O}_{W,\xi_0}}$, the fractions $\{\frac{f}{y^{j-1}}\}$ belong to $\Delta^{b-(j-1)}(J_1)$ in
$\widehat{{\Cal O}_{W_1,\xi_1}}$.

We take the elements $f$ of $\Delta^{b-j}(J_0)$ and elements of the form $f = D(g)$, where $g \in
\Delta^{b-j}(J_0)$ with $D$ a $k$-derivation on $\widehat{{\Cal O}_{W,\xi_0}}$, as generators
of
$\Delta^{b-(j-1)}(J_0)$.

\vskip.1in

Case: $f \in \Delta^{b-j}(J_0)$.

\vskip.1in

By induction we have
$$\frac{f}{y^j} \in \Delta^{b-j}(J_1) \subset \Delta^{b-(j-1)}(J_1)$$
and hence
$$\frac{f}{y^{j-1}} = y \cdot \frac{f}{y^j} \in \Delta^{b-(j-1)}(J_1).$$

\vskip.1in

Case: $f = D(g)$ where $g \in \Delta^{b-j}(J_0)$ and $D$ is a $k$-derivation on
\linebreak
$\widehat{{\Cal O}_{W,\xi_0}} \cong k[[y = x_1, x_2, \cdot\cdot\cdot, x_s, x_{s+1},
\cdot\cdot\cdot, x_d]]$

\vskip.1in

Note first that $D$ can be extended to a $k$-derivation
$$D:k[[y,\frac{x_2}{y}, \cdot\cdot\cdot, \frac{x_s}{y}, x_{s+1}, \cdot\cdot\cdot, x_d]] \rightarrow
k[[y, \frac{x_2}{y}, \cdot\cdot\cdot, \frac{x_s}{y}, x_{s+1}, \cdot\cdot\cdot, x_d, \frac{1}{y}]]$$
in the obvious way (by Leibniz rule).

We claim that $D' = yD$ is a $k$-derivation on $\widehat{{\Cal O}_{W_1,\xi_1}}$.  In fact, it
suffices to check
$$D'(y), D'(\frac{x_2}{y}), ... , D'(\frac{x_s}{y}), D'(x_{s+1}), ... ,
D'(x_d) \in \widehat{{\Cal O}_{W_1,\xi_1}} \cong k[[y, \frac{x_2}{y}, \cdot\cdot\cdot,
\frac{x_s}{y}, x_{s+1}, \cdot\cdot\cdot, x_d]].$$ 
We clearly see
$$D'(y), D'(x_{s+1}), ... , D'(x_d) \in \widehat{{\Cal O}_{W_1,\xi_1}},$$
while for $m = 2, ... , s$ we have
$$D'(\frac{x_m}{y}) = \frac{yD(x_m) \cdot y - x_m \cdot yD(y)}{y^2} = D(x_m) - \frac{x_m}{y}D(y)
\in \widehat{{\Cal O}_{W_1,\xi_1}}.$$
Now since by induction
$$\frac{g}{y^j} \in \frac{1}{I(H_{r+1})^j}\Delta^{b-j}(J_0){\Cal O}_{W_1} \subset \Delta^{b-j}(J_1)
(\subset \Delta^{b-(j-1)}(J_1)),$$
we have by the claim
$$D'(\frac{g}{y^j}) = \frac{D(g)}{y^{j-1}} - jD(y) \frac{g}{y^j} \in \Delta^{b-(j-1)}(J_1).$$
Therefore, we conclude finally
$$\frac{f}{y^{j-1}} = \frac{D(g)}{y^{j-1}} = D'(\frac{g}{y^j}) + jD(y)\frac{g}{y^j} \in
\Delta^{b-(j-1)}(J_1),$$
completing the proof of the lemma. 

\vskip.2in

We go back to the proof of Proposition 1-12 (ii).

\vskip.1in

When 
$$(W_{i-1}, (J_{i-1},b), E_{i-1}) \overset{\pi_i}\to{\leftarrow} (W_i, (J_i,b), E_i)$$
is a smooth morphism, since $\pi$ is \'etale equivalent to the projection \linebreak
$W_{i-1}
\leftarrow W_{i-1} \times {\Bbb A}^n$ for some $n$, it is obvious that
$$\frac{\nu_{\xi_{i-1}}(\overline{J_{i-1}})}{b} = \frac{\nu_{\xi}(\overline{J_i})}{b}$$
for a point $\xi_i \in W_i$ and its image $\xi_{i-1} = \pi_i(\xi_i) \in W_{i-1}$, and hence
$$\align
\roman{Sing}(J_{i-1},b) &\supset \pi_i(\roman{Sing}(J_i,b)) \\
w\text{-}\roman{ord}_{i-1}(\xi_{i-1}) &= w\text{-}\roman{ord}_i(\xi_i) \text{\ for\ }\xi_i \in
\roman{Sing}(J_i,b).\\
\endalign$$
Since $\max\ w\text{-}\roman{ord}_{i-1} \geq \max\ w\text{-}\roman{ord}_i$ and since $E_i$ is the
inverse image of
$E_{i-1}$ by the smooth morphism, we also have
$$t_{i-1}(\xi_{i-1}) = t_i(\xi_i) \text{\ for\ }\xi_i \in
\roman{Sing}(J_i,b).$$

Thus we have only to deal with the case where $\pi_i$ is a transformation.

\vskip.1in

Since $Y_{i-1} \subset \roman{Sing}(J_{i-1},b)$ as a part of the requirement for the transformation
of basic objects by definition, we clearly have
$$\roman{Sing}(J_{i-1},b) \supset \pi_i(\roman{Sing}(J_i,b)).$$

Condition $(\heartsuit) \hskip.1in Y_{i-1} \subset \underline{\roman{Max}}\
w\text{-}\roman{ord}_{i-1}$ implies that the transformation of the basic objects
$$(W_{i-1}, (J_{i-1},b), E_{i-1}) \overset{\pi_i}\to{\leftarrow} (W_i, (J_i,b), E_i)$$
induces another
$$(W_{i-1}, (\overline{J_{i-1}},c), E_{i-1}) \overset{\pi_i}\to{\leftarrow}
(W_i, ((\overline{J_{i-1}})_1,c), E_i),$$
where $c = b \cdot \max w\text{-}\roman{ord}_{i-1}$.

Let $\xi_i \in \roman{Sing}(J_i,b)$ be a point and $\xi_{i-1} = \pi_i(\xi_i) \in
\roman{Sing}(J_{i-1},b)$.  We have by Lemma 1-13
$$\align
\Delta^{\nu_{\xi_{i-1}}(\overline{J_{i-1}})}(\overline{J_{i-1}})
{\Cal O}_{W_i} &= I(H_{r+i})^{c -
\nu_{\xi_{i-1}}(\overline{J_{i-1}})} \cdot \frac{1}{I(H_{r+i})^{c -
\nu_{\xi_{i-1}}(\overline{J_{i-1}})}}\Delta^{\nu_{\xi_{i-1}}(\overline{J_{i-1}})}(\overline{J_{i-1}})
{\Cal O}_{W_i}
\\ 
&\subset
I(H_{r+i})^{c -
\nu_{\xi_{i-1}}(\overline{J_{i-1}})} \cdot
\Delta^{\nu_{\xi_{i-1}}(\overline{J_{i-1}})}((\overline{J_{i-1}})_1)
\\
&\subset
\Delta^{\nu_{\xi_{i-1}}(\overline{J_{i-1}})}((\overline{J_{i-1}})_1).\\
\endalign$$ 
Since 
$$\Delta^{\nu_{\xi_{i-1}}(\overline{J_{i-1}})}(\overline{J_{i-1}})_{\xi_{i-1}} = {\Cal
O}_{W_{i-1},\xi_{i-1}},$$
we have
$$\Delta^{\nu_{\xi_{i-1}}(\overline{J_{i-1}})}((\overline{J_{i-1}})_1)_{\xi_i} = {\Cal
O}_{W_1,\xi_i},$$ which implies
$$\nu_{\xi_{i-1}}(\overline{J_{i-1}}) \geq \nu_{\xi_i}((\overline{J_{i-1}})_1).$$
Observing $\overline{J_i} = (\overline{J_{i-1}})_1$, we conclude
$$w\text{-}\roman{ord}_{i-1}(\xi_{i-1}) = \frac{\nu_{\xi_{i-1}}(\overline{J_{i-1}})}{b} \geq
\frac{\nu_{\xi_i}(\overline{J_i})}{b} = w\text{-}\roman{ord}_i(\xi_i).$$

\vskip.1in

We want to show that the invariant $t_k$ is an upper semi-continuous function, i.e.,
$$F_{(\alpha,\beta)} = \{p \in \roman{Sing}(J_k,b);t_k(p) =
(w\text{-}\roman{ord}_k(p),n_k(p)) \geq (\alpha,\beta)\}$$ is closed for any $(\alpha,\beta)
\in
\frac{1}{b}{\Bbb Z}_{\geq 0} \times {\Bbb Z}_{\geq 0}$, where the set $\frac{1}{b}{\Bbb
Z}_{\geq 0} \times {\Bbb Z}_{\geq 0}$ is given lexicographical order.

Note that the sets
$$\align
G_{\alpha} &= \{p \in \roman{Sing}(J_k,b); w\text{-}\roman{ord}_k(p) \geq \alpha\} \\
G_{\alpha}^+ &= \{p \in \roman{Sing}(J_k,b); w\text{-}\roman{ord}_k(p) > \alpha\} \\
\endalign$$
are closed, since $w\text{-}\roman{ord}_k$ is an upper semi-continuous function with images in
$\frac{1}{b}{\Bbb Z}_{\geq 0}$.

Therefore, if $\alpha < \max\ w\text{-}\roman{ord}_k$, then
$$F_{(\alpha,\beta)} = G_{\alpha}^+ \cup \cup_{H_{i_1}, \cdot\cdot\cdot, H_{i_{\beta}} \in
E_k}(G_{\alpha} \cap H_{i_1} \cap \cdot\cdot\cdot \cap H_{i_{\beta}}),$$
and if $\alpha = \max\ w\text{-}\roman{ord}_k$, then
$$F_{(\alpha,\beta)} = G_{\alpha}^+ \cup \cup_{H_{i_1}, \cdot\cdot\cdot, H_{i_{\beta}} \in
E_k^-}(G_{\alpha} \cap H_{i_1} \cap \cdot\cdot\cdot \cap H_{i_{\beta}}).$$
In both cases, $F_{(\alpha,\beta)}$ is a closed subset.

\vskip.1in

Finally we show the inequality
$$t_{i-1}(\xi_{i-1}) \geq t_i(\xi_i)$$
for $\xi_i \in \roman{Sing}(J_i,b)$ and its image $\xi_{i-1} = \pi_i(\xi_i) \in
\roman{Sing}(J_{i-1},b)$.

\vskip.1in

From the first part we have
$$w\text{-}\roman{ord}_{i-1}(\xi_{i-1}) \geq w\text{-}\roman{ord}_i(\xi_i).$$

Suppose $w\text{-}\roman{ord}_{i-1}(\xi_{i-1}) > w\text{-}\roman{ord}_i(\xi_i)$.  Then we obviously
have  
$$t_{i-1}(\xi_{i-1}) = (w\text{-}\roman{ord}_{i-1}(\xi_{i-1}),n_{i-1}(\xi_{i-1}) >
(w\text{-}\roman{ord}_i(\xi_i),n_i(\xi_i)) = t_i(\xi_i).$$

Suppose $w\text{-}\roman{ord}_{i-1}(\xi_{i-1}) = w\text{-}\roman{ord}_i(\xi_i)$. 

In case $\max\ w\text{-}\roman{ord}_{i-1} > \max\ w\text{-}\roman{ord}_i$, we have by definition $i_o
= i$ and hence $E_i = E_i^-$.  (See Definition 1-10 (iii) for the meaning of the number $i_o$.) 
Moreover, by condition $(\heartsuit)$ the center $Y_{i-1} \subset \underline{\roman{Max}}\
w\text{-}\roman{ord}_{i-1}$ is disjoint from $\xi_{i-1}$.  Thus $E_{i-1}$ and $E_i$ are identical in
a neighborhood of $\xi_{i-1} = \xi_i$.  Therefore, we conclude
$$\align
n_{i-1}(\xi_{i-1}) &= \#\{H \in E_{i-1};\xi_{i-1} \in H\} \\
&= \#\{H \in E_i;\xi_i \in H\} = \#\{H \in
E_i^-;\xi_i \in H\} = n_i(\xi_i)\\
\endalign$$ 
and hence
$$t_{i-1}(\xi_{i-1}) = t_i(\xi_i).$$ 

In case $\max\ w\text{-}\roman{ord}_{i-1} = \max\ w\text{-}\roman{ord}_i >
w\text{-}\roman{ord}_i(\xi_i) = w\text{-}\roman{ord}_{i-1}(\xi_{i-1})$, by condition $(\heartsuit)$
the center
$Y_{i-1} \subset
\underline{\roman{Max}}\ w\text{-}\roman{ord}_{i-1}$ is disjoint from $\xi_{i-1}$.  Thus $E_{i-1}$
and
$E_i$ are identical in a neighborhood of
$\xi_{i-1} = \xi_i$.  Therefore, we conclude
$$n_{i-1}(\xi_{i-1}) = \#\{H \in E_{i-1};\xi_{i-1} \in H\} = \#\{H \in E_i;\xi_i \in H\} =
n_i(\xi_i)$$ and hence
$$t_{i-1}(\xi_{i-1}) = t_i(\xi_i).$$

In case $\max\ w\text{-}\roman{ord}_{i-1} = \max\ w\text{-}\roman{ord}_i =
w\text{-}\roman{ord}_i(\xi_i) = w\text{-}\roman{ord}_{i-1}(\xi_{i-1})$, we have by definition
$(i-1)_o = i_o$.  Thus the strict transforms of the divisors in $E_{i-1}^-$ of $E_{i-1}$ contain the
divisors in $E_i^-$ of
$E_i$.  Therefore, we conclude
$$n_{i-1}(\xi_{i-1}) = \#\{H \in E_{i-1}^-;\xi_{i-1} \in H\} \geq \#\{H \in E_i^-;\xi_i \in H\} =
n_i(\xi_i)$$ and hence
$$t_{i-1}(\xi_{i-1}) \geq t_i(\xi_i).$$
This completes the proof of Proposition 1-12.

\vskip.1in

\proclaim{Remark 1-14}\endproclaim

(i) The use of Lemma 1-13, whose proof may look non-trivial, if not tricky at first, is
something of an overkill just for the purpose of verifying Proposition 1-12.  

For example, in order to see the inequality 
$$w\text{-}\roman{ord}_{i-1}(\xi_{i-1}) \geq w\text{-}\roman{ord}_i(\xi_i) \text{\ for\ }\xi_i \in
\roman{Sing}(J_i,b)
\text{\ and\ }\xi_{i-1} = \pi_i(\xi_i) \in \roman{Sing}(J_{i-1},b)$$
under the condition $Y_{i-1} \subset \underline{\roman{Max}}\ w\text{-}\roman{ord}_{i-1}$, we have
only to prove the inequality
$$\nu_q(J_0) \geq \nu_p(J_1) \text{\ for\ }p \in W_1 \text{\ and\ }q =
\pi_1(p)$$
for a transformation of basic objects
$$(W_0, (J_0,b), E_0) \overset{\pi_1}\to{\leftarrow} (W_1, (J_1,b), E_1)$$
with a permissible center $Y_0 \subset \roman{Sing}(J_0,b)$ under the condition $b = \max
\nu_{J_0}$.  (See the argument in the proof above for Proposition 1-12 (ii).  We may have to
shrink $W_{i-1}$ when we consider the basic object $(W_{i-1}, (\overline{J_{i-1}},c), E_{i-1})$)  
This can be seen easily using the Taylor expansion expressions of the completions of the local rings
as follows:

We take the completions $\widehat{{\Cal O}_{W_0,q}}$ and $\widehat{{\Cal O}_{W_1,p}}$ with systems
of regular parameters $(y = x_1, ... , x_s, x_{s+1}, ... , x_d)$ and $(y,
\frac{x_2}{y}, ... , \frac{x_s}{y}, x_{s+1}, ... , x_d)$, where $x_1,
... , x_s$ in the former ring define the center $Y_0$ and $y$ in the latter ring defines
the (exceptional) divisor $H_{r+1}$.  (Again we may assume $p \in H_{r+1}$, as the
assertion is obvious otherwise.)  Then 
$$\align
\widehat{{\Cal O}_{W_0,q}} &\cong k[[y = x_1, ... , x_s, x_{s+1}, ... ,
x_d]],\\
\widehat{{\Cal O}_{W_1,p}} &\cong k[[y,
\frac{x_2}{y}, ... , \frac{x_s}{y}, x_{s+1}, ... , x_d]].\\
\endalign$$
where the homomorphism $\pi_1^*:\widehat{{\Cal O}_{W_0,q}} \rightarrow \widehat{{\Cal O}_{W_1,p}}$
corresponds to the obvious inclusion $k[[y = x_1, ... , x_s, x_{s+1}, ... ,
x_d]] \hookrightarrow k[[y,
\frac{x_2}{y}, ... , \frac{x_s}{y}, x_{s+1}, ... , x_d]]$.  The condition \linebreak
$Y_0
\subset
\roman{Sing}(J_0,b)$ translates into the statement that for any $f
\in J_0$ all the monomials appearing in the Taylor expansion of $f$ should contain $x_1, ...
, x_s$ total of degree at least $b$.  From this it follows immediately $\frac{f}{y^b} \in
\widehat{{\Cal O}_{W_1,p}}$ and $J_1$ is indeed well-defined.

For $f \in J_0$ with $\nu_q(f) = \nu_q(J_0) = b$, there exists in the Taylor expansion of $f$ a
monomial containing $x_1, ... , x_s$ precisely of total degree $b$.  Then it is clear
that in the Taylor expansion of $\frac{f}{y^b}$ there appears a monomial of degree $\leq b$,
and hence
$\nu_q(f) = b \geq \nu_p(\frac{f}{y^b}) \geq \nu_p(J_1)$. 

\vskip.1in

Lemma 1-13, however, will be crucial when we define certain ideals in terms of the extensions and
analyze their behavior under transformations (cf. Lemma 3-1 and Claim 3-5).

(ii) If we use the local description of the $t$-invariant (cf. Remark 1-11 (v)), then the
proof for the statement $t_{i-1}(\xi_{i-1}) \geq t_i(\xi_i)$ becomes simpler: 

As before, we
have only to deal with the case where $\pi_i$ is a transformation.  

When
$w\text{-}\roman{ord}_{i-1}(\xi_{i-1}) > w\text{-}\roman{ord}_i(\xi_i)$, we obviously have
$$t_{i-1}(\xi_{i-1}) = (w\text{-}\roman{ord}_{i-1}(\xi_{i-1}), n_{i-1}(\xi_{i-1})) >
(w\text{-}\roman{ord}_i(\xi_i), n_i(\xi_i)) = t_i(\xi_i).$$
When $w\text{-}\roman{ord}_{i-1}(\xi_{i-1}) = w\text{-}\roman{ord}_i(\xi_i)$, we have
$(i-1)_{o\xi_{i-1}} = i_{o\xi_i}$.  (See Remark 1-11 (v) for the meaning of the numbers
$(i-1)_{o\xi_{i-1}}$ and $i_{o\xi_i}$.)  Thus the strict transforms of the divisors in
$E_{i-1,\xi_{i-1}}^-$ of $E_{i-1}$ contain the divisors in $E_{i,\xi_i}^-$ of $E_i$. 
Therefore, we conclude
$$n_{i-1}(\xi_{i-1}) = \#\{H \in E_{i-1,\xi_{i-1}}^-;\xi_{i-1} \in H\} \geq \#\{H \in
E_{i,\xi_i}^-;\xi_i
\in H\} = n_i(\xi_i)$$ and hence
$$t_{i-1}(\xi_{i-1}) = (w\text{-}\roman{ord}_{i-1}(\xi_{i-1}), n_{i-1}(\xi_{i-1})) \geq
(w\text{-}\roman{ord}_i(\xi_i), n_i(\xi_i)) = t_i(\xi_i).$$

\newpage

\hskip.3in $\bold{CHAPTER\ 2.\ RESOLUTION\ OF\ SINGULARITIES}$ \linebreak
${}\hskip.9in$ $\bold{OF\ MONOMIAL\ BASIC\ OBJECTS}$

\vskip.1in

In this chapter, we present an algorithm for resolution of singularities of the ``monomial" basic
objects.  This settles the problem of resolution of singularities for any basic object (in a sequence of transformations and smooth morphisms of basic
objects) with (the
maximum of) the invariant
$w\text{-}\roman{ord}$ (cf. Definition 1-10) being equal to 0, since it is easily reduced to
the problem of resolution of singularities for some monomial basic object.

\vskip.1in

\proclaim{Definition 2-1 (Monomial basic object)} Let $B = (W, (J,b), E)$ be a basic object
of dimension $\dim W = d$ with
$E = \{H_1, ... , H_r\}$ (cf. Definition 1-6 and Note 1-7).  We say $B$ is a
monomial basic object if
$$J = I(H_1)^{a_1} \cdot\cdot\cdot I(H_r)^{a_r}.$$
(See Definition 1-10 (ii) for the meaning of the multi-index notation $I(H_j)^{a_j}$.) 
\endproclaim

\vskip.1in

\proclaim{Remark 2-2}\endproclaim

If $\max\ w\text{-}\roman{ord}_k = 0$ for a basic object $(W_k,
(J_k,b), E_k)$ (in a sequence of transformations and smooth morphisms described as in Definition
1-8), then \linebreak
$J_k = I(H_{r+1})^{a_{r+1}} \cdot\cdot\cdot I(H_{r+k})^{a_{r+k}}$, i.e., $\overline{J_k} = {\Cal
O}_{W_k}$ in a neighborhood of the singular locus $\roman{Sing}(J_k,b)$.  Recall that the
definition of $\overline{J_k}$ only involves $I(H_{r+1}), ... , I(H_{r+k})$ but not $I(H_1),
... , I(H_r)$ (cf. Remark 1-11 (ii)).  

However, once $J_k = I(H_{r+1})^{a_{r+1}} \cdot\cdot\cdot
I(H_{r+k})^{a_{r+k}}$ and hence $(W_k, (J_k,b), E_k)$ is a monomial basic object (in a
neighborhood of $\roman{Sing}(J_k,b)$), then our algorithm of resolution of singularities
depends only on $(W_k, (J_k,b), E_k)$ or only on $(W_k, (J_k, b), \{H_{r+1}, ... , H_{r+k}\})$,
and is independent of the sequence.  This is why we characterize a monomial basic object $(W,
(J,b), E)$ by the condition \linebreak
$J = I(H_1)^{a_1} \cdot\cdot\cdot I(H_r)^{a_r}$, involving
all
$I(H_1), ... , I(H_r)$.

\vskip.1in

\proclaim{Definition 2-3 (The $\Gamma$-invariant on a monomial basic object)} Let $(W, (J,b),
E)$ be a monomial basic object of dimension $\dim W = d$ with $J = I(H_1)^{a_1} \cdot\cdot\cdot
I(H_r)^{a_r}$ where $E = \{H_1, ... , H_r\}$.  The
invariant $\Gamma:\roman{Sing}(J,b) \rightarrow {\Bbb Z}_{\geq - d} \times \frac{1}{b}{\Bbb
Z}_{\geq 0}
\times {{\Bbb Z}_{\geq 0}}^d$ is a function defined over $\roman{Sing}(J,b)$ such that
$$\Gamma(p) = (\Gamma_1(p),\Gamma_2(p),\Gamma_3(p)) \text{\ for\ }p \in \roman{Sing}(J,b)$$
where
$$\align
- \Gamma_1(p) &= \min \{n; \exists j_1, ... , j_n \text{\ s.t.\ }a_{j_1}(p) +
\cdot\cdot\cdot + a_{j_n}(p) \geq b, p \in H_{j_1} \cap \cdot\cdot\cdot \cap H_{j_n}\} \\
\Gamma_2(p) &= \max \{\frac{a_{j_1}(p) + \cdot\cdot\cdot + a_{j_n}(p)}{b};n = - \Gamma_1(p),
a_{j_1}(p) +
\cdot\cdot\cdot + a_{j_n}(p) \geq b, \\
&p \in H_{j_1} \cap \cdot\cdot\cdot \cap H_{j_n}\} \\
\Gamma_3(p) &= \max \{(j_1, ... , j_n);n = -
\Gamma_1(p), \Gamma_2(p) =
\frac{a_{j_1}(p) + \cdot\cdot\cdot + a_{j_n}(p)}{b}, \\
&p \in H_{j_1} \cap \cdot\cdot\cdot \cap
H_{j_n}, j_1 \geq \cdot\cdot\cdot \geq j_n\}\\
&\text{with\ the maximum\ taken\ with\ respect\ to\ the\
lexicographical\ order\ given\ to\ }{{\Bbb Z}_{\geq 0}}^d.\\
&\text{We\ identify\ }(j_1, ... , j_n) \text{\ with\ }(j_1, ... , j_n, 0,
... , 0) \in {{\Bbb Z}_{\geq 0}}^d.\\
\endalign$$  

We order the values of $\Gamma$ according to the lexicographical order given to \linebreak
${\Bbb Z}_{\geq - d}
\times
\frac{1}{b}{\Bbb Z}_{\geq 0}
\times {{\Bbb Z}_{\geq 0}}^d$.
\endproclaim

\vskip.1in

\proclaim{Remark 2-4}\endproclaim

The number $- \Gamma_1(p)$ is the minimum of the codimensions of the components
given as the intersections of the hypersurfaces (in $E$) in
$\roman{Sing}(J,b)$ with order (at the generic points of the components) $\geq b$ and containing
the point
$p$.  We take its negaive $\Gamma_1(p) = - (- \Gamma_1(p))$ for the first factor of the invariant
$\Gamma$.  The moral here is: 
 
The less the codimension is, the worse the locus is (and hence to
be blown up earlier).

The number $\Gamma_2(p)$ is the maximum of the orders (devided by $b$) along (the generic
points of) the components containing $p$ of codimension $- \Gamma_1(p)$.  The moral here is: 
 
The more the order is, the worse the locus is (and hence to be blown up earlier). 

The third factor $\Gamma_3(p)$ is the ``tie breaker" given by the indices of the divisors in $E$.  Without this third factor, two irreducible components of the maximum locus of the pair
$(\Gamma_1, \Gamma_2)$ may meet at a point $p$ and hence not be smooth.  This third factor
guarantees that the maximum locus of the invariant $\Gamma$ is smooth.

\vskip.1in

\proclaim{Proposition 2-5 (Canonical center for a monomial basic object)} Let \linebreak
$B_0 =
(W_0, (J_0,b), E)$ be a monomial basic object and
$$\align
\underline{\roman{Max}}\ \Gamma_{B_0} &= \{p \in \roman{Sing}(J_0,b); \Gamma(p) = \max\ \Gamma_{B_0}\}
\\
\max\ \Gamma_{B_0} &= \max \{\Gamma(p); p \in \roman{Sing}(J_0,b)\}. \\
\endalign$$
Observe that $Y_0 = \underline{\roman{Max}}\ \Gamma _{B_0} \subset \roman{Sing}(J_0,b)$ is a
smooth closed subset permissible with respect to
$E_0$.  Take the transformation of basic objects with center $Y_0$
$$B_0 = (W_0, (J_0,b), E_0) \overset{\pi_1}\to{\leftarrow} B_1 = (W_1, (J_1,b), E_1).$$
Then $B_1$ is a monomial basic object and we have
$$\max\ \Gamma_{B_0} > \max\ \Gamma_{B_1}.$$
\endproclaim

\demo{proof}\enddemo  The proof is straightforward and left to the reader as an exercise.

\vskip.1in

\proclaim{Corollary 2-6 (Resolution of singularities of a monomial basic object)} Let $(W,
(J,b), E)$ be a monomial basic object.  Then there exists a sequence of transformations of
monomial basic objects
$$\align
B_0 = (W, (J,b), E) = (W_0, (J_0,b), E_0) &\overset{\pi_1}\to{\leftarrow} B_1 = (W_1, (J_1,b),
E_1)
\overset{\pi_2}\to{\leftarrow} \cdot\cdot\cdot \\
B_{i-1} = (W_{i-1}, (J_{i-1},b), E_{i-1}) &\overset{\pi_i}\to{\leftarrow}
B_i = (W_i, (J_i,b), E_i) \\
\cdot\cdot\cdot \overset{\pi_{k-1}}\to{\leftarrow} B_{k-1} = (W_{k-1}, (J_{k-1},b),
E_{k-1}) &\overset{\pi_k}\to{\leftarrow} B_k = (W_k, (J_k,b), E_k) \\
\endalign$$
which represents resolution of singularities, i.e., 
$$\roman{Sing}(J_k,b) = \emptyset,$$ 
where
$$B_{i-1} = (W_{i-1}, (J_{i-1},b), E_{i-1}) \overset{\pi_i}\to{\leftarrow}
B_i = (W_i, (J_i,b), E_i) \hskip.1in i = 1, ... , k$$ 
are the transformations with centers $Y_{i-1} = \underline{\roman{Max}}\ \Gamma_{B_{i-1}}$.
\endproclaim

\demo{Proof}\enddemo It follows immediately from Proposition 2-5 and the observation that
the set ${\Bbb Z}_{\geq - d} \times \frac{1}{b}{\Bbb Z}_{\geq 0}
\times {{\Bbb Z}_{\geq 0}}^d$ satisfies the descending chain condition (i.e., it admits no
infinite strictly decreasing sequence).

\vskip.1in

\proclaim{Corollary 2-7 (Resolution of singularities of a basic object with \linebreak
$\roman{max}\
w\text{-}\roman{ord} = 0$)} Let 
$$\align
B_0 = (W, (J,b),E) = (W_0, (J_0,b), E_0) &\overset{\pi_1}\to{\leftarrow} B_1 = (W_1, (J_1,b),
E_1)
\overset{\pi_2}\to{\leftarrow} \cdot\cdot\cdot \\
B_{i-1} = (W_{i-1}, (J_{i-1},b), E_{i-1}) &\overset{\pi_i}\to{\leftarrow}
B_i = (W_i, (J_i,b), E_i) \\
\cdot\cdot\cdot \overset{\pi_{k-1}}\to{\leftarrow} B_{k-1} = (W_{k-1}, (J_{k-1},b),
E_{k-1}) &\overset{\pi_k}\to{\leftarrow} B_k = (W_k, (J_k,b), E_k) \\
\endalign$$
be a sequence of transformations and smooth morphisms of basic objects.  

Suppose $\max w\text{-}\roman{ord}_k = 0.$

Then there exists an open
neighborhood $\roman{Sing}(J_k,b) \subset V_k \subset W_k$ of $\roman{Sing}(J_k,b)$ such that
$(V_k, (J_k|_{V_k},b), E_k|_{V_k})$ is a monomial basic object.

Take the sequence of
transformations of monomial basic objects as described in Corollary 2-6  
$$B_k|_{V_k} \overset{\pi_{k+1}|_{V_{k+1}}}\to{\leftarrow} B_{k+1}|_{V_{k+1}} \cdot\cdot\cdot
\overset{\pi_{k+N-1}|_{V_{k+N-1}}}\to{\leftarrow} B_{k+N-1}|_{V_{k+N-1}}
\overset{\pi_{k+N}|_{V_{k+N}}}\to{\leftarrow} B_{k+N}|_{V_{k+N}}$$ with centers $Y_{i-1} =
\underline{\roman{Max}}\ \Gamma_{B_{i-1}|_{V_{i-1}}}$ for $i = k+1,
... , k+N$, which represents resolution of singularities of the monomial basic object
$B_k|_{V_k}$.

The sequence can naturally be expanded
to a sequence of transformations of the original basic object $B_k = (W_k, (J_k,b), E_k)$
$$B_k = (W_k, (J_k,b), E_k) \overset{\pi_{k+1}}\to{\leftarrow} \cdot\cdot\cdot
\overset{\pi_{k+N}}\to{\leftarrow} B_{k+N} = (W_{k+N}, (J_{k+N},b), E_{k+N})$$
with the same centers $Y_{i-1}$, which
repesents resolution of singularities of the basic object $B_k$, i.e.,
$$\roman{Sing}(J_{k+N},b) = \emptyset.$$ 
Moreover, the expanded sequence is independent of the choice of the open neighborhood $V_k$.
\endproclaim

\demo{Proof}\enddemo Since $\max\ w\text{-}\roman{ord}_k = 0$ and since the order function
$\nu_{\overline{J_k}}$ is upper semi-continuous, writing (cf. Definition 1-10 (ii))
$$J_k = I(H_{r+1})^{a_{r+1}} \cdot\cdot\cdot I(H_{r+k})^{a_{r+k}} \cdot \overline{J_k}$$
we conclude that $S = \roman{Supp}({\Cal O}_{W_k}/\overline{J_k})$ is a closed subset disjoint
from
$\roman{Sing}(J_k,b)$.  Take $V_k$ to be any open subset such that $\roman{Sing}(J_k,b)
\subset V_k \subset W_k \setminus S$.  Then $B_k|_{V_k}$ is a monomial basic object.

Remark that the centers $\underline{\roman{Max}}\ \Gamma_{B_i|_{V_i}}$ for the sequence of
transformations decribed as in Corollary 2-6 are all over $\roman{Sing}(J_k,b)$ and hence
that the sequence can be expanded as claimed.  

As the centers chosen according to Proposition 2-5 are easily seen to be independent of the choice of
$V_k$, so is the sequence. 

We note that under the specified transformations $\roman{max}\ w\text{-}\roman{ord}$ remains zero, i.e.,
$$\max\ w\text{-}\roman{ord}_{k+i} = 0 \text{\ for\ }i = 0, ... , k-1$$
and hence that the invariant $\Gamma_{k+i}$ is well-defined.

\newpage

$$\bold{CHAPTER\ 3.\ KEY\ INDUCTIVE\ LEMMA}$$

\vskip.1in

In this chapter, we prove the key inductive lemma , which reduces the problem of resolution
of singularities of a ``simple" basic object (See Remark 3-2.) of dimension $d$ to that of
``charts" consisting of basic objects of dimension $d-1$, provided the existence of smooth
hypersurfaces (in the open subsets which give rise to the charts) which cover the singularities of
the simple basic object and which cross transversally with the specified boundary divisors
of the original simple basic object of dimension $d$.  This lemma will become the
prototype of the inductive argument which follows, leading to the notion of general basic
objects.  The shortcomings of the key inductive lemma, namely the requirements for the
original basic object to be simple (while the resulting basic objects of dimension
$(d-1)$ in charts may not be simple) and for the existence of certain smooth
hypersurfaces with the transversality condition, will be overcome in the ultimate
inductive algorithm toward the general solution of resolution of singularities in Chapter
5 via the brilliant use of the
$t$-invariant. 

The underlying idea of the key inductive lemma may be most transparent when we consider
an ideal $\langle f \rangle \subset k[x_1, ... , x_{d-1}, x_d]$ generated by an element
$f$ of the form
$$f = x_d^n + c_{n-2}x_d^{n-2} + \cdot\cdot\cdot + c_1x_d + c_0$$
where the coefficient $c_{n-1}$ of the term $x_d^{n-1}$ is zero after the Tshirnhausen
transformation and where the coefficients $c_{n-2}, ... , c_0$ depend only on the
variables $x_1, ... , x_{d-1}$.  We reduce the problem of resolution of singularities of
$f$ on ${\Bbb A}^d$ of dimension $d$ (around the origin) to that of the coefficients
$c_{n-2}, ... , c_1, c_0$ on $\{x_d = 0\}$ of dimension $d-1$, which is a hypersurafce of
maximal contact.  (See Remark 1-5 for the notion of a hypersurface of maximal contact. 
Check that the condition $c_{d-1} = 0$ immediately implies that $\{x_d = 0\}$ indeed
is.  See Remark 3-7 for more details.)

\proclaim{Lemma 3-1 (Key inductive lemma)} Let $B = (W, (J,b), E)$ be a basic object
with an open covering $\{W^{\lambda}\}_{\lambda \in \Lambda}$ satisfying the
following conditions: for each $\lambda \in \Lambda$, there exists a smooth hypersurface $W_h^{\lambda} \subset
W^{\lambda}$, embedded as a closed subscheme, such that

1. $I(W_h^{\lambda}) \subset \Delta^{b-1}(J)|_{W^{\lambda}} \hskip.1in (\text{and\ hence\
}W_h^{\lambda} \supset \roman{Sing}(J,b) \cap {W^{\lambda}})$, and  

2. $W_h^{\lambda}$ is permissible with respect to
$E \cap W^{\lambda}$, and $W_h^{\lambda}$ is not contained in  $E$, i.e., $W_h^{\lambda}
\not\subset E$.

\vskip.1in

Then, $R(1)(\roman{Sing}(J,b))$ denoting the union of irreducible components in
$\roman{Sing}(J,b)$ of codimension one (i.e., of dimension $\dim W - 1$), we have the
following:

\vskip.1in

$\bold{Case\ A}$: $R(1)(\roman{Sing}(J,b)) \neq \emptyset$.

\vskip.1in

In this case, the set $R(1)(\roman{Sing}(J,b))$ is smooth, open and closed in
$\roman{Sing}(J,b)$ (i.e., a union of smooth connected components of $\roman{Sing}(J,b)$
disjoint from each other).  Condition 2 guarantees that it is also permissible with
respect to $E$.

Take
$$(W, (J,b), E) = (W_0, (J_0,b), E_0) \leftarrow (W_1, (J_1,b), E_1)$$
to be the transformation with center $Y_0 = R(1)(\roman{Sing}(J,b))$.

Then
$$R(1)(\roman{Sing}(J_1,b)) = \emptyset.$$

\vskip.1in

$\bold{Case\ B}$: $R(1)(\roman{Sing}(J,b)) = \emptyset$.

\vskip.1in

In this case, let ${\goth C}$ be the collection of all the sequences of transformations
and smooth morphisms of pairs with specified closed subsets
$$(F_0, (W_0,E_0)) \leftarrow \cdot\cdot\cdot \leftarrow (F_k, (W_k,E_k))$$
induced from the sequences of transformations and smooth morphisms of basic objects
$$(W, (J,b), E) = (W_0, (J_0,b), E_0) \leftarrow \cdot\cdot\cdot \leftarrow (W_k, (J_k,b),
E_k)$$
where the specified closed subsets are the singular loci of the corresponding basic
objects, i.e.,
$$F_i = \roman{Sing}(J_i,b) \hskip.1in \text{\ for\ }i = 0, 1, ... , k.$$

Then with respect to the open covering $\{W^{\lambda}\}_{\lambda \in \Lambda}$ we can
construct the following data ${\Cal D}_{\lambda}$ for each ${\lambda}$:

\vskip.1in

(i) $j_0^{\lambda}:(\widetilde{W_0^{\lambda}},\widetilde{E_0^{\lambda}}) \hookrightarrow
(W_0^{\lambda},E_0^{\lambda}) = (W^{\lambda},E_0 \cap W^{\lambda})$ is an immersion of
pairs, that is to say, $\widetilde{W_0^{\lambda}} = W_h^{\lambda} \hookrightarrow W^{\lambda}$ is a
closed immersion of a $(\dim W - 1)$-dimensional smooth variety $\widetilde{W_0^{\lambda}}$ into
$W^{\lambda}$,
$\widetilde{W_0^{\lambda}}$ is permissible with respect to
$E_0^{\lambda}$, $\widetilde{W_0^{\lambda}}$ is not contained in $E_0^{\lambda}$, i.e.,
$\widetilde{W_0^{\lambda}}
\not\subset E_0^{\lambda}$, and $\widetilde{E_0^{\lambda}} = E_0^{\lambda} \cap
\widetilde{W_0^{\lambda}}$, 

(ii) a basic object $(\widetilde{W_0^{\lambda}}, ({\goth
a}_0^{\lambda},b^{\lambda}), \widetilde{E_0^{\lambda}}) = (W_h^{\lambda}, (C(J^{\lambda}),b!),
E_h^{\lambda})$ where $E_h^{\lambda} = E \cap W_h^{\lambda}$

(See the proof below for the definition of the ideal $C(J^{\lambda})$.),

\vskip.1in

satisfying the following
conditions (GB-0,1,2,3):

\vskip.1in

(GB-0) The trivial sequence consisting only of $(F_0, (W_0,E_0))$ is in the collection ${\goth
C}$, i.e.,
$$(F_0, (W_0,E_0)) \in {\goth C}$$
and
$$F_0 = \roman{Sing}(J_0,b) = \cup \roman{Sing}({\goth a}_0^{\lambda},b^{\lambda}) \text{\
with\ }F_0
\cap W_0^{\lambda} = \roman{Sing}({\goth a}_0^{\lambda},b^{\lambda}).$$

\vskip.1in

(GB-1) For any sequence of transformations and smooth morphisms in the collection
${\goth C}$
$$(F_0, (W_0,E_0)) \leftarrow \cdot\cdot\cdot \leftarrow (F_k, (W_k,E_k))$$
there corresponds for each $\lambda$ a sequence of transformations (with the same centers) and (the
same) smooth morphisms (obtained by taking the Cartesian products)
$$(\widetilde{W_0^{\lambda}}, ({\goth a}_0^{\lambda},b^{\lambda}),
\widetilde{E_0^{\lambda}}) \leftarrow \cdot\cdot\cdot \leftarrow (\widetilde{W_k^{\lambda}}, ({\goth
a}_k^{\lambda},b^{\lambda}),
\widetilde{E_k^{\lambda}})$$
with the natural immersions
$$\CD
(W_0^{\lambda},E_0^{\lambda}) @<<< \cdot\cdot\cdot @<<< (W_k^{\lambda},E_k^{\lambda}) \\
@AAA @. @AAA \\
(\widetilde{W_0^{\lambda}},\widetilde{E_0^{\lambda}}) @<<< \cdot\cdot\cdot @<<<
(\widetilde{W_k^{\lambda}},\widetilde{E_k^{\lambda}})
\endCD$$
and we have
$$F_i = \roman{Sing}(J_i,b) = \cup \roman{Sing}({\goth a}_i^{\lambda},b^{\lambda}) \text{\
with\ }F_i
\cap W_i^{\lambda} = \roman{Sing}({\goth a}_i^{\lambda},b^{\lambda})$$
for $i = 0, 1, ... , k$.

\vskip.1in

(We note here that in the above clause ``there corresponds ...", it is required that whenever
$(W_{i-1},E_{i-1}) \overset{\pi_i}\to{\leftarrow} (W_i,E_i)$ is a transformation with center
\linebreak
$Y_{i-1} \subset W_{i-1}$, the center $Y_{i-1}$ is permissible for each
$(\widetilde{W_{i-1}^{\lambda}}, ({\goth a}_{i-1}^{\lambda},b^{\lambda}),
\widetilde{E_{i-1}^{\lambda}})$, i.e., \linebreak
$Y_{i-1} \cap W_{i-1}^{\lambda} \subset
\widetilde{W_{i-1}^{\lambda}}$,
$Y_{i-1} \cap W_{i-1}^{\lambda}$ is permissible with respect to
$\widetilde{E_{i-1}^{\lambda}}$, and \linebreak
$Y_{i-1}
\cap W_{i-1}^{\lambda} \subset \roman{Sing}({\goth a}_{i-1}^{\lambda},b^{\lambda})$.)

\vskip.1in
 
(GB-2) Let
$$(F_0, (W_0,E_0)) \leftarrow \cdot\cdot\cdot \leftarrow (F_k, (W_k,E_k))$$
be a sequence of transformations and smooth morphisms in ${\goth C}$ and
$$\{(\widetilde{W_0^{\lambda}}, ({\goth a}_0^{\lambda},b^{\lambda}),
\widetilde{E_0^{\lambda}}) \leftarrow \cdot\cdot\cdot \leftarrow (\widetilde{W_k^{\lambda}}, ({\goth
a}_k^{\lambda},b^{\lambda}),
\widetilde{E_k^{\lambda}})\}$$
the corresponding sequences (indexed by $\lambda \in \Lambda$) of transformations and smooth
morphisms as in (GB-1).

\vskip.1in

We take a morphism of pairs $(W_k,E_k) \overset{\pi_{k+1}}\to{\leftarrow} (W_{k+1},E_{k+1})$ which is either in
Case T or Case S.

\vskip.1in

Case T: $(W_k,E_k) \overset{\pi_{k+1}}\to{\leftarrow} (W_{k+1},E_{k+1})$ is a transformation
with center $Y_k \subset W_k$, satisfying the condition that $Y_k$ is permissible for each
$(\widetilde{W_k^{\lambda}}, ({\goth a}_k^{\lambda},b^{\lambda}),
\widetilde{E_k^{\lambda}})$, i.e., \linebreak
$Y_k \cap W_k^{\lambda} \subset
\widetilde{W_k^{\lambda}}$,
$Y_k \cap W_k^{\lambda}$ is permissible with respect to $\widetilde{E_k^{\lambda}}$, and
\linebreak
$Y_k
\cap W_k^{\lambda} \subset \roman{Sing}({\goth a}_k^{\lambda},b^{\lambda})$.

\vskip.1in

Case S: $(W_k,E_k) \overset{\pi_{k+1}}\to{\leftarrow} (W_{k+1},E_{k+1})$ is a smooth morphism.

\vskip.2in

Then we have the following assertions on the extension of the original sequence:

\vskip.1in

Case T: Take for each $\lambda$ the corresponding transformation of basic objects with center $Y_k \cap
W_k^{\lambda}$
$$(\widetilde{W_k^{\lambda}}, ({\goth a}_k^{\lambda},b^{\lambda}),
\widetilde{E_k^{\lambda}}) \overset{\pi_{k+1}^{\lambda}}\to{\leftarrow}
(\widetilde{W_{k+1}^{\lambda}}, ({\goth a}_{k+1}^{\lambda},b^{\lambda}),
\widetilde{E_{k+1}^{\lambda}}).$$

In this case, $Y_k$ is permissible with respect to $(W_k, (J_k,b), E_k)$, i.e., $Y_k$ is
permissible with respect to $E_k$ and $Y_k \subset \roman{Sing}(J_k,b)$, and we have the
induced transformation of basic objects
$$(W_k, (J_k,b), E_k) \overset{\pi_{k+1}}\to{\leftarrow} (W_{k+1}, (J_{k+1},b), E_{k+1}).$$

Then
$$F_{k+1} := \cup \roman{Sing}({\goth
a}_{k+1}^{\lambda},b^{\lambda})$$ is a closed subset of $W_{k+1}$ with
$$F_{k+1} \cap W_{k+1}^{\lambda} = \roman{Sing}({\goth a}_{k+1}^{\lambda},b^{\lambda}),$$
and the extended sequence belongs to ${\goth C}$, i.e., $F_{k+1} = \roman{Sing}(J_{k+1},b)$ and
$$(F_0, (W_0,E_0)) \leftarrow \cdot\cdot\cdot \leftarrow (F_k, (W_k,E_k)) \leftarrow (F_{k+1},
(W_{k+1},E_{k+1})) \in {\goth C}.$$

\vskip.1in

Case S: Take for each $\lambda$ the corresponding smooth morphism of basic objects
$$(\widetilde{W_k^{\lambda}}, ({\goth a}_k^{\lambda},b^{\lambda}),
\widetilde{E_k^{\lambda}}) \overset{\pi_{k+1}^{\lambda}}\to{\leftarrow}
(\widetilde{W_{k+1}^{\lambda}}, ({\goth a}_{k+1}^{\lambda},b^{\lambda}),
\widetilde{E_{k+1}^{\lambda}}).$$
where
$$\widetilde{W_{k+1}^{\lambda}} = \widetilde{W_k^{\lambda}} \times_{W_k}W_{k+1}, {\goth
a}_{k+1}^{\lambda} = {\goth a}_k^{\lambda}{\Cal O}_{\widetilde{W_{k+1}^{\lambda}}},
\widetilde{E_{k+1}^{\lambda}} = {\pi_{k+1}^{\lambda}}^{-1}(\widetilde{E_k^{\lambda}})$$ 
and where $\pi_{k+1}^{\lambda}$ is the
projection onto the first factor.

Take the smooth morphism of basic objects
$$(W_k, (J_k,b), E_k) \overset{\pi_{k+1}}\to{\leftarrow} (W_{k+1}, (J_{k+1},b), E_{k+1}).$$

Then
$$F_{k+1} := \cup \roman{Sing}({\goth
a}_{k+1}^{\lambda},b^{\lambda})$$ is a closed subset of $W_{k+1}$ with
$$F_{k+1} \cap W_{k+1}^{\lambda} = \roman{Sing}({\goth a}_{k+1}^{\lambda},b^{\lambda}),$$
and the extended sequence belongs to ${\goth C}$, i.e., $F_{k+1} =
\roman{Sing}(J_{k+1},b)$ and
$$(F_0, (W_0,E_0)) \leftarrow \cdot\cdot\cdot \leftarrow (F_k, (W_k,E_k))
\leftarrow (F_{k+1}, (W_{k+1},E_{k+1})) \in {\goth C}.$$

\vskip.1in

(GB-3) There exists $c\ (= b!) \in {\Bbb N}$ such that $c \geq b^{\lambda} \hskip.1in \forall
\lambda$.
\endproclaim

\vskip.1in

\proclaim{Remark 3-2}\endproclaim

(i) The main point of the key inductive lemma is that the
problem of resolution of singularities for a basic object $(W,(J,b),E) = (W_0,(J_0,b),E_0)$ (under conditions 1
and 2) of dimension
$d$, i.e., the problem of finding a sequence of transformations
$$(W_0,(J_0,b),E_0) \leftarrow \cdot\cdot\cdot (W_k,(J_k,b),E_k) \text{\ with\
}\roman{Sing}(J_k,b) = \emptyset,$$
which gives rise to a sequence of transformations in ${\goth C}$
$$(F_0,(W_0,E_0)) \leftarrow \cdot\cdot\cdot \leftarrow (F_k,(W_k,E_k)) \text{\ with\ }F_k =
\emptyset,$$
can be reduced to
that of charts of basic objects $\{(\widetilde{W_0^{\lambda}},({\goth
a}_0^{\lambda},b^{\lambda}),\widetilde{E_0^{\lambda}})\}$ of dimension $d-1$, if we could find the
global centers which are permissible with respect to all the local charts.  

\vskip.1in

That is to say, starting with the trivial sequence (condition (GB-0)), we build up a resolution sequence of a
(simple) basic object of dimension $d$ by extending the trivial one via the repeated use of condition (GB-2),
based upon the resolution sequence of the charts of basic objects in dimension $d-1$, which is obtained by
induction.

\vskip.1in

Together with condition (GB-2) which characterizes the sequences in the collection ${\goth C}$ and with
condition (GB-3) which ensures the boundedness of the integers $b^{\lambda}$ in order to gurantee the
descending chain condition of our invariants, (GB-0,1,2,3) will be used as
the defining conditions for general basic objects in Definition 4-1.

\vskip.1in

We discuss the shortcomings of the key inductive lemma toward a complete inductive algorithm in (ii), (iii), and
(iv) below.

\vskip.1in

(ii) ($\bold{Simple\ basic\ object}$) Condition 1 implies (cf. Lemma 1-4) that 
$$\nu_p(\Delta^{b-1}(J)) = 1 \hskip.1in \forall p \in V(\Delta^{b-1}(J)),$$
which is equivalent to the characterization of what we call a $\bold{simple\ basic\
object}$: 

\vskip.1in

A basic object $(W,(J,b),E)$ is called simple if
$$\nu_p(J) = b = b_{\max} = \max\{\nu_q(J);q \in W\} \hskip.1in \forall p \in
\roman{Sing}(J,b) (= V(\Delta^{b-1}(J))).$$
Conversely, for a simple basic object $(W, (J,b), E)$ it is easy to find an open covering
$\{W^{\lambda}\}$ and smooth hypersurfaces $W_h^{\lambda} \subset W^{\lambda}$ with the
property $I(W_h^{\lambda}) \subset \Delta^{b-1}(J)$ and hence satisfying condition 1. 
(See Remark 1-5 where we discuss the primitive but fundamental idea toward the inductional
argument via the notion of a hypersurface of maximal contact.) 

Thus the key inductive lemma restrictively applies only to simple basic objects, for which an open covering
satisfying condition 1 comes almost for free.

(iii) Even if we start from a simple basic object of dimension $d$, the resulting basic
objects $\{(\widetilde{W_0^{\lambda}}, ({\goth
a}_0^{\lambda},b^{\lambda}), \widetilde{E_0^{\lambda}})\}$ of dimension $(d-1)$, in Case B, are
almost never simple.  Therefore, even though we call this lemma with the adjective
``inductive", it is not clear at this point (until Chapter 5) how this induction would
actually work. 

(iv) Condition 2 is more problematic, if one tries to see the inductive structure of a
possible proof for resolution of singularities in a naive way suggested by the above lemma: 

It
is NOT true that we can find an open covering $\{W^{\lambda}\}$ which satisfies conditions
1 and 2 for an arbitrary simple basic object. 

\vskip.1in

Take $(W, (J,b), E) = ({\Bbb A}^2 = \roman{Spec}[x,y], (\langle x^2\rangle),2), \{y - x^2 = 0\})$.  It is
easy to check that
$\nu_p(J) = 2 = b \hskip.1in \forall p \in \roman{Sing}(J,b) = \{x = 0\}$ and hence that it is
a simple basic object.  Observe $\Delta^{b-1}(J) = (x)$ and hence by condition 2
the smooth hypersurface $W_h^{\lambda}$ has to contain $\{x = 0\} \cap W^{\lambda}$, no matter
how we choose an open covering
$\{W^{\lambda}\}$.  However, in any open subset (in the covering) containing $(0,0)$ the hypersurface $\{x = 0\}$
does not cross $E =
\{y - x^2\}$ transversally, failing to satisfy condition 2.

\vskip.1in

(v) The shortcomings of the key inductive lemma expressed in the above (ii), (iii), and
(iv) and the problem of how to find the global centers permissible with respect the local charts expressed in
(i), will be so elegantly and beautifully resolved in Chapter 5 via the power of the
$t$-invariant.  See also Chapter 6 for a more-down-to-earth interpretation of this inductive procedure.

(vi) Since $E = \{H_1, ... , H_r\}$ is a collection of hypersurfaces and $E^{\lambda}$ is a
collection of the restrictions of the hypersurfcaes to the open subset $W^{\lambda}$, it is more
appropriate logically to write $E^{\lambda} = \{H_1 \cap W^{\lambda}, ... , H_r \cap W^{\lambda}\}$
or
$E^{\lambda} = \{H_1|_{W^{\lambda}}, ... , H_r|_{W^{\lambda}}\}$ than to write $E^{\lambda} = E
\cap W^{\lambda}$ or $E^{\lambda} = E|_{W^{\lambda}}$, which we use, however, by abuse of and
for simplicity of notation.  We also write $\widetilde{E_0^{\lambda}} = E_0^{\lambda} \cap
\widetilde{W_0^{\lambda}}$ instead of \linebreak
$\widetilde{E_0^{\lambda}} = \{H_1 \cap
\widetilde{W_0^{\lambda}}, ... , H_r \cap \widetilde{W_0^{\lambda}}\}$.

\vskip.1in

\demo{Proof of Lemma 3-1}\enddemo

$\bold{Case\ A}$: $R(1)(\roman{Sing}(J,b)) \neq \emptyset$.

\vskip.1in

Since condition 1 implies $\roman{Sing}(J,b) \cap W^{\lambda} \subset
W_h^{\lambda}$, we have
$$\align
&\text{either\ }R(1)(\roman{Sing}(J,b)) \cap W^{\lambda} = \emptyset \\
&\text{\ or\
}R(1)(\roman{Sing}(J,b)) \cap W^{\lambda} \text{\ open\ and\ closed\ in\ } W_h^{\lambda}.\\
\endalign$$
Therefore, we conclude that $R(1)(\roman{Sing}(J,b))$ is smooth, open and closed in
$\roman{Sing}(J,b)$ and that it is permissible with respect to $E$, since $W_h^{\lambda}$ is smooth
and since $W_h^{\lambda}$ is permissible with respect to $E \cap W^{\lambda}$ for each $\lambda \in
\Lambda$.

Take
$$(W, (J,b), E) = (W_0, (J_0,b), E_0) \leftarrow (W_1, (J_1,b), E_1)$$
to be the transformation with center $Y_0 = R(1)(\roman{Sing}(J,b))$.

Since $Y_0$ is a smooth divisor of codimension one, the ambient space remains unchanged, i.e.,
$W_0 \overset{\sim}\to{\leftarrow} W_1$.  What changes is the ideal, from $J_0$ to $J_1$.  By
definition, we have
$$J_0{\Cal O}_{W_1} = I(H_{r+1})^b \cdot J_1$$
and hence
$$\nu_p(J_1) = b - b = 0 < b$$
for any codimension one point $p$ which is the generic point of an irreducible component in $H_{r+1} =
Y_0 = R(1)(\roman{Sing}(J,b))$.  Therefore, $\roman{Sing}(J_1,b) \subset \roman{Sing}(J,b)$
has no codimension one point, i.e.,
$$R(1)(\roman{Sing}(J_1,b)) = \emptyset.$$

\vskip.1in

$\bold{Case\ B}$: $R(1)(\roman{Sing}(J,b)) = \emptyset$.

\vskip.1in

We take with respect to the given open covering $\{W^{\lambda}\}$ a collection of basic objects
$$\{(\widetilde{W_0^{\lambda}}, ({\goth a}^{\lambda},b^{\lambda}), \widetilde{E_0^{\lambda}}) =
(W_h^{\lambda}, (C(J^{\lambda}),b!), W_h^{\lambda} \cap E)\},$$
where $J^{\lambda} = J|_{W^{\lambda}}$ and where the $\bold{coefficient\ ideal}$
$C(J^{\lambda})$ is defined to be

$$C(J^{\lambda}) := \Sigma_{i=0}^{b-1}\Delta^i(J^{\lambda})^{\frac{b!}{b-i}}{\Cal
O}_{W_h^{\lambda}}.$$

(The definition of the coefficient ideal forms the technical core of the proof of the key inductive
lemma.  For a background motivation,  the reader is invited to look at Remark 3-7 regarding the Tschirmhausen
transformation.)  

Note that in order to verify conditions (GB-0,1,2,3) for $(W, (J,b), E)$ with \linebreak
$\{(\widetilde{W_0^{\lambda}}, ({\goth a}^{\lambda},b^{\lambda}),
\widetilde{E_0^{\lambda}})\}$, it suffices to verify conditions (GB-0,1,2,3) for
$(W^{\lambda}, (J|_{W^{\lambda}},b), E|_{W^{\lambda}})$ with
$(\widetilde{W_0^{\lambda}}, ({\goth a}^{\lambda},b^{\lambda}), \widetilde{E_0^{\lambda}})$ for
each
$\lambda$.
 
Therefore, we drop the superscript $\lambda$ from the following argument.  That is to say,
we assume $(W, (J,b), E)$ is a basic object with a smooth hypersurface $W_h \subset W$, embedded as a closed
subscheme, such that

1. $I(W_h) \subset \Delta^{b-1}(J)$,

2. $W_h$ is permissible with respect to $E$, and $W_h$ is not contained in $E$, i.e., $W_h
\not\subset E$.

\vskip.1in

(Note, however, that we do keep the superscript $\lambda$ in $b^{\lambda}$, since it may actually
be different from $b$.)

\vskip.1in

\proclaim{Claim 3-3}: $\roman{Sing}(J,b) = \roman{Sing}(C(J),b!)$.\endproclaim

\demo{Proof}\enddemo We observe that
$$\align
& p \in \roman{Sing}(J,b)\ (\subset W_h) \\
\Longleftrightarrow \hskip.1in & \nu_p(\Delta^i(J)) \geq b - i \text{\ for\ }i = 0,
... , b - 1\\
\Longleftrightarrow \hskip.1in & \nu_p({\Delta^i(J)}^{\frac{b!}{b - i}}) \geq b! \text{\ for\
}i = 0,
... , b - 1\\
\Longrightarrow \hskip.1in & \nu_p(\Sigma_{i = 0}^{b - 1}{\Delta^i(J)}^{\frac{b!}{b - i}}{\Cal
O}_{W_h}) \geq b! \\
\Longleftrightarrow \hskip.1in & p \in \roman{Sing}(C(J),b!)\\
\endalign$$
and that
$$\align
& p \not\in \roman{Sing}(J,b) \\
\Longleftrightarrow \hskip.1in & \Delta^i(J)_p = {\Cal O}_{W,p} \text{\ for\ some\ }i = 0,
... , b - 1 \\
\Longrightarrow \hskip.1in & \Sigma_{i = 0}^{b - 1}\Delta^i(J)^{\frac{b!}{b - i}}{\Cal
O}_{W_h,p} = {\Cal O}_{W_h,p} \\
\Longrightarrow \hskip.1in & p \not\in \roman{Sing}(C(J),b!).\\
\endalign$$

This proves the assertion.

\vskip.1in

\proclaim{Claim 3-4 (Giraud's Lemma)} Let $(W, (J,b), E)$ be a basic object and $W_h \subset
W$ be a smooth hypersurface, embedded as a closed subscheme, satisfying conditions 1 and 2 as above.

\vskip.1in
 
Case T: $(W, (J,b), E) \overset{\pi_1}\to{\leftarrow} (W_1, (J_1,b), E_1)$ is a transformation
with permissible center $Y \subset W$ for $(W, (J,b), E)$, i.e., $Y$ is permissible with
respect to
$E$ and \linebreak
$Y
\subset \roman{Sing}(J,b)$.

In this case, $(W_h)_1 \subset W_1$, the strict transform of $W_h$, is a smooth
hypersurface satisfying conditions 1 and 2.

\vskip.1in
 
Case S: $(W, (J,b), E) \overset{\pi_1}\to{\leftarrow} (W_1, (J_1,b), E_1)$ is a smooth morphism.

In this case, $(W_h)_1 = \pi_1^{-1}(W_h) \subset W_1$ is a
smooth hypersurface satisfying conditions 1 and 2.

\vskip.1in

\endproclaim

\demo{Proof}\enddemo The assertions in Case S are obvious, since $\pi_1$ is \'etale equivalent to
the projection $W \leftarrow W \times {\Bbb A}^n$ for some $n$.  In Case T, since
$Y$ is permissible with respect to $E$ and $Y \subset \roman{Sing}(J,b) \subset W_h$, where $W_h$ is
permissible with respect to $E$, we see that $Y$ is permissible with respect to $\{H_1,
... , H_r, W_h\}$.  Therefore, we conclude that \linebreak
$H_{r+1} = \pi_1^{-1}(Y)$ is
permissible with respect to
$\{H_1, ... , H_r, (W_h)_1\}$, where
$H_i$ denotes the strict transform of $H_i$
(which is denoted by the same letter by abuse of notation) and hence that
$(W_h)_1$ is permissible with respect to $E_1 = \{H_1, ... , H_r, H_{r+1}\}$.  This shows
condition 2 for $(W_h)_1$.

Moreover, Lemma 1-13 implies
$$I((W_h)_1) = \frac{1}{I(H_{r+1})}I(W_h){\Cal O}_{W_1}
\subset \frac{1}{I(H_{r+1})}\Delta^{b-1}(J) \subset \Delta^{b-1}(J_1),$$
and hence that
$$(W_h)_1 \supset \roman{Sing}(J_1,b).$$
This shows condition 1 for $(W_h)_1$.

\vskip.1in

We go back to the analysis of $\bold{Case\ B}$.

\vskip.1in

Let
$$(W, (J,b), E) = (W_0, (J_0,b), E_0) \leftarrow \cdot\cdot\cdot \leftarrow (W_k, (J_k,b),
E_k)$$
be a sequence of transformations and smooth morphisms of basic objects.

Claim 3-4 implies that there corresponds a sequence of transformations and smooth morphisms of pairs
with the natural immersions
$$\CD
(W_0,E_0) @<<< \cdot\cdot\cdot @<<< (W_k,E_k) \\
@AAA @. @AAA \\
(\widetilde{W_0},\widetilde{E_0}) @<<< \cdot\cdot\cdot @<<<
(\widetilde{W_k},\widetilde{E_k})
\endCD$$
where
$$\align
&\widetilde{W_0} = (W_0)_h = W_h, \text{and\ inductively}\\
&\widetilde{W_{i-1}} = (W_{i-1})_h \subset W_{i-1} \text{\ is\ a\ smooth\ hypersurface\
satisfying\ conditions\ 1\ and\ 2,\ and\ }\\
&\widetilde{W_i} = (W_i)_h = ((W_{i-1})_h)_1 \subset W_i \\
\endalign$$
such that
$$\roman{Sing}(J_i,b) \subset \widetilde{W_i} \text{\ for\ }i = 0, 1, ... , k.$$

\proclaim{Claim 3-5} The following assertions hold inductively for $i = 0, 1, ... , k$:

\vskip.1in

(i) $\roman{Sing}(J_i,b) = \roman{Sing}(C(J_i),b!) = \roman{Sing}(C(J_0)_i,b!).$

(ii) (The assertion (ii) is void when $i = 0$.)

Case T: $(W_{i-1}, (J_{i-1},b), E_{i-1}) \overset{\pi_i}\to{\leftarrow} (W_i, (J_i,b), E_i)$ is
a transformation with center $Y_{i-1} \subset \roman{Sing}(J_{i-1},b) \subset W_{i-1}$ and $\pi_i^{-1}(Y_{i-1}) =
H_i$.  

(Note that by induction $Y_{i-1} \subset \roman{Sing}(J_{i-1},b) = \roman{Sing}(C(J_0)_{i-1},b!)$ and
hence that we have an induced transformation of basic objects 
$$(\widetilde{W_{i-1}},{\goth a}_{i-1} =
C(J_0)_{i-1},\widetilde{E_{i-1}}) \leftarrow (\widetilde{W_i},{\goth a}_i =
C(J_0)_i,\widetilde{E_i}).)$$

\vskip.1in

The ideal $[\Delta^{b-j}(J_0)]_{i-1}{\Cal O}_{W_i}$ is divisible by $I(H_i)^j$, and we
set	 
$$[\Delta^{b-j}(J_0)]_i = \frac{1}{I(H_i)^j}[\Delta^{b-j}(J_0)]_{i-1}{\Cal O}_{W_i} \text{\
for\ }j = 1, ... , b$$
and we have an equality
$$C(J_0)_i = \Sigma_{j = 1}^b[\Delta^{b-j}(J_0)]_i^{\frac{b!}{j}}{\Cal O}_{\widetilde{W_i}}.$$

\vskip.1in

Case S: $(W_{i-1}, (J_{i-1},b), E_{i-1}) \overset{\pi_i}\to{\leftarrow} (W_i, (J_i,b), E_i)$ is
a smooth morphism.  

(We have an induced smooth morphism of basic objects 
$$(\widetilde{W_{i-1}},{\goth
a}_{i-1} = C(J_0)_{i-1},\widetilde{E_{i-1}}) \leftarrow (\widetilde{W_i},{\goth a}_i =
C(J_0)_i,\widetilde{E_i}).)$$

\vskip.1in

In this case, we set 
$$[\Delta^{b-j}(J_0)]_i = [\Delta^{b-j}(J_0)]_{i-1}{\Cal O}_{W_i} =
\frac{1}{I(H_i)^j}[\Delta^{b-j}(J_0)]_{i-1}{\Cal O}_{W_i} \text{\
for\ }j = 1, ... , b$$ 
where we use the convention $H_i = \emptyset$ and $I(H_i) = {\Cal O}_{W_i}$, and we have an
equality
$$C(J_0)_i = \Sigma_{j = 1}^b[\Delta^{b-j}(J_0)]_i^{\frac{b!}{j}}{\Cal O}_{\widetilde{W_i}} =
C(J_0)_{i-1}{\Cal O}_{\widetilde{W_i}}.$$

\vskip.1in

Moreover, both in Case T and in Case S, we have an inclusion 
$$[\Delta^{b-j}(J_0)]_i \subset \Delta^{b-j}(J_i) \text{\ for\ } j = 1, ... , b.$$

\vskip.1in

We understand by definition
$$[\Delta^{b-j}(J_0)]_0 = \Delta^{b-j}(J_0) \text{\ and\ }C(J_0)_0 = C(J_0).$$

\vskip.1in

(iii) At any closed point $\xi_i \in \roman{Sing}({\goth
a}_i,b^{\lambda}) = \roman{Sing}(C(J_0)_i,b!)$, considered as a closed point in $W_i$ ($\xi_i \in
\roman{Sing}({\goth
a}_i,b^{\lambda}) \subset \widetilde{W_i} \subset W_i$), there exists a system of regular parameters $z_i,
x_{i,1}, ... , x_{i,d-1}$, where
$d =
\dim W_i$, such that

\vskip.1in

\hskip.2in (a) $I(\widetilde{W_i}) = \langle z_i \rangle$,

\hskip.2in (b) setting $R_i = \widehat{{\Cal O}_{W_i,\xi_i}}, \overline{R_i} = \widehat{{\Cal
O}_{\widetilde{W_i},\xi_i}}$ there is a set of generators $\{f_i^{(\sigma)}\}$ for the ideal
$J_iR_i$
$$f_i^{(\sigma)} = \Sigma_{\alpha}a_{i,\alpha}^{(\sigma)}z_i^{\alpha} \hskip.1in \text{where}
\hskip.1in a_{i,\alpha}^{(\sigma)}
\in k[[x_{i,1}, ... , x_{i,d-1}]] \subset k[[z_i,x_{i,1}, ... , x_{i,d-1}]] = R_i$$
so that
$$(a_{i,\alpha}^{(\sigma)})^{\frac{b!}{b-\alpha}} \in C(J_0)_i\overline{R_i} \text{\ for\ all\
}\alpha
\text{\ with\ }\alpha < b.$$
\endproclaim

\vskip.1in

\proclaim{Remark 3-6}\endproclaim

The claim may look technical at first sight.  The calamity is
that the operation of taking extensions and that of taking transformations of an ideal (i.e.,
taking the pull-back divided by some suitable multiple of the ideal defining the exceptional divisor)
do not commute under a transformation of basic objects
$(W_0,(J_0,b),E)
\overset{\pi_1}\to{\leftarrow} (W_1,(J_1,b),E_1)$, that is to say, in general
$$\frac{1}{I(H_{r+1})^j}\Delta^{b-j}(J_0){\Cal O}_{W_1}
\neq
\Delta^{b-j}(J_1),$$  
where $I(H_{r+1})$ is the ideal defining the exceptional divisor $H_{r+1}$ for $\pi_1$.  (If the
equality were to hold, we would have $C(J_0)_1 = C(J_1)$, from which the assertions in
$\bold{Case\ B}$ would have followed directly.)  The purpose
of the claim is to overcome this calamity using Lemma 1-13, which analyzes the behavior of
extensions under transformations.

\vskip.1in

\demo{Proof of Claim 3-5}\enddemo

(i) The assertion (i) is an immediate consequence of the assertions (ii) and (iii) as follows:

\vskip.1in

Observe first that Claim 3-3 and Claim 3-4 imply
$$\roman{Sing}(J_i,b) = \roman{Sing}(C(J_i),b!).$$
The ``moreover" part of the assertion (ii) implies $C(J_0)_i \subset C(J_i)$ (Note that the inclusion
obviously holds even when $i = 0$.) and hence
$$\roman{Sing}(C(J_i),b!) \subset \roman{Sing}(C(J_0)_i,b!).$$ 
On the other hand, the assertion (iii) implies that for $p \in \roman{Sing}(C(J_0)_i,b!)$ we have
$$\nu_p(a_{i,\alpha}^{(\sigma)}) \geq b - \alpha \text{\ and\ hence\ }\nu_p(f_i^{(\sigma)})
\geq b,$$
and therefore we have
$$\roman{Sing}(C(J_0)_i,b!) \subset \roman{Sing}(J_i,b),$$
completing the argument for the assertion (i).

\vskip.1in

We prove the assertions (ii) and (iii) by induction on $i$.

\vskip.1in

For $i = 0$, the assertion (iii) is obvious from the definition, while the assertion (ii) is void.  (Note
that $a_{0,\alpha}^{(\sigma)} =\frac{1}{\alpha!}(\frac{\partial}{\partial z_0})^{\alpha}f_0^{(\sigma)}|_{z_0 = 0}
\in
\Delta^{\alpha}(J_0)\overline{R_0}$ and hence \linebreak
$(a_{0,\alpha}^{(\sigma)})^{\frac{b!}{b-\alpha}} \in
C(J_0)\overline{R_0}$.)

\vskip.1in

Assuming that the assertions hold for $i \leq s$,  we prove the assertions for $i = s + 1$ by
induction.

\vskip.1in

(ii) In Case S, since $\pi_{s+1}$ is \'etale equivalent to the projection $W_s \leftarrow W_s \times
{\Bbb A}^n$ for some $n$, all the assertions immediately hold by induction.  So we will concentrate
our consideration on Case T where
$(W_s, (J_s,b), E_s)
\overset{\pi_{s+1}}\to{\leftarrow} (W_{s+1},(J_{s+1},b), E_{s+1})$ is a transformation with center
$Y_s \subset W_s$ and $H_{s+1} =
\pi_{s+1}^{-1}(Y_s)$.

\vskip.1in

By the inductional hypothesis, we have
$$[\Delta^{b-j}(J_0)]_s \subset \Delta^{b-j}(J_s) \text{\ for\ }j = 1, ... , b.$$ 
Lemma 1-13 implies that $\Delta^{b-j}(J_s){\Cal O}_{W_{s+1}}$ is divisible by $I(H_{s+1})^j$
and that
$$\frac{1}{I(H_{s+1})^j}\Delta^{b-j}(J_s){\Cal O}_{W_{s+1}} \subset \Delta^{b-j}(J_{s+1}).$$
Therefore, we conclude that $[\Delta^{b-j}(J_0)]_s{\Cal O}_{W_{s+1}}$ is also divisible by
$I(H_{s+1})^j$ and that
$$\align
[\Delta^{b-j}(J_0)]_{s+1} &= \frac{1}{I(H_{s+1})^j}[\Delta^{b-j}(J_0)]_s \\
&\subset \frac{1}{I(H_{s+1})^j}\Delta^{b-j}(J_s) \subset \Delta^{b-j}(J_{s+1}).\\
\endalign$$

\vskip.1in

(iii) Let $\xi_{s+1} \in \roman{Sing}(C(J_0)_{s+1},b!) \subset \widetilde{W_{s+1}}$ be a closed
point and
$\xi_s =
\pi_{s+1}(\xi_s)$ be its image in $\widetilde{W_s} \subset W_s$.

We may only consider the case where $\xi_{s+1} \in H_{s+1}$, since otherwise the assertions are
automatic by the inductional
hypothesis as $\pi_{s+1}$ is isomorphic in a neighborhood of $\xi_{s+1}$.

By taking a change of variables involving only $x_{s,l}$'s ($l = 1, ... , c = \dim_{\xi_s} Y_s$) in the system
of regular parameters for $R_s$ given by the inductional hypothesis for (iii) where $Y_s$ is defined by the
ideal $\langle z_s, x_{s,1}, ... , x_{s,c-1} \rangle$, we may assume that
$R_{s+1} = \widehat{{\Cal O}_{W_{s+1},\xi_{s+1}}}$ has a system of regular parameters \hfill
\linebreak
$z_{s+1}, x_{s+1,1}, ... , x_{s+1,c-1}, x_{s+1,c}, ... , x_{s+1,d-1}$ with $d =
\dim W_{s+1} =
\dim W_s$ such that
$$\align
I(H_{s+1})_{\xi_{s+1}} &= \langle x_{s+1,1} \rangle \\
x_{s+1,1} &= x_{s,1} \\
I(\widetilde{W_{s+1}})_{\xi_{s+1}} &= \langle z_{s+1} \rangle \\
z_{s+1} &= \frac{z_s}{x_{s,1}} \\
x_{s+1,l} &= \frac{x_{s,l}}{x_{s,1}} \text{\ for\ }l = 2, ... , c-1\\
x_{s+1,l} &= x_{s,l} \text{\ for\ } l = c, ... , d-1.\\
\endalign$$

Take a set of generators $\{f_{s+1}^{(\sigma)}\}$ for $J_{s+1}R_{s+1}$ where
$$f_{s+1}^{(\sigma)} = \frac{f_s^{(\sigma)}}{x_{s,1}^b} =
\Sigma_{\alpha}a_{s+1,\alpha}^{(\sigma)}z_{s+1}^{\alpha}$$
so that
$$a_{s+1,\alpha}^{(\sigma)} = \frac{a_{s,\alpha}^{(\sigma)}}{x_{s,1}^{b-\alpha}}.$$
Therefore, we conclude
$$\align
(a_{s+1,\alpha}^{(\sigma)})^{\frac{b!}{b-\alpha}} &=
(\frac{a_{s,\alpha}^{(\sigma)}}{x_{s,1}^{b-\alpha}})^{\frac{b!}{b-\alpha}} \\
&= (a_{s,\alpha}^{(\sigma)})^{\frac{b!}{b-\alpha}} \cdot (\frac{1}{x_{s,1}})^{b!} \\
&\in C(J_0)_s\overline{R_s} \cdot (\frac{1}{x_{s,1}})^{b!} \cdot \overline{R_{s+1}}\
(\text{by\ inductional\ hypothesis})\\
&= \Sigma_{j = 1}^b[\Delta^{b-j}(J_0)]_s^{\frac{b!}{j}}{\Cal O}_{\widetilde{W_s}} \cdot
\overline{R_s} \cdot (\frac{1}{x_{s,1}})^{b!} \cdot \overline{R_{s+1}}\\
&= \Sigma_{j = 1}^b\{\frac{1}{I(H_{s+1})^j}[\Delta^{b-j}(J_0)]_s\}^{\frac{b!}{j}}{\Cal
O}_{\widetilde{W_{s+1}}} \cdot \overline{R_{s+1}} \\
&= C(J_0)_{s+1}\overline{R_{s+1}}.\\
\endalign$$
This completes the proof of Claim 3-5.

\vskip.1in

We go back to the proof of Lemma 3-1.

\vskip.1in

Condition (GB-0) is obvious.  Condition (GB-1) is an immediate consequence of Claim 3-5 (i),
and so is condition (GB-2).  Condition (GB-3) is obvious from the construction.

\vskip.1in

This completes the proof of Lemma 3-1 (the key inductive lemma).

\vskip.1in

\proclaim{Remark 3-7}\endproclaim

(i) Consider a basic object $(W, (J,b), E)$ where
$$\left\{\aligned
W &= {\Bbb A}^d = \roman{Spec}\ k[x_1, ... , x_{d-1}, x_d] \\
J &= \langle f\rangle \text{\ with\ }\\
f &= x_d^n + c_{n-2}x_d^{n-2} + \cdot\cdot\cdot + c_1x_d + c_0 \text{\ with\
}c_i
\in k[x_1, ... , x_{d-1}] \\
b &= n\\
E &= \emptyset.\\
\endaligned\right.$$

We would like to emphasize that $f$ is already in the form after a \linebreak
$\bold{Tschirnhausen\ transformation}$, which can be always carried out over a field of characteristic zero,
i.e., 
$$\text{the\ coefficient\ $c_{n-1}$\ of\ the\ term\ $x_d^{n-1}$\ is\ equal\ to\ $0$\ 
in\  $f$.}$$ 
This implies that (Recall $b = n$.)
$$\frac{1}{(b-1)!}\frac{\partial^{b-1} f}{\partial x_d^{b-1}} = x_d \in \Delta^{b-1}(J).$$

Set
$$W_h = V(x_d).$$

Then conditions 1 and 2 of Lemma 3-1 are clearly satisfied.

\vskip.1in

Observe
$$\align
p \in \roman{Sing}(J,b) &\Longleftrightarrow x_d(p) = 0\quad \&\quad \nu_p(c_i) \geq b - i
\text{\ for\ }i = 0, ... , b-1\\ 
&\Longleftrightarrow x_d(p) = 0\quad \&\quad \text{\ for\ }i = 0, ... , b-1\\
&\hskip.3in \nu_p(\langle c_i\rangle) \geq b - i, \nu_p(\Delta^1\langle c_{i-1}\rangle)) \geq b
- i, ... , \nu_p(\Delta^i\langle c_0\rangle) \geq b - i\\ 
&\Longleftrightarrow p \in W_h\quad
\&\quad
\nu_p(\Sigma_{i=0}^{b-1}\{\langle c_i\rangle +
\Delta^1\langle c_{i-1}\rangle + \cdot\cdot\cdot + \Delta^i\langle
c_0\rangle\}^{\frac{b!}{b-i}}) \geq b!.\\
\endalign$$

This observation might give a justification for calling
$$\Sigma_{i = 0}^{b-1} \Delta^i(J)^{\frac{b!}{b-i}}{\Cal O}_{W_h} = \Sigma_{i=0}^{b-1}\{\langle
c_i\rangle +
\Delta^1\langle c_{i-1}\rangle + \cdot\cdot\cdot + \Delta^i\langle
c_0\rangle\}^{\frac{b!}{b-i}}$$
the coefficient ideal.

\vskip.1in

(ii) Let $(W,(J,b,E)$ be a (simple) basic object with a smooth hypersurface \linebreak
$W_h
\subset W$ satisfying conditions 1 and 2:

1. $I(W_h) \subset \Delta^{b-1}(J)$,

2. $W_h$ is permissible with respect to $E$ and $W_h \not\subset E$.

\vskip.1in

Since
$$\roman{Sing}(J,b) = V(\Delta^{b-1}(J)) = V(\Delta^{b-1}(J))|_{W_h},$$
it might look plausible (and simpler) to consider the basic object
$$(W_h,(D(J),1),E \cap W_h) \text{\ with\ }D(J) = \Delta^{b-1}(J){\Cal O}_{W_h},$$
instead of
$$(W_h,(C(J),b!),E \cap W_h) \text{\ with\ }C(J) = \Sigma_{i=0}^{b-1}\Delta^i(J){\Cal
O}_{W_h},$$
as a candidate for the basic object of dimension one less in order for the key
inductive lemma to work.

This alternative definition, however, does NOT work.

\vskip.1in

Look at a simple basic object $(W,(J,b),E)$ where
$$\left\{\aligned
W &= {\Bbb A}^2 = \roman{Spec}\ k[x,y] \\
J &= \langle x^2 - y^3 \rangle \\
b &= 2 \\
E &= \emptyset \\
\endaligned\right.$$
with a smooth hypersurface $W_h = \{x = 0\}$, clearly satisfying conditions 1 and 2 of Lemma 3-1.

If we take the transformation of basic objects
$$(W,(J,b),E) = (W_0,(J_0,b),E_0) \leftarrow (W_1,(J_1,b),E_1),$$
which is the blowup of the origin, then in the chart with coordinate system \linebreak
$(t
= \frac{x}{y},y)$ we have $J_1 = \langle t^2 - y \rangle$ and hence
$$\roman{Sing}(J_1,b) = \roman{Sing}(J_1,b)|_{(W_h)_1} = \emptyset.$$
On the other hand, if we look at the corresponding transformation of basic objects
$$\align
(W_h,(D(J) = \Delta^{b-1}(J){\Cal O}_{W_h},1),E \cap W_h) &=
((W_h)_0,(D(J)_0,b),(E_h)_0) \\
&\leftarrow ((W_h)_1,(D(J)_1,1),(E_h)_1),\\
\endalign$$
then we have $W_h = (W_h)_0 = (W_h)_1 = {\Bbb A}^1 = \roman{Spec}\ k[y]$ with $D(J)_0 =
\langle y^2\rangle$ and $D(J)_1 = \langle y \rangle$, and hence
$$\roman{Sing}(D(J)_1,1) \neq \emptyset.$$
Therefore, we have
$$\roman{Sing}(J_1,b) \neq \roman{Sing}(D(J)_1,1),$$
failing to have the desired equality between the singular loci.

\vskip.1in

It is worthwhile to note that if the equality in Remark
3-6, which was remarked there not to hold, were true, then the above simpler
candidate would work. 

\vskip.1in

It is our intention to emphasize the subtlety involving the definition of the coefficient
ideal and Claim 3-5.

\newpage

$$\bold{CHAPTER\ 4.\ GENERAL\ BASIC\ OBJECTS\ AND\ INVARIANTS}$$

\vskip.1in

In this chapter, we introduce the notion of a general basic object, which turns out to be the
right framework, in the solution by Encinas and Villamayor, to extract the inductive nature of
the problem of resolution of singularities.  We also show that the invariants defined on the
singular loci of the individual basic objects, as in Definition 1-10, in the charts of a general
basic object, patch together to define well-defined invariants on the singular locus of
the general basic object.

\proclaim{Definition 4-1 (General basic object)} A general basic object over $(F_0, (W_0,
E_0))$, where $(W_0,E_0)$ is a pair (cf. Definition 1-6) and $F_0 \subset W_0$ is a closed
subset, with a $d$-dimensional structure $(d \leq \dim W_0)$, is a collection ${\goth C}$ of
sequences of transformations and smooth morphisms of pairs with specified closed
subsets starting with $(F_0,(W_0,E_0))$
$$(F_0, (W_0,E_0)) \leftarrow \cdot\cdot\cdot \leftarrow (F_k, (W_k,E_k))$$
and an open covering $\{W_0^{\lambda}\}_{\lambda \in \Lambda}$ with the following data ${\Cal
D}_{\lambda}$ for each $\lambda \in \Lambda$:

(i) $j_0^{\lambda}:(\widetilde{W_0^{\lambda}},\widetilde{E_0^{\lambda}}) \hookrightarrow
(W^{\lambda}_0,E^{\lambda}_0)$ is an immersion of pairs where $\dim \widetilde{W_0^{\lambda}} =
d$, that is to say, $\widetilde{W_0^{\lambda}} \hookrightarrow W_0^{\lambda}$ is a closed
immersion of a
$d$-dimensional smooth variety $\widetilde{W_0^{\lambda}}$ into $W_0^{\lambda}$,
$\widetilde{W_0^{\lambda}}$ is permissible with respect to $E_0^{\lambda} = E_0 \cap
W_0^{\lambda}$ and
$\widetilde{W_0^{\lambda}} \not\subset E_0^{\lambda}$, and $\widetilde{E_0^{\lambda}} =
E_0^{\lambda} \cap \widetilde{W_0^{\lambda}}$,

(ii) a basic object $(\widetilde{W_0^{\lambda}}, ({\goth a}_0^{\lambda},b^{\lambda}),
\widetilde{E_0^{\lambda}})$,

\vskip.1in

satisfying the following conditions (GB-0,1,2,3):

\vskip.1in

(GB-0) The trivial sequence consisting only of $(F_0, (W_0,E_0))$ is in the collection ${\goth
C}$, i.e.,
$$(F_0, (W_0,E_0)) \in {\goth C}$$
and
$$F_0 = \cup \roman{Sing}({\goth a}_0^{\lambda},b^{\lambda}) \text{\ with\ }F_0 \cap
W_0^{\lambda} = \roman{Sing}({\goth a}_0^{\lambda},b^{\lambda}).$$

\vskip.1in
 
(GB-1) With any sequence of transformations and smooth morphisms in the collection
${\goth C}$
$$(F_0, (W_0,E_0)) \leftarrow \cdot\cdot\cdot \leftarrow (F_k, (W_k,E_k))$$
there corresponds for each $\lambda$ a sequence of transformations (with the same centers) and
(the same) smooth morphisms (obtained by taking the Cartesian products)
$$(\widetilde{W_0^{\lambda}}, ({\goth a}_0^{\lambda},b^{\lambda}),
\widetilde{E_0^{\lambda}}) \leftarrow \cdot\cdot\cdot \leftarrow (\widetilde{W_k^{\lambda}},
({\goth a}_k^{\lambda},b^{\lambda}),
\widetilde{E_k^{\lambda}})$$
with the natural immersions
$$\CD
(W_0^{\lambda},E_0^{\lambda}) @<<< \cdot\cdot\cdot @<<< (W_k^{\lambda},E_k^{\lambda}) \\
@AAA @. @AAA \\
(\widetilde{W_0^{\lambda}},\widetilde{E_0^{\lambda}}) @<<< \cdot\cdot\cdot @<<<
(\widetilde{W_k^{\lambda}},\widetilde{E_k^{\lambda}})
\endCD$$
and we have
$$F_i = \cup \roman{Sing}({\goth a}_i^{\lambda},b^{\lambda}) \text{\ with\ }F_i \cap
W_i^{\lambda} = \roman{Sing}({\goth a}_i^{\lambda},b^{\lambda})$$
for $i = 0, 1, ... , k$.

\vskip.1in
 
(We note here that in the above clause ``there corresponds ...", it is required that whenever
$(W_{i-1},E_{i-1}) \overset{\pi_i}\to{\leftarrow} (W_i,E_i)$ is a transformation with center
$Y_{i-1} \subset W_{i-1}$, the center $Y_{i-1}$ is permissible for
$(\widetilde{W_{i-1}^{\lambda}}, ({\goth a}_{i-1}^{\lambda},b^{\lambda}),
\widetilde{E_{i-1}^{\lambda}})$, i.e., $Y_{i-1} \cap W_{i-1}^{\lambda} \subset
\widetilde{W_{i-1}^{\lambda}}$,
$Y_{i-1} \cap W_{i-1}^{\lambda}$ is permissible with respect to
$\widetilde{E_{i-1}^{\lambda}}$, and
$Y_{i-1}
\cap W_{i-1}^{\lambda} \subset \roman{Sing}({\goth a}_{i-1}^{\lambda},b^{\lambda})$.)

\vskip.1in

(GB-2) Let
$$(F_0, (W_0,E_0)) \leftarrow \cdot\cdot\cdot \leftarrow (F_k, (W_k,E_k))$$
be a sequence of transformations and smooth morphisms in ${\goth C}$ and
$$\{(\widetilde{W_0^{\lambda}}, ({\goth a}_0^{\lambda},b^{\lambda}),
\widetilde{E_0^{\lambda}}) \leftarrow \cdot\cdot\cdot \leftarrow (\widetilde{W_k^{\lambda}},
({\goth a}_k^{\lambda},b^{\lambda}),
\widetilde{E_k^{\lambda}})\}$$
the corresponding sequences (indexed by $\lambda \in \Lambda$) of transformations and smooth
morphisms as in (GB-1).

\vskip.1in

We take a morphism of pairs $(W_k,E_k) \overset{\pi_{k+1}}\to{\leftarrow} (W_{k+1},E_{k+1})$ which is
either in Case T or Case S.

\vskip.1in

Case T: $(W_k,E_k) \overset{\pi_{k+1}}\to{\leftarrow} (W_{k+1},E_{k+1})$ is a transformation
with center $Y_k \subset W_k$, satisfying the condition that $Y_k$ is permissible for each
$(\widetilde{W_k^{\lambda}}, ({\goth a}_k^{\lambda},b^{\lambda}),
\widetilde{E_k^{\lambda}})$, i.e., \linebreak
$Y_k \cap W_k^{\lambda} \subset
\widetilde{W_k^{\lambda}}$,
$Y_k \cap W_k^{\lambda}$ is permissible with respect to $\widetilde{E_k^{\lambda}}$, and
\linebreak
$Y_k \cap W_k^{\lambda} \subset \roman{Sing}({\goth a}_k^{\lambda},b^{\lambda})$.

\vskip.1in

Case S: $(W_k,E_k) \overset{\pi_{k+1}}\to{\leftarrow} (W_{k+1},E_{k+1})$ is a smooth morphism.

\vskip.2in

Then we have the following assertions on the extension of the original sequence:

\vskip.1in

Case T: Take for each $\lambda$ the corresponding transformation of basic objects with center $Y_k \cap
W_k^{\lambda}$
$$(\widetilde{W_k^{\lambda}}, ({\goth a}_k^{\lambda},b^{\lambda}),
\widetilde{E_k^{\lambda}}) \overset{\pi_{k+1}^{\lambda}}\to{\leftarrow}
(\widetilde{W_{k+1}^{\lambda}}, ({\goth a}_{k+1}^{\lambda},b^{\lambda}),
\widetilde{E_{k+1}^{\lambda}}).$$

Then
$$F_{k+1} := \cup \roman{Sing}({\goth a}_{k+1}^{\lambda},b^{\lambda})$$
is a closed subset of $W_{k+1}$ with
$$F_{k+1} \cap W_{k+1}^{\lambda} = \roman{Sing}({\goth a}_{k+1}^{\lambda},b^{\lambda}),$$
and the extended sequence belongs to ${\goth C}$, i.e.,
$$(F_0, (W_0,E_0)) \leftarrow \cdot\cdot\cdot \leftarrow (F_k, (W_k,E_k)) \leftarrow (F_{k+1},
(W_{k+1},E_{k+1})) \in {\goth C}.$$

\vskip.1in

Case S: Take for each $\lambda$ the corresponding morphism of basic objects
$$(\widetilde{W_k^{\lambda}}, ({\goth a}_k^{\lambda},b^{\lambda}),
\widetilde{E_k^{\lambda}}) \overset{\pi_{k+1}^{\lambda}}\to{\leftarrow}
(\widetilde{W_{k+1}^{\lambda}}, ({\goth a}_{k+1}^{\lambda},b^{\lambda}),
\widetilde{E_{k+1}^{\lambda}})$$
where
$$\widetilde{W_{k+1}^{\lambda}} = \widetilde{W_k^{\lambda}} \times_{W_k}W_{k+1}, {\goth
a}_{k+1}^{\lambda} = {\goth a}_k^{\lambda}{\Cal O}_{\widetilde{W_{k+1}}},
\widetilde{E_{k+1}^{\lambda}} = {\pi_{k+1}^{\lambda}}^{-1}(\widetilde{E_k^{\lambda}})$$ 
and where $\pi_{k+1}^{\lambda}$ is the
projection onto the first factor.

Then
$$F_{k+1} := \cup \roman{Sing}({\goth a}_{k+1}^{\lambda},b^{\lambda})$$
is a closed subset of $W_{k+1}$ with
$$F_{k+1} \cap W_{k+1}^{\lambda} = \roman{Sing}({\goth a}_{k+1}^{\lambda},b^{\lambda}),$$
and the extended sequence belongs to ${\goth C}$, i.e.,
$$(F_0, (W_0,E_0)) \leftarrow \cdot\cdot\cdot \leftarrow (F_k, (W_k,E_k)) \leftarrow (F_{k+1},
(W_{k+1},E_{k+1})) \in {\goth C}.$$

\vskip.1in

(GB-3) There exists $c \in {\Bbb N}$ such that $c \geq b^{\lambda} \hskip.1in \forall \lambda$.

\vskip.1in

We denote by $({\Cal F}_0, (W_0,E_0))$ a general basic object over $(F_0, (W_0,E_0))$.

\vskip.1in

We say that the $d$-dimensional structure of $({\Cal F}_0, (W_0,E_0))$ is given by the charts
$\{(\widetilde{W_0^{\lambda}}, ({\goth a}_0^{\lambda},b^{\lambda}),\widetilde{E_0^{\lambda}})\}$
of basic objects of dimension $d$.

\vskip.1in 

We also say by abuse of language that the collection ${\goth C}$ is represented by the general
basic object $({\Cal F}_0, (W_0,E_0))$.

\vskip.1in 

We identify two general basic objects $({\Cal F}_0, (W_0,E_0))$ and
$({\Cal F}_0', (W_0,E_0))$ if and only if the collections
${\goth C}$ and
${\goth C}'$, represented by the general basic objects, coincide.
\endproclaim

\vskip.1in

\proclaim{Remark 4-2}\endproclaim

(i) In the previous chapters, the letter $d$ was used to denote the dimension of the ambient
space $W$ of a basic object $(W, (J,b), E)$.  When we say a general basic object $({\Cal F}_0,
(W_0, E_0))$ with a $d$-dimensional structure, the letter $d$ refers to the dimension of the
basic objects $\{(\widetilde{W_0}, ({\goth a}_{\lambda},b^{\lambda}),
\widetilde{E_0^{\lambda}}\}$ in the charts, and not to the dimension of $W_0$.  In general, $d
\leq \dim W_0$.

(ii) ($\bold{The\ general\ basic\ object\ defined\ by\ a\ basic\ object)}$ Let $(W, (J,b), E)$
be a basic object with
$d =
\dim W$.  Set
$(F_0, (W_0,E_0)) = (\roman{Sing}(J,b), (W,E))$ and take
${\goth C}$ to be the collection of all the sequences of transformations and smooth morphisms
$$(F_0, (W_0,E_0)) \leftarrow \cdot\cdot\cdot \leftarrow (F_r, (W_r,E_r))$$
induced by the sequences of transformations and smooth morphisms of basic objects
$$(W, (J,b), E) = (W_0, (J_0,b), E_0) \leftarrow \cdot\cdot\cdot \leftarrow (W_k, (J_k,b), E_k)$$
with
$$F_i = \roman{Sing}(J_i,b) \text{\ for\ }i = 0, 1, ... , k.$$
This defines a general basic object over $(F_0, (W_0,E_0))$ with a $d$-dimensional structure.

\vskip.1in

Note that two different basic objects, e.g. $(W,(J,b),E)$ and $(W,(J^2,2b),E)$, may define the
same general basic object, as they represent the same collection ${\goth C}$.

\vskip.1in

(iii) ($\bold{The\ meaning\ of\ the\ key\ inductive\ lemma}$) Let $(W, (J,b), E)$ be a
simple basic object with an open covering
$\{W^{\lambda}\}_{\lambda \in \Lambda}$ satisfying conditions 1 and 2 of the key inductive
lemma (Lemma 3-1).  Then, in case $R(1)(\roman{Sing}(J,b)) = \emptyset$, the general basic object
$({\Cal F}_0, (W_0,E_0))$ with a
$d$-dimensional structure defined as above by the basic object $(W, (J,b), E)$, also has a
$(d-1)$-dimensional structure given by $\{(\widetilde{W_0^{\lambda}}, ({\goth
a}_0^{\lambda},b^{\lambda}),
\widetilde{E_0^{\lambda}}) = (W_h^{\lambda}, (C(J^{\lambda}), b!), E_h^{\lambda})\}$.

This drop in the dimension of the structure is the very essence of the key inductive lemma in
terms of the notion of general basic objects.

\vskip.1in

(iv) ($\bold{Why\ smooth\ morphisms ?}$) In the process of resolution of singularities, we only consider a
sequence of transformations.  So why do we have to consider smooth morphisms in the definition of a general
basic object ?  One main reason is that we would like to use Hironaka's trick (cf. the proof of
Definition-Proposition 4-5) to guarantee that the invariants defined on the individual charts
patch together to give well-defined invariants on the general basic object.  Another reason is
that, by including open immersions, we also want to guarantee that the general basic objects
behave well under localization.

\vskip.1in

(v) In the original papers by Encinas and Villamayor, they do not include general smooth
morphisms in the sequences to consider for the collection to characterize a general basic
object, but include only special smooth morphisms, namely, the projections of the form $W_{i-1}
\leftarrow W_i = W_{i-1}
\times {\Bbb A}^1$, which they call the restrictions.  However, their definition causes a few
problems:  

\ \ $\circ$ It is not
clear by their definition whether their general basic objects behave well under localization. 
That is to say, it is not clear a priori whether the charts $\{(\widetilde{W_0^{\lambda}} \cap V,
({\goth a}_0^{\lambda}|_V,b^{\lambda}),
\widetilde{E_0^{\lambda}} \cap V)\}$ would define a general basic object for an open subset $V
\subset W_0$, since the permissiblity of the centers is a global condition.  (Though this can be
proved using, e.g., the embedded resolution of singularities of the closure of the center taken
in the open subset without affecting the open subset itself.)

\ \ $\circ$ When we want to discuss the stability of the process of resolution of singularities
under smooth morphisms, we would like to have the definition of a smooth morphism between
general basic objects, which we would lack under their definition.

\vskip.1in

Including (general) smooth morphisms into the definition of a general basic object requires no
change in the structure of the argument and brings some theoretical clarity avoiding the
problems as above.

\vskip.1in

(vi) ($\bold{An\ alternative\ way\ of\ defining\ the\ notion\ of\ a\ general\ basic\ object}$)
Recall that a differentiable (resp. complex) manifold $W$ is defined to consist of a
topological space $W$ and an open covering $\{W^{\lambda}\}_{\lambda \in \Lambda}$ with charts
$h^{\lambda}:W^{\lambda}
\rightarrow U^{\lambda} \subset {\Bbb R}^n$ (resp. $\subset {\Bbb C}^n$) so that the
$W^{\lambda}$ patch up in the sense that $h^{\mu} \circ {h^{\lambda}}^{-1}$ are
invertible $C^{\infty}$-functions (resp. holomorphic functions).

\vskip.1in

We can give an alternative definition of a general basic object in a similar manner: 

\vskip.1in

A general
basic object consists of a pair $(W,E)$ and an open covering $\{W^{\lambda}\}_{\lambda \in
\Lambda}$ with charts $(\widetilde{W^{\lambda}},({\goth
a}^{\lambda},b^{\lambda}),\widetilde{E^{\lambda}})$, i.e., basic objects where
$\widetilde{W^{\lambda}} \hookrightarrow W^{\lambda}$ is a closed immersion of a $d$-dimensional
smooth variety $\widetilde{W^{\lambda}}$ into $W^{\lambda}$, $\widetilde{W^{\lambda}}$ is
permissible with respect to $E^{\lambda} = E \cap W^{\lambda}$, $\widetilde{W^{\lambda}}
\not\subset E^{\lambda}$, and where
$\widetilde{E^{\lambda}} = E^{\lambda} \cap \widetilde{W^{\lambda}}$.  We require that there
exists $c \in {\Bbb N}$ such that $c \geq b^{\lambda} \hskip.1in \forall \lambda$.  
 
We also require the following
patching condition among the charts.

Let
$$(W^{\lambda} \cap W^{\mu}, E \cap W^{\lambda} \cap W^{\mu}) =
(W_0^{\lambda\mu},E_0^{\lambda\mu}) \leftarrow \cdot\cdot\cdot \leftarrow (W_k^{\lambda\mu},
E_k^{\lambda\mu})$$ be a sequence of transformations and smooth morphisms of pairs, starting with
the intersection of
$(W^{\lambda},E^{\lambda})$ and $(W^{\mu},E^{\mu})$, such that there corresponds
a sequence of transformations (with the same centers) and (the same) smooth morphisms (obtained by
taking the Cartesian products) of basic objects 
$$(\widetilde{W^{\lambda}},({\goth
a}^{\lambda},b^{\lambda}),\widetilde{E^{\lambda}}) \cap W^{\mu} =
(\widetilde{W_0^{\lambda\mu}},({\goth
a}_0^{\lambda},b^{\lambda}),\widetilde{E_0^{\lambda\mu}}) 
\leftarrow \cdot\cdot\cdot \leftarrow (\widetilde{W_k^{\lambda\mu}},({\goth
a}_k^{\lambda},b^{\lambda}),\widetilde{E_k^{\lambda\mu}}).$$ 
Then there should
correspond a sequence of transformations (with the same centers) and (the same) smooth morphisms
(obtained by taking the Cartesian products) of basic objects 
$$(\widetilde{W^{\mu}},({\goth
a}^{\mu},b^{\mu}),\widetilde{E^{\mu}}) \cap W^{\lambda} =
(\widetilde{W_0^{\mu\lambda}},({\goth a}_0^{\mu},b^{\mu}),\widetilde{E_0^{\mu\lambda}})
\leftarrow \cdot\cdot\cdot \leftarrow (\widetilde{W_k^{\mu\lambda}},({\goth
a}_k^{\mu},b^{\mu}),\widetilde{E_k^{\mu\lambda}})$$ satisfying
$$\roman{Sing}({\goth a}_i^{\lambda},b^{\lambda}) = \roman{Sing}({\goth
a}_i^{\mu},b^{\mu}).$$

\vskip.1in
 
The equivalence of this alternative definition and Definition 4-1 is straightforward and its
verification is left to the reader as an exercise.

\vskip.1in

\proclaim{Note 4-3 (Sequence of transformations and smooth morphisms of general basic
objects)}\endproclaim

Let $({\Cal F}_0, (W_0,E_0))$ be a general basic object with a $d$-dimensional structure given
by the charts $\{(\widetilde{W_0^{\lambda}}, ({\goth a}_0^{\lambda},b^{\lambda}),
\widetilde{E_0^{\lambda}})\}$.

Let
$$(F_0, (W_0,E_0)) \leftarrow \cdot\cdot\cdot \leftarrow (F_k, (W_k,E_k))$$
be a sequence of transformations and smooth morphisms in ${\goth C}$.

Then the seqeunce induces a general basic object over $(F_i, (W_i,E_i))$ for \linebreak
$i = 0, 1,
... , k$, with a $\{\dim W_i - (\dim W_0 - d)\}$-dimensional structure, in the
following way: 
 
We take ${\goth C}_i$ to be the collection of those sequences which are the truncations of the
sequences in
${\goth C}$ whose first $(i+1)$-terms coincide with
$$(F_0, (W_0,E_0)) \leftarrow \cdot\cdot\cdot \leftarrow (F_i, (W_i,E_i))$$
of the given sequence. 

\vskip.1in

We take the data ${\Cal D}_i^{\lambda}$ for each $\lambda$ to consist of:

(i) the induced immersion $j_i^{\lambda}:(\widetilde{W_i^{\lambda}}, \widetilde{E_i^{\lambda}})
\hookrightarrow (W_i^{\lambda},E_i^{\lambda})$, where \hfill\hfill\linebreak
$\dim \widetilde{W_i^{\lambda}} =
\dim W_i - (\dim W_0 - d)$,

(ii) the induced basic object $(\widetilde{W_i^{\lambda}}, ({\goth a}_i^{\lambda},b^{\lambda}),
\widetilde{E_i^{\lambda}})$.

\vskip.1in

We denote the general basic object as above $({\Cal F}_i, (W_i,E_i))$.

\vskip.1in

Therefore, by abuse of notation, we write the sequence
$$({\Cal F}_0, (W_0,E_0)) \leftarrow \cdot\cdot\cdot \leftarrow ({\Cal F}_k, (W_k,E_k))$$
and call it a sequence of transformations and smooth morphisms of general basic objects.

\vskip.1in

\proclaim{Definition 4-4 (Resolution of singularities of a general basic object)} We call a
sequence of transformations only of general basic objects
$$({\Cal F}_0, (W_0,E_0)) \leftarrow \cdot\cdot\cdot \leftarrow ({\Cal F}_k, (W_k,E_k))$$
resolution of singularities of a general basic object $({\Cal F}_0, (W_0,E_0))$ if
$$F_k = \emptyset.$$
\endproclaim 

\proclaim{Definition-Proposition 4-5 (Key invariants of general basic objects)} Let
$$({\Cal F}_0, (W_0,E_0)) \leftarrow \cdot\cdot\cdot \leftarrow ({\Cal F}_k, (W_k,E_k))$$
be a sequence of transformations and smooth morphisms of general basic objecs (cf. Note 4-3).

(i) The invariant $\roman{ord}_k:F_k \rightarrow \frac{1}{c!}{\Bbb Z}{\geq 0}$ is a function
defined over $F_k$ such that for each $\lambda$ it restricts to the invariant
$\roman{ord}_k^{\lambda}$ of the basic object $(\widetilde{W_k^{\lambda}}, ({\goth
a}_k^{\lambda},b^{\lambda}),
\widetilde{E_k^{\lambda}})$, i.e., we have the following commutative diagram
$$\CD
F_k \cap W_k^{\lambda} @. \hookrightarrow @. F_k @. \overset{\roman{ord}_k}\to{\rightarrow}
@. \frac{1}{c!}{\Bbb Z}_{\geq 0} \\
@| @. @. @. \cup \\
\roman{Sing}({\goth a}_k^{\lambda},b^{\lambda}) @.\hskip.2in @.
\overset{\roman{ord}_k^{\lambda}}\to{\rightarrow} @.@. \frac{1}{b^{\lambda}}{\Bbb Z}_{\geq 0}. \\
\endCD$$

(ii) The invariant $w\text{-}\roman{ord}_k:F_k \rightarrow \frac{1}{c!}{\Bbb Z}{\geq 0}$ is a
function defined over $F_k$ such that for each $\lambda$ it restricts to the invariant
$w\text{-}\roman{ord}_k^{\lambda}$ of the basic object $(\widetilde{W_k^{\lambda}}, ({\goth
a}_k^{\lambda},b^{\lambda}),
\widetilde{E_k^{\lambda}})$, i.e., we have the following commutative diagram
$$\CD
F_k \cap W_k^{\lambda} @. \hookrightarrow @. F_k @.
\overset{w\text{-}\roman{ord}_k}\to{\rightarrow} @. \frac{1}{c!}{\Bbb Z}_{\geq 0} \\
@| @. @. @. \cup \\
\roman{Sing}({\goth a}_k^{\lambda},b^{\lambda}) @.\hskip.2in @.
\overset{w\text{-}\roman{ord}_k^{\lambda}}\to{\rightarrow} @.@. \frac{1}{b^{\lambda}}{\Bbb
Z}_{\geq 0}.
\\
\endCD$$

(iii) First note that in order to define the invariant $t_k$ we require the following extra
condition $(\heartsuit)$ on the sequence of transformations and smooth morphisms of general basic
objects
$$(\heartsuit) \hskip.1in \left\{\aligned &Y_{i-1} \subset
\underline{\roman{Max}}\ w\text{-}\roman{ord}_{i-1} (\subset F_{i-1}) \\
&\text{whenever\ } \pi_i \text{\ is\ a\ transformation\ with\ center\
}Y_{i-1}\\
\endaligned\right\}$$  
where 
$$\align
\underline{\roman{Max}}\ w\text{-}\roman{ord}_{i-1} &= \{p \in F_{i-1};w\text{-}\roman{ord}_{i-1}(p) =
\max w\text{-}\roman{ord}_{i-1}\}\\
\max w\text{-}\roman{ord}_{i-1} &= \max\{w\text{-}\roman{ord}_{i-1}(p);p \in
F_{i-1}\}.\\
\endalign$$
Under condition $(\heartsuit)$ it follows that we have inequalities (See Proposition 1-12.)
$$\align
\max w\text{-}\roman{ord}_0 &\geq \max w\text{-}\roman{ord}_1 \geq \cdot\cdot\cdot \\
\max w\text{-}\roman{ord}_{i-1}
&\geq \max w\text{-}\roman{ord}_i \\
\cdot\cdot\cdot \geq \max w\text{-}\roman{ord}_{k-1} &\geq \max
w\text{-}\roman{ord}_k.\\
\endalign$$

Let $k_o$ be the index so that
$$\max w\text{-}\roman{ord}_{k_o-1} > \max w\text{-}\roman{ord}_{k_o} = \cdot\cdot\cdot = \max
w\text{-}\roman{ord}_k.$$ 
(We let $k_o = 0$ if $\max w\text{-}\roman{ord}_0 = \cdot\cdot\cdot = \max
w\text{-}\roman{ord}_k$.)  Set $E_k = E_k^{-} \cup E_k^{+}$ where $E_k^{-} = \{H_1, ... , H_r,
... , H_{r + k_o}\}$ as a subset of $E_k = \{H_1, ... , H_r, ... , H_{r + k_o},
... , H_{r + k}\}$ and where
$E_k^+$ is the complement of
$E_k^-$ in
$E_k$.  (Look also at the convention explained in Definition 1-8 (iii).)  

\vskip.1in

The invariant $t_k:F_k \rightarrow \frac{1}{c!}{\Bbb Z}_{\geq 0} \times {\Bbb Z}_{\geq
0}$ is a function defined over $F_k$ such that
$$t_k(p) = (w\text{-}\roman{ord}_k(p),n_k(p)) \text{\ for\ }p \in F_k$$
where 
$$n_k(p) = \left\{\aligned
\#\{H_i \in E_k;p \in H_i\} &\text{\ if\ } w\text{-}\roman{ord}_k(p) < \max
w\text{-}\roman{ord}_k
\\
\#\{H_i \in E_k^{-}; p \in H_i\} &\text{\ if\ } w\text{-}\roman{ord}_k(p) = \max
w\text{-}\roman{ord}_k.
\\
\endaligned
\right.$$

Moreover, for each $\lambda$ it restricts to the invariant
$t_k^{\lambda}$ of the basic object $(\widetilde{W_k^{\lambda}}, ({\goth
a}_k^{\lambda},b^{\lambda}),
\widetilde{E_k^{\lambda}})$, i.e., we have the following commutative diagram
$$\CD
F_k \cap W_k^{\lambda} @. \hookrightarrow \hskip.1in @. F_0 @. \overset{t_k}\to{\rightarrow}
@. \frac{1}{c!}{\Bbb Z}_{\geq 0} \times {\Bbb Z}_{\geq 0}\\
@| @. @. @. \cup \\
\roman{Sing}({\goth a}_k^{\lambda},b^{\lambda}) @.\hskip.2in @.
\overset{t_k^{\lambda}}\to{\rightarrow} @.@. \frac{1}{b^{\lambda}}{\Bbb Z}_{\geq 0} \times {\Bbb
Z}_{\geq 0}.
\\
\endCD$$

\vskip.1in

(iv) Suppose
$$\max\ w\text{-}\roman{ord}_k = 0.$$

Then the invariant $\Gamma_k:F_k \rightarrow {\Bbb Z}_{\geq - \dim W_k} \times \frac{1}{c!}{\Bbb
Z}_{\geq 0}
\times {\Bbb Z}_{\geq 0}^{\dim W_k}$ is a function defined over $F_k$ such that for each $\lambda$
it restricts to the invariant
$\Gamma_k^{\lambda}$ of the monomial basic object $(\widetilde{W_k^{\lambda}}, ({\goth a}_k^{\lambda},b^{\lambda}),
\widetilde{E_k^{\lambda}})$ (shrinking $\widetilde{W_k^{\lambda}}$ to an open neighborhood of
\linebreak
$\roman{Sing}({\goth a}_k^{\lambda},b^{\lambda}) = F_k \cap W_k^{\lambda}$ if necessary (cf.
Corollary 2-7)), i.e., we have the following commutative diagram
$$\CD
F_k \cap W_k^{\lambda} @. \hookrightarrow \hskip.1in @. F_k @. \overset{\Gamma_k}\to{\rightarrow}
@. {\Bbb Z}_{\geq - \dim W_k} \times \frac{1}{c!}{\Bbb Z}_{\geq 0} \times
{\Bbb Z}_{\geq 0}^{\dim W_k} \\
@| @. @. @. \cup \\
\roman{Sing}({\goth a}_k^{\lambda},b^{\lambda}) @.\hskip.2in @.
\overset{\Gamma_k^{\lambda}}\to{\rightarrow} @.@. {\Bbb Z}_{\geq - \dim \widetilde{W_k^{\lambda}}}
\times
\frac{1}{b^{\lambda}}{\Bbb Z}_{\geq 0}
\times {\Bbb Z}_{\geq 0}^{\dim \widetilde{W_k^{\lambda}}}. \\
\endCD$$ 

Finally, the invariants $\roman{ord}_k, w\text{-}\roman{ord}_k, t_k$ and $\Gamma_k$ are
determined purely in terms of the collection ${\goth C}_k$ of sequences of transformations and
smooth morphisms represented by the general basic object $({\Cal F}_k,(W_k,E_k))$, in terms
of the original sequence
$$(F_0,(W_0,E_0)) \leftarrow \cdot\cdot\cdot \leftarrow (F_k,(W_k,E_k)),$$
and in terms of the specified dimension $d_k = \dim W_k - \dim W_0 + d_0$ of the structure of
the general basic object $({\Cal F}_k,(W_k,E_k))$, but free of the presentation using charts. 
\endproclaim

\demo{Proof}\enddemo (i) Since the invariant $\roman{ord}_k$ depends only on the general basic
object
$({\Cal F}_k, (W_k,E_k))$ and not on the sequence, in order to avoid the complication which may
be caused by the subscripts, we prove that the invariant $\roman{ord}_0$ exists for the general
basic object $({\Cal F}_0, (W_0,E_0))$ with the required property.

It suffices to show that the functions $\{\roman{ord}_0^{\lambda}\}$, defined on the individual charts as
in Definition 1-10 (i), patch up.  That is to say, it suffices to show that for any closed point $x_0 \in
F_0$ and commutative diagrams of the form
$$\CD
x_0^{\lambda} @. \in @. \roman{Sing}({\goth a}_0^{\lambda},b^{\lambda}) @. \subset @.
\widetilde{W_0^{\lambda}} \\
@.@.@.@. \cap \\
@VV j_0^{\lambda} V @.@.@. W_0^{\lambda} \\
@.@.@.@. \cup \\
x_0 @.@. \in @.@. W_0^{\lambda} \cap W_0^{\lambda'}\\
@.@.@.@. \cap \\
@AA j_0^{\lambda'} A @.@.@. W_0^{\lambda'} \\
@.@.@.@. \cup \\
x_0^{\lambda'} @. \in @. \roman{Sing}({\goth a}_0^{\lambda'},b^{\lambda'}) @. \subset @.
\widetilde{W_0^{\lambda'}} \\
\endCD$$
we have
$$\roman{ord}_0^{\lambda}(x_0^{\lambda}) = \frac{\nu_{x_0^{\lambda}}({\goth
a}_0^{\lambda})}{b^{\lambda}} =
\frac{\nu_{x_0^{\lambda'}}({\goth a}_0^{\lambda'})}{b^{\lambda'}} =
\roman{ord}_0^{\lambda'}(x_0^{\lambda'}).$$
We will show that the number
$$\frac{\nu_{x_0^{\lambda}}({\goth a}_0^{\lambda})}{b^{\lambda}}$$
can be determined purely in terms of the collection ${\goth C}$ of transformations and
smooth morphisms representing the general basic object $({\Cal F}_0, (W_0,E_0))$ and hence
is independent of $\lambda$.

\vskip.1in

The method of the proof below is what Encinas and Villamayor call \linebreak
``$\bold{Hironaka's\
trick}$" and this is the only place in the paper where we make use of smooth morphisms (other than
open immersions) in the sequences in the collection represented by a general basic object.

\vskip.1in

We construct a sequence of transformations and smooth morphisms in the following manner.  (A warning to the
reader: Though we use the same notation, the sequence we construct below has nothing to do with the
original sequence described in the statement of Definition-Proposition 4-5.) 

\vskip.1in

Step 1. First we consider the following smooth morphism, which is nothing but the projection onto
the first factor from the product with ${\Bbb A}^1$
$$W_0 \overset{\pi_1}\to{\leftarrow} W_1 = W_0 \times {\Bbb A}^1,$$
where we denote $L_1 = \pi_1^{-1}(x_0)$ and choose a point \linebreak
$x_1 = (x_0,0) \in L_1 \subset
W_0
\times {\Bbb A}^1 = W_1$.  For each $\lambda$, we have the corresponding smooth morphisms
$$(\widetilde{W_0^{\lambda}},({\goth a}_0^{\lambda},b^{\lambda}), \widetilde{E_0^{\lambda}})
\overset{\pi_1^{\lambda}}\to{\leftarrow} (\widetilde{W_1^{\lambda}},({\goth
a}_1^{\lambda},b^{\lambda}),
\widetilde{E_0^{\lambda}}) =
(\widetilde{W_0^{\lambda}} \times {\Bbb A}^1, ({\goth
a}_0^{\lambda}{\Cal O}_{\widetilde{W_1^{\lambda}}},b^{\lambda}), \widetilde{E_0^{\lambda}} \times
{\Bbb A}^1)$$
where, for the indices $\lambda$ with $x_0 = x_0^{\lambda} \in F_0 \cap W_0^{\lambda} \subset
\widetilde{W_0^{\lambda}}
\subset W_0^{\lambda}$, we denote \linebreak
$L_1^{\lambda} =
{\pi_1^{\lambda}}^{-1}(x_0^{\lambda})$ and the corresponding point of choice by $x_1^{\lambda} =
(x_0^{\lambda},0) \in L_1^{\lambda} \subset \widetilde{W_0^{\lambda}}
\times {\Bbb A}^1 =
\widetilde{W_1^{\lambda}}$.

\vskip.1in

Step 2. Secondly we consider the sequence of transformations of pairs
$$(W_1,E_1) \overset{\pi_2}\to{\leftarrow} \cdot\cdot\cdot \overset{\pi_N}\to{\leftarrow}
(W_N,E_N)$$
where inductively

$\pi_2$ is the blowup at $x_1 \subset L_1 \subset F_1$, 

and for $i \geq 3$

$\pi_i$ is the blowup at $x_{i-1} = L_{i-1} \cap H_{r + i-1} \subset L_{i-1} \subset F_{i-1}$ with
$H_{r + i-1}$ being the exceptional divisor for $\pi_{i-1}$ and $L_{i-1}$ being the strict
transform of $L_1$.

For any $\lambda$ with $W_1^{\lambda} \ni x_1$ (i.e., $W_0^{\lambda} \ni x_0$), we have the
corresponding sequence of transformations of basic objects in the charts
$$(\widetilde{W_1^{\lambda}},({\goth a}_1^{\lambda},b^{\lambda}),\widetilde{E_1^{\lambda}})
\overset{\pi_2^{\lambda}}\to{\leftarrow} \cdot\cdot\cdot
\overset{\pi_N^{\lambda}}\to{\leftarrow} (\widetilde{W_N^{\lambda}},({\goth
a}_N^{\lambda},b^{\lambda}),\widetilde{E_N^{\lambda}})$$
where

$\pi_2^{\lambda}$ is the blowup at $x_1^{\lambda} \subset L_1^{\lambda} \subset
\roman{Sing}({\goth a}_1^{\lambda},b^{\lambda})$, 

and for $i \geq 3$

$\pi_i^{\lambda}$ is the blowup at $x_{i-1}^{\lambda} = L_{i-1}^{\lambda} \cap
\widetilde{H_{r+i-1}^{\lambda}}
\subset L_{i-1}^{\lambda}
\subset
\roman{Sing}({\goth a}_{i-1}^{\lambda},b^{\lambda})$ with
$\widetilde{H_{r+i-1}^{\lambda}}$ being the exceptional divisor for $\pi_{i-1}^{\lambda}$ and
$L_{i-1}^{\lambda}$ being the strict transform of $L_1^{\lambda}$.

Note that under the inclusion $\widetilde{W_{i-1}^{\lambda}} \subset W_{i-1}^{\lambda}$ we
identify
$$x_{i-1}^{\lambda} = x_{i-1}, L_{i-1}^{\lambda} = L_{i-1}, \widetilde{H_{r+i-1}^{\lambda}} =
H_{r+i-1}
\cap
\widetilde{W_{i-1}^{\lambda}}.$$

By conditions (GB-0) and (GB-2) this gives rise to a
sequence in the collection
${\goth C}$ of the general basic object
$({\Cal F}_0, (W_0,E_0))$
$$(F_0, (W_0,E_0)) \overset{\pi_1}\to{\leftarrow} (F_1, (W_1,E_1))
\overset{\pi_2}\to{\leftarrow} \cdot\cdot\cdot \overset{\pi_N}\to{\leftarrow} (F_N, (W_N,E_N))$$
where
$$F_i = \cup \roman{Sing}({\goth
a}_i^{\lambda},b^{\lambda})$$
and for each $\lambda$ we have
$$F_i \cap W_i^{\lambda} = F_i \cap \widetilde{W_i^{\lambda}} = \roman{Sing}({\goth
a}_i^{\lambda},b^{\lambda}).$$

\vskip.1in

We compute the transformations of ideals
$${\goth a}_2^{\lambda} =
I(\widetilde{H_{r+2}^{\lambda}})^{(\beta^{\lambda}-b^{\lambda})}\overline{{\goth a}_2^{\lambda}}
\text{\ where\ }\beta^{\lambda} = \nu_{x_0^{\lambda}}({\goth a}_0^{\lambda}),$$
since
$$\nu_{\xi_1}(\overline{{\goth a}_1^{\lambda}}) = \nu_{\xi_1}({\goth a}_1^{\lambda}) =
\beta^{\lambda}
\hskip.1in
\forall \xi_1 \in L_1^{\lambda}.$$ 

Note that
$$\nu_{\xi_2}(\overline{{\goth a}_2^{\lambda}}) = \beta^{\lambda} \hskip.1in \forall \xi_2 \in
L_2^{\lambda}.$$  

In fact, it is clear that
$$\nu_{\xi_2}(\overline{{\goth a}_2^{\lambda}}) = \beta^{\lambda} \hskip.1in \forall \xi_2 \in
L_2^{\lambda}
\setminus x_2^{\lambda}$$
and hence by the upper semi-continuity
$$\nu_{\xi_2}(\overline{{\goth a}_2^{\lambda}}) \geq \beta^{\lambda} \text{\ for\ }\xi_2 =
x_2^{\lambda}.$$ 

On the other hand, applying Proposition 1-12 (ii) locally, we conclude
$$\nu_{x_2^{\lambda}}(\overline{{\goth a}_2^{\lambda}}) \leq
\nu_{x_1^{\lambda}}(\overline{{\goth a}_1^{\lambda}}) = \nu_{x_1^{\lambda}}({\goth
a}_1^{\lambda}) =
\beta^{\lambda}.$$ 

Therefore, inductively locally around $x_i^{\lambda}$ for $i = 2,
... , N$, we compute
$${\goth a}_i^{\lambda} = I(\widetilde{H_{r+i}^{\lambda}})^{(i-1)(\beta^{\lambda} -
b^{\lambda})}\overline{{\goth a}_i^{\lambda}}$$
with
$$\nu_{\xi_i}(\overline{{\goth a}_i^{\lambda}}) = \beta^{\lambda} \hskip.1in \forall \xi_i \in
L_i^{\lambda}.$$

Therefore, we conclude that
$$\align
&\dim F_N \cap H_{r+N} = d = (d + 1) - 1 \\
\Longleftrightarrow \hskip.2in & \widetilde{H_{r+N}^{\lambda}} = H_{r+N} \cap
\widetilde{W_N^{\lambda}}
\subset F_N^{\lambda} \\
\Longleftrightarrow \hskip.2in &(N - 1)(\beta^{\lambda} - b^{\lambda}) \geq b^{\lambda}.\\
\endalign$$

Remark that the condition on the first line is determined purely in terms of the collection
${\goth C}$ represented by the general basic object $({\Cal F}_0, (W_0,E_0))$ and the condition
on the third line is a numerical one about the number $b^{\lambda}$. 

\vskip.1in

Step 3. Therefore we conclude that
$$\align
& \beta^{\lambda} = b^{\lambda} \\
\Longleftrightarrow \hskip.2in & \widetilde{H_{r+N}^{\lambda}} = H_{r+N} \cap
\widetilde{W_N^{\lambda}}
\not\subset F_N^{\lambda} \text{\ for\ all\ }N \in {\Bbb N}\\
\Longleftrightarrow \hskip.2in & \dim F_N \cap H_{r+N} < d = (d + 1) - 1 \text{\ for\ all\ }N \in
{\Bbb N}\\
\endalign$$
and that
$$\align
& \beta^{\lambda} > b^{\lambda} \\
\Longleftrightarrow \hskip.2in & \widetilde{H_{r+N}^{\lambda}} = H_{r+N} \cap
\widetilde{W_N^{\lambda}} \subset F_N^{\lambda} \text{\ for\ all\ sufficiently\ large\ }N \in
{\Bbb N}\\
\Longleftrightarrow \hskip.2in & \dim F_N \cap H_{r+N} = d = (d + 1) - 1 \text{\ for\ all\
sufficiently\ large\ }N
\in {\Bbb N}.\\
\endalign$$

In the latter case, we consider the further extension of the sequence of transformations of pairs
$$(W_N, E_N) \overset{\pi_{N+1}}\to{\leftarrow} \cdot\cdot\cdot
\overset{\pi_{N+S}}\to{\leftarrow} (W_{N+S}, E_{N+S})$$ where $\pi_{N+i}$ is the blowup with
center
$Y_{N+i-1} = F_{N+i-1} \cap H_{r+N+i-1}.$

For each $\lambda$ there corresponds the sequence of transformations of basic objects
$$(\widetilde{W_N^{\lambda}}, ({\goth a}_N^{\lambda},b^{\lambda}), \widetilde{E_N^{\lambda}})
\overset{\pi_{N+1}^{\lambda}}\to{\leftarrow} \cdot\cdot\cdot
\overset{\pi_{N+S}}\to{\leftarrow} (\widetilde{W_{N+S}^{\lambda}}, ({\goth
a}_{N+S}^{\lambda},b^{\lambda}), \widetilde{E_{N+S}^{\lambda}})$$ 
where $\pi_{N+i}^{\lambda}$ is the blowup with center
$$\align
\widetilde{H_{r+N+i-1}^{\lambda}} &= Y_{N+i-1} \cap \widetilde{W_{N+i-1}^{\lambda}} \\
&= F_{N+i-1}
\cap H_{r+N+i-1} \cap \widetilde{W_{N+i-1}^{\lambda}} \\
&= F_{N+i-1}
\cap H_{r+N+i-1} \cap W_{N+i-1}^{\lambda}.\\
\endalign$$

Note that these transformations are set-theoretically nothing but identities with
$$\widetilde{H_{r+N}^{\lambda}} = \widetilde{H_{r+N+1}^{\lambda}} = \cdot\cdot\cdot =
\widetilde{H_{r+N+S-1}^{\lambda}}.$$

Therefore, by condition (GB-2) this gives rise to a sequence of transformations in the
collection
${\goth C}$ of the general basic object $({\Cal F}_0, (W_0,E_0))$ as long as \linebreak
$\widetilde{H_{r+N+i-1}^{\lambda}} \subset \roman{Sing}({\goth
a}_{N+i-1}^{\lambda},b^{\lambda})$ for $i = 1, ... , S$, the condition which
translates into the following equivalent conditions:
$$\align
& (N - 1)(\beta^{\lambda} - b^{\lambda}) - (i - 1)b^{\lambda} \geq b^{\lambda}\ \text{for}\ i = 1, ... , S
\\ 
\Longleftrightarrow \hskip.2in & (N - 1)(\beta^{\lambda} - b^{\lambda}) - (S - 1)b^{\lambda} \geq
b^{\lambda} \\
\Longleftrightarrow \hskip.2in & (N - 1)(\beta^{\lambda} - b^{\lambda}) \geq S b^{\lambda} \\
\Longleftrightarrow \hskip.2in & \left[\frac{(N - 1)(\beta^{\lambda} -
b^{\lambda})}{b^{\lambda}}\right] \geq S,\\
\endalign$$ 
where $\left[\quad\right]$ is the Gauss symbol, representing the integer $\left[x\right] = \alpha
\in {\Bbb N}$ such that
$\alpha
\leq x <
\alpha + 1$.

Therefore, we finally conclude by conditions (GB-0) (GB-1) and (GB-2) that for a fixed
sufficiently large
$N
\in {\Bbb N}$,
$$\left[\frac{(N - 1)(\beta^{\lambda} -
b^{\lambda})}{b^{\lambda}}\right]$$
is characterized as the largest integer $S_N$ such that the sequence described as above of one
smooth morphism followed by the transformations
$$\align
(F_0, (W_0,E_0)) &\leftarrow (F_1, (W_1,E_1)) \leftarrow \cdot\cdot\cdot \leftarrow (F_N, (W_N,E_N)) \\
&\leftarrow (F_{N+1}, (W_{N+1},E_{N+1})) \leftarrow \cdot\cdot\cdot \leftarrow
(F_{N+S_N},(W_{N+S_N},E_{N+S_N}))\\
\endalign$$
is in the collection ${\goth C}$ represented by the general basic object $({\Cal F}_0,
(W_0,E_0))$, and hence that the number
$$\frac{\nu_{x_0^{\lambda}}({\goth a}_0^{\lambda})}{b^{\lambda}} - 1 =
\frac{\beta^{\lambda}}{b^{\lambda}} - 1 = \lim_{N \rightarrow \infty}\frac{1}{N -
1}\left[\frac{(N - 1)(\beta^{\lambda} - b^{\lambda})}{b^{\lambda}}\right] = \lim_{N \rightarrow \infty}\frac{1}{N -
1}S_N$$
is characterized purely in terms of the collection ${\goth C}$ represented by the general basic object $({\Cal F}_0,
(W_0,E_0))$ and hence of independent of $\lambda$.

\vskip.1in

This completes the proof of (i), verifying that $\roman{ord}_0$ is a well-defined function on the
singular locus $F_0$ of the general basic object $({\Cal F}_0, (W_0,E_0))$.  (Therefore,
bringing back the subscripts right, we complete the proof that $\roman{ord}_k$ is a well-defined
function on the singular locus $F_k$ of the general basic object $({\Cal F}_k, (W_k,E_k))$.

\vskip.1in

(ii) As in the proof of (i), it suffices to show that the functions
$\{w\text{-}\roman{ord}_k^{\lambda}\}$, defined on the individual charts as in Definition 1-10 (ii), patch
up.  That is to say, it suffices to show that for any closed point $x_k \in F_k$ and commutative diagram of
the form
$$\CD
x_k^{\lambda} @. \in @. \roman{Sing}({\goth a}_k^{\lambda},b^{\lambda}) @. \subset @.
\widetilde{W_k^{\lambda}} \\
@.@.@.@. \cap \\
@VV j_k^{\lambda} V @.@.@. W_k^{\lambda} \\
@.@.@.@. \cup \\
x_k @.@. \in @.@. W_k^{\lambda} \cap W_k^{\lambda'}\\
@.@.@.@. \cap \\
@AA j_k^{\lambda'} A @.@.@. W_k^{\lambda'} \\
@.@.@.@. \cup \\
x_k^{\lambda'} @. \in @. \roman{Sing}({\goth a}_k^{\lambda'},b^{\lambda'}) @. \subset @.
\widetilde{W_k^{\lambda'}} \\
\endCD$$
we have
$$w\text{-}\roman{ord}_k^{\lambda}(x_k^{\lambda}) = \frac{\nu_{x_k^{\lambda}}(\overline{{\goth
a}_k^{\lambda}})}{b^{\lambda}} =
\frac{\nu_{x_k^{\lambda'}}(\overline{{\goth a}_k^{\lambda'}})}{b^{\lambda'}} =
w\text{-}\roman{ord}_k^{\lambda'}(x_k^{\lambda'}).$$

We will show that the number
$$w\text{-}\roman{ord}_k^{\lambda}(x_k^{\lambda})$$
can be determined purely in terms of the invariants $\roman{ord}_i$ for the general basic objects
$({\Cal F}_i, (W_i,E_i))$ for $i = 0, ... , k$ and in terms of the sequence
$$(F_0, (W_0,E_0)) \leftarrow \cdot\cdot\cdot \leftarrow (F_k, (W_k,E_k)),$$
and hence independent of $\lambda$.
 
\vskip.1in

\proclaim{Claim 4-6} We have the formula
$$(*)\ w\text{-}\roman{ord}_k^{\lambda}(x_k) = \roman{ord}_k(x_k) - \Sigma_{j = 1}^k
\Sigma_{H_{r+j,l} \subset H_{r+j}}\{\roman{ord}_{i_{H_{r+j,l}}}(\eta_{Y_{H_{r+j,l}}}) - 1\}
\cdot \epsilon_{H_{r+j,l},x_k}$$ 
where 

$E_k = \{H_1, ... , H_r, H_{r+1}, ... , H_{r+k}\}$ (cf. the convention in
Definition 1-8 (iii)),

$H_{r+j,l}$ are the irreducible components of $H_{r+j}$ with the generic points
$\eta_{H_{r+j,l}}$,

the number $i_{H_{r+j,l}}$ is the maximum of such $i$ that $\overline{\pi_{i+1} \circ
\cdot\cdot\cdot \circ \pi_k(\eta_{H_{r+j,l}})}$ is an irreducible component of the center $Y_i$
(See the convention explained in Definition 1-8 (iii).),

$Y_{H_{r+j,l}} = \pi_{i_{H_{r+j,l}} + 1} \circ \cdot\cdot\cdot \circ \pi_k(H_{r+j,l})$ with the
generic point $\eta_{Y_{H_{r+j,l}}}$, and

$$\epsilon_{H_{r+j,l},x_k} = \left\{\aligned
0 & \text{\ if\ }x_k \not\in H_{r+j,l} \\
1 & \text{\ if\ }x_k \in H_{r+j,l}.\\
\endaligned\right.$$

In particular, since the right hand side does not depend on $\lambda$,
$w\text{-}\roman{ord}_k^{\lambda}(x_k)$ is also independent of $\lambda$ as desired.
\endproclaim

\demo{Proof}\enddemo Observe
$${\goth a}_k^{\lambda} = I(\widetilde{H_{r+1}^{\lambda}})^{a_{r+1}^{\lambda}} \cdot\cdot\cdot
I(\widetilde{H_{r+k}^{\lambda}})^{a_{r+k}}\overline{{\goth a}_k^{\lambda}}$$
where $\widetilde{H_{r+j}^{\lambda}} = H_{r+j} \cap \widetilde{W_k^{\lambda}}$ and
$a_{r+j}^{\lambda}$ are multi-indices so that
$$I(\widetilde{H_{r+j}^{\lambda}})^{a_{r+j}^{\lambda}} = \prod
I(\widetilde{H_{r+j,l}^{\lambda}})^{a_{r+j,l}^{\lambda}}$$
where $\widetilde{H_{r+j,l}^{\lambda}} = H_{r+j,l} \cap \widetilde{W_j^{\lambda}}$ are
the irreducible components of $\widetilde{H_{r+j}^{\lambda}}$.  (Note that the irreducible
components $H_{r+j,l}$ of $H_{r+j}$ are in one-to-one correspondence with the irreducible
components $\widetilde{H_{r+j,l}}$ of $\widetilde{H_{r+j}}$.)  

Therefore, we conclude
$$(*)_{\lambda}\hskip.1in \roman{ord}_k^{\lambda}(x_k^{\lambda}) =
w\text{-}\roman{ord}_k(x_k^{\lambda}) +
\Sigma_{j = 1}^k
\Sigma_{\widetilde{H_{r+j,l}}
\subset
\widetilde{H_{r+j}}}\{\roman{ord}_{i_{\widetilde{H_{r+j,l}}}}(\eta_{Y_{\widetilde{H_{r+j,l}}}}) -
1\}
\cdot \epsilon_{\widetilde{H_{r+j,l}},x_k^{\lambda}}$$
where

\vskip.1in

$\widetilde{H_{r+j,l}}$ are the irreducible components of $\widetilde{H_{r+j}}$ with the generic
points
$\eta_{\widetilde{H_{r+j,l}}}$,

the number $i_{\widetilde{H_{r+j,l}}}$ is the maximum of such $i$ that
$\overline{\pi_{i+1}^{\lambda}
\circ
\cdot\cdot\cdot \pi_k^{\lambda}(\eta_{\widetilde{H_{r+j,l}}})}$ is an irreducible component of
$Y_i$,

$Y_{\widetilde{H_{r+j,l}}} = \pi_{i_{\widetilde{H_{r+j,l}}} + 1}^{\lambda} \circ \cdot\cdot\cdot
\circ
\pi_k^{\lambda}(\widetilde{H_{r+j,l}})$ with the generic point
$\eta_{Y_{\widetilde{H_{r+j,l}}}}$, and
$$\epsilon_{\widetilde{H_{r+j,l}},x_k^{\lambda}} = \left\{\aligned
0 & \text{\ if\ }x_k^{\lambda} \not\in \widetilde{H_{r+j,l}} \\
1 & \text{\ if\ }x_k^{\lambda} \in \widetilde{H_{r+j,l}}.\\
\endaligned\right.$$

Note that, denoting $\pi_{i_{\widetilde{H_{r+j,l}}}+1}^{-1}(Y_{\widetilde{H_{r+j,l}}})$
by $\widetilde{H_{r+j,l}}$ and its generic point $\eta_{\widetilde{H_{r+j,l}}}$ by abuse of
notation, we compute
$$\align
\frac{a_{r+j,l}^{\lambda}}{b^{\lambda}} =
\frac{\nu_{\eta_{\widetilde{H_{r+j,l}}}}({\goth a}_k^{\lambda})}{b^{\lambda}} &=
\frac{\nu_{\eta_{\widetilde{H_{r+j,l}}}}({\goth
a}_{i_{\widetilde{H_{r+j,l}}}+1}^{\lambda})}{b^{\lambda}} \\ 
&= \frac{\nu_{\eta_{Y_{\widetilde{H_{r+j,l}}}}}({\goth a}_{i_{\widetilde{H_{r+j,l}}}}^{\lambda})
- b^{\lambda}}{b^{\lambda}} \\
&= \roman{ord}_{i_{H_{r+j,l}}}(\eta_{Y_{\widetilde{H_{r+j,l}}}}) - 1.\\
\endalign$$

Remark that for the corresponding irreducible components $\widetilde{H_{r+j,l}} \subset
H_{r+j,l}$, the numbers and the generic points of the center coincide
$$\align
i_{\widetilde{H_{r+j,l}}} &= i_{H_{r+j,l}} \\
\eta_{Y_{\widetilde{H_{r+j,l}}}} &= \eta_{Y_{H_{r+j,l}}}.\\
\endalign$$

The formula $(*)$ in the claim now follows from this remark and the formula $(*)_{\lambda}$. 
\vskip.1in

(iii) Since $w\text{-}\roman{ord}$ is a well-defined invariant on a general basic object by (ii),
and since the inequalities
$$\align
\max w\text{-}\roman{ord}_0 &\geq \max w\text{-}\roman{ord}_1 \geq \cdot\cdot\cdot \\
\max w\text{-}\roman{ord}_{i-1}
&\geq \max w\text{-}\roman{ord}_i \\
\cdot\cdot\cdot \geq \max w\text{-}\roman{ord}_{k-1} &\geq \max
w\text{-}\roman{ord}_k.\\
\endalign$$
follow from those for the basic objects under condition $(\heartsuit)$, the invariant
\linebreak
$t_k:F_k \rightarrow \frac{1}{c!}{\Bbb Z}_{\geq 0} \times {\Bbb Z}_{\geq 0}$ is a
function well-defined globally on $F_k$.

In order to verify that the restriction $t_k|_{F_k \cap W_k^{\lambda}}$ coincides with
$t_k^{\lambda}$ defined on $\roman{Sing}({\goth a}_k^{\lambda},b^{\lambda}) = F_k \cap
W_k^{\lambda}$, one has only to observe that $t_k$ has the local description identical to the one
given in Remark 1-11 (v), which is easily seen to coincide with the local description of
$t_k^{\lambda}$ also given in Remark 1-11 (v).

\vskip.1in

(iv) Let $E_k = \{H_1, ... , H_r, H_{r+1}, ... , H_{r+k}\}$.  To a point
$p \in H_{r+j}$ we assign the
following number $\alpha_{r+j}(p)$
$$\alpha_{r+j}(p) = \roman{ord}_{i_{H_{r+j,l}}}(\eta_{Y_{\widetilde{H_{r+j,l}}}}) - 1,$$
where $H_{r+j,l}$ is the irreducible component of $H_{r+j}$ containing $p \in H_{r+j,l}$ (cf.
the formula for $\frac{a_{r+j,l}^{\lambda}}{b^{\lambda}}$ at the end of the proof for (ii)).

We define the
invariant 
$$\Gamma_k:F_k \rightarrow {\Bbb Z}_{\geq - d_k} \times \frac{1}{c!}{\Bbb Z}_{\geq
0}
\times {{\Bbb Z}_{\geq 0}}^{d_k},$$
where
$d_k = \dim \widetilde{W_k^{\lambda}}$
(independent of
$\lambda$), to be a function defined over
$F_k$ such that
$$\Gamma_k(p) = (\Gamma_{k1}(p),\Gamma_{k2}(p),\Gamma_{k3}(p)) \text{\ for\ }p \in F_k$$
where
$$\align
- \Gamma_{k1}(p) &= \min \{n; \exists r+j_1, ... , r+j_n \text{\ s.t.\
}{\alpha}_{r+j_1}(p) +
\cdot\cdot\cdot + {\alpha}_{r+j_n}(p) \geq 1, \\
&\hskip.1in p \in H_{r+j_1} \cap \cdot\cdot\cdot \cap
H_{r+j_n}\}
\\
\Gamma_{k2}(p) &= \max \{{\alpha}_{r+j_1}(p) + \cdot\cdot\cdot + {\alpha}_{r+j_n}(p);n = -
\Gamma_{k1}(p), {\alpha}_{r+j_1}(p) +
\cdot\cdot\cdot + {\alpha}_{r+j_n}(p) \geq 1, \\
&\hskip.1in p \in H_{r+j_1} \cap \cdot\cdot\cdot \cap
H_{r+j_n}\}
\\
\Gamma_{k3}(p) &= \max \{(r+j_1, ... , r+j_n);n = -
\Gamma_{k1}(p), \Gamma_{k2}(p) = {\alpha}_{r+j_1}(p) + \cdot\cdot\cdot + {\alpha}_{r+j_n}(p),\\
&\hskip.1in p
\in H_{r+j_1}
\cap
\cdot\cdot\cdot \cap H_{r+j_n}, r + j_1 \geq \cdot\cdot\cdot \geq r + j_n\}\\
&\text{with\ the maximum\ taken\ with\ respect\ to\ the\
lexicographical\ order.}\\
&\text{We\ identify\ }(r+j_1, ... , r+j_n) \text{\ with\ }(r+j_1, ... ,
r+j_n, 0, ... , 0) \in {{\Bbb Z}_{\geq 0}}^{d_k}.\\
\endalign$$  
(We order the values of $\Gamma_k$ according to the lexicographical order given to \linebreak
${\Bbb
Z}_{\geq - d_k}
\times
\frac{1}{c!}{\Bbb Z}_{\geq 0}
\times {{\Bbb Z}_{\geq 0}}^{d_k}$.)

\vskip.1in

It is clear from the above description that for each $\lambda$ the invariant $\Gamma_k$
restricts to the invariant $\Gamma_k^{\lambda}$ of the monomial basic object
$(\widetilde{W_k^{\lambda}}, ({\goth a}_k^{\lambda},b^{\lambda}), \widetilde{E_k^{\lambda}})$
(cf. Definition 2-3) (shrinking $\widetilde{W_k^{\lambda}}$ to an open neighborhood of
$\roman{Sing}({\goth a}_k^{\lambda},b^{\lambda}) = F_k \cap W_k^{\lambda}$ if necessary (cf.
Corollary 2-7)). 

\vskip.1in

Finally, from the proof of (i), it is clear that the invariant $\roman{ord}_k$ is purely determined
in terms of the collection ${\goth C}_k$ of sequences of transformations and smooth morphisms
represented by the general basic object $({\Cal F}_k,(W_k,E_k))$, once the dimension $d_k =
\dim W_k - \dim W_0 + d_0$ of the structure is specified.  From Claim 4-6, it follows that the
invariant
$w\text{-}\roman{ord}_k$ is purely determined in terms of
$\roman{ord}_0, ... ,
\roman{ord}_k$ and by looking at the set-theoretical behavior of $E_1, ... , E_k$ in the
original sequence.  The invariant $t_k$ is purely detrmined in terms of $w\text{-}\roman{ord}_0,
... , w\text{-}\roman{ord}_k$ and by looking at the set-theoretical behavior of $E_1, ... ,
E_k$ in the original sequence.  The invariant $\Gamma_k$ is purely determined in terms of
$\roman{ord}_k$ and by looking at the set-theoretical behavior of $E_1, ... , E_k$ in the original
sequence.  This verifies the ``Finally" part of definition-Proposition 4-5.  

\vskip.1in

This completes the proof of Definition-Proposition 4-5.

\vskip.1in

\proclaim{Remark 4-7}\endproclaim

A general basic object $({\Cal F},(W,E))$ sometimes can have a $d$-dimensional structure as
well as a $d'$-dimensional structure for two different numbers $d \neq d'$.  That is to say, we
 can have two different sets of charts
$\{(\widetilde{W^{\lambda}},({\goth
a}^{\lambda},b^{\lambda}),\widetilde{E^{\lambda}})\}_{\lambda \in \Lambda}$ and
$\{(\widetilde{W^{\mu}},({\goth b}^{\mu},c^{\mu}),\widetilde{E^{\mu}})\}_{\mu \in M}$ being of
different dimensions $d$ and $d'$, i.e., $\dim \widetilde{W^{\lambda}} = d \neq d' =
\widetilde{W^{\mu}}$, but giving rise to the same collection ${\goth C}$ of sequences of smooth
morphisms and transformations of pairs with specified closed subsets, represented by the
general basic object $({\Cal F},
(W,E))$.

The invariants $\roman{ord}$ and $w\text{-}\roman{ord}$ DO depend on the specification of the
dimension of the structure of your choice, as the following example demonstrates (cf. Remark
3-7 (ii)).

\vskip.1in

Take a basic object 
$$(W,(J,b),E) = ({\Bbb A}^2,(\langle x^2 - y^3\rangle,2),\emptyset),$$ 
which
defines a general basic object with a $2$-dimensional structure via Remark 4-2 (ii).

By the key inductive lemma, the same general basic object has a \linebreak
$2 - 1 = 1$-dimensional
structure with a (global) chart 
$$({\Bbb A}^1 = \{x = 0\} = \roman{Spec}\ k[y],(\langle
(y^3)^{\frac{2!}{2-0}},(y^2)^{\frac{2!}{2-1}}\rangle = \langle y^3\rangle,2!),\emptyset).$$

Denoting by $\roman{ord}_0^{(2)}, w\text{-}\roman{ord}_0^{(2)}$ the $\roman{ord}-$ and
$w\text{-}\roman{ord}-$invarinats with respect to the 2-dimensional structure of the general
basic object (considered to form a trivial sequence by itself), we have
$$\roman{ord}_0^{(2)}(0) = w\text{-}\roman{ord}_0^{(2)}(0) = \frac{2}{2} = 1.$$
On the other hand, denoting by $\roman{ord}_0^{(1)}, w\text{-}\roman{ord}_0^{(1)}$ the
$\roman{ord}-$ and
$w\text{-}\roman{ord}-$invarinats with respect to the 1-dimensional structure of the general
basic object (considered to form a trivial sequence by itself), we have
$$\roman{ord}_0^{(1)}(0) = w\text{-}\roman{ord}_0^{(1)}(0) = \frac{3}{2!} = \frac{3}{2}.$$

Therefore, theoretically and strictly speaking, it might be more appropriate to put the
superscript
$d_k$, such as $\roman{ord}_k^{(d_k)}$ or $w\text{-}\roman{ord}_k^{(d_k)}$ to indicate the
dependence of the invariants $\roman{ord}_k$ or $w\text{-}\roman{ord}_k$ upon the specified
dimension $d_k$ of the structure.  However, we will omit the superscript for notational
simplicity, since little confusion is likely to occur.  

\vskip.2in

\proclaim{Proposition 4-8 (Properties of key invariants for general basic objects)} Let 
$$\align
({\Cal F}_0, (W_0,E_0)) &\overset{\pi_1}\to{\leftarrow} ({\Cal F}_1, (W_1,E_1))
\overset{\pi_2}\to{\leftarrow} \cdot\cdot\cdot \\
({\Cal F}_{i-1}, (W_{i-1}, E_{i-1})) &\overset{\pi_i}\to{\leftarrow}
({\Cal F}_i, (W_i,E_i)) \\
\cdot\cdot\cdot \overset{\pi_{k-1}}\to{\leftarrow} ({\Cal F}_{k-1}, (W_{k-1},E_{k-1})
&\overset{\pi_k}\to{\leftarrow} ({\Cal F}_k, (W_k,E_k)) \\
\endalign$$
be a sequence of transformations and smooth morphisms of general basic objects.

\vskip.1in

(i) The invariants $\roman{ord}_k$ and $w\text{-}\roman{ord}_k$ are upper semi-continuous
functions.

(ii) Note first that
$$F_{i-1} \supset \pi_i(F_i) \text{\ for\ } i = 1, ... , k.$$

Suppose that the sequence satisfies condition $(\heartsuit)$ (See Definition-Proposition 4-5
(iii).).  

Then for $i = 1, ... , k$ we have inequalities
$$w\text{-}\roman{ord}_{i-1}(\xi_{i-1})
\geq w\text{-}\roman{ord}_i(\xi_i)$$
where $\xi_i \in F_i$ and $\xi_{i-1} =
\pi_i(\xi_i) \in F_{i-1}$, which imply
$$\max\ w\text{-}\roman{ord}_{i-1}
\geq \max\ w\text{-}\roman{ord}_i.$$
That is to say, we have
$$\align
\max\ w\text{-}\roman{ord}_0 &\geq \max\ w\text{-}\roman{ord}_1 \geq \cdot\cdot\cdot \\
\max\ w\text{-}\roman{ord}_{i-1} &\geq \max\ w\text{-}\roman{ord}_i \\
\cdot\cdot\cdot \geq \max\ w\text{-}\roman{ord}_{k-1} &\geq \max\ w\text{-}\roman{ord}_k. \\
\endalign$$ 

The invariant $t_k$ is an upper semi-continuous function and we have
inequalities
$$t_{i-1}(\xi_{i-1}) \geq t_i(\xi_i)$$
where $\xi_i \in F_i$ and $\xi_{i-1} =
\pi_i(\xi_i) \in F_{i-1}$, which imply
$$\max\ t_{i-1}
\geq \max\ t_i.$$
That is to say, we have
$$\align
\max\ t_0 &\geq \max\ t_1 \geq \cdot\cdot\cdot \\
\max\ t_{i-1} &\geq \max\ t_i \\
\cdot\cdot\cdot \geq \max\ t_{k-1} &\geq \max\ t_k. \\
\endalign$$ 
\endproclaim

\demo{Proof}\enddemo The proof is identical to the one for basic objects (cf. Proposition
1-12) and left to the reader as an exercise.

\vskip.1in

\proclaim{Corollary 4-9 (Resolution of singularities of a general basic object with
$\max\ w\text{-}\roman{ord} = 0$)} Let
$$({\Cal F}_0, (W_0,E_0)) \leftarrow \cdot\cdot\cdot \leftarrow ({\Cal F}_k, (W_k,E_k))$$
be a sequence of transformations and smooth morphisms of general basic objecs.  

Suppose $\max\
w\text{-}\roman{ord}_k = 0$.  
 
Then there exists a sequence of transformations only of general basic objects
$$\align
({\Cal F}_k, (W_k,E_k)) &\overset{\pi_{k+1}}\to{\leftarrow} ({\Cal F}_{k+1},
(W_{k+1},E_{k+1})) \overset{\pi_{k+2}}\to{\leftarrow}
\cdot\cdot\cdot \\ 
({\Cal F}_{k+i-1},
(W_{k+i-1},E_{k+i-1})) &\overset{\pi_{k+i}}\to{\leftarrow} ({\Cal F}_{k+i},
(W_{k+i},E_{k+i})) \\
\cdot\cdot\cdot \overset{\pi_{k+N-1}}\to{\leftarrow} ({\Cal F}_{k+N-1}, (W_{k+N-1},E_{k+N-1}))
&\overset{\pi_{k+N}}\to{\leftarrow} ({\Cal F}_{k+N}, (W_{k+N},E_{k+N}))\\
\endalign$$
which represents resolution of singularities, i.e.,
$$F_{k+N} = \emptyset,$$
where $\pi_{k+i}\hskip.1in \text{for}\hskip.1in i = 1, ... , N$ are the transformations with
centers
$$Y_{k+i-1} = \underline{\roman{Max}}\ \Gamma_{k+i-1} \subset F_{k+i-1}.$$
\endproclaim

\demo{Proof}\enddemo We note that under the specified transformations $\max\ w\text{-}\roman{ord}$
remains zero, i.e.,
$$\max\ w\text{-}\roman{ord}_{k+i} = 0 \text{\ for\ }i = 0, ... , N-1$$
and hence that the invariant $\Gamma_{k+i}$ is well-defined.  The rest of the proof is identical to
the one for resolution of singularities for basic objects with $\max\ w\text{-}\roman{ord} = 0$ and
left to the reader as an exercise (cf. Corollary 2-7).

\newpage

$\bold{CHAPTER\ 5.\ INDUCTIVE\ ALGORITHM\ FOR\ RESOLUTION}$ \linebreak
${}\hskip.3in$ $\bold{OF\ SINGULARITIES\ OF\ GENERAL\ BASIC\ OBJECTS}$

\vskip.1in

This chapter is the culmination of the ideas of Encinas and Villamayor, seeing how we overcome
the shortcomings (cf. Remark 3-2) of the key inductive lemma (Lemma 3-1) via the use of the
$t$-invariant and transform the lemma into a genuine inductive algorithm of resolution of
singularities of general basic objects.

\proclaim{Theorem 5-1 (Inductive algorithm for resolution of singularities of general basic
objects)} Let $({\Cal F}_0, (W_0,E_0))$ be a general basic object over $(F_0, (W_0,E_0))$ with a
$d$-dimensional structure given by an open covering $\{W^{\lambda}\}_{\lambda \in \Lambda}$ and
charts
$\{(\widetilde{W_0^{\lambda}}, ({\goth a}_0^{\lambda},b^{\lambda}),
\widetilde{E_0^{\lambda}})\}$ of $d$-dimensional basic objects.  Let ${\goth C}$ denote the
collection of sequences of transformations and smooth morphisms represented by $({\Cal F}_0,
(W_0,E_0))$ (cf. Definition 4-1).

Then there exists an inductive algorithm which provides a seqeuence of trasformations of
pairs with specified closed subsets in the collection ${\goth C}$, representing resolution of
singularities of the general basic object
$({\Cal F}_0,(W_0,E_0))$
$$(F_0,(W_0,E_0)) \leftarrow \cdot\cdot\cdot \leftarrow (F_l,(W_l,E_l)) \text{\
with\ }F_l = \emptyset,$$
by uniquely specifying the centers satisfying the following condition $(\heartsuit')$
\linebreak
for
$i = 1, ... , l$
$$(\heartsuit')\ Y_{i-1} \subset \underline{\roman{Max}}\ t_{i-1} \subset
\underline{\roman{Max}}\ w\text{-}\roman{ord}_{i-1} \subset F_{i-1} \text{\ if\ }\max\
w\text{-}\roman{ord}_{i-1} > 0.$$

\vskip.1in

The process of the inductive algorithm is described below:

\vskip.1in
 
Suppose we have constructed the sequence of transformations up to the $k$-th stage via the
inductive algorithm 
$$(F_0, (W_0,E_0)) \leftarrow \cdot\cdot\cdot \leftarrow (F_k, (W_k,E_k))$$
with the centers satisfying condition $(\heartsuit')$ for $i = 1, ... , k$
$$(\heartsuit')\ Y_{i-1} \subset \underline{\roman{Max}}\ t_{i-1} \subset
\underline{\roman{Max}}\ w\text{-}\roman{ord}_{i-1} \subset F_{i-1} \text{\ if\ }\max\
w\text{-}\roman{ord}_{i-1} > 0.$$
(Recall that condition $(\heartsuit)$ 
$$(\heartsuit)\ Y_{i-1} \subset \underline{\roman{Max}}\ w\text{-}\roman{ord}_{i-1} \subset
F_{i-1}
\text{\ for\ }i = 1, ... , k,$$
which obviously follows from condition $(\heartsuit')$, implies inequalities
$$\align
\max\ w\text{-}\roman{ord}_0 \geq \cdot\cdot\cdot \geq \max\ w\text{-}\roman{ord}_{i-1} &\geq
\max\ w\text{-}\roman{ord}_i
\geq \cdot\cdot\cdot \geq \max\ w\text{-}\roman{ord}_k \\
\max\ t_0 \geq \cdot\cdot\cdot \geq \max\ t_{i-1} &\geq \max\ t_i
\geq \cdot\cdot\cdot \geq \max\ t_k. \\
\endalign$$
Recall also that the sequence induces an sequence of general basic objects
$$({\Cal F}_0, (W_0,E_0)) \leftarrow \cdot\cdot\cdot \leftarrow ({\Cal F}_k, (W_k,E_k))$$
as explained in Note 4-3.)

\vskip.1in

Then we have the following three possibilities:

\vskip.1in

$\bold{\underline{P1}}: F_k = \emptyset.$ 

\vskip.1in

Under this possibility, the sequence represents resolution of singularities of the general basic
object
$({\Cal F}_0, (W_0,E_0))$.

\vskip.1in

$\bold{\underline{P2}}: F_k \neq \emptyset$ and $\max\ w\text{-}\roman{ord}_k = 0.$

\vskip.1in

Under this possibility, we can apply Corollary 4-9 to $({\Cal F}_k, (W_k,E_k))$ and create a
sequence of transformations representing resolution of singularities of $({\Cal F}_k, (W_k,E_k))$. 
Attaching this to the original sequence, we obtain a sequence representing resolution of
singularities of
$({\Cal F}_0, (W_0,E_0))$.

\vskip.1in

$\bold{\underline{P3}}: F_k \neq \emptyset$ and $\max\ w\text{-}\roman{ord}_k > 0.$

\vskip.1in

Under this possibility $\bold{\underline{P3}}$, there are two cases $\bold{Case\ A}$ and
$\bold{Case\ B}$.  

\vskip.1in

We denote by
$R(1)(\underline{\roman{Max}}\ t_k)$ the $(d-1)$-dimensional part (i.e., codimension one with
respect to the $d$-dimensional $\widetilde{W_k^{\lambda}}$'s) of the locus
$\underline{\roman{Max}}\ t_k \subset F_k$.

\vskip.1in

$\bold{Case\ A}$: $R(1)(\underline{\roman{Max}}\ t_k) \neq \emptyset$.

\vskip.1in

In this case, $R(1)(\underline{\roman{Max}}\ t_k) (\subset \underline{\roman{Max}}\ t_k
\subset
\underline{\roman{Max}}\ w\text{-}\roman{ord}_k \subset F_k)$ is smooth, open and closed in
$\underline{\roman{Max}}\ w\text{-}\roman{ord}_k$ (i.e., a union of smooth connected components of
$\underline{\roman{Max}}\ w\text{-}\roman{ord}_k$ disjoint from each other).

The locus $R(1)(\underline{\roman{Max}}\ t_k) \cap W_k^{\lambda}$ is permissible for each
$(\widetilde{W_k^{\lambda}}, ({\goth a}_k^{\lambda},b^{\lambda}),
\widetilde{E_k^{\lambda}})$, i.e.,\linebreak 
$R(1)(\underline{\roman{Max}}\ t_k) \cap
W_k^{\lambda}
\subset
\widetilde{W_k^{\lambda}}$, $R(1)(\underline{\roman{Max}}\ t_k) \cap W_k^{\lambda}$ is
permissible with respect to $\widetilde{E_k^{\lambda}}$, and $R(1)(\underline{\roman{Max}}\ t_k)
\cap W_k^{\lambda} \subset \roman{Sing}({\goth a}_k^{\lambda},b^{\lambda})$. 

\vskip.1in

Take the transformation of pairs
$$(W_k,E_k) \overset{\pi_{k+1}}\to{\leftarrow} (W_{k+1},E_{k+1})$$
with center $Y_k = R(1)(\underline{\roman{Max}}\ t_k)$.

Take for each $\lambda$ the corresponding transformation of basic objects with center $Y_k \cap
W_k^{\lambda}$
$$(\widetilde{W_k^{\lambda}}, ({\goth a}_k^{\lambda},b^{\lambda}),
\widetilde{E_k^{\lambda}}) \overset{\pi_{k+1}^{\lambda}}\to{\leftarrow}
(\widetilde{W_{k+1}^{\lambda}}, ({\goth a}_{k+1}^{\lambda},b^{\lambda}),
\widetilde{E_{k+1}^{\lambda}}).$$

Then
$$F_{k+1} := \cup \roman{Sing}({\goth a}_{k+1}^{\lambda},b^{\lambda})$$
is a closed subset of $W_{k+1}$ with
$$F_{k+1} \cap W_{k+1}^{\lambda} = \roman{Sing}({\goth a}_{k+1}^{\lambda},b^{\lambda}),$$
and the extended sequence belongs to ${\goth C}$, i.e.,
$$(F_0, (W_0,E_0)) \leftarrow \cdot\cdot\cdot \leftarrow (F_k, (W_k,E_k)) \leftarrow (F_{k+1},
(W_{k+1},E_{k+1})) \in {\goth C}.$$

\vskip.1in

We have one of the following four cases:

\vskip.1in

$\bold{A\text{-}1}$: $F_{k+1} = \emptyset$.

\vskip.1in

In this case, the extended sequence represents resolution of singularities of the general basic
object $({\Cal F}_0, (W_0,E_0))$.

\vskip.1in

$\bold{A\text{-}2}$: $F_{k+1} \neq \emptyset$ and $\max\ w\text{-}\roman{ord}_{k+1} =
0.$

\vskip.1in

In this case, we can apply Corollary 4-9 to $({\Cal F}_{k+1}, (W_{k+1},E_{k+1}))$ and create a
sequence of transformations representing resolution of singularities of $({\Cal F}_{k+1},
(W_{k+1},E_{k+1}))$.  Attaching this to the extended sequence, we obtain a sequence representing
resolution of singularities of
$({\Cal F}_0, (W_0,E_0))$.

\vskip.1in

$\bold{A\text{-}3}$: $F_{k+1} \neq \emptyset$, $\max\ w\text{-}\roman{ord}_{k+1} > 0$, and $\max\
t_k >
\max\ t_{k+1}.$

\vskip.1in

In this case, obviously the maximum of the $t$-invariant drops.

\vskip.1in

$\bold{A\text{-}4}$: $F_{k+1} \neq \emptyset$, $\max\ w\text{-}\roman{ord}_{k+1} > 0$, $\max\
t_k =
\max\ t_{k+1}$, and \hfill\hfill\linebreak
$R(1)(\underline{\roman{Max}}\ t_{k+1}) = \emptyset.$

\vskip.1in

In this case, we go to Case B for the general basic object 
$({\Cal F}_{k+1}, (W_{k+1},E_{k+1}))$.  

\vskip.2in

$\bold{Case\ B}$: $R(1)(\underline{\roman{Max}}\ t_k) = \emptyset$.

\vskip.1in

In this case, we construct a general basic object over $(G_k = \underline{\roman{Max}}\ t_k,
(W_k,E_k''))$ with a $(d-1)$-dimensional structure with the following property:

We construct a sequence of transformations representing resolution of singularities of the
general basic object $({\Cal G}_k, (W_k,E_k''))$ by induction on the dimension of the structure
$$(G_k, (W_k,E_k'')) \leftarrow \cdot\cdot\cdot \leftarrow (G_{k+N} = \emptyset,
(W_{k+N},E_{k+N}'')).$$
We have the corresponding sequence of transformations (with the same centers)
$$(F_k, (W_k,E_k)) \leftarrow \cdot\cdot\cdot \leftarrow (F_{k+N}, (W_{k+N},E_{k+N}))$$
which belongs to the collection ${\goth C}$ represented by the original general basic object
$({\Cal F}_0, (W_0,E_0))$ when attached to the original sequence of transformations
$$(F_0, (W_0,E_0)) \leftarrow \cdot\cdot\cdot \leftarrow (F_k, (W_k,E_k)) \leftarrow
\cdot\cdot\cdot \leftarrow (F_{k+N}, (W_{k+N},E_{k+N}))$$
such that it satisfies the conditions
$$\align
&(i)\ \text{\ for\ } j = 1, ... , N, \text{\ we\ have}\\
&(\heartsuit')\ Y_{k+j-1} \subset \underline{\roman{Max}}\ t_{k+j-1} \subset
\underline{\roman{Max}}\ w\text{-}\roman{ord}_{k+j-1} \subset F_{k+j-1} \text{\ where\ }\max\
w\text{-}\roman{ord}_{k+j-1} > 0,\\
&\hskip.5in \text{(Recall\ that\ we\ have\ condition}\\
& \hskip1.2in (\heartsuit')\ Y_{i-1} \subset \underline{\roman{Max}}\ t_{i-1} \subset
\underline{\roman{Max}}\ w\text{-}\roman{ord}_{i-1}\text{\
for\ }i = 1, ... , k\\ 
&\hskip.5in \text{where\ }\max\ w\text{-}\roman{ord}_{i-1} \geq \max\ w\text{-}\roman{ord}_k > 0
\text{\ by\ the\ case\ assumption.}\\
&\hskip.5in \text{Note\ also\ that\ the\ condition\ (iii)\ below\
implies}\\   
&\hskip1.2in \max\ w\text{-}\roman{ord}_k = \cdot\cdot\cdot = \max\ w\text{-}\roman{ord}_{k+N-1}
> 0.)\\ 
&(ii)\ \underline{\roman{Max}}\ t_{k+j-1} = G_{k+j-1}\text{\ for\ }j = 1, ... ,
N,\\ &(iii)\ \max\ t_k = \cdot\cdot\cdot = \max\ t_{k+N-1}\\
\endalign$$
(Note that the invariants $t_{k+j-1}$ and $w\text{-}\roman{ord}_{k+j-1}$ for $j =
1, ... , N$ are the ones defined for the general basic objects $({\Cal
F}_{k+j-1},(W_{k+j-1},E_{k+j-1}))$ with the $d$-dimensional structures.)

and that we have one of the following three cases:

\vskip.1in

$\bold{B\text{-}1}$: $F_{k+N} = \emptyset$.

\vskip.1in

In this case, the extended sequence represents resolution of
singularities of the general basic object $({\Cal F}_0, (W_0,E_0))$.

\vskip.1in

$\bold{B\text{-}2}$: $F_{k+N} \neq \emptyset$ and $\max\ w\text{-}\roman{ord}_{k+N} =
0.$

\vskip.1in

In this case, we can apply Corollary 4-9 to $({\Cal F}_{k+N}, (W_{k+N},E_{k+N}))$ and create a
sequence of transformations representing resolution of singularities of \linebreak
$({\Cal
F}_{k+N}, (W_{k+N},E_{k+N}))$.  Attaching this to the extended sequence, we
obtain a sequence representing resolution of singularities of
$({\Cal F}_0, (W_0,E_0))$.

\vskip.1in

$\bold{B\text{-}3}$: $F_{k+N} \neq \emptyset$, $\max\ w\text{-}\roman{ord}_{k+N} > 0$, and $\max\
t_k >
\max\ t_{k+N}.$

\vskip.1in

In this case, obviously the maximum of the $t$-invariant drops.

\vskip.1in

Since the set of values of the $t$-invariant satisfies the descending chain condition, after
finitely many executions of the process described as above, we obtain the uniquely determined
sequence of transformations representing resolution of singularities of
$({\Cal F}_0,(W_0,E_0))$ via the inductive algorithm. 
\endproclaim

\vskip.1in

\proclaim{Remark 5-2}\endproclaim

(i) In $\bold{Case\ B}$, we first construct such a general basic object $({\Cal G}_k, (W_k,
E_k''))$ over
$(G_k, (W_k,E_k''))$ with a
$d$-dimensional structure $\{(\widetilde{W_k^{\lambda}}, ({\goth b}_k^{\lambda},
e^{\lambda}),\widetilde{D^{\lambda}})\}$, where the basic objects in the charts are ``simple",
that satisfies the requirements specified in $\bold{Case\ B}$.  Then we find an open covering
of each
$\widetilde{W_k^{\lambda}}$ together with smooth hypersurfaces satisfying
conditions 1 and 2 as described in the key inductive lemma (Lemma 3-1). 
This is done via the power of the
$t$-invariant.  Now the key inductive lemma implies that the general basic object $({\Cal
G}_k, (W_k,E_k''))$ has a
$(d-1)$-dimensional structure, and hence completes the inductive step of the algorithm above.

(ii) In $\bold{Case\ B}$, the extended sequence
$$\align
&(F_0,(W_0,E_0)) \leftarrow \cdot\cdot\cdot \leftarrow (F_k,(W_k,E_k)) \leftarrow
\cdot\cdot\cdot \leftarrow (F_{k+j-1},(W_{k+j-1},E_{k+j-1})), \\
\endalign$$
for $j = 1, ... , N$, remains in $\bold{Case\ B}$ (cf. the construction described
in (i), Giraud's Lemma (Claim 3-4), and requirements (i) (ii) for
$\bold{Case\ B}$).  The general basic object over
$(\underline{\roman{Max}}\ t_{k+j},(W_{k+j},E_{k+j}''))$ we construct for
$({\Cal F}_{k+j},(W_{k+j},E_{k+j}))$ under the prescription of $\bold{Case\ B}$ coincides
with the general basic object $({\Cal G}_{k+j},(W_{k+j},E_{k+j}''))$ induced from $({\Cal
G}_k,(W_k,E_k''))$ via the sequence
$$(G_k,(W_k,E_k'')) \leftarrow \cdot\cdot\cdot \leftarrow (G_{k+j},(W_{k+j},E_{k+j}'')).$$

(iii) The process of extending the sequence
$$(F_0,(W_0,E_0)) \leftarrow \cdot\cdot\cdot \leftarrow (F_k,(W_k,E_k))$$
prescribed by the inductive algorithm, in order to obtain resolution of singularities of $({\Cal
F}_0,(W_0,E_0))$, may seem depend on the number $d$, which specifies the dimension of the
structure.  However, it is not difficult to see that the process is actually independent of $d$ and
that it is purely determined by the sequence of general basic objects
$$({\Cal F}_0,(W_0,E_0)) \leftarrow \cdot\cdot\cdot \leftarrow ({\Cal F}_k,(W_k,E_k)).$$   
In fact, if we are in possibility $\bold{\underline{P1}}$, then we are done and have nothing more
to do.  The process in possibility $\bold{\underline{P2}}$ is independent of $d$, since it only
depends on the invariant $\Gamma_k$, which is determined purely in terms of the collection ${\goth
C}_k$ represented by the general basic object $({\Cal F}_k,(W_k,E_k))$ (cf. Definition-Proposition
4-5).  So suppose we are in possibility
$\bold{\underline{P3}}$.  Suppose that the general basic objects have $d$-dimensional
structures and that at the $k$-th stage after $l$-repetitions of the processes in $\bold{Case\
B}$ we reach a general basic object with a $(d-l)$-dimensional structure, for which we are
either in possibility $\bold{\underline{P2}}$ or in $\bold{Case\ A}$ with center $Y_k
\subset F_k$.  If the same general basic objects have $d'$-dimensional structures (with $d'
\geq d$), then what happens at the $k$-th stage is that after $(d' - d) + l$-repetitions of
the processes in $\bold{Case\ B}$ we reach a general basic object with a $(d-l)
= (d' - (d' - d + l))$-dimensional structure, for which we are either in possibility
$\bold{\underline{P2}}$ or in
$\bold{Case\ A}$ with the same center $Y_k
\subset F_k$.  Therefore, we conclude that the dimension of the structures of the general basic
objects has no effect on the process.  

\vskip.1in

\demo{Proof of Theorem 5-1}\enddemo First we check the assertions for
$\bold{\underline{P1}}$, $\bold{\underline{P2}}$, and $\bold{Case\ A}$
under $\bold{\underline{P3}}$, in the process prescribed by the inductive algorithm.

\vskip.1in
 
$\bold{\underline{P1}}$: Under this possibility, the sequence represents
resolution of singularities of the general basic object
$({\Cal F}_0, (W_0,E_0))$ and we are done.

\vskip.1in

$\bold{\underline{P2}}$: Under this possibility, we can apply Corollary 4-9 to $({\Cal F}_k,
(W_k,E_k))$ and create a sequence of transformations representing resolution of singularities of
$({\Cal F}_k, (W_k,E_k))$.  Attaching this to the original sequence, we obtain a sequence
representing resolution of singularities of
$({\Cal F}_0, (W_0,E_0))$ and we are done.  

\vskip.1in

$\bold{\underline{P3}}$: So we may assume in the following that we are under possibility
$\bold{\underline{P3}}$.

\vskip.1in

$\bold{Case\ A}$: $R(1)(\underline{\roman{Max}}\ t_k) \neq \emptyset$. 

\vskip.1in

Let $p \in R(1)(\underline{\roman{Max}}\ t_k)$ be an arbitrary point.  

Then there exists
$W_k^{\lambda}$ such that 
$$p \in R(1)(\underline{\roman{Max}}\ t_k) \cap W_k^{\lambda} \subset
\underline{\roman{Max}}\ t_k \cap W_k^{\lambda} \subset \underline{\roman{Max}}\
w\text{-}\roman{ord}_k \cap W_k^{\lambda} = \underline{\roman{Max}}\
w\text{-}\roman{ord}_k^{\lambda}.$$  

Set $c_k^{\lambda} = b^{\lambda} \cdot \max\ w\text{-}\roman{ord}_k^{\lambda}$.  
 
Then since $\underline{\roman{Max}}\
w\text{-}\roman{ord}_k^{\lambda} = V(\Delta^{c_k^{\lambda}-1}(\overline{{\goth a}_k^{\lambda}}))
\subset
\widetilde{W_k^{\lambda}}$ and since
$\nu_p(\Delta^{c_k^{\lambda}-1}(\overline{{\goth a}_k^{\lambda}})) = 1$ (cf. Lemma 1-4), there
exists an open subset
$p
\in U_p \subset \widetilde{W_k^{\lambda}}$ and a regular parameter $f_p$ defined over $U_p$
such that $\underline{\roman{Max}}\
w\text{-}\roman{ord}_k^{\lambda} \cap U_p \subset \{f_p = 0\} \subset U_p$ where $\{f_p = 0\}$ is
a nonsingular closed subvariety of codimension one in $U_p$.

We see by shrinking $U_p$ if necessary that this implies
$$p \in R(1)(\underline{\roman{Max}}\ t_k) \cap U_p = \{f_p = 0\} =
\underline{\roman{Max}}\ w\text{-}\roman{ord}_k \cap U_p.$$

Since $p \in R(1)(\underline{\roman{Max}}\ t_k)$ is arbitrary, we conclude that
$R(1)(\underline{\roman{Max}}\ t_k)$ \linebreak
$(\subset \underline{\roman{Max}}\ t_k \subset
\underline{\roman{Max}}\ w\text{-}\roman{ord}_k \subset F_k)$ is smooth, open and closed in
$\underline{\roman{Max}}\ w\text{-}\roman{ord}_k$ (i.e., a union of smooth connected components of
$\underline{\roman{Max}}\ w\text{-}\roman{ord}_k$ disjoint from each other).

\vskip.1in

Since obviously $R(1)(\underline{\roman{Max}}\ t_k) \cap W_k^{\lambda} \subset
\widetilde{W_k^{\lambda}}$ and $R(1)(\underline{\roman{Max}}\ t_k) \cap W_k^{\lambda} \subset F_k
\cap W_k^{\lambda} = \roman{Sing}({\goth a}_k^{\lambda},b^{\lambda})$, we have only to show
that $R(1)(\underline{\roman{Max}}\ t_k) \cap W_k^{\lambda}$ is permissible with respect to
$\widetilde{E_k^{\lambda}} = E_k \cap \widetilde{W_k^{\lambda}}$.  

\vskip.1in

Observe first that
$$\max\ t_k = (\max\ w\text{-}\roman{ord}_k, 0) \text{\ or\ }(\max\ w\text{-}\roman{ord}_k, 1),$$ 
since if the second factor is $\geq 2$, then $\underline{\roman{Max}}\ t_k$ has codimension at
least two in $\widetilde{W_k^{\lambda}}$, which is against the case assumption of
$R(1)(\underline{\roman{Max}}\ t_k) \neq \emptyset$.

\vskip.1in

Suppose $\max\ t_k = (\max\ w\text{-}\roman{ord}_k, 1)$.  Then for any point $p \in
R(1)(\underline{\roman{Max}}\ t_k) \cap W_k^{\lambda}$ there exists an open neighborhood $U_p
\subset \widetilde{W_k^{\lambda}}$ and $p \in \widetilde{H_j} \in E_k^- \cap
\widetilde{W_k^{\lambda}} \subset \widetilde{E_k^{\lambda}}$ (See
Definition-Proposition 4-5 for the definition of
$E_k^-$.) such that
$$R(1)(\underline{\roman{Max}}\ t_k) \cap U_p = \widetilde{H_j} \cap U_p.$$
Since $\widetilde{H_j}$ is clearly permissible with respect to $\widetilde{E_k^{\lambda}}$,
we verify the permissibility in this case.

\vskip.1in

Suppose $\max\ t_k = (\max\ w\text{-}\roman{ord}_k, 0)$.  Then, by definition, no point \linebreak
$p
\in R(1)(\underline{\roman{Max}}\ t_k) \cap W_k^{\lambda}$ is contained in
$E_k^- \cap W_k^{\lambda}$.  (Hence $R(1)(\underline{\roman{Max}}\ t_k) \cap
W_k^{\lambda}$ is disjoint from
$E_k^- \cap \widetilde{W_k^{\lambda}} = \widetilde{E_k^{\lambda}}^{-}$.)  Thus we have only to
check
$R(1)(\underline{\roman{Max}}\ t_k)
\cap W_k^{\lambda}$ is permissible with respect to $E_k^+ \cap
\widetilde{W_k^{\lambda}} = \widetilde{E_k^{\lambda}}^{+}$.

Let $D \subset R(1)(\underline{\roman{Max}}\ t_k) \cap W_k^{\lambda}$ be an irreducible
(and hence connected) component.  Let $k_o$ be the index such that
$$\max\ w\text{-}\roman{ord}_{k_o - 1} > \max\ w\text{-}\roman{ord}_{k_o} = \cdot\cdot\cdot =
\max\ w\text{-}\roman{ord}_k$$ 
as described in Definition-Proposition 4-5 (iii).

If $D$ is one of the exceptional divisors for the morphism $\widetilde{W_{k_o}^{\lambda}}
\leftarrow \widetilde{W_k^{\lambda}}$ (That is to say, if the strict transform of $D$ is the
pull-back of the irreducible component of the center for some transformation.  Note that even in
the case where we take the center to be a divisor (a subvariety of codimension one) and hence
where the transformation is an isomorphism set-theoretically, we call $D$ one of the \it
exceptional \rm divisors.  See the convention of Definition 1-8 (iii).), then
$D
\subset
\widetilde{E_k^{\lambda}}$ and hence is permissible with respect to
$\widetilde{E_k^{\lambda}}^{+}$.  (In fact, $D \subset
\widetilde{E_k^{\lambda}}^{+}$ would imply that
$w\text{-}\roman{ord}_k(\eta_D) = 0$, where $\eta_D$ is the generic point of $D$, which is
against the case assumption $\max\ w\text{-}\roman{ord}_k > 0$ of
$\bold{\underline{P3}}$.  So this case does not happen.)

If $D$ is not any one of the exceptional divisors for the morphism $\widetilde{W_{k_o}^{\lambda}}
\leftarrow \widetilde{W_k^{\lambda}}$, then $D_i \subset \underline{\roman{Max}}\
w\text{-}\roman{ord}_i$, for $i = k_o, ... , k$, where $D_i$ is the image of $D$ on
$\widetilde{W_i^{\lambda}}$.  (See the above note for the meaning of the \it exceptional
\rm divisors.)  Therefore,
$D_i$ is a smooth connected component of
$\underline{\roman{Max}}\ w\text{-}\roman{ord}_i$, by the same argument at the beginning of the
proof for $\bold{Case\ A}$.  Therefore, any irreducible component of the center $Y_i \subset
\underline{\roman{Max}}\ w\text{-}\roman{ord}_i$ is either contained in $D_i$ or disjoint from
$D_i$, while $Y_i$ is permissible with respect to $\widetilde{E_i^{\lambda}}$.  Observe that
$D_{k_o}$ is clearly permissible with respect to $\widetilde{E_{k_o}^{\lambda}}^{+}$, since
$\widetilde{E_{k_o}^{\lambda}}^{+} =
\emptyset$.  The above analysis of the centers inductively implies that $D_i$ is permissible with respect to
$\widetilde{E_i^{\lambda}}^{+}$ for $i = k_o, ... , k$.

\vskip.1in

This completes the proof that $R(1)(\underline{\roman{Max}}\ t_k) \cap W_k^{\lambda}$ is permissible with respect to
$\widetilde{E_k^{\lambda}} = E_k \cap \widetilde{W_k^{\lambda}}$.

\vskip.1in

Now for the transformation
$$(F_k, (W_k,E_k)) \overset{\pi_{k+1}}\to{\leftarrow} (F_{k+1}, (W_{k+1},E_{k+1}))$$
with center $Y_k = R(1)(\underline{\roman{Max}}\ t_k)$, it is obvious, since $\max\ t_k \geq \max\
t_{k+1}$ (cf. Proposition 4-8), that we have the following four cases:

\vskip.1in

$\bold{A\text{-}1}$: $F_{k+1} = \emptyset$.

$\bold{A\text{-}2}$: $F_{k+1} \neq \emptyset$ and $\max\ w\text{-}\roman{ord}_{k+1}
= 0$.

$\bold{A\text{-}3}$: $F_{k+1} \neq \emptyset$, $\max\ w\text{-}\roman{ord}_{k+1} > 0$, and $\max\
t_k >
\max\ t_{k+1}$.

$\bold{A\text{-}4}$: $F_{k+1} \neq \emptyset$, $\max\ w\text{-}\roman{ord}_{k+1} > 0$, and $\max\
t_k =
\max\ t_{k+1}$.

\vskip.1in

The assertions for cases $\bold{A\text{-}1, A\text{-}2, A\text{-}3}$ are obvious.

\vskip.1in

So we have only to prove that in case $\bold{A\text{-}4}$ we have
$R(1)(\underline{\roman{Max}}\ t_{k+1}) = \emptyset.$  In that case, for each $\lambda$, the
corresponding transformation of basic objects
$$(\widetilde{W_k^{\lambda}}, ({\goth a}_k^{\lambda},b^{\lambda}), \widetilde{E_k^{\lambda}})
\overset{\pi_{k+1}^{\lambda}}\to{\leftarrow} (\widetilde{W_{k+1}^{\lambda}}, ({\goth
a}_{k+1}^{\lambda},b^{\lambda}), \widetilde{E_{k+1}^{\lambda}})$$
is an isomorphism between $\widetilde{W_k^{\lambda}}$ and $\widetilde{W_{k+1}^{\lambda}}$. 
Moreover, for each irreducible component $D \subset R(1)(\underline{\roman{Max}}\ t_k) \cap
W_k^{\lambda}$ we have
$w\text{-}\roman{ord}_{k+1}(\eta_D) = w\text{-}\roman{ord}_{k+1}^{\lambda}(\eta_D) = 0$, since $D
\subset
\widetilde{E_{k+1}^{\lambda}}$, where $\eta_D$ is the generic point of
$D$.  Since $\max\ w\text{-}\roman{ord}_{k+1} > 0$
and since $\max\ t_k = \max\ t_{k+1}$, we have
$\eta_D \not\in \underline{\roman{Max}}\ t_{k+1}$ and we also have \linebreak
$(\underline{\roman{Max}}\ t_k) \cap W_k^{\lambda} 
\supset
(\underline{\roman{Max}}\ t_{k+1}) \cap W_{k+1}^{\lambda}$ (cf. Proposition 4-8).  Therefore,
we finally conclude that $R(1)(\underline{\roman{Max}}\ t_{k+1}) = \emptyset$.  

\vskip.1in

This completes the proof of the assertions in $\bold{Case\ A}$ under possibility
$\bold{\underline{P3}}$.

\vskip.2in

The rest of the proof is dedicated to verifying the assertions in $\bold{Case\ B}$ under
possibility $\bold{\underline{P3}}$.

\vskip.1in

$\bold{Case\ B}$: $R(1)(\underline{\roman{Max}}\ t_k) = \emptyset$.

\vskip.1in

In order to construct a general basic object over $(G_k = \underline{\roman{Max}}\ t_k, (W_k,
E_k''))$ with the specified properties, we prove the following two lemmas, which,
given a sequence of transformations and smooth morphisms of basic objects, construct basic
objects whose singular loci coincide with
$\underline{\roman{Max}}\ w\text{-}\roman{ord}$ and
$\underline{\roman{Max}}\ t$, respectively.

\proclaim{Lemma 5-3} Let
$$\align
(W_0, (J_0,b), E_0) &\overset{\pi_1}\to{\leftarrow} (W_1, (J_1,b), E_1)
\overset{\pi_2}\to{\leftarrow} \cdot\cdot\cdot \\
(W_{i-1}, (J_{i-1},b), E_{i-1}) &\overset{\pi_i}\to{\leftarrow}
(W_i, (J_i,b), E_i) \\
\cdot\cdot\cdot \overset{\pi_{k-1}}\to{\leftarrow} (W_{k-1}, (J_{k-1},b),
E_{k-1}) &\overset{\pi_k}\to{\leftarrow} (W_k, (J_k,b), E_k) \\
\endalign$$
be a sequence of transformations and smooth morphisms of basic objects such that
$$\max\ w\text{-}\roman{ord}_k > 0.$$

Then there exists a simple basic object $(W_k' = W_k, (J_k',b'),E_k')$ whose singular
locus (as well as the singular loci of the basic objects in the sequences of
transformations and smooth morphisms starting from it) coincides with the locus
$\underline{\roman{Max}}\ w\text{-}\roman{ord}_k$ of $(W_k,(J_k,b),E_k)$ (as well as
the loci $\underline{\roman{Max}}\ w\text{-}\roman{ord}$ of the basic objects in the
extended sequences of transformations and smooth morphisms satisfying condition
$(\heartsuit)$) in the sense precisely formulated as follows:

\vskip.1in

$(\alpha)$ With each sequence of transformations and smooth morphisms of basic
objects starting from $(W_k', (J_k',b'),E_k')$
$$\align
(W_k', (J_k',b'), E_k') = (W_k', ((J_k')_0,b'), E_k') &\overset{\pi_{k+1}'}\to{\leftarrow}
(W_{k+1}', ((J_k')_1,b'), E_{k+1}') \overset{\pi_{k+2}'}\to{\leftarrow} \cdot\cdot\cdot \\
\cdot\cdot\cdot \overset{\pi_{N-1}'}\to{\leftarrow} (W_{N-1}', ((J_k')_{N-1},b'),
E_{k+N-1}') &\overset{\pi_{k+N}'}\to{\leftarrow} (W_{k+N}', ((J_k')_N,b'), E_{k+N}')\\
\endalign$$
satisfying the condition
$$\roman{Sing}((J_k')_j,b') \neq \emptyset \text{\ for\ }j = 0, ... , N-1,$$
there corresponds an extension of the original sequence of transformations and smooth
morphisms, satisfying condition $(\heartsuit)$ for $i = k + j + 1 \hskip.1in (j = 0, ...
, N-1)$,
$$\align
&(W_0, (J_0,b), E_0) \overset{\pi_1}\to{\leftarrow}
\cdot\cdot\cdot \overset{\pi_{k-1}}\to{\leftarrow} (W_{k-1}, (J_{k-1},b), E_{k-1})
\overset{\pi_k}\to{\leftarrow}
\\  
&(W_k, (J_k,b), E_k) \overset{\pi_{k+1}}\to{\leftarrow} (W_{k+1},
(J_{k+1},b), E_{k+1})
\overset{\pi_{k+2}}\to{\leftarrow} \cdot\cdot\cdot \\
&\cdot\cdot\cdot \overset{\pi_{k+N-1}}\to{\leftarrow} (W_{k+N-1}, (J_{k+N-1},b),
E_{k+N-1}) \overset{\pi_{k+N}}\to{\leftarrow} (W_{k+N}, (J_{k+N},b), E_{k+N})\\
\endalign$$
such that the following conditions are satisfied:

\vskip.1in

(i) $\pi_{k+j+1}'$ and $\pi_{k+j+1}$ are the transformations with the same
centers or the same smooth morphisms (as abstract varieties) for $j = 0, ... , N-1$
with
$W_{k+j+1}' = W_{k+j+1}$ (which means, in particular, that if
$\pi_{k+j+1}'$ is the transformation with center $Y_{k+j}' \subset W_{k+j}'$ which is
permissible for
$(W_{k+j}', ((J_k')_j,b'), E_{k+j}')$, then $\pi_{k+j+1}$ is the transformation with the same
center $Y_{k+j}' \subset W_{k+j}' =
W_{k+j}$ which is also permissible for $(W_{k+j}, (J_{k+j},b), E_{k+j})$), 

\vskip.1in

(ii) we have 

\vskip.1in

\noindent either
$$\left\{\aligned
&\max\ w\text{-}\roman{ord}_k = \max\ w\text{-}\roman{ord}_{k+1} = \cdot\cdot\cdot = \max\
w\text{-}\roman{ord}_{k+N},
\text{\ and\ }\\
&\roman{Sing}((J_k')_j,b') = \underline{\roman{Max}}\ w\text{-}\roman{ord}_{k+j} \text{\ for\ }j =
0, ... , N\\
\endaligned\right.$$
or
$$\left\{\aligned
&\max\ w\text{-}\roman{ord}_k = \max\ w\text{-}\roman{ord}_{k+1} = \cdot\cdot\cdot = \max\
w\text{-}\roman{ord}_{k+N-1} >
\max\ w\text{-}\roman{ord}_{k+N}\\
&(\text{or\ }\max\ w\text{-}\roman{ord}_k = \max\ w\text{-}\roman{ord}_{k+1} = \cdot\cdot\cdot = \max\
w\text{-}\roman{ord}_{k+N-1} \hskip.1in \& \hskip.1in \roman{Sing}(J_{k+N},b) = \emptyset), \text{\
and}\\ &\roman{Sing}((J_k')_j,b') = \underline{\roman{Max}}\ w\text{-}\roman{ord}_{k+j} \text{\
for\ }j = 0, ... , N-1 \hskip.1in \& \hskip.1in \roman{Sing}((J_k')_N,b') = \emptyset.\\
\endaligned\right.$$

$(\beta)$ Conversely, with each extension of the original sequence of transformations and
smooth morphisms
$$\align
&(W_0, (J_0,b), E_0) \overset{\pi_1}\to{\leftarrow}
\cdot\cdot\cdot \overset{\pi_{k-1}}\to{\leftarrow} (W_{k-1}, (J_{k-1},b), E_{k-1})
\overset{\pi_k}\to{\leftarrow}
\\  
&(W_k, (J_k,b), E_k) \overset{\pi_{k+1}}\to{\leftarrow} (W_{k+1},
(J_{k+1},b), E_{k+1})
\overset{\pi_{k+2}}\to{\leftarrow} \cdot\cdot\cdot \\
&\cdot\cdot\cdot \overset{\pi_{k+N-1}}\to{\leftarrow} (W_{k+N-1}, (J_{k+N-1},b),
E_{k+N-1}) \overset{\pi_{k+N}}\to{\leftarrow} (W_{k+N}, (J_{k+N},b), E_{k+N}),\\
\endalign$$
satisfying condition $(\heartsuit)$ for $i = k + j + 1 \hskip.1in (j = 0, ... , N-1)$ and
the condition
$$\max\ w\text{-}\roman{ord}_k = \max\ w\text{-}\roman{ord}_{k+1} = \cdot\cdot\cdot =
\max\ w\text{-}\roman{ord}_{k+N-1},$$
there corresponds a sequence of transformations and smooth morphisms of basic
objects starting from $(W_k', (J_k',b'),E_k')$
$$\align
(W_k', (J_k',b'), E_k') = (W_k', ((J_k')_0,b'), E_k') &\overset{\pi_{k+1}'}\to{\leftarrow}
(W_{k+1}', ((J_k')_1,b'), E_{k+1}') \overset{\pi_{k+2}'}\to{\leftarrow} \cdot\cdot\cdot \\
\cdot\cdot\cdot \overset{\pi_{N-1}'}\to{\leftarrow} (W_{N-1}', ((J_k')_{N-1},b'),
E_{k+N-1}') &\overset{\pi_{k+N}'}\to{\leftarrow} (W_{k+N}', ((J_k')_N,b'), E_{k+N}')\\
\endalign$$
satisfying the condition
$$\roman{Sing}((J_k')_j,b') \neq \emptyset \text{\ for\ }j = 0, ... , N-1$$
and conditions (i) and (ii) as in $(\alpha)$.
\endproclaim

\demo{Proof}\enddemo Recall the characterization of the ideal $\overline{J_k}$ (cf. Definition
1-10 (ii))
$$J_k = I(H_{r+1})^{a_{r+1}} \cdot\cdot\cdot I(H_{r+k})^{a_{r+k}} \cdot \overline{J_k}.$$
We set
$$b_k = b \cdot \max\ w\text{-}\roman{ord}_k.$$
We define the basic object $(W_k', (J_k',b'), E_k')$ in the following way:

$$\left\{\aligned
W_k' &= W_k\\
J_k' &= \left\{\aligned
&\overline{J_k} \hskip2.5in \text{\ if\ }b_k \geq b \\
&\overline{J_k}^{b - b_k} + \{I(H_{r+1})^{a_{r+1}} \cdot\cdot\cdot
I(H_{r+k})^{a_{r+k}}\}^{b_k}
\ \ \text{\ if\ }b_k < b\\
\endaligned\right.\\
b' &= \left\{\aligned
&b_k \hskip2.52in \text{\ if\ }b_k \geq b\\
&b_k(b - b_k) \hskip2.05in \text{\ if\ }b_k < b\\
\endaligned\right.\\
E_k' &= E_k.
\endaligned\right.$$

\vskip.1in

We check that the basic object $(W_k', (J_k',b'), E_k')$ has the required properties.

\newpage

Property $(\alpha)$ 

\vskip.1in

$\boxed{\roman{Case}: b_k \geq b}$

\vskip.1in

In this case, we observe
$$\align
\xi_k \in \roman{Sing}(J_k',b') &\Longleftrightarrow \xi_k \in
\roman{Sing}(\overline{J_k},b_k)\\
&\Longleftrightarrow \nu_{\xi_k}(\overline{J_k}) \geq b_k \\
&\Longleftrightarrow \nu_{\xi_k}(\overline{J_k}) \geq b_k \hskip.1in \& \hskip.1in
\nu_{\xi_k}(J_k) \geq b\\ 
&\Longleftrightarrow \nu_{\xi_k}(\overline{J_k}) = b_k \hskip.1in \& \hskip.1in
\nu_{\xi_k}(J_k) \geq b\\ 
&\Longleftrightarrow \xi_k \in \underline{\roman{Max}}\ w\text{-}\roman{ord}_k,\\
\endalign$$
which implies
$$\roman{Sing}(J_k',b') = \underline{\roman{Max}}\ w\text{-}\roman{ord}_k.$$
Moreover, the equivalence of the conditions above also shows
$$\xi_k \in \roman{Sing}(J_k',b') \Longrightarrow \nu_{\xi_k}(J_k') =
\nu_{\xi_k}(\overline{J_k}) = b_k = b',$$
verifying that $(W_k', (J_k',b'), E_k')$ is a simple basic object. 

Inductively, we can construct an extension of the original sequence (of the
transformations with the same centers and the same smooth morphisms) such that for
$j = 0, ... , N-1$, since
$\roman{Sing}((J_k')_j,b') \neq \emptyset$, we have
$$\left\{\aligned
&(J_k')_j = \overline{J_{k+j}}\\
&\max\ w\text{-}\roman{ord}_k = \cdot\cdot\cdot = \max\
w\text{-}\roman{ord}_{k+j} \text{\ i.e.,\ }b_k = \cdot\cdot\cdot = b_{k+j} \geq b\\
&\text{and\ hence\ }(J_k')_j = J_{k+j}'.\\
\endaligned\right.$$
Therefore, by the argument at the beginning applied to $(W_{k+j}', ((J_k')_j,b'), E_{k+j}')
= (W_{k+j}', (J_{k+j}',b'), E_{k+j}')$, we conclude
$$\roman{Sing}((J_k')_j,b') = \roman{Sing}(J_{k+j}',b') = \underline{\roman{Max}}\
w\text{-}\roman{ord}_{k+j}
\text{\ for\ }j = 0, ... , N-1,$$
which also implies condition (i) as $E_{k+j}' = E_{k+j}$.

\vskip.1in

In the case $\roman{Sing}((J_k')_N,b') \neq \emptyset$, we have
$$\left\{\aligned
&(J_k')_N = \overline{J_{k+N}}\\
&\max\ w\text{-}\roman{ord}_k = \cdot\cdot\cdot = \max\
w\text{-}\roman{ord}_{k+N} \text{\ i.e.,\ }b_k = \cdot\cdot\cdot = b_{k+N} \geq b\\
&\text{and\ hence\ }(J_k')_N = J_{k+N}',\\
\endaligned\right.$$
which implies
$$\roman{Sing}((J_k')_N,b') = \underline{\roman{Max}}\ w\text{-}\roman{ord}_{k+N}.$$
Thus we are in the former case stated in condition (ii). 

\vskip.1in

In the case $\roman{Sing}((J_k')_N,b') = \emptyset$, we have
$$\left\{\aligned
&(J_k')_N = \overline{J_{k+N}}\\
&\max\ w\text{-}\roman{ord}_k = \cdot\cdot\cdot = \max\
w\text{-}\roman{ord}_{k+N-1} > \max\ w\text{-}\roman{ord}_{k+N} \text{\ i.e.,\ }b_k =
\cdot\cdot\cdot = b_{k+N-1} > b_{k+N}\\
&(\text{or\ }\max\ w\text{-}\roman{ord}_k = \cdot\cdot\cdot = \max\
w\text{-}\roman{ord}_{k+N-1} \text{\ i.e.,\ }b_k =
\cdot\cdot\cdot = b_{k+N-1} \hskip.1in \& \hskip.1in \roman{Sing}(J_{k+N},b) = \emptyset).\\
\endaligned\right.$$
Thus we are in the latter case stated in condition (ii).

\vskip.1in

$\boxed{\roman{Case}: b_k < b}$

\vskip.1in

In this case, we observe (cf. Remark 1-2 (ii))
$$\align
\xi_k &\in \roman{Sing}(J_k',b') \\
&\Longleftrightarrow \xi_k \in
\roman{Sing}(\overline{J_k}^{b-b_k} + \{I(H_{r+1})^{a_{r+1}} \cdot\cdot\cdot
I(H_{r+k})^{a_{r+k}}\}^{b_k},b_k(b-b_k))\\ 
&\Longleftrightarrow \nu_{\xi_k}(\overline{J_k}^{b-b_k}) \geq b_k(b-b_k) \hskip.1in \&
\hskip.1in
\nu_{\xi_k}(\{I(H_{r+1})^{a_{r+1}} \cdot\cdot\cdot I(H_{r+k})^{a_{r+k}}\}^{b_k}) \geq
b_k(b-b_k)\\  
&\Longleftrightarrow \nu_{\xi_k}(\overline{J_k}) \geq b_k \hskip.1in \& \hskip.1in
\nu_{\xi_k}(I(H_{r+1})^{a_{r+1}} \cdot\cdot\cdot I(H_{r+k})^{a_{r+k}}) \geq
b-b_k\\  
&\Longleftrightarrow \nu_{\xi_k}(I(H_{r+1})^{a_{r+1}} \cdot\cdot\cdot
I(H_{r+k})^{a_{r+k}} \cdot \overline{J_k}) \geq b_k + (b - b_k) = b, \\
&\hskip.34in \nu_{\xi_k}(\overline{J_k}) = b_k \hskip.1in \& \hskip.1in
\nu_{\xi_k}(I(H_{r+1})^{a_{r+1}} \cdot\cdot\cdot I(H_{r+k})^{a_{r+k}}) \geq b - b_k\\
&\Longleftrightarrow \xi_k \in \underline{\roman{Max}}\ w\text{-}\roman{ord}_k,\\
\endalign$$
which implies
$$\roman{Sing}(J_k',b') = \underline{\roman{Max}}\ w\text{-}\roman{ord}_k.$$
Moreover, the equivalence of the conditions above also shows
$$\xi_k \in \roman{Sing}(J_k',b') \Longrightarrow \nu_{\xi_k}(J_k') =
\nu_{\xi_k}(\overline{J_k}^{b-b_k}) = b_k(b - b_k) = b',$$
verifying that $(W_k', (J_k',b'), E_k')$ is a simple basic object. 

Inductively, we can construct an extension of the original sequence (of the
transformations with the same centers and the same smooth morphisms) such that, since
$\roman{Sing}((J_k')_j,b') \neq \emptyset$, for
$j = 0, ... , N-1$
$$\left\{\aligned
&(J_k')_j = \overline{J_{k+j}}^{b-b_{k+j}} + \{I(H_{r+1})^{a_{r+1}} \cdot\cdot\cdot
I(H_{r+k+j})^{a_{r+k+j}}\}^{b_{k+j}}\\ 
&\max\ w\text{-}\roman{ord}_k = \cdot\cdot\cdot = \max\
w\text{-}\roman{ord}_{k+j} \text{\ i.e.,\ }b_k = \cdot\cdot\cdot = b_{k+j} < b\\
&\text{and\ hence\ }(J_k')_j = J_{k+j}'.\\
\endaligned\right.$$

Therefore, by the argument at the beginning applied to $(W_{k+j}', ((J_k')_j,b'), E_{k+j}')
= (W_{k+j}', (J_{k+j}',b'), E_{k+j}')$, we conclude
$$\roman{Sing}((J_k')_j,b') = \roman{Sing}(J_{k+j}',b') = \underline{\roman{Max}}\
w\text{-}\roman{ord}_{k+j}
\text{\ for\ }j = 0, ... , N-1,$$
which also implies condition (i) as $E_{k+j}' = E_{k+j}$.

\vskip.1in

In the case $\roman{Sing}((J_k')_N,b') \neq \emptyset$, we have
$$\left\{\aligned
&(J_k')_N = \overline{J_{k+N}}^{b-b_{k+N}} + \{I(H_{r+1})^{a_{r+1}} \cdot\cdot\cdot
I(H_{r+k+N})^{a_{r+k+N}}\}^{b_{k+N}}\\
&\max\ w\text{-}\roman{ord}_k = \cdot\cdot\cdot = \max\
w\text{-}\roman{ord}_{k+N} \text{\ i.e.,\ }b_k = \cdot\cdot\cdot = b_{k+N} < b\\
&\text{and\ hence\ }(J_k')_N = J_{k+N}',\\
\endaligned\right.$$
which implies
$$\roman{Sing}((J_k')_N,b') = \underline{\roman{Max}}\ w\text{-}\roman{ord}_{k+N}.$$
Thus we are in the former case stated in condition (ii). 

\vskip.1in

In the case $\roman{Sing}((J_k')_N,b') = \emptyset$, we have
$$\left\{\aligned
&(J_k')_N = \overline{J_{k+N}}^{b-b_{k+N-1}} + \{I(H_{r+1})^{a_{r+1}} \cdot\cdot\cdot
I(H_{r+k+N})^{a_{r+k+N}}\}^{b_{k+N-1}}\\
&\max\ w\text{-}\roman{ord}_k = \cdot\cdot\cdot = \max\
w\text{-}\roman{ord}_{k+N-1} > \max\ w\text{-}\roman{ord}_{k+N} \text{\ i.e.,\ }b_k =
\cdot\cdot\cdot = b_{k+N-1} > b_{k+N}\\
&(\text{or\ }\max\ w\text{-}\roman{ord}_k = \cdot\cdot\cdot = \max\
w\text{-}\roman{ord}_{k+N-1} \text{\ i.e.,\ }b_k =
\cdot\cdot\cdot = b_{k+N-1} \hskip.1in \& \hskip.1in \roman{Sing}(J_{k+N},b) = \emptyset).\\
\endaligned\right.$$
Thus we are in the latter case stated in condition (ii).

\vskip.1in

Property $(\beta)$

\vskip.1in 

The argument for verification of property $(\beta)$ is identical to the one for
verification of property $(\alpha)$, and is left to the reader as an exercise.

\vskip.1in

This completes the proof for Lemma 5-3.

\vskip.2in

\proclaim{Lemma 5-4} Let
$$\align
(W_0, (J_0,b), E_0) &\overset{\pi_1}\to{\leftarrow} (W_1, (J_1,b), E_1)
\overset{\pi_2}\to{\leftarrow} \cdot\cdot\cdot \\
(W_{i-1}, (J_{i-1},b), E_{i-1}) &\overset{\pi_i}\to{\leftarrow}
(W_i, (J_i,b), E_i) \\
\cdot\cdot\cdot \overset{\pi_{k-1}}\to{\leftarrow} (W_{k-1}, (J_{k-1},b),
E_{k-1}) &\overset{\pi_k}\to{\leftarrow} (W_k, (J_k,b), E_k) \\
\endalign$$
be a sequence of transformations and smooth morphisms of basic objects with condition
$(\heartsuit)$
$$(\heartsuit) \hskip.1in \left\{\aligned &Y_{i-1} \subset
\underline{\roman{Max}}\ w\text{-}\roman{ord}_{i-1} (\subset
\roman{Sing}(J_{i-1},b)) \\
&\text{whenever\ } \pi_i \text{\ is\ a\ transformation\ with\ center\
}Y_{i-1}\\
\endaligned\right\} \text{\ for\ }i = 1, ... , k$$
such that
$$\max\ w\text{-}\roman{ord}_k > 0.$$

Then there exists a simple basic object $(W_k'' = W_k, (J_k'',b''),E_k'')$ whose
singular locus (as well as the singular loci of the basic objects in the sequences of
transformations and smooth morphisms starting from it) coincides with the locus
$\underline{\roman{Max}}\ t_k$ of $(W_k,(J_k,b),E_k)$ (as well as the loci
$\underline{\roman{Max}}\ t$ of the basic objects in the extended sequences of
transformations and smooth morphisms satisfying condition $(\heartsuit')$) in the sense
precisely formulated as follows:

\vskip.1in

$(\alpha)$ With each
sequence of transformations and smooth morphisms of basic objects starting with $(W_k'',
(J_k'',b''),E_k'')$
$$\align
(W_k'', (J_k'',b''), E_k'') = (W_k'', ((J_k'')_0,b''), E_k'')
&\overset{\pi_{k+1}''}\to{\leftarrow} (W_{k+1}'', ((J_k'')_1,b'), E_{k+1}'')
\overset{\pi_{k+2}''}\to{\leftarrow}
\cdot\cdot\cdot \\
\cdot\cdot\cdot \overset{\pi_{k+N-1}''}\to{\leftarrow} (W_{k+N-1}'', ((J_k'')_{N-1},b''),
E_{k+N-1}'') &\overset{\pi_{k+N}''}\to{\leftarrow} (W_{k+N}'', ((J_k'')_N,b''), E_{k+N}'')\\
\endalign$$
with the condition
$$\roman{Sing}((J_k'')_j,b'') \neq \emptyset \text{\ for\ }j = 0, ... , N-1,$$
there corresponds an extension of the original sequence of transformations and smooth
morphisms
$$\align
&(W_0, (J_0,b), E_0) \overset{\pi_1}\to{\leftarrow}
\cdot\cdot\cdot \overset{\pi_{k-1}}\to{\leftarrow} (W_{k-1}, (J_{k-1},b), E_{k-1})
\overset{\pi_k}\to{\leftarrow}\\ 
&(W_k, (J_k,b), E_k) \overset{\pi_{k+1}}\to{\leftarrow} (W_{k+1},
(J_{k+1},b), E_{k+1})
\overset{\pi_{k+2}}\to{\leftarrow} \cdot\cdot\cdot \\
&\cdot\cdot\cdot \overset{\pi_{k+N-1}}\to{\leftarrow} (W_{k+N-1}, (J_{k+N-1},b), E_{k+N-1})
\overset{\pi_{k+N}}\to{\leftarrow} (W_{k+N}, (J_{k+N},b), E_{k+N})\\
\endalign$$
with condition
$$(\heartsuit') \hskip.1in \left\{\aligned &Y_{i-1} \subset \underline{\roman{Max}}\
t_{i-1} \subset
\underline{\roman{Max}}\ w\text{-}\roman{ord}_{i-1} (\subset
\roman{Sing}(J_{i-1},b)) \\
&\text{whenever\ } \pi_i \text{\ is\ a\ transformation\ with\ center\
}Y_{i-1}\\
\endaligned\right\} \text{\ for\ }i = k+1 , ... , N$$
satisfying the following conditions:

\vskip.1in

(i) $\pi_{k+j+1}''$ and $\pi_{k+j+1}$ are the transformations with the same
centers or the same smooth morphisms (as abstract varieties) for $j = 0, ... , N-1$
with
$W_{k+j+1}'' = W_{k+j+1}$ (which means, in particular, if
$\pi_{k+j+1}''$ is the transformation with center $Y_{k+j}'' \subset W_{k+j}''$ which is
permissible for
$(W_{k+j}'', ((J_k'')_j,b''), E_{k+j}'')$, then $\pi_{k+j+1}$ is the
transformation with the same center $Y_{k+j}''
\subset W_{k+j}'' = W_{k+j}$ which is also permissible for $(W_{k+j}, (J_{k+j},b),
E_{k+j})$), 

\vskip.1in

(ii) we have 

\vskip.1in

\noindent either
$$\left\{\aligned
&\max\ t_k = \max\ t_{k+1} = \cdot\cdot\cdot = \max\ t_{k+N},
\text{\ and\ }\\
&\roman{Sing}((J_k'')_j,b'') = \underline{\roman{Max}}\ t_{k+j} \text{\ for\ }j = 0, ... , N\\
\endaligned\right.$$
or
$$\left\{\aligned
&\max\ t_k = \max\ t_{k+1} = \cdot\cdot\cdot = \max\ t_{k+N-1} >
\max\ t_{k+N}\\
&(\text{or\ }\max\ t_k = \max\ t_{k+1} = \cdot\cdot\cdot = \max\ t_{k+N-1}
\hskip.1in \& \hskip.1in \roman{Sing}(J_{k+N},b) =
\emptyset), \text{\ and}\\ &\roman{Sing}((J_k'')_j,b'') = \underline{\roman{Max}}\ t_{k+j} \text{\
for\ }j = 0, ... , N-1
\hskip.1in \& \hskip.1in \roman{Sing}((J_k'')_N,b'') = \emptyset.\\
\endaligned\right.$$

$(\beta)$ Conversely, with each extension of the original sequence of transformations
and smooth morphisms
$$\align
&(W_0, (J_0,b), E_0) \overset{\pi_1}\to{\leftarrow}
\cdot\cdot\cdot \overset{\pi_{k-1}}\to{\leftarrow} (W_{k-1}, (J_{k-1},b), E_{k-1})
\overset{\pi_k}\to{\leftarrow}\\ 
&(W_k, (J_k,b), E_k) \overset{\pi_{k+1}}\to{\leftarrow} (W_{k+1},
(J_{k+1},b), E_{k+1})
\overset{\pi_{k+2}}\to{\leftarrow} \cdot\cdot\cdot \\
&\cdot\cdot\cdot \overset{\pi_{k+N-1}}\to{\leftarrow} (W_{k+N-1}, (J_{k+N-1},b), E_{k+N-1})
\overset{\pi_{k+N}}\to{\leftarrow} (W_{k+N}, (J_{k+N},b), E_{k+N})\\
\endalign$$
with condition
$$(\heartsuit') \hskip.1in \left\{\aligned &Y_{i-1} \subset \underline{\roman{Max}}\
t_{i-1} \subset
\underline{\roman{Max}}\ w\text{-}\roman{ord}_{i-1} (\subset
\roman{Sing}(J_{i-1},b)) \\
&\text{whenever\ } \pi_i \text{\ is\ a\ transformation\ with\ center\
}Y_{i-1}\\
\endaligned\right\} \text{\ for\ }i = k+1 , ... , N$$
and the condition
$$\max\ t_k = \max\ t_{k+1} = \cdot\cdot\cdot = \max\ t_{k+N-1},$$
there corresponds a sequence of transformations and smooth morphisms of basic objects
starting from $(W_k'', (J_k'',b''),E_k'')$
$$\align
(W_k'', (J_k'',b''), E_k'') = (W_k'', ((J_k'')_0,b''), E_k'')
&\overset{\pi_{k+1}''}\to{\leftarrow} (W_{k+1}'', ((J_k'')_1,b'), E_{k+1}'')
\overset{\pi_{k+2}''}\to{\leftarrow}
\cdot\cdot\cdot \\
\cdot\cdot\cdot \overset{\pi_{k+N-1}''}\to{\leftarrow} (W_{k+N-1}'', ((J_k'')_{N-1},b''),
E_{k+N-1}'') &\overset{\pi_{k+N}''}\to{\leftarrow} (W_{k+N}'', ((J_k'')_N,b''), E_{k+N}'')\\
\endalign$$
satisfying the condition
$$\roman{Sing}((J_k'')_j,b'') \neq \emptyset \text{\ for\ }j = 0, ... , N-1$$
and conditions (i) and (ii) as in $(\alpha)$.

\vskip.2in

Moreover, the basic object $(W_k'' = W_k, (J_k'',b''),E_k'')$ has an open covering
$\{(W_k'')^{\gamma}\}_{\gamma \in \Gamma}$ satisfying conditions 1 and 2 of the key
inductive lemma (Lemma 3-1): for each $\gamma \in \Gamma$, there exists a smooth hypersurface $(W_k'')_h^{\gamma}
\subset (W_k'')^{\gamma}$, embedded as a closed subscheme, such that

1. $I((W_k'')_h^{\gamma}) \subset \Delta^{b''-1}(J_k'')|_{(W_k'')^{\gamma}} \hskip.1in
(\text{and\ hence\ }(W_k'')_h^{\gamma} \supset \roman{Sing}(J_k'',b'') \cap
(W_k'')^{\gamma})$, and

2. $(W_k'')_h^{\gamma}$ is permissible with respect to $E_k''
\cap (W_k'')^{\gamma}$, and $(W_k'')_h^{\gamma}$ is not contained in $E_k''$,
i.e., $(W_k'')_h^{\gamma}
\not\subset E_k''$.
\endproclaim

\demo{Proof}\enddemo Let $(W_k', (J_k',b'), E_k')$ be the basic object we constructed as
in Lemma 5.3.  

Let
$$\max\ t_k = (\max\ w\text{-}\roman{ord}_k,n).$$

We define the basic object $(W_k'', (J_k'',b'), E_k'')$ in the following
way:
$$\left\{\aligned
W_k'' &= W_k\\
J_k'' &= J_k' + \prod_{\{H_{s_1}, ... , H_{s_n}\} \subset E_k^- \text{\ with\ } H_{s_1},
... , H_{s_n} \text{\ distinct}}\{I(H_{s_1}) +
\cdot\cdot\cdot + I(H_{s_n})\}^{b'}
\\ b'' &= b'\\
E_k'' &= E_k^+.\\
\endaligned\right.$$
Recall that $E_k^- = \{H_1, ... , H_r, ... , H_{r+k_o}\}$ as a subset of \linebreak
$E_k =
\{H_1,
... , H_r, ... , H_{r+k_o}, H_{r+k_o+1}... , H_{r+k}\}$ and that $E_k^+$ is the
complement of $E_k^-$ in $E_k$, where $k_o$ is the index (See Definition 1-10 (iii).) so that
$$\max\ w\text{-}\roman{ord}_{k_o-1} > \max\ w\text{-}\roman{ord}_{k_o} = \cdot\cdot\cdot = \max\
w\text{-}\roman{ord}_k.$$  
Note that the order of the ideal $I(H_{s_1}) +
\cdot\cdot\cdot + I(H_{s_n})$ is the characteristic function of the set $H_{s_1} \cap
\cdot\cdot\cdot \cap H_{s_n}$, i.e.,
$$\nu_{\xi_k}(I(H_{s_1}) +
\cdot\cdot\cdot + I(H_{s_n})) = \left\{\aligned
0 &\hskip.1in \text{\ if\ }\hskip.1in \xi_k \not\in H_{s_1} \cap \cdot\cdot\cdot \cap H_{s_n}\\
1 &\hskip.1in \text{\ if\ }\hskip.1in \xi_k \in H_{s_1} \cap \cdot\cdot\cdot \cap H_{s_n}.\\
\endaligned\right.$$  
Note also that over a point $\xi_k \in \underline{\roman{Max}}\ w\text{-}\roman{ord}_k$ NO two
distinct intersections \linebreak
$H_{s_1} \cap \cdot\cdot\cdot \cap H_{s_n}$ (with $H_{s_1}, ... , H_{s_n}$
distinct) and
$H_{s_1'} \cap
\cdot\cdot\cdot \cap H_{s_n'}$ (with $H_{s_1'}, ... , H_{s_n'}$
distinct) can meet, since $n$ is the maximum of the number of such
divisors in $E_k^-$ that intersect at a point in $\underline{\roman{Max}}\
w\text{-}\roman{ord}_k$.  

\vskip.1in

We check that the basic object $(W_k'', (J_k'',b'), E_k'')$ has the required properties.

\newpage

Property $(\alpha)$

\vskip.1in

We observe (cf. Remark 1-2 (ii))
$$\align
&\xi_k \in \roman{Sing}(J_k'',b'') \\
&\Longleftrightarrow \xi_k \in
\roman{Sing}(J_k' + \prod_{\{H_{s_1}, \cdot\cdot\cdot, H_{s_n}\} \subset E_k^- \text{\ with\ } H_{s_1},
... , H_{s_n} \text{\ distinct}}\{I(H_{s_1}) +
\cdot\cdot\cdot + I(H_{s_n})\}^{b'},b')\\
&\Longleftrightarrow \nu_{\xi_k}(J_k') \geq b' \hskip.1in \&
\hskip.1in \nu_{\xi_k}(\prod_{\{H_{s_1},
\cdot\cdot\cdot, H_{s_n}\} \subset E_k^- \text{\ with\ } H_{s_1},
... , H_{s_n} \text{\ distinct}}\{I(H_{s_1}) +
\cdot\cdot\cdot + I(H_{s_n})\}^{b'}) \geq b'\\ 
&\Longleftrightarrow \nu_{\xi_k}(J_k') = b' \\
&\hskip.35in \& \hskip.1in \xi_k \in H_{s_1}
\cap \cdot\cdot\cdot \cap H_{s_n} \text{\ for\ some\ }\{H_{s_1}, \cdot\cdot\cdot,
H_{s_n}\} \subset E_k^- \text{\ with\ } H_{s_1},
... , H_{s_n} \text{\ distinct}\\  
&\Longleftrightarrow
\xi_k
\in
\underline{\roman{Max}}\ t_k,\\
\endalign$$
which implies
$$\roman{Sing}(J_k'',b'') = \underline{\roman{Max}}\ t_k.$$
Moreover, the equivalence of the conditions above also implies
$$\xi_k \in \roman{Sing}(J_k'',b'') \Longrightarrow \nu_{\xi_k}(J_k'') = \nu_{\xi_k}(J_k') = b'
= b'',$$ 
verifying that $(W_k'', (J_k'',b''), E_k'')$ is a simple basic object. 

Inductively, we can construct an extension of the original sequence (of the
transformations with the same centers and the same smooth morphisms) such that for
$j = 0,
\cdot\cdot\cdot, N-1$, since
$\roman{Sing}((J_k'')_j,b'') \neq \emptyset$, we have
$$\left\{\aligned
&(J_k'')_j = J_{k+j}' + \prod_{\{H_{s_1}, \cdot\cdot\cdot, H_{s_n}\} \subset
E_{k+j}^- \text{\ with\ } H_{s_1},
... , H_{s_n} \text{\ distinct}}\{I(H_{s_1}) +
\cdot\cdot\cdot + I(H_{s_n})\}^{b'}\\
&\max\ t_k = \cdot\cdot\cdot = \max\ t_{k+j}\\ 
&\text{and\ hence\ }(J_k'')_j = J_{k+j}''.\\
\endaligned\right.$$
Therefore, by the argument at the beginning applied to $(W_{k+j}'', ((J_k'')_j,b''),
E_{k+j}'') = (W_{k+j}'', (J_{k+j}'',b''), E_{k+j}'')$, we conclude
$$\roman{Sing}((J_k'')_j,b'') = \roman{Sing}(J_{k+j}'',b'') = \underline{\roman{Max}}\
t_{k+j}
\text{\ for\ }j = 0, \cdot\cdot\cdot, N-1.$$

Suppose that $\pi_{k+j+1}''$ is the transformation with the center \linebreak
$Y_{k+j}'' \subset
W_{k+j}'' = W_{k+j}$ permissible for
$(W_{k+j}'', ((J_k'')_j = J_{k+j}'',b''), E_{k+j}'')$ where $j = 0, \cdot\cdot\cdot,
N-1$.  Then remark that
$Y_{k+j}''$ is contained in $H_s \in E_{k+j}^-$ if $Y_{k+j}'' \cap H_s \neq
\emptyset$, since
$Y_{k+j}'' \subset
\roman{Sing}((J_k'')_j,b'') = \underline{\roman{Max}}\ t_{k+j}$, and that $Y_{k+j}''$
is permissible with respect to $E_{k+j}'' = E_{k+j}^+$ by definition.  This implies
that the center $Y_{k+j}'' \subset W_{k+j}'' = W_{k+j}$ is permissible for
$(W_{k+j}, (J_{k+j},b), E_{k+j})$, verifying condition (i). 

\vskip.1in

In the case $\roman{Sing}((J_k'')_N,b'') \neq \emptyset$, we have
$$\left\{\aligned
&(J_k'')_N = J_{k+N}' + \prod_{\{H_{s_1}, \cdot\cdot\cdot, H_{s_n}\} \subset
E_{k+N}^- \text{\ with\ } H_{s_1},
... , H_{s_n} \text{\ distinct}}\{I(H_{s_1}) +
\cdot\cdot\cdot + I(H_{s_n})\}^{b'}\\
&\max\ t_k = \cdot\cdot\cdot = \max\
t_{k+N} \\
&\text{and\ hence\ }(J_k'')_N = J_{k+N}'',\\
\endaligned\right.$$
which implies
$$\roman{Sing}((J_k'')_N,b'') = \underline{\roman{Max}}\ t_{k+N}.$$
Thus we are in the former case stated in condition (ii). 

In the case $\roman{Sing}((J_k'')_N,b'') = \emptyset$, we have
$$\left\{\aligned
&(J_k'')_N = (J_k')_N + \prod_{\{H_{s_1}, \cdot\cdot\cdot, H_{s_n}\} \subset
E_{k+N-1}^- \subset E_{k+N} \text{\ with\ } H_{s_1},
... , H_{s_n} \text{\ distinct}}\{I(H_{s_1}) +
\cdot\cdot\cdot + I(H_{s_n})\}^{b'}\\
&\max\ t_k = \cdot\cdot\cdot = \max\ t_{k+N-1} > \max\ t_{k+N}\\
&(\text{or\ } \max\ t_k = \cdot\cdot\cdot = \max\ t_{k+N-1}
\hskip.1in \& \hskip.1in \roman{Sing}(J_{k+N},b) =
\emptyset).\\
\endaligned\right.$$
Thus we are in the latter case stated in condition (ii).  (Note in the last case
of
\linebreak
$\roman{Sing}((J_k'')_N,b'') = \emptyset$ that when the morphism $\pi_{k+N}'' = \pi_{k+N}$ is
a transformation the above notation $E_{k+N-1}^- \subset E_{k+N}$ denotes the set
of the strict transforms in $E_{k+N}$ of the divisors in $E_{k+N-1}^-$ and that when 
the morphism $\pi_{k+N}'' = \pi_{k+N}$ is a
smooth morphism the above notation $E_{k+N-1}^- \subset E_{k+N}$ denotes the set
of the inverse images in $E_{k+N}$ of the divisors in $E_{k+N-1}^-$.)  

\vskip.1in

Property $(\beta)$

\vskip.1in 

The argument for verification of property $(\beta)$ is identical to the one for
verification of property $(\alpha)$, and is left to the reader as an exercise.

\vskip.2in

Now we show that $(W_k'', (J_k'',b''), E_k'')$ has an open covering
$\{(W_k'')^{\gamma}\}_{\gamma \in \Gamma}$ with smooth hypersurfaces
$(W_k'')_h^{\gamma} \subset (W_k'')^{\gamma}$, embedded as closed
subschemes, satisfying conditions 1 and 2 of the key inductive lemma (cf. Lemma
3-1).

\vskip.1in

Let $k_o$ be the index as before (See Definition 1-10 (iii).) so that
$$\max\ w\text{-}\roman{ord}_{k_o-1} > \max\ w\text{-}\roman{ord}_{k_o} = \cdot\cdot\cdot = \max\
w\text{-}\roman{ord}_k.$$

Looking at the basic object $(W_{k_o}, (J_{k_o},b), E_{k_o})$, we consider the
basic object
\linebreak
$(W_{k_o}' = W_{k_o}, (J_{k_o}',b'), E_{k_o}' = E_{k_o})$ as constructed in Lemma 5-3.

First remark that by Lemma 5-3 the part after the index $k_o$ of the original
sequence of transformations and smooth morphisms of basic objects (satisfying
condition $(\heartsuit)$)
$$(W_{k_o}, (J_{k_o},b), E_{k_o}) \leftarrow \cdot\cdot\cdot \leftarrow (W_k, (J_k,b), E_k)$$
gives rise to a sequence of transformations with the same centers and the same smooth
morphisms (as abstract varieties) of basic objects starting from \linebreak
$(W_{k_o}',(J_{k_o}',b'), E_{k_o}')$
$$(W_{k_o}' = W_{k_o}, (J_{k_o}',b'), E_{k_o}' = E_{k_o}) \leftarrow \cdot\cdot\cdot
\leftarrow (W_k' = W_k, ((J_{k_o}')_{k-k_o},b'), E_k' = E_k).$$

Secondly, since $(W_{k_o}',(J_{k_o}',b'), E_{k_o}')$ is simple (cf. Remark 1-5),
there exists an open covering $\{(W_{k_o}')^{\gamma}\}_{\gamma \in \Gamma}$
satisfying the following condition: for each $\gamma \in \Gamma$, there exists
a smooth hypersurface $(W_{k_o}')_h^{\gamma}$, embedded as a closed subscheme,
such that

1. $I((W_{k_o}')_h^{\gamma}) \subset
\Delta^{b'-1}(J_{k_o}')|_{(W_{k_o}')^{\gamma}}.$

Remark that we do not (and
can not) require any transversality condition on
$(W_{k_o}')_h^{\gamma}$ with respect to $E_{k_o}' = E_{k_o}$.

\vskip.1in

We claim that the open covering $\{(W_k'')^{\gamma}\}_{\gamma \in \Gamma}$ where
$(W_k'')^{\gamma}$ are the inverse images $(W_k')^{\gamma}$ of $(W_{k_o}')^{\gamma}$, together with
$(W_k'')_h^{\gamma}
\subset (W_k'')^{\gamma}$ where $(W_k'')_h^{\gamma}$ are the strict transforms
$(W_k')_h^{\gamma}$ of
$(W_{k_o}')_h^{\gamma}$, satisfies conditions 1 and 2 of the key inductive
lemma as required.  (We note that when we have a smooth morphism we call the
inverse image of the hypersurface ``the strict transform" by abuse of language.)

\vskip.1in

In order to check condition 1, by the repetitive use of Giraud's Lemma (cf. Claim 3-4),
we see that
$$I((W_k'')_h^{\gamma}) \subset \Delta^{b'-1}((J_{k_o}')_{k-k_o}).$$
Then noting that
$$\align
&(J_{k_o}')_{k-k_o} = J_{k_o + (k - k_o)}' = J_k'\\
&J_k' \subset J_k'' \text{\ and\ hence\ }\Delta^{b'-1}(J_k') \subset \Delta^{b'-1}(J_k''), \\
\endalign$$ 
we finally conclude that
$$I((W_k'')_h^{\gamma}) \subset \Delta^{b'-1}(J_k'').$$

In order to check condition 2, $(W_j')_h^{\gamma}$ being the strict transform of
$(W_{k_o}')_h^{\gamma}$, inductively for $j = k_o+1, ... , k$, we see, whenever
$\pi_j'$ is the transformation with center $Y_{j-1}' \subset
\roman{Sing}((J_{k_o}')_{j-1-k_o},b') = \roman{Sing}(J_{j-1}',b') =
\underline{\roman{Max}}\ w\text{-}\roman{ord}_{j-1}$, that 
$$\align
&Y_{j-1}' \cap
(W_{j-1}')^{\gamma} \subset (W_{j-1}')_h^{\gamma} \hskip.1in (\text{cf.\ Claim\ 3-4})\\
&Y_{j-1}' \text{\ permissible\ with\ respect\ to\ }E_{j-1}' \cap
(W_{j-1}')^{\gamma} = E_{j-1} \cap (W_{j-1}')^{\gamma},\\ &(\text{and\ hence\
}Y_{j-1}' \text{\ permissible\ with\ respect\ to\ }E_{j-1}^+ \cap
(W_{j-1}')^{\gamma}),\\
\endalign$$
which implies inductively (See the argument in $\bold{Case\ A}$ for permissibility of an
irreducible component $D \subset R(1)(\underline{\roman{Max}}\ t_k \cap W_k^{\lambda})$ with
respect to
$E_k^+ \cap W_k^{\lambda}$ in the case $D$ is not any one of the exceptional divisor
$\widetilde{W_{k_o}^{\lambda}} \leftarrow \widetilde{W_k^{\lambda}}$.) that

for each $\gamma \in \Gamma$, $(W_j')_h^{\gamma}$ is permissible with respect to $E_j^+ \cap
(W_j')^{\gamma}$ and $(W_j')_h^{\gamma} \not\subset E_j^+$.

\vskip.1in

In particular, for $j = k$, we have condition 2: 

2. $(W_k'')_h^{\gamma} = (W_k')_h^{\gamma}$ is permissible with
respect to
$E_k^+ \cap (W_k'')^{\gamma} = E_k''
\cap (W_k'')^{\gamma}$ and $(W_k'')_h^{\gamma} \not \subset E_k''$.
 
\vskip.1in

This completes the proof of Lemma 5-4.

\vskip.2in

\noindent $\italic{Conclusion\ of\ the\ proof\ for\ the\ assertions\ in} \bold{\ Case\ B}
\italic{\ under\ possibility\ }\bold{\underline{P3}}$

\vskip.2in

Now we go back to the proof of the assertions in $\bold{Case\ B}$: 
$R(1)(\underline{\roman{Max}}\
t_k) =
\emptyset$.

\vskip.1in

We construct a general basic object over $(G_k = \underline{\roman{Max}}\ t_k,
(W_k,E_k'' = E_k^+))$, with a $d$-dimensional structure first, by specifying its charts of
basic objects
$\{(\widetilde{{W_k''}^{\lambda}}, ({{\goth a}_k''}^{\lambda},{b''}^{\lambda}),
\widetilde{{E_k''}^{\lambda}})\}$ of dimension $d$ in the following way:

Let $\{(\widetilde{W_k^{\lambda}}, (\goth a_k^{\lambda},b^{\lambda}),
\widetilde{E_k^{\lambda}})$ be the charts for the general basic objects $({\Cal F}_k,
(W_k,E_k))$ arising from the sequence (cf. Note 4-3)
$$(F_0,(W_0,E_0)) \leftarrow \cdot\cdot\cdot \leftarrow (F_k,(W_k,E_k)).$$ 

We take $\widetilde{{W_k''}^{\lambda}} = \widetilde{W_k^{\lambda}}$.

If $\widetilde{{W_k''}^{\lambda}} \cap \underline{\roman{Max}}\ t_k = \emptyset$, then we
take the basic object $(\widetilde{{W_k''}^{\lambda}}, ({{\goth
a}_k''}^{\lambda},{b''}^{\lambda}),
\widetilde{{E_k''}^{\lambda}})$ to be
$$\left\{\aligned
\widetilde{{W_k''}^{\lambda}} &= \widetilde{W_k^{\lambda}} \\
{{\goth a}_k''}^{\lambda} &= {\Cal O}_{\widetilde{{W_k''}^{\lambda}}} \\
{b''}^{\lambda} &= 1 \\
\widetilde{{E_k''}^{\lambda}} &= E_k^+ \cap \widetilde{{W_k''}^{\lambda}}.\\
\endaligned\right.$$

If $\widetilde{{W_k''}^{\lambda}} \cap \underline{\roman{Max}}\ t_k \neq \emptyset$, then we
take the basic object $(\widetilde{{W_k''}^{\lambda}}, ({{\goth
a}_k''}^{\lambda},{b''}^{\lambda}),
\widetilde{{E_k''}^{\lambda}})$ to be
$$\left\{\aligned
\widetilde{{W_k''}^{\lambda}} &= \widetilde{W_k^{\lambda}} \\
{{\goth a}_k''}^{\lambda} &= ({\goth a}_k^{\lambda})'' \text{\ as\ constructed\ in\ Lemma\
5-4}\\ 
{b''}^{\lambda} &= (b^{\lambda})'' \text{\ as\ constructed\ in\ Lemma\
5-4}\\
\widetilde{{E_k''}^{\lambda}} &= \widetilde{E_k^{\lambda}}^+ = E_k^+ \cap
\widetilde{{W_k''}^{\lambda}}.\\
\endaligned\right.$$

Let ${\goth C}_G$ be the collection of all the sequences of transformations and smooth
morphisms of pairs with specified closed subsets, starting with $(G_k,(W_k,E_k''))$, which
satisfy condition (GB-1) with respect to the charts $\{(\widetilde{{W_k''}^{\lambda}},
({{\goth a}_k''}^{\lambda},{b''}^{\lambda}),
\widetilde{{E_k''}^{\lambda}})\}$.  Then condition (GB-3) is trivially satisfied by the
construction, whereas condition (GB-0) is a consequence of the statements of Lemma 5-4
for $N = 0$ and condition (GB-2) a consequence of the statements of Lemma 5-4 for $N$
general.

Therefore, the collection ${\goth C}_G$ is represented by a general basic object $({\Cal
G}_k,(W_k,E_k''))$ over
$(G_k,(W_k,E_k''))$ with a $d$-dimensional structure, having charts
$\{(\widetilde{{W_k''}^{\lambda}}, ({{\goth a}_k''}^{\lambda},{b''}^{\lambda}),
\widetilde{{E_k''}^{\lambda}})\}$ of dimension $d$.  (In short and roughly
speaking, the general basic object
$({\Cal G}_k,(W_k,E_k''))$ is the one whose specified closed subsets coincides with the
loci $\underline{\roman{Max}}\ t$ of the sequences of transformations and smooth morphisms
of general basic objects satisfying condition $(\heartsuit')$ and extending the original
sequence
$$({\Cal F}_0, (W_0,E_0)) \leftarrow \cdot\cdot\cdot \leftarrow ({\Cal F}_k, (W_k,E_k).)$$

\vskip.1in

Now the ``Moreover" part of Lemma 5-4, applied to the charts $\{(\widetilde{{W_k''}^{\lambda}},
({{\goth a}_k''}^{\lambda},{b''}^{\lambda}),
\widetilde{{E_k''}^{\lambda}})\}$, and the key inducive lemma imply
that the general basic object $({\Cal G}_k,(W_k,E_k''))$, which represents the collection
${\goth C}_G$, has a
$(d-1)$-dimensional structure.

\vskip.1in

It also follows from Lemma 5-4 that the general basic object $({\Cal G}_k,(W_k,E_k''))$ has
the following two properties $(\alpha)$ and $(\beta)$:

\vskip.1in

$(\alpha)$ With each sequence in ${\goth C}_G$
$$(G_k,(W_k,E_k'')) \overset{\pi_{k+1}''}\to{\leftarrow} \cdot\cdot\cdot
\overset{\pi_{k+N}''}\to{\leftarrow} (G_{k+N}, (W_{k+N},E_{k+N}''))$$
satisfying the condition
$$G_{k+j} \neq \emptyset \text{\ for\ }j = 0, ... , N - 1,$$
there corresponds an extension of the original sequence of transformations and smooth
morphisms of general basic objects
$$\align
&({\Cal F}_0,(W_0,E_0)) \overset{\pi_1}\to{\leftarrow} \cdot\cdot\cdot
\overset{\pi_{k-1}}\to{\leftarrow} ({\Cal F}_{k-1},(W_{k-1},E_{k-1}))
\overset{\pi_k}\to{\leftarrow} \\  
&({\Cal F}_k,(W_k,E_k)) \overset{\pi_{k+1}}\to{\leftarrow} \cdot\cdot\cdot
\overset{\pi_{k+N}}\to{\leftarrow} ({\Cal F}_{k+N},(W_{k+N},E_{k+N}))\\
\endalign$$
with condition
$$(\heartsuit') \hskip.1in \left\{\aligned &Y_{i-1} \subset \underline{\roman{Max}}\
t_{i-1} \subset
\underline{\roman{Max}}\ w\text{-}\roman{ord}_{i-1} (\subset F_{i-1})\\
&\text{whenever\ } \pi_i \text{\ is\ a\ transformation\ with\ center\
}Y_{i-1}\\
\endaligned\right\} \text{\ for\ }i = k+1 , ... , N$$
satisfying the following conditions:

\vskip.1in

(i) $\pi_{k+j+1}''$ and $\pi_{k+j+1}$ are the transformations with the same
centers or the same smooth morphisms (as abstract varieties) for $j = 0, ... , N - 1$
(which means, in particular, if
$\pi_{k+j+1}''$ is the transformation with center $Y_{k+j}'' \subset W_{k+j}$ which is
permissible for
$(G_{k+j}, (W_{k+j}, E_{k+j}''))$, then $\pi_{k+j+1}$ is the transformation with the
same center $Y_{k+j}'' \subset W_{k+j} = W_{k+j}''$ which is also permissible for
$(F_{k+j}, (W_{k+j}, E_{k+j})$), 

\vskip.1in
 
(ii) we have 

\vskip.1in

\noindent either
$$\left\{\aligned
&\max\ t_k = \max\ t_{k+1} = \cdot\cdot\cdot = \max\ t_{k+N},
\text{\ and\ }\\
&G_{k+j} = \underline{\roman{Max}}\ t_{k+j} \text{\ for\ }j = 0,
\cdot\cdot\cdot, N\\
\endaligned\right.$$
or
$$\left\{\aligned
&\max\ t_k = \max\ t_{k+1} = \cdot\cdot\cdot = \max\ t_{k+N-1} >
\max\ t_{k+N}\\
&(\text{or\ }\max\ t_k = \max\ t_{k+1} = \cdot\cdot\cdot = \max\ t_{k+N-1} \hskip.1in \&
\hskip.1in F_{k+N} =
\emptyset),
\text{\ and}\\ 
&G_{k+j} = \underline{\roman{Max}}\ t_{k+j} \text{\ for\ }j = 0,
\cdot\cdot\cdot, N-1 \hskip.1in \& \hskip.1in G_{k+N} = \emptyset.\\
\endaligned\right.$$

$(\beta)$ Conversely, with each extension of the original sequence of transformations and
smooth morphisms of general basic objects 
$$\align
&({\Cal F}_0,(W_0,E_0)) \overset{\pi_1}\to{\leftarrow} \cdot\cdot\cdot
\overset{\pi_{k-1}}\to{\leftarrow} ({\Cal F}_{k-1},(W_{k-1},E_{k-1}))
\overset{\pi_k}\to{\leftarrow} \\  
&({\Cal F}_k,(W_k,E_k)) \overset{\pi_{k+1}}\to{\leftarrow} \cdot\cdot\cdot
\overset{\pi_{k+N}}\to{\leftarrow} ({\Cal F}_{k+N},(W_{k+N},E_{k+N}))\\
\endalign$$
with condition
$$(\heartsuit') \hskip.1in \left\{\aligned &Y_{i-1} \subset \underline{\roman{Max}}\
t_{i-1} \subset
\underline{\roman{Max}}\ w\text{-}\roman{ord}_{i-1} (\subset
F_{i-1}) \\
&\text{whenever\ } \pi_i \text{\ is\ a\ transformation\ with\ center\
}Y_{i-1}\\
\endaligned\right\} \text{\ for\ }i = k+1 , ... , N$$

and the condition
$$\max\ t_k = \max\ t_{k+1} = \cdot\cdot\cdot = \max\ t_{k+N-1},$$
there corresponds a sequence in ${\goth C}_G$
$$(G_k,(W_k,E_k'')) \overset{\pi_{k+1}''}\to{\leftarrow} \cdot\cdot\cdot
\overset{\pi_{k+N}''}\to{\leftarrow} (G_{k+N}, (W_{k+N},E_{k+N}''))$$
satisfying the condition
$$G_{k+j} \neq \emptyset \text{\ for\ }j = 0, ... , N - 1$$
and conditions (i) and (ii) as in $(\alpha)$.

(Note that the above properties $(\alpha)$ and $(\beta)$ provide a characterization of the
general basic object
$({\Cal G}_k,(W_k,E_k''))$ free of the description using the charts we construct, and that it is via this property
of the collections ${\goth C}$ and ${\goth C}_G$ that we verify ${\goth C}_G$ satisfies
conditions (GB-0) and (GB-2).)

\vskip.2in

The assertions in
$\bold{B\text{-}1,\ B\text{-}2,\ B\text{-}3}$ follow immediately from this.

\vskip.1in

Starting from $(F_0,(W_0,E_0))$, the inductive algorithm allows one by
induction on the dimension $d$ of the structure to construct a unique extension of the
sequence of transformations (constructed up to the $k$-th stage) satisfying condition
$(\heartsuit')$ such that
$$\align
\text{either\ } &F_{k+N} = \emptyset \\
&\text{\ where\ resolution\ of\ singularities\ is\
already\ achieved},\\ 
\text{or\ }&F_{k+N} \neq \emptyset \ \&\ \max\ w\text{-}\roman{ord}_{k+N} = 0 \\
&\text{\ where\ reslotuion\ of\ singularities\ is\
reduced\ to that\ of\ a\ monomial\ case},\\
\text{or\ }&F_{k+N} \neq \emptyset, \max\ w\text{-}\roman{ord}_{k+N} > 0\ \&\ \emptyset\max t_k
>
\max\ t_{k+N}.\\
\endalign$$

Thanks to condition (GB-3), the values of the $t$-inavariant are in $\frac{1}{c!}{\Bbb
Z}_{\geq 0} \times {\Bbb Z}_{\geq 0}$ and hence satisfy the descending chain condition. 
Therefore, the third possibility can not happen infinitely many times.  Thus after finitely
many executions of the inductive algorithm, we obtain the sequence representing resolution
of singularities as asserted.

\vskip.1in

We note that resolution of singularities of a general basic object with a $1$-dimensional
structure is obvious.  In fact, at the stage $i = 0$, we are always in $\bold{Case\ A}$ only
with the possibilities $\bold{A\text{-}1,\ A\text{-}2}$ at the stage $i = 1$.  If
$\bold{A\text{-}1}$ is the case, then resolution of singularities is already achieved.  If
$\bold{A\text{-}2}$ is the case, then we are reduced to the monomial case, where resolution
of singularities is achieved by nothing but a sequence of consecutive point blowups.  This
supports the tower of the inductional algorithm at the bottom $d = 1$.

\vskip.1in

This completes the proof of Theorem 5-1.

\vskip.2in

\proclaim{Remark 5-5}\endproclaim

The role of the invariant $t$ is absolutely crucial at a couple of places, in the inductive
algorithm presented as Theorem 5-1, e.g.:

\hskip.1in $\circ$ permissiblity of the center $R(1)(\underline{\roman{Max}}\ t_k)$ in
$\bold{Case\ A}$.

\hskip.1in $\circ$ the proof that $(W_k'',(J_k'',b''),E_k'')$ satisfies conditions 1
and 2 as stated in the 

\hskip.1in key inductive lemma. 

We refer the reader to Chapter 6 for a more and natural explanation of the origin of the
$t$-invariant by breaking up the inductive algorithm into a couple of natural reduction steps.

\vskip.1in

\proclaim{Corollary 5-6 (Resolution of singularities of a basic object)} Let
$(W_0,(J_0,b),E_0)$ be a basic object with $d = \dim W_0$.  Then there exists a sequence of
transformations only
$$(W_0,(J_0,b),E_0) \leftarrow \cdot\cdot\cdot \leftarrow (W_k,(J_k,b),E_k)$$
satisfying the following condition $(\heartsuit')$ for $i = 1, ... , k$
$$(\heartsuit')\hskip.1in Y_{i-1} \subset \underline{\roman{Max}}\ t_{i-1} \subset
\underline{\roman{Max}}\ w\text{-}\roman{ord}_{i-1} \text{\ if\ }\max\ w\text{-}\roman{ord}_{i-1}
> 0$$ 
such that
$$\max\ w\text{-}\roman{ord}_0 \geq \max\ w\text{-}\roman{ord}_1 \geq \cdot\cdot\cdot \geq
\max\ w\text{-}\roman{ord}_{k-1} \geq
\max\ w\text{-}\roman{ord}_k$$
and that
$$\roman{Sing}(J_k,b) = \emptyset,$$
i.e., the sequence represents resolution of singularities of the basic object
$(W_0,(J_0,b),E_0)$.
\endproclaim

\demo{Proof}\enddemo A basic object $(W_0,(J_0,b),E_0)$ defines a general basic object
$({\Cal F}_0,(W_0,E_0))$ with a $d$-dimensional structure, as explained in Remark 4-2 (i).  A
sequence representing resolution of singularities of $({\Cal F}_0,(W_0,E_0))$, obtained via
the inductive algorithm of Theorem 5-1, provides that of $(W_0,(J_0,b),E_0)$ with the
required properties. 

\newpage

${}\hskip.4in$$\bold{CHAPTER\ 6.\ A\ MORE\ DOWN\text{-}TO\text{-}EARTH}$ \linebreak
${}\hskip.4in$$\bold{APPROACH\ TO\ THE\ INDUCTIVE\ ALGORITHM}$

\vskip.1in

In Chapter 5 we saw how the inductive algorithm for resolution of singularities of a (general)
basic object works.  However, in the untrained eyes (e.g. those of the author), its
mechanism may look more like a miracle than a natural process.  Especially the
ingenious $t$-invariant seems to have come ``out of blue".  The purpose of this chapter
is to explain the mathematical origin of the
$t$-invariant and show how natural the inductive algorithm is.

For this purpose, firstly we break up the problem of resolution of singularities of a (general)
basic object into the following three stages, depending upon the restrictions on the basic
objects to deal with (and hence with increasing difficulties, going from $A_d$ through
$B_d$ to
$C_d$).

The descriptions of the restrictions we put at these stages on the basic objects are:

\vskip.1in

$A_d$: simple, with empty boundary, $d$-dimensional,

$B_d$: simple, with possibly non-empty boundary, $d$-dimensional, and

$C_d$: the general $d$-dimensional case with no restrictions, i.e., not necessarily 

\hskip.2in simple, with
possibly non-empty boundary, $d$-dimensional.

\vskip.1in

Secondly we establish the reduction steps.

\vskip.1in

Reduction $A_d$ to $C_{d-1}$: via the key inductive lemma,

Reduction $B_d$ to $C_{d-1}\ (+\ A_d)$: via the key inductive lemma and
 
\hskip1.4in the introduction of the
$n$-invariant, 

Reduction $C_d$ to $B_d$: via the introduction of the invariant $w\text{-}\roman{ord}$,

\vskip.1in

and hence establish the inductive algorithm for resolution of singularities of
a (general) basic object.  

\vskip.1in

We denote the three reduction steps as above figuratively by

\vskip.07in

$\boxed{\text{Reduction}\ A_d \leftarrow C_{d-1}}$

$\boxed{\text{Reduction}\ B_d \leftarrow C_{d-1}}$
 
$\boxed{\text{Reduction}\ C_d \leftarrow B_d\hskip.15in}$

\vskip.07in

\noindent where, e.g., $A_d \leftarrow C_{d-1}$ indicates that a solution to the problem
of resolution of singularities at stage $C_{d-1}$ implies a solution to the problem of
resolution of singularities at stage $A_d$, with the arrow $\leftarrow$ representing the
implication.
 
(We note and emphasize that logically we only need the last two reduction steps $C_d
\leftarrow B_d, B_d \leftarrow C_{d-1}$ to complete the inductive algorithm and not 
the first \linebreak
$A_d \leftarrow C_{d-1}$, whose construction in showing the reduction step, but not whose
consequence, is used in the argument to show the reduction step
$B_d
\leftarrow C_{d-1}$.  The first reduction step is included solely to demonstrate how
the key inductive lemma directly and easily solves the problem of resolution of singularities
for a simple basic object with an empty boundary by induction.) 

We observe that the presentation of the inductive algorithm in this chapter is a
decomposition of the one in Chapter 5 into a few reduction steps as above and that the
$t$-invariant is the natural combination of the $n$-invariant and the invariant
$w\text{-}\roman{ord}$ when one wants to put these reduction steps together back into one
nice package as presented in Chapter 5.

This chapter is based upon one of the lectures delivered by Prof. Villamayor under the title
``Constructive Desingularization" at Purdue University.

\vskip.1in

We start with giving a more precise description of the three stages $A_d$, $B_d$, and $C_d$ of the
problem of resolution of singularities of a (general) basic object.

\vskip.1in

\noindent $\bold{Description\ of\ \text{$A_d$},\ \text{$B_d$},\ and\ \text{$C_d$}}$

\vskip.1in

$\underline{A_d: \text{\ simple,\ with\ empty\ boundary,\ $d$-dimensional}}$

\vskip.1in

Resolution of singularities of a simple (general) basic object
$$(W, (J,b), E)$$
with an empty boundary divisor of dimension $d$, i.e.,
$$\left\{\aligned
\dim W &= d, \\
E &= \emptyset, \\
\nu_p(J) &= b \hskip.1in \forall p \in \roman{Sing}(J,b). \\
\endaligned
\right.$$

(In case of a general basic object given by the charts $\{(\widetilde{W_0^{\lambda}},
({\goth a}_0^{\lambda},b^{\lambda}), \widetilde{E_0^{\lambda}})\}$ we require that all
the charts are simple basic objects with empty boundary divisors of dimension
$d$,
i.e.,
$$\forall \lambda \in \Lambda \hskip.1in
\left\{\aligned
\dim \widetilde{W_0^{\lambda}} &= d, \\
\widetilde{E_0^{\lambda}} &= \emptyset, \\
\nu_p({\goth a}_0^{\lambda}) &= b^{\lambda} \hskip.1in \forall p \in \roman{Sing}({\goth
a}_0^{\lambda},b^{\lambda}).)
\\
\endaligned
\right.$$

\vskip.1in

$\underline{B_d: \text{\ simple,\ with\ possibly\ non-empty\ boundary,\ $d$-dimensional}}$

\vskip.1in

Resolution of singularities of a simple (general) basic object
$$(W, (J,b), E)$$
with a possibly non-empty boundary divisor of dimension $d$, i.e.,
$$\left\{\aligned
\dim W &= d, \\
E & \text{\ an\ arbitrary\ divisor\ with\ simple\ normal\ crossings}, \\
\nu_p(J) &= b \hskip.1in \forall p \in \roman{Sing}(J,b). \\
\endaligned
\right.$$

(In case of a general basic object given by the charts $\{(\widetilde{W_0^{\lambda}},
({\goth a}_0^{\lambda},b^{\lambda}), \widetilde{E_0^{\lambda}})\}$ we require that all
the charts are simple basic objects with possibly non-empty boundary divisors of
dimension
$d$, i.e.,
$$\forall \lambda \in \Lambda \hskip.1in
\left\{\aligned
\dim \widetilde{W_0^{\lambda}} &= d, \\
\widetilde{E_0^{\lambda}} &= E \cap \widetilde{W_0^{\lambda}}, \\
&\text{where\ }E \text{\ is\ an\ arbitrary\ divisor\ with\ simple\ normal\ crossings,}\\
&\text{intersecting\ }\widetilde{W_0^{\lambda}} \text{\ transversally},\\
\nu_p({\goth a}_0^{\lambda}) &= b^{\lambda} \hskip.1in \forall p \in \roman{Sing}({\goth
a}_0^{\lambda},b^{\lambda}).)
\\
\endaligned
\right.$$

\vskip.1in

$\underline{C_d: \text{the\ general\ $d$-dimensional\ case\ with\ no\ restrictions,\ i.e.,}}$ 
 
\hskip.2in $\underline{\text{not\ necessarily\ simple,\ with\ possibly\ non-empty\ boundary,\
$d$-dimensional}}$

\vskip.1in

Resolution of singularities of a (general) basic object
$$(W, (J,b), E)$$
without any restrictions of dimension $d$, i.e., it may not be simple and with a possibly non-empty
boundary divisor.

(In case of a general basic object given by the charts $\{(\widetilde{W_0^{\lambda}},
({\goth a}_0^{\lambda},b^{\lambda}),
\widetilde{E_0^{\lambda}})\}$ there is no restriction on the basic objects
$(\widetilde{W_0^{\lambda}}, ({\goth a}_0^{\lambda},b^{\lambda}),
\widetilde{E_0^{\lambda}})$ other than $\dim \widetilde{W_0^{\lambda}} = d$, i.e., they
may not be simple but with possibly non-empty boundary divisors, and of dimension $d$.)

\vskip.1in

\noindent $\bold{Arguments\ for\ the\ reduction\ steps}$

\vskip.1in

Now we are ready to provide the arguments for the following three reduction steps:

\vskip.07in
 
$\boxed{\text{Reduction}\ A_d \leftarrow C_{d-1}}$

$\boxed{\text{Reduction}\ B_d \leftarrow C_{d-1}}$
 
$\boxed{\text{Reduction}\ C_d \leftarrow B_d\hskip.15in}$.

\vskip.2in

\noindent $\boxed{\text{Reduction}\ A_d \leftarrow C_{d-1}}$

\vskip.1in

Let $(W, (J,b), E)$ be a simple basic object with $E = \emptyset$.  

\vskip.1in

Then we can find an open covering
$\{W^{\lambda}\}_{\lambda \in \Lambda}$ of $W$ with smooth hypersurfaces $W_h^{\lambda} \subset
W^{\lambda}$, embedded as closed subschemes, satisfying conditions 1 and 2 of
the key inductive lemma.  In fact, by Remark 1-5 with
$$S = \{p \in W;\nu_p(J) = b_{\max} = b\} = \roman{Sing}(J,b)$$
where $b_{\max} = b$ since $(W, (J,b), E)$ is simple, we can find an open covering
$\{W^{\lambda}\}_{\lambda
\in
\Lambda}$ of $W$ with smooth hypersurfaces $W_h^{\lambda} \subset W^{\lambda}$ satisfying
condition 1.  Now condition 2 is void and hence automatically satisfied, since $E$ is empty.

\vskip.1in

$\bold{Case\ A}$: $R(1)(\roman{Sing}(J,b)) \neq \emptyset.$

\vskip.1in

As in the proof of the key inductive lemma, it is easy to see in this case that
$R(1)(\roman{Sing}(J,b))$ is smooth and open in $\roman{Sing}(J,b)$ and that it is permissible with respect to $E =
\emptyset$ automatically.  After taking the transformation with center $R(1)(\roman{Sing}(J,b))$,
we are in
$\bold{Case\ B}$.

\vskip.1in

$\bold{Case\ B}$: $R(1)(\roman{Sing}(J,b)) = \emptyset.$

\vskip.1in

Having an open covering $\{W^{\lambda}\}_{\lambda \in \Lambda}$ of $W$ with smooth hypersurfaces
\linebreak
$W_h^{\lambda}
\subset W^{\lambda}$ satisfying conditions 1 and 2, we are in a position to apply the key
inductive lemma to conclude that resolution of singularities of $(W, (J,b), E)$ is reduced
to resolution of singularities of the general basic objects whose charts are given by
$\{(\widetilde{W_0^{\lambda}}, ({\goth a}_0^{\lambda},b^{\lambda}), \widetilde{E_0^{\lambda}})\}$
as constructed in the key inductive lemma.

\vskip.1in

Thus we see that $A_d$ is reduced to $C_{d-1}$.

\vskip.1in

\noindent $\bold{Crucial\ Remark\ in\ \text{$A_d$}}$

\vskip.1in

Remark that if we have a sequence of transformations of basic objects
$$(W, (J,b), E = \emptyset) = (W_0, (J,b), E_0) \leftarrow \cdot\cdot\cdot \leftarrow (W_k,
(J_k,b), E_k)$$
and if we have an open covering $\{W^{\lambda}\}_{\lambda \in \Lambda}$ with
smooth hypersurfaces $W_h^{\lambda} \subset W^{\lambda}$ satisfying condition 1
(and 2), then the open covering $\{W_k^{\lambda}\}_{\lambda \in \Lambda}$ consisting of the
inverse images $W_k^{\lambda}$ of
$W^{\lambda} = W_0^{\lambda}$ and the smooth hypersurfaces consisting of the strict transforms
$(W_k^{\lambda})_h
\subset W_k^{\lambda}$ of
$(W_0^{\lambda})_h = W_h^{\lambda}
\subset W^{\lambda} = W_0^{\lambda}$, satisfy conditions 1 and 2 of the key inductive lemma for
the simple basic object $(W_k, (J_k,b), E_k)$.  This is a consequence of Giraud's Lemma (cf.
Claim 3-4), which is an essential part of the construction in the key inductive lemma.

\vskip.2in

\noindent $\boxed{\text{Reduction}\ B_d \leftarrow C_{d-1}}$

\vskip.1in

Let $(W, (J,b), E)$ be a simple basic object with $E = \{H_1, ... , H_r\}$.

\vskip.1in

Let
$$(W, (J,b), E) = (W_0, (J_0,b), E_0) \leftarrow (W_1, (J_1,b), E_1) \leftarrow \cdot\cdot\cdot
\leftarrow (W_k, (J_k,b), E_k)$$
be a sequence of transformations of basic objects.  We decompose $E_k = E_k^- \cup E_k^+$ into two
disjoint subsets $E_k^- = \{H_1, ... , H_r\}$ and its complement $E_k^+ = \{H_{r+1}, ... ,
H_{r+k}\}$ \linebreak
in
$E_k =
\{H_1, ... , H_r, H_{r+1}, ... , H_{r+k}\}$.  (Look at the convention in Definition 1-8 (iii)
and see also Definition 1-10 (iii).)  Note that the assumption $(W, (J,b), E)$ being simple
implies
$(W_i, (J_i,b), E_i)$ also being simple and $w\text{-}\roman{ord}_i = \roman{ord}_i$ for $i =
0, ... , k$.)  Remark that this decomposition is
motivated by the crucial remark at the end of the discussion of $A_d$.

\vskip.1in

Define the function (cf. Definition 1-10 (iii))
$$n_k:\roman{Sing}(J_k,b) \rightarrow {\Bbb Z}_{\geq 0}$$
by
$$n_k(p) = \#\{H_i \in E_k^-; p \in H_i\} \text{\ for\ }p \in \roman{Sing}(J_k,b).$$

\vskip.1in

\noindent Case: $\max\ n_k = 0$.

\vskip.1in

Observe that corresponding to the original sequence of transformations of basic objects
$$(W, (J,b), E) = (W_0, (J_0,b), E_0) \leftarrow (W_1, (J_1,b), E_1) \leftarrow \cdot\cdot\cdot
\leftarrow (W_k, (J_k,b), E_k)$$
we have another sequence of transformations of basic objects
$$\align
(W, (J,b), \emptyset) = (W_0, (J_0,b), E_0^+) &\leftarrow (W_1, (J_1,b), (E_0^+)_1 = E_1^+)
\leftarrow
\cdot\cdot\cdot \\
&\hskip.8in \cdot\cdot\cdot \leftarrow (W_k, (J_k,b), (E_0^+)_k = E_k^+).\\
\endalign$$

Since $\roman{Sing}(J_k,b) \cap E_k^- = \emptyset$ under the case assumption $\max\ n_k = 0$, we
conclude that a sequence representing resolution of singularities of $(W_k, (J_k,b), E_k^+)$ is
also a sequence representing resolution of singularities of $(W_k, (J_k,b), E_k)$.  Thus we only
have to find a sequence representing resolution of singularities of $(W_k, (J_k,b), E_k^+)$.

We now observe, by the crucial remark at the end of the reduction step \linebreak
$A_d
\leftarrow C_{d-1}$, that we can find an open covering
$\{W_k^{\lambda}\}_{\lambda \in
\Lambda}$ of
$W_k$ with smooth hypersurfaces $(W_k^{\lambda})_h \subset W_k^{\lambda}$, embedded as closed
subschemes, satisfying conditions 1 and 2 of the key inductive lemma.

\vskip.1in

$\bold{Case\ A}$: $R(1)(\roman{Sing}(J_k,b)) \neq \emptyset$.

\vskip.1in

Again as in the proof of the key inductive lemma, it is easy to see in this case that
$R(1)(\roman{Sing}(J,b))$ is smooth and open in $\roman{Sing}(J_k,b)$ and that it is
permissible with  respect to $E_k^+$, as so is $(W_k^{\lambda})_h$.  After taking the
transformation with center
$R(1)(\roman{Sing}(J_k,b))$, we are in
$\bold{Case\ B}$.

\vskip.1in

$\bold{Case\ B}$: $R(1)(\roman{Sing}(J_k,b)) = \emptyset$.

\vskip.1in

Having an open covering $\{W_k^{\lambda}\}_{\lambda \in \Lambda}$ of $W_k$ with smooth
hypersurfaces
\linebreak
$(W_k^{\lambda})_h
\subset W_k^{\lambda}$ satisfying conditions 1 and 2, we are in a position to apply the key
inductive lemma to conclude that resolution of singularities of $(W_k, (J_k,b), E_k^+)$ is reduced
to resolution of singularities of the general basic objects whose charts are given by
$\{(\widetilde{W_k^{\lambda}}, ({\goth a}_k^{\lambda},b^{\lambda}),
\widetilde{(E_k^+)^{\lambda}})\}$ as constructed in the key inductive lemma.

\vskip.1in

Thus we see that $B_d$ is reduced to $C_{d-1}$ via the key inductive lemma in this case of
$\max\ n_k = 0$.
 
\vskip.1in

\noindent Case: $\max\ n_k = l_k > 0$.

\vskip.1in

Consider the following basic object $(V, (I,c), G)$ \hfill\hfill\linebreak
where $(V, (I,c), G) = ({W_k}'',
({J_k}'',b''), {E_k}'')$ in the notation of Lemma 5-4, i.e.,
$$\left\{\aligned
V &= W_k \\
I &= (J_k) + \prod_{\{H_{s_1}, ... , H_{s_{l_k}}\} \subset E_k^- \text{\ with\ }H_{s_1}, ... ,
H_{s_{l_k}} \text{\ distinct}}\{I(H_{s_1}) +
\cdot\cdot\cdot I(H_{s_{l_k}})\}^b \\
c &= b \\
G &= E_k^+. \\
\endaligned\right.$$

(Remark that, since $(W,(J,b),E) = (W_0,(J_0,b),E_0)$ is simple, so is $(W_i,(J_i,b),E_i)$ and
$\roman{ord}_i = w\text{-}\roman{ord}_i\ \&\ J_i =
\overline{J_i}$ for $i = 0, ... , k$.  This implies \hfill\hfill \linebreak$b_k = b \cdot \max\
w\text{-}\roman{ord}_k = b \cdot \roman{ord}_k = b = b' = b''$ in the notation of Lemma 5-3 and
Lemma 5-4).

\vskip.1in

First it is easy to see that $(V, (I,c), G)$ is a simple basic object (as shown in Lemma 5-4). 
Second we claim that for the simple basic object
$(V, (I,c), G)$ we can find an open covering
$\{V^{\lambda}\}_{\lambda
\in
\Lambda}$ of
$V$ with smooth hypersurfaces $V_h^{\lambda} \subset V^{\lambda}$ satisfying conditions 1 and 2
of the key inductive lemma.  In fact, by the crucial reamrk at the end of the reduction step $A_d
\leftarrow C_{d-1}$, for the simple basic object $(W_k, (J_k,b), E_k^+)$ we can find an open
covering
$\{W_k^{\lambda}\}_{\lambda \in \Lambda}$ of $W_k$ with smooth hypersurfaces $(W_k^{\lambda})_h
\subset W_k^{\lambda}$ satisfying conditions 1 and 2 of the key inductive lemma.  Set
$$V^{\lambda} = W_k^{\lambda}, V_h^{\lambda} = (W_k^{\lambda})_h.$$
Then we have for the simple basic object $(V, (I,c), G)$
$$I(V_h^{\lambda}) = I((W_k^{\lambda})_h) \subset \Delta^{b-1}(J_k) \subset \Delta^{c-1}(I),$$
and hence satisfying condition 1.  Condition 2 is identical and satisfied both for the simple
basic object
$(V, (I,c), G)$ and for $((W_k, (J_k,b), E_k^+)$, since $G = E_k^+$.

Therefore, by the key inductive lemma and $C_{d-1}$, possibly after going through $\bold{Case\ A}$
first, we find a sequence representing resolution of singularities of $(V, (I,c), G)$
$$(V, (I,c), G) = (V_0, (I_0,c), G_0) \leftarrow \cdot\cdot\cdot \leftarrow (V_N, (I_N,c), G_N)$$
where $\roman{Sing}(J_N,c) = \emptyset$.  Then we observe that there corresponds an extension of
the original sequence of transformations with the same centers
$$(W_k, (J_k,b), E_k) \leftarrow \cdot\cdot\cdot \leftarrow (W_{k+N}, (J_{k+N},b), E_{k+N}),$$
where inductively for $j = 0, ... , N$ we see that
$$I_j = (J_{k+j}) + \prod_{\{H_{s_1}, ... , H_{s_{l_k}}\} \subset E_{k+j}^- \text{\
with\ }H_{s_1}, ... , H_{s_{l_k}} \text{\ distinct}}\{I(H_{s_1}) +
\cdot\cdot\cdot I(H_{s_{l_k}})\}^b$$
and that for $j = 1, ... , N$ the center $Y_{k + j - 1} = Y_{G,j-1}$ is permissible with respect
to $(W_{k + j - 1}, (J_{k + j - 1},b), E_k)$, since
$$Y_{G,j-1} \subset \roman{Sing}(I_{j-1},c) \subset \roman{Sing}(J_{k + j - 1},b),$$
since $Y_{G,j-1}$ is contained in $H_s \in E_{k + j - 1}^-$ if $Y_{G,j-1} \cap H_s \neq
\emptyset$
\linebreak
as $\max\ n_k = l_k = l_{k + j - 1} = \max\ n_{k + j - 1}$, and since $Y_{G,j-1}$ is
permissible with respect to $G_{j-1} = E_{k + j - 1}^+$.

Finally, since
$$I_N = (J_{k+N}) + \prod_{\{H_{s_1}, ... , H_{s_{l_k}}\} \subset E_{k+N}^- \text{\
with\ }H_{s_1}, ... , H_{s_{l_k}} \text{\ distinct}}\{I(H_{s_1}) +
\cdot\cdot\cdot I(H_{s_{l_k}})\}^b$$
and since
$$\roman{Sing}(I_N,c) = \emptyset,$$
we conclude that either
$$\roman{Sing}(J_{k+N},b) = \emptyset,$$
in which case the extension realizes a sequence representing resolution of singularities of $(W,
(J,b), E)$, or
$$\max\ n_k > \max\ n_{k+N},$$
in which case by induction on the maximum of the invariant $n$ we also obtain a sequence representing resolution of singularities of $(W,
(J,b), E)$.

\vskip.1in

This completes the proof of the reduction step $B_d \leftarrow C_{d-1}$ via the key inductive
lemma and the introduction of the invariant $n$.

\newpage

\noindent $\boxed{\text{Reduction}\ C_d \leftarrow B_d}$ 

\vskip.1in

Let $(W, (J,b), E)$ be a basic object without any restrictions except that \linebreak
$\dim W = d$,
i.e., it may not be simple but with a possibly non-empty boundary divisor, and of dimension $d$.

\vskip.1in

Let
$$(W, (J,b), E) = (W_0, (J_0,b), E_0) \leftarrow (W_1, (J_1,b), E_1) \leftarrow \cdot\cdot\cdot
\leftarrow (W_k, (J_k,b), E_k)$$
be a sequence of transformations of basic objects.

Define the function
$$w\text{-}\roman{ord}_k:\roman{Sing}(J_k,b) \rightarrow \frac{1}{b}{\Bbb Z}_{\geq 0}$$
by
$$w\text{-}\roman{ord}_k(p) = \frac{\nu_p(\overline{J_k})}{b} \text{\ for\ }p \in
\roman{Sing}(J_k,b)$$ 
as in Definition 1-10 (ii).

\vskip.1in

\noindent Case: $\max\ w\text{-}\roman{ord}_k = 0.$

\vskip.1in

In this case, the problem of resolution of singularities of $(W_k, (J_k,b), E_k)$ is reduced to
that of monomial basic objects, which is settled in Chapter 2.

\vskip.1in

\noindent Case: $\max\ w\text{-}\roman{ord}_k > 0.$

\vskip.1in

In this case, consider the following basic object $(V, (I,c), G)$ where \linebreak
$(V, (I,c), G) =
(W_k', (J_k',b'), E_k')$ in the notation of Lemma 5-3, i.e.,
$$\left\{\aligned
V &= W_k' = W_k\\
I &= J_k' = \left\{\aligned
&\overline{J_k} \hskip2.5in \text{\ if\ }b_k \geq b \\
&\overline{J_k}^{b - b_k} + \{I(H_{r+1})^{a_{r+1}} \cdot\cdot\cdot
I(H_{r+k})^{a_{r+k}}\}^{b_k}
\ \ \text{\ if\ }b_k < b\\
\endaligned\right.\\
c &= b' = \left\{\aligned
&b_k \hskip2.52in \text{\ if\ }b_k \geq b\\
&b_k(b - b_k) \hskip2.05in \text{\ if\ }b_k < b\\
\endaligned\right.\\
G &= E_k' = E_k.
\endaligned\right.$$

Recall that
$$\align
J_k &= I(H_{r+1})^{a_{r+1}} \cdot\cdot\cdot I(H_{r+k})^{a_{r+k}} \cdot \overline{J_k} \\
b_k &= b \cdot (\max\ w\text{-}\roman{ord}_k).\\
\endalign$$

Then $(V, (I,c), G)$ is a simple basic object of dimension $d$ (as shown in Lemma 5-3).

Therefore, by $B_d$, we find a sequence representing resolution
of singularities of $(V, (I,c), G)$
$$(V, (I,c), G) = (V_0, (I_0,c), G_0) \leftarrow \cdot\cdot\cdot \leftarrow (V_N, (I_N,c), G_N)$$
where $\roman{Sing}(J_N,c) = \emptyset$.  Then we observe that there corresponds an extension of
the original sequence of transformations with the same centers
$$(W_k, (J_k,b), E_k) \leftarrow \cdot\cdot\cdot \leftarrow (W_{k+N}, (J_{k+N},b), E_{k+N}),$$
where inductively for $j = 0, ... , N-1$ we see that 
$$b_k = b_{k+j}, \text{\ i.e.\ }, \max\ w\text{-}\roman{ord}_k = \max\
w\text{-}\roman{ord}_{k+j} \hskip.1in \& \hskip.1in I_j = J_{k + j}'$$
and for $j = N$
$$I_N = \left\{\aligned
&\overline{J_{k+N}} \hskip3.3in\text{\ if\ }b_{k+N-1} = b_k \geq b \\
&\overline{J_{k+N}}^{b - b_{k+N-1}} + \{I(H_{r+1})^{a_{r+1}} \cdot\cdot\cdot
I(H_{r+k+N})^{a_{r+k+N}}\}^{b_{k+N-1}} \text{\ if\ }b_{k+N-1} = b_k < b,\\
\endaligned
\right.$$
and that for $j = 1, ... , N$ the center 
$$Y_{k+j-1} = Y_{G,j-1} \subset
\roman{Sing}(I_{j-1},c) = \roman{Sing}(J_{k+j-1}',c) \subset \roman{Sing}(J_{k+j-1},b)$$ 
is permissible with respect to $E_{k+j-1} = G_{j-1}$.

Finally since
$$\roman{Sing}(I_N,c) = \emptyset,$$
we conclude that either
$$\roman{Sing}(J_{k+N},b) = \emptyset,$$
in which case the extension realizes a sequence representing resolution of singularities of $(W,
(J,b), E)$, or
$$\max\ w\text{-}\roman{ord}_k > \max\ w\text{-}\roman{ord}_{k+N},$$
in which case by induction on the maximum of the invariant $w\text{-}\roman{ord}$ we also
obtain a sequence representing resolution of singularities of $(W, (J,b), E)$.

\vskip.1in

This completes the proof of the reduction step $C_d \leftarrow B_d$.

\vskip.1in

(The argument for resolution of singualrities of a general basic object is identical and
left to the reader as an exercise.)

\vskip.1in

\proclaim{Remark 6-1}\endproclaim

Though suppressed in the above argument of the reduction steps, it is absolutely necessary and
crucial to argue and verify that the processes of resolution of singularities of charts in the
general basic object patch up, via the observation that our choice of the centers only
depend on the invariants we set up and that those invariants are independent of charts, which
is one of the essential points in Chapter 4 via Hironaka's trick.  It is also necessary and
crucial to generalize the notion of a basic object to that of a general basic object to
complete the inductive step.

\vskip.1in

\proclaim{Exercise 6-2}\endproclaim

Check that the inductive algorithm for resolution of
singularities given by the reduction steps $C_d \leftarrow B_d$, $B_d \leftarrow C_{d-1}$
described as above actually coincides with the one described in Chapter 5 using the
$t$-invariant.

\newpage

$\bold{CHAPTER\ 7.\ EMBEDDED\ RESOLUTION\ OF\ SINGULARITIES}$

\vskip.1in

In this chapter, we present a proof for (embedded resoluion of singularities) stated
in Main Theme 0-2, as an easy consequence of (resolution of singularities of a basic
object) proved in Corollary 5-6.

\proclaim{Theorem 7-1 (Embedded resolution of singularities)} Let $X \subset W$ be a
variety, embedded as a closed subscheme of another variety $W$ smooth over a field $k$ of
charactersitic zero.  

We can construct a sequence of blowups
$$X = X_0 \subset W = W_0 \overset{\pi_1}\to{\leftarrow} X_1 \subset
W_1 \overset{\pi_2}\to{\leftarrow} \cdot\cdot\cdot \overset{\pi_{l-1}}\to{\leftarrow} X_{l-1}
\subset W_{l-1} \overset{\pi_l}\to{\leftarrow} X_l
\subset W_l$$
so that

(i) the centers $Y_{i-1} \subset W_{i-1}$ of the blowups $\pi_i \hskip.1in (i = 1,
... , l)$ are over $\roman{Sing}(X) = X \setminus \roman{Reg}(X)$, 

(ii) the centers $Y_{i-1} \subset W_{i-1}$ are permissible with respect to the
exceptional divisors
$E_{i-1} \subset W_{i-1}$ for the morphisms $\psi_{i-1} = \pi_1 \circ \pi_2 \circ \cdot\cdot\cdot
\circ \pi_{i-2} \circ \pi_{i-1}$ (which are simple normal crossing divisors),

(iii) the strict transform $X_l$ (of $X_0$) $\subset W_l$ is a variety smooth over $k$, permissibe
with respect to $E_l$, and the induced morphism $X = X_0 \overset{\pi}\to{\leftarrow} X_l$,
where $\pi = \psi_l = \pi_1
\circ
\pi_2 \circ \cdot\cdot\cdot \circ \pi_{l-1} \circ \pi_l$, is a projective
birational morphism isomorphic over $\roman{Reg}(X)$.
\endproclaim

\proclaim{Remark 7-2} \endproclaim

(i) We want to emphasize that we do NOT
require that $Y_{i-1} \subset X_{i-1}$, i.e., the center $Y_{i-1}$ be contained in the strict
transform
$X_{i-1}$ of $X = X_0$, or that it be contained in its singular locus,
i.e., $Y_{i-1} \subset \roman{Sing}(X_{i-1})$, as Hironaka or Bierstone-Milman does in their
formulation of embedded resolution of singularities.  Therefore, though the centers
$Y_{i-1}$ are smooth in the ambient varieties $W_{i-1}$ and permissible with respect to
$E_{i-1}$, their restrictions
$Y_{i-1}
\cap X_{i-1}$ to the strict transforms may not be smooth or reduced in general.  See Chapter 11
for some examples.

\vskip.1in

(ii) When $X$ is a hypersurface in $W$, i.e., $\dim X = \dim W - 1$, our algorithm
provides a sequence satisfying

\hskip.1in (i') the centers $Y_{i-1} \subset W_{i-1}$ of the blowups $\pi_i \hskip.1in
(i = 1, ... , l)$ are contained in $\roman{Sing}(X_{i-1})$,

\vskip.1in

which is stronger than condition (i) and coincides with the requirement that Hironaka or
Bierstone-Milman makes.

This is because, in the case of $X$ being a hypersurface, the weak transform coincides
with the strict transform and because our centers are contained in the maximum loci of
the invariant $w\text{-}\roman{ord}$ (by construction), which are necessarily sitting
inside of the singular loci of $X_{i-1} \hskip.1in (i = 1, ... , l)$.

\newpage

\demo{Proof}\enddemo We consider the following basic object $(W_0,(J_0,b),E_0)$ where
$$\left\{\aligned
W_0 &= W \\
J_0 &= {\Cal I}_X \hskip.1in (\text{the\ defining\ ideal\ of\ }X \text{\ in\ }W)\\
b &= 1 \\
E_0 &= \emptyset \\
\endaligned\right.$$
and take the sequence of transformations of basic objects, which represents resolution of
singularities of $(W_0,(J_0,b),E_0)$, constructed by the inductive algorithm of Theorem
5-1 and with the properties as specified in Corollary 5-6
$$(W_0,(J_0,b),E_0) \leftarrow \cdot\cdot\cdot \leftarrow (W_k,(J_k,b),E_k).$$
Observe that if for $i = 1, ... , l$ with $l \leq k$ the centers $Y_{i-1} \subset W_{i-1}$ do
not contain the strict transforms $X_{i-1}$, then $X_l$ is an irreducible component of
$\roman{Sing}(J_l,b)$.  Since $\roman{Sing}(J_k,b) = \emptyset$, we conclude that there exists
$1 \leq l \leq k - 1$ such that $Y_l$ contains $X_l$, while $Y_{i-1}$ does not contain
$X_{i-1}$ for $i = 1, ... , l$.  Moreover, since $Y_l \subset \roman{Sing}(J_l,b)$, we see that
$X_l$ is an irreducible component of $Y_l$.  Since $Y_l$ is smooth over $k$, a property which
is guaranteed by the inductive algorithm, we conclude that $X_l$ is smooth over $k$.

If we look at the sequence up to the $l$-th stage, it is immediate that it satisfies conditions
(ii) and (iii) (except for the claim that $\pi$ is isomorphic over $\roman{Reg}(X)$, which
follows once we check condition (i)). 

\vskip.1in

Condition (i) is a consequence of the process prescribed by the inductive algorithm of
Theorem 5-1.  In fact, let $p \in \roman{Reg}(X)$ be an arbitrary point of $\roman{Reg}(X)$. 
Let $l_p$ be the smallest number so that $p \in Y_{l_p}$.  The condition
$(\heartsuit')$ 
$$(\heartsuit')\hskip.1in
Y_{i-1} \subset \underline{\roman{Max}}\ t_{i-1} \subset \underline{\roman{Max}}\
w\text{-}\roman{ord}_{i-1} \text{\ if\ }w\text{-}\roman{ord}_{i-1} > 0 \text{\ for\ }i =
1, ... , k$$ implies that
$$\max\ w\text{-}\roman{ord}_{l_p} = 1 = w\text{-}\roman{ord}_{l_p}(p) \hskip.1in \&
\hskip.1in
\max\ t_{l_p} = t_{l_p}(p) = (1,0).$$
There exists an open neighborhood $p
\in U_p = W_{l_p}^{\lambda} \subset W_{l_p}$ such that $\roman{Reg}(X_{l_p}) \cap U_p = V(x_1,
\cdot\cdot\cdot, x_r)$ where
$x_1,
\cdot\cdot\cdot, x_r$ are regular parameters with $r = \dim W_{l_p} - \dim X_{l_p}$. 
According to the inductive algorithm described in Theorem 5-1, after $(r-1)$-repetitions
of $\bold{Case\ B}$, we reach a
$(\dim X_{l_p} + 1)$-dimensional basic object
$(\widetilde{W_{l_p}^{\lambda}},({\goth a}_{l_p}^{\lambda},b^{\lambda}),
\widetilde{E_{l_p}^{\lambda}}) = (V(x_1,
\cdot\cdot\cdot, x_{r-1}), ((x_r),1), \emptyset)$ where $R(1)(\underline{\roman{Max}}\
t_{l_p})
\cap W_{l_p}^{\lambda} = V(x_1, \cdot\cdot\cdot, x_r) = \roman{Reg}(X_{l_p}) \cap U_p$
and hence we are in $\bold{Case\ A}$.  This implies that $X_{l_p}$ is contained in the
center $Y_{l_P}$ and hence that $l_p = l$.  

Since $l_p$ is the smallest number so that $p
\in Y_{l_p}$ and since $p \in \roman{Reg}(X)$ is arbitrary, we conclude that the centers $Y_{i-1} \subset W_{i-1}$ of the blowups $\pi_i \hskip.1in (i = 1,
... , l)$ are over $\roman{Sing}(X) = X \setminus \roman{Reg}(X)$, verifying condition
(i). 

\vskip.1in

This completes the proof of Theorem 7-1.

\vskip.1in

\proclaim{Remark 7-3}\endproclaim

\vskip.1in

(i) Resolution of singularities of a basic object $(W, ({\Cal I}_X,1), \emptyset)$ is called
\linebreak
the ``$\bold{principalization}$" of the ideal ${\Cal I}_X$, since as a consequence we
obtain
$$\pi^{-1}{\Cal I}_X \cdot {\Cal O}_{W_k} = I(H_{r+1})^{a_{r+1}} \cdot\cdot\cdot
I(H_{r+k})^{a_{r+k}}$$
making the total transform of the ideal ``principal" (where actually $r = 0$ as \linebreak
$E_0 =
\{H_1, ... , H_r\} = \emptyset$).  (See the remark right after Main Theme 0-3
(Principalization of ideals for our restrictive use of the word ``principal".)

\vskip.1in
 
(ii) In the paper ``A course on constructive desingularization and equivariance", embedded
resolution of \it hypersurface \rm singularities is proved, starting with resolution of
singularities of a basic object $(W_0,(J_0,b),E_0)$ where
$$\left\{\aligned
W_0 &= W \\
J_0 &= {\Cal I}_X \\
b &= b_{\max} = \max\{\nu_p({\Cal I}_X);p \in W_0\} \\
E_0 &= \emptyset \\
\endaligned\right.$$ 
and then continuing with the descending induction on $b_{\max}$.  Embedded resolution of
\it non-hypersurface \rm singularities quotes the results of Hironaka and others (without
proof) which use the Hilbert-Samuel function.  It was a big discouragement for those of
us in the seminar who were hoping to have a self-contained course with complete proofs
for the entire picture of the process of resolution of (\it hypersurface \rm and
\it non-hypersurface \rm) singularities.  

But then after moments of thoughts, we the students realize, as the teachers Encinas and
Villamayor reveal to us \footnote"${}^1$"{Prof. Bierstone informed us in an informal way
that this idea of achieving resolution of singularities, \it hypersurface
\rm and
\it non-hypersurafce \rm, as a consequence of principalization imposing the
permissibility condition on the centers, could also be traced back to Hironaka.} in the
second paper ``A new theorem of desingularization over fields of characteristic zero",
that they have already told us ALL in the first paper, i.e., the inductive algorithm of
resolution of singularities of general basic objects applied to the principalization of
the ideal
${\Cal I}_X$ gives us reslotuion of singularities, \it hypersurface \rm and
\it non-hypersurface \rm, even $\bold{without}$ the use of Hilbert-Samuel function.  (Of
course, however, a price has to be paid in not being able to claim the stronger
properties on the centers in the process of resolution of singularities, as mentioned in
Remark 7-2.) 

\newpage

${}\hskip.6in$$\bold{CHAPTER\ 8.\ EQUIVARIANCE}$ \linebreak
${}\hskip1.8in$$\bold{AND}$ \linebreak
${}\hskip.6in$$\bold{RESOLUTION\ OF\ SINGULARITIES}$ \linebreak
$\bold{OVER\ BASE\ FIELDS\ (OF\ CHARACTERISTIC\ ZERO})$ \linebreak
${}\hskip1.2in$$\bold{WHICH\ ARE\ POSSIBLY}$ \linebreak
${}\hskip.8in$$\bold{NOT\ ALGEBRAICALLY\ CLOSED}$

\vskip.1in

In this chapter, we prove that the inductive algorithm for resolution of singularities of
(general) basic objects presented in Chapter 5 is equivariant under any ``action" (of a group)
and hence that all the centers are invariant under the action.  This implies, as an easy
corollary, that given a basic object $(W,(J,b),E)$ defined over a field $k$ which is
of characteristic zero but which may not be algebraically closed, the inductive algorithm
of resolution of singularities for the basic object $(W,(J,b),E) \times
\roman{Spec}\ \overline{k}$ is equivariant under the action of the Galois group
$\roman{Gal}(\overline{k}/k)$, all the centers are invariant under the action of
$\roman{Gal}(\overline{k}/k)$ and hence defined over
$k$ and that it induces resolution of singularities of $(W,(J,b),E)$ over $k$.

\proclaim{Definition 8-1 (``Action" on a basic object)} Let $(W,E = \{H_1, ... ,
H_r\})$ and $(W',E' = \{H_1', ... , H_{r'}'\})$ be pairs (cf. Definition 1-6).  An
isomorphism of pairs $\theta:(W,E) \overset{\sim}\to{\rightarrow} (W',E')$ is an isomorphism
$\theta:W \overset{\sim}\to{\rightarrow} W'$ as abstract varieties (not necessarily over the
base field $k$) such that $r = r'$ and that
$$\theta(H_i) = H_i' \text{\ for\ }i = 1, ... , r.$$  
Let $(W,(J,b),E)$ and $(W',(J',b'),E')$ be basic objects.  An isomorphism of basic objects
$\theta:(W,(J,b),E) \overset{\sim}\to{\rightarrow} (W',(J',b'),E')$ is an isomorphism of pairs
\linebreak
$\theta:(W,E) \overset{\sim}\to{\rightarrow} (W',E')$ such that $b' = b$ and that it induces an
isomorphism of ideals
$$J' = \theta_*(J) \subset \theta_*({\Cal O}_W) = {\Cal O}_{W'}.$$
An action on a pair $(W,E)$, by definition, is an isomorphism of pairs of
$(W,E)$ onto itself.
An action on a basic object $(W,(J,b),E)$, by definition, is an isomorphism of basic objects of
$(W,(J,b),E)$ onto itself.
\endproclaim

\proclaim{Remark 8-2}\endproclaim

(i) Recall that $H_i$ actually consists of smooth irreducible components \linebreak
$H_{i,1}, ... , H_{i,l_i}$ in our notation (cf. Note 1-7) and that so does $H_i'$ of
$H_{i,1}', ... , H_{i,{l_i}'}'$.  Therefore, when we state the condition 
$$\theta(H_i) = H_i' \text{\ for\ }i = 1, ... , r,$$ 
what we really mean is that for each $i = 1, \cdot\cdot\cdot, r$ we have $l_i = l_i'$ and
that there is a permutation of
$\{1,
\cdot\cdot\cdot, l_i\}$ (which we denote by the same letter $\theta$ by abuse of notation) with
$$\theta(H_{i,j}) = H_{i,\theta(j)}' \text{\ for\ }j = 1, ... , l_i.$$

(ii) As stated above, in order for $\theta$ to be an action on a pair, we do require not
only $E$ to be preserved as a whole but also each
$H_i$ to be preserved by $\theta$, fixing each index $i$.  Since our algorithm for resolution
of singularities of (monomial and hence all the general) basic objects depends on the indexing of
the divisors in
$E$ (cf. Definition 2-3, Proposition 2-5, Corollary 2-6, Definition-Proposition 4-5 and Theorem 5-1),
this requirement is necessary for us to claim that our algorithm for resolution of
singularities of (general) basic objects is equivariant under any action.  However, this
requirement makes little difference when we consider the equivariance of (embedded or non-embedded)
resolution of singulaities of a variety $X$, since the basic object of concern $(W,({\Cal I}_X,1),
\emptyset)$ that we start with (cf. Chapter 7) has empty boundary divisor $E = \emptyset$ and
since the indexing of the subsequent exceptional divisors are determined by the resolution
process itself.

(iii) We do NOT require in the definition of an action for an isomorphism $\theta$ to be over the base field
$k$.

(iv) Two non-isomorphic basic objects may define isomorphic general basic objects: Take
$(W,(J,b),E) = ({\Bbb A}^2 = \roman{Spec}\ k[x,y],((x),1),\emptyset)$ and $(W',(J',b'),E') =
({\Bbb A}^2,((x^2),2),\emptyset)$.  Then $(W,(J,b),E)$ and $(W',(J',b'),E')$ are non-isomorphic as
basic objects, though they define the same general basic object representing the same collection
${\goth C}$ of sequences of transformations and smooth morphisms of pairs with specified closed
subsets according to Remark 4-2 (ii) (cf. Definition 8-5).

\vskip.1in

\proclaim{Definition 8-3 (Equivariant sequence of basic objects)} Let
$$\align
(W, (J,b),E) = (W_0, (J_0,b), E_0) &\overset{\pi_1}\to{\leftarrow} (W_1, (J_1,b), E_1)
\overset{\pi_2}\to{\leftarrow} \cdot\cdot\cdot \\
(W_{i-1}, (J_{i-1},b), E_{i-1}) &\overset{\pi_i}\to{\leftarrow}
(W_i, (J_i,b), E_i) \\
\cdot\cdot\cdot \overset{\pi_{k-1}}\to{\leftarrow} (W_{k-1}, (J_{k-1},b),
E_{k-1}) &\overset{\pi_k}\to{\leftarrow} (W_k, (J_k,b), E_k) \\
\endalign$$
be a sequence of transformations and smooth morphisms of basic objects.  

We say that the sequence is
$\theta$-equivariant, given an action $\theta$ on the basic object $(W,(J,b),E)$, if inductively
for $i = 1, ... , k$ we
have a commutative diagram
$$\CD
(W_{i-1},(J_{i-1},b),E_{i-1}) \hskip.1in @. \overset{\pi_i}\to{\leftarrow} @.
\hskip.1in (W_i,(J_i,b),E_i)
\\ @V \theta VV @. @V \theta VV \\
(W_{i-1},(J_{i-1},b),E_{i-1}) \hskip.1in @. \overset{\pi_i}\to{\leftarrow} @.
\hskip.1in (W_i,(J_i,b),E_i)
\\
\endCD$$
i.e., the action $\theta:(W_{i-1},(J_{i-1},b),E_{i-1}) \overset{\sim}\to{\rightarrow}
(W_{i-1},(J_{i-1},b),E_{i-1})$ induced by
$\theta:(W_0,(J_0,b),E_0) \overset{\sim}\to{\rightarrow} (W_0,(J_0,b),E_0)$ lifts to
an action $\theta:(W_i,(J_i,b),E_i) \overset{\sim}\to{\rightarrow} (W_i,(J_i,b),E_i)$.

\vskip.1in

We note that if the sequence is $\theta$-equivariant, then
$$\left\{\aligned
&\theta(Y_{i-1}) = Y_{i-1} \\
&\text{whenever\ }\pi_i \text{\ is\ a\ transformation\ with\
center\ }Y_{i-1} \subset \roman{Sing}(J_{i-1},b) \subset W_{i-1}.\\
\endaligned\right.$$

We say that the sequence is equivariant if it is $\theta$-equivariant for \rm any \it action
$\theta$ on
$(W,(J,b),E)$. 
\endproclaim

\proclaim{Definition 8-4 (Equivariant resolution of singularities of a basic object)} Let
$$(W,(J,b),E) = (W_0,(J_0,b),E_0) \leftarrow \cdot\cdot\cdot \leftarrow (W_k,(J_k,b),E_k)$$
be a sequence of transformations representing resolution of singularities of the
basic object $(W,(J,b),E)$ (i.e., $\roman{Sing}(J_k,b) = \emptyset$).  

We say that the resolution of singularities is
$\theta$-equivariant, given an action $\theta$ on $(W,(J,b),E)$, if the sequence is
$\theta$-equivariant.

We say that the resolution of singularities is equivariant if the sequence is equivariant.
\endproclaim 

\proclaim{Definition 8-5 (``Action" on a general basic object)} Let $({\Cal F}_0,(W_0,E_0))$ be
a general basic object over $(F_0,(W_0,E_0))$ with a $d$-dimensional structure,
representing the collection
${\goth C}$ of sequences of transformations and smooth morphisms of pairs with specified closed
subsets.  Let $({\Cal F}_0',(W_0',E_0'))$ be another over $(F_0',(W_0',E_0'))$
with a $d' = d$-dimensional structure, representing the collection ${\goth C}'$.  An
isomorphism of general basic objects
$\theta:({\Cal F}_0,(W_0,E_0))
\overset{\sim}\to{\rightarrow} ({\Cal F}_0',(W_0',E_0'))$ is an isomorphism of pairs
$\theta:(W_0,E_0) \overset{\sim}\to{\rightarrow} (W_0',E_0')$ with $\theta(F_0) = F_0'$ which
satisfies the following condition: 

For each commutative diagram of sequences of transformations
and smooth morphisms of pairs with specified closed subsets
$$\CD
(F_0,(W_0,E_0)) @. \leftarrow @. \hskip.1in \cdot\cdot\cdot \hskip.1in @. \leftarrow @. (F_k,(W_k,E_k)) \\ 
@V \theta VV @. @. @. @V \theta VV \\
(F_0',(W_0',E_0')) @. \leftarrow @. \hskip.1in \cdot\cdot\cdot \hskip.1in @. \leftarrow @.
(F_k',(W_k',E_k')),
\\ 
\endCD$$
where the vertical arrows are all isomorphisms of pairs $\theta:(W_i,E_i)
\overset{\sim}\to{\rightarrow} (W_i',E_i')$ with $\theta(F_i) = F_i'$, induced by the
original isomorphism of pairs $\theta:(W_0,E_0) \overset{\sim}\to{\rightarrow} (W_0',E_0')$ with $\theta(F_0)
= F_0'$, the sequence 
$$(F_0,(W_0,E_0)) \leftarrow \cdot\cdot\cdot \leftarrow (F_k,(W_k,E_k))$$
is in the collection ${\goth C}$ if and only if the sequence
$$(F_0',(W_0',E_0')) \leftarrow \cdot\cdot\cdot \leftarrow (F_k',(W_k',E_k'))$$
is in the collection ${\goth C}'$.

(By abuse of notation, we could express the last condition as requiring $\theta({\goth C}) = {\goth C}'$.) 

\vskip.1in

An action on a general basic object $({\Cal F}_0,(W_0,E_0))$ is an isomorphism of general basic
objects of $({\Cal F}_0,(W_0,E_0))$ onto itself.

\vskip.1in

Or equivalently, we can define an action on a general basic object $({\Cal F}_0,(W_0,E_0))$ in the
following way.

Let $({\Cal F}_0,(W_0,E_0))$ be
a general basic object over $(F_0,(W_0,E_0))$, with a $d$-dimensional structure given by the
charts $\{(\widetilde{W_0^{\lambda}},({\goth
a}_0^{\lambda},b^{\lambda}),\widetilde{E_0^{\lambda}})\}$ of basic objects of dimension $d$,
with the collection ${\goth C}$ of transformations and smooth morphisms represented by $({\Cal
F}_0,(W_0,E_0))$. 

Let $\theta$ be an action on the pair $(W_0,E_0)$.  

Remark that the new charts $\{({}^\theta\widetilde{W_0^{\lambda}},({}^{\theta}{\goth
a}_0^{\lambda},{}^{\theta}b^{\lambda}),{}^{\theta}\widetilde{E_0^{\lambda}})\}$, where
$$\left\{\aligned
{}^{\theta}\widetilde{W_0^{\lambda}} &= \theta(\widetilde{W_0^{\lambda}}) \\
{}^{\theta}{\goth a}_0^{\lambda} &= \theta_*({\goth a}_0^{\lambda}) \\
{}^{\theta}b^{\lambda} &= b^{\lambda} \\
{}^{\theta}\widetilde{E_0^{\lambda}} &= \theta(\widetilde{E_0^{\lambda}}), \\
\endaligned\right.$$
defines a general basic object $({}^{\theta}{\Cal F}_0,({}^{\theta}W_0 = W_0,{}^{\theta}E_0
= E_0))$ over
$({}^{\theta}F_0 =
\theta(F_0),(W_0,E_0))$ with the collection ${}^{\theta}{\goth C}$ of transformations and
restrictions represented by $({}^{\theta}{\Cal F}_0,(W_0,E_0))$. 

We say that $\theta$ is an action on the general basic object $({\Cal F}_0,(W_0,E_0))$ if
$$F_0 = {}^{\theta}F_0 \hskip.1in \& \hskip.1in {\goth C} = {}^{\theta}{\goth C}.$$
That is to say, $\theta$ is an action if $F_0 = {}^{\theta}F_0$ and if a necessary and
sufficient condition for a sequence of pairs with specified basic object
$$(F_0,(W_0,E_0)) \overset{\pi_1}\to{\leftarrow} (F_1,(W_1,E_1)) \overset{\pi_2}\to{\leftarrow}
\cdot\cdot\cdot
\overset{\pi_{k-1}}\to{\leftarrow} (F_{k-1},(W_{k-1},E_{k-1})) \overset{\pi_k}\to{\leftarrow}
(F_k,(W_k,E_k))$$ 
to be in the collection
${\goth C}$ is for its
$\theta$-counterpart (the sequence below which makes an obvious commutative diagram with the sequence above
having the vertical arrows being $\theta$)
$$\CD
({}^{\theta}F_0,({}^{\theta}W_0,{}^{\theta}E_0)) @.\overset{{}^{\theta}\pi_1}\to{\leftarrow}
({}^{\theta}F_1,({}^{\theta}W_1,{}^{\theta}E_1))
\overset{{}^{\theta}\pi_2}\to{\leftarrow}
\cdot\cdot\cdot
\overset{{}^{\theta}\pi_{k-1}}\to{\leftarrow}
({}^{\theta}F_{k-1},({}^{\theta}W_{k-1},{}^{\theta}E_{k-1}))
\overset{{}^{\theta}\pi_k}\to{\leftarrow} ({}^{\theta}F_k,({}^{\theta}W_k,{}^{\theta}E_k))\\
@| @. \\
(F_0,(W_0,E_0)) @. @. \\
\endCD$$
to be in the collection ${\goth C}$.

\vskip.1in

We emphasize that, in order for $\theta$ to be an action on the general basic object, we are NOT requiring
the charts
$\{(\widetilde{W_0^{\lambda}},({\goth a}_0^{\lambda},b^{\lambda}),\widetilde{E_0^{\lambda}})\}$ to coincide
with the new charts $\{({}^{\theta}\widetilde{W_0^{\lambda}},({}^{\theta}{\goth
a}_0^{\lambda},{}^{\theta}b^{\lambda}),{}^{\theta}\widetilde{E_0^{\lambda}})\}$ but
that we are requiring the collection ${\goth C}$ to coincide with the new collection
${}^{\theta}{\goth C}$. 
\endproclaim

\proclaim{Definition 8-6 (Equivariant sequence of general basic objects)} Let
$$({\Cal F}_0,(W_0,E_0)) \leftarrow \cdot\cdot\cdot \leftarrow ({\Cal F}_k,(W_k,E_k))$$
be a sequence of transformations and smooth morphisms of general basic objects.  (Remark that the above is
nothing but a notational convention (cf. Note 4-3) expressing a sequence of transformations and
smooth morphisms
$$(F_0,(W_0,E_0)) \leftarrow \cdot\cdot\cdot \leftarrow (F_k,(W_k,E_k))$$
in the collection ${\goth C}$ represented by $({\Cal F}_0,(W_0,E_0))$.  We say that the sequence is
$\theta$-equivariant, given an action $\theta$ on the general basic object $({\Cal
F}_0,(W_0,E_0))$, if inductively for $i = 1, ... , k$
we have a commutative diagram
$$\CD
({\Cal F}_{i-1},(W_{i-1},E_{i-1})) \hskip.1in @. \overset{\pi_i}\to{\leftarrow} @.
\hskip.1in ({\Cal F}_i,(W_i,E_i)) \\ 
@V \theta VV @. @V \theta VV \\
({\Cal F}_{i-1},(W_{i-1},E_{i-1})) \hskip.1in @. \overset{\pi_i}\to{\leftarrow} @.
\hskip.1in ({\Cal F}_i,(W_i,E_i))
\\
\endCD$$
i.e., the action on
$({\Cal F}_{i-1},(W_{i-1},E_{i-1}))$ induced by $\theta$ on $({\Cal F}_0,(W_0,E_0))$ lifts to the action on
$({\Cal F}_i,(W_i,E_i))$.  (Remark that it is simply equivalent to requiring that the action on the pair
$(W_{i-1},E_{i-1})$ induced by $\theta$ on $(W_0,E_0)$ lifts to the action on
$(W_i,E_i)$ and $\theta(F_i) = F_i$.)  

\vskip.1in

We note that if the sequence is $\theta$-equivariant, then
$$\left\{\aligned
&\theta(Y_{i-1}) = Y_{i-1}^{\theta} = Y_{i-1} \\
&\text{whenever\ }\pi_i \text{\ is\ a\ transformation\ with\ center\ }Y_{i-1} \subset F_{i-1}
\subset W_{i-1}.\\
\endaligned\right.$$

\vskip.1in

We say that the sequence is equivariant if it is $\theta$-equivariant for \rm any \it action $\theta$ on the
general basic object $({\Cal F}_0,(W_0,E_0))$. 
\endproclaim

\vskip.1in

\proclaim{Definition 8-7 (Equivariant resolution of singularities of a general basic object)} Let
$$({\Cal F}_0,(W_0,E_0)) \leftarrow \cdot\cdot\cdot \leftarrow ({\Cal F}_k,(W_k,E_k))$$
be a sequence of transformations representing resolution of singularities of the
general basic object $({\Cal F}_0,(W_0,E_0))$ (i.e., $F_k = \emptyset$).  

We say that the resolution of singularities is
$\theta$-equivariant, given an action $\theta$ on $({\Cal F}_0,(W_0,E_0))$, if the sequence is
$\theta$-equivariant.

We say that the resolution of singularities is equivariant if the sequence is equivariant.
\endproclaim 

\proclaim{Remark 8-8}\endproclaim

There is some confusion concerning the definitions of an action and an equivariant sequence
of general basic objects in the original paper ``A course on constructive desingularization
and equivariance" by Encinas and Villamayor.  

In Definition 6.20 (of the paper) they define:

We say an automorphism $\theta:W_0 \overset{\sim}\to{\rightarrow} W_0$ acts on the general
basic object $({\Cal F}_0,(W_0,E_0))$ if:

\ \ (a) $\theta$ acts on the pair $(W_0,E_0)$ and $\theta(F_0) = F_0$, and

\ \ (b) for any sequence in ${\goth C}$
$$(F_0,(W_0,E_0)) \leftarrow \cdot\cdot\cdot \leftarrow (F_k,(W_k,E_k))$$
which is $\theta$-equivariant in the sense that
$$\left\{\aligned
&\theta(Y_{i-1}) = Y_{i-1} \\
&\text{whenever\ }\pi_i \text{\ is\ a\ transformation\ with\
center\ }Y_{i-1} \subset F_{i-1} \subset W_{i-1},\\
\endaligned\right.$$
we have
$$\theta(F_i) = F_i \text{\ for\ }i = 0, ... , k.$$

Their definition is clearly different from ours.  With their definition of an ``action" on a
general basic object, one has trouble, e.g., in proving $\roman{ord}_0(x_0) = \roman{ord}_0(\theta(x_0))$
for an action on a general basic object $({\Cal F}_0,(W_0,E_0))$ and $x_0 \in F_0 \subset W_0$.  The
sequence we construct in Hironaka's trick is NOT
$\theta$-equivariant (in their sense as above) and their definition of the action does not provide any
information.  They tried to compensate for this calamity of their definition by looking at $x_0$ and
$\theta(x_0)$ simultaneously in Hironaka's trick (in a vain attempt to make the sequence equivariant) and
claiming in their proof of Proposition 7.4 that ``$\theta
\times Id$ acts on $(W_1,E_1)$ interchanging the components of $x_0 \times {\Bbb A}^1$ and
$\theta(x_0) \times {\Bbb A}^1$ and interchanging $(x_0,0)$ and $(\theta(x_0),0)$", which
is sheer nonsense.

\vskip.1in

It is very clear, however, from the context of the paper that what Encinas and Villamayor really
mean is the definition(s) that we give here in these notes and that the discrepancies mentioned
above should be considered mere ``typos" in the paper. 

\vskip.1in

\proclaim{Lemma 8-9 (Invariance of key invariants for basic objects under action)} Let
$$\align
(W, (J,b),E) = (W_0, (J_0,b), E_0) &\overset{\pi_1}\to{\leftarrow} (W_1, (J_1,b), E_1)
\overset{\pi_2}\to{\leftarrow} \cdot\cdot\cdot \\
(W_{i-1}, (J_{i-1},b), E_{i-1}) &\overset{\pi_i}\to{\leftarrow}
(W_i, (J_i,b), E_i) \\
\cdot\cdot\cdot \overset{\pi_{k-1}}\to{\leftarrow} (W_{k-1}, (J_{k-1},b),
E_{k-1}) &\overset{\pi_k}\to{\leftarrow} (W_k, (J_k,b), E_k) \\
\endalign$$
be a sequence of transformations and smooth morphisms of basic objects.

Let $\theta$ be a action on the basic object $(W,(J,b),E)$.

Let 
$$\align
(W_0,(J_0,b),E_0) = ({}^{\theta}W_0, ({}^{\theta}J_0,b), {}^{\theta}E_0)
&\overset{{}^{\theta}\pi_1}\to{\leftarrow} ({}^{\theta}W_1, ({}^{\theta}J_1,b),
{}^{\theta}E_1)
\overset{{}^{\theta}\pi_2}\to{\leftarrow} \cdot\cdot\cdot \\
({}^{\theta}W_{i-1}, ({}^{\theta}J_{i-1},b), {}^{\theta}E_{i-1})
&\overset{{}^{\theta}\pi_i}\to{\leftarrow} ({}^{\theta}W_i, ({}^{\theta}J_i,b),
{}^{\theta}E_i)
\\
\cdot\cdot\cdot \overset{{}^{\theta}\pi_{k-1}}\to{\leftarrow} ({}^{\theta}W_{k-1},
({}^{\theta}J_{k-1},b), {}^{\theta}E_{k-1}) &\overset{{}^{\theta}\pi_k}\to{\leftarrow}
({}^{\theta}W_k, ({}^{\theta}J_k,b), {}^{\theta}E_k)
\\
\endalign$$
be the $\theta$-counterpart of the sequence. 

Then for $\xi_k \in \roman{Sing}(J_k,b)$ and $\theta(\xi_k) = {}^{\theta}\xi_k \in
\theta(\roman{Sing}(J_k,b)) = \roman{Sing}({}^{\theta}J_k,b)$ we have
$$\left\{\aligned
\roman{ord}_k(\xi_k) & = {}^{\theta}\roman{ord}_k({}^{\theta}\xi_k), \\
w\text{-}\roman{ord}_k(\xi_k) & = {}^{\theta}w\text{-}\roman{ord}_k({}^{\theta}\xi_k)
\hskip.1in
\&
\hskip.1in \theta(\underline{\roman{Max}}\ w\text{-}\roman{ord}_k) =
\underline{\roman{Max}}\ {}^{\theta}w\text{-}\roman{ord}_k,\\
\Gamma_k(\xi_k) & = {}^{\theta}\Gamma_k({}^{\theta}\xi_k) \hskip.36in \& \hskip.1in
\theta(\underline{\roman{Max}}\ \Gamma_k) = \underline{\roman{Max}}\ {}^{\theta}\Gamma_k\\
&\text{if\ }\max\ w\text{-}\roman{ord}_k = 0
\text{\ and\ hence\ we\ may\ regard\ }\\ &(W_k,(J_k,b),E_k) \text{\ and\
}({}^{\theta}W_k,({}^{\theta}J_k,b),{}^{\theta}E_k)
\text{\ as\ monomial\ basic\ objects\ }\\
&\text{in\ open\ neighborhoods\ of\ }\roman{Sing}(J_k,b) \text{\ and\
}\theta(\roman{Sing}(J_k,b)) = \roman{Sing}({}^{\theta}J_k,b) \text{\ respectively,}\\
t_k(\xi_k) &= {}^{\theta}t_k({}^{\theta}\xi_k) \hskip.38in \& \hskip.1in
\theta(\underline{\roman{Max}}\ t_k) = \underline{\roman{Max}}\ {}^{\theta}t_k\\ 
&\text{if\ the\
sequence\ satisfies\ condition\ }(\heartsuit),\\
\endaligned\right.$$
where ${}^{\theta}\roman{ord}_k, {}^{\theta}w\text{-}\roman{ord}_k, {}^{\theta}\Gamma_k,
{}^{\theta}t_k$ are the
$ord, w\text{-}ord, \Gamma, t$-invariants on
$({}^{\theta}W_k,({}^{\theta}J_k,b),{}^{\theta}E_k)$ (defined with respect to the
$\theta$-counterpart of the sequence).

If the sequence is $\theta$-equivariant, then $\theta$ lifts to an action on $(W_k,(J_k,b),E_k)$
and we have
$$\left\{\aligned
\roman{ord}_k(\xi_k) & = \roman{ord}_k({}^{\theta}\xi_k), \\
w\text{-}\roman{ord}_k(\xi_k) & = w\text{-}\roman{ord}_k({}^{\theta}\xi_k) \hskip.1in \&
\hskip.1in \theta(\underline{\roman{Max}}\ w\text{-}\roman{ord}_k) =
\underline{\roman{Max}}\ w\text{-}\roman{ord}_k,\\
\Gamma_k(\xi_k) & = \Gamma_k({}^{\theta}\xi_k) \hskip.36in \& \hskip.1in
\theta(\underline{\roman{Max}}\ \Gamma_k) = \underline{\roman{Max}}\ \Gamma_k\\
&\text{if\ }\max\ w\text{-}\roman{ord}_k = 0
\text{\ and\ hence\ we\ may\ regard\ }\\ &(W_k,(J_k,b),E_k) \text{\ as\ a\ monomial\ basic\ object\
}\\ 
&\text{in\ an\ open\ neighborhood\ of\ }\roman{Sing}(J_k,b)\\
t_k(\xi_k) &= t_k({}^{\theta}\xi_k) \hskip.38in \& \hskip.1in
\theta(\underline{\roman{Max}}\ t_k) = \underline{\roman{Max}}\ t_k\\ 
&\text{if\ the\
sequence\ satisfies\ condition\ }(\heartsuit).\\
\endaligned\right.$$
\endproclaim

\demo{Proof}\enddemo The proof is obvious from the definition.  Note that stability (equivariance)
of the key invariants under isomorphism of basic objects can be stated and proved in a similar
manner.
   
\proclaim{Lemma 8-10 (Invariance of key invariants for general basic objects under
action)} Let
$$({\Cal F}_0,(W_0,E_0)) \leftarrow \cdot\cdot\cdot \leftarrow ({\Cal F}_k,(W_k,E_k))$$
be a sequence of transformations and smooth morphisms of general basic objects.

Let $\theta$ be a action on the general basic object $({\Cal F}_0,(W_0,E_0))$ with a
$d$-dimensional structure.

Let
$$({}^{\theta}{\Cal F}_0,({}^{\theta}W_0,{}^{\theta}E_0)) \leftarrow \cdot\cdot\cdot
\leftarrow ({}^{\theta}{\Cal F}_k,({}^{\theta}W_k,{}^{\theta}E_k))$$
be the $\theta$-counterpart of the sequence.

Then for $\xi_k \in F_k$ and $\theta(\xi_k) = {}^{\theta}\xi_k \in
\theta(F_k) = {}^{\theta}F_k$ we have
$$\left\{\aligned
\roman{ord}_k(\xi_k) & = {}^{\theta}\roman{ord}_k({}^{\theta}\xi_k), \\
w\text{-}\roman{ord}_k(\xi_k) & = {}^{\theta}w\text{-}\roman{ord}_k({}^{\theta}\xi_k)
\hskip.1in
\&
\hskip.1in \theta(\underline{\roman{Max}}\ w\text{-}\roman{ord}_k) =
\underline{\roman{Max}}\ {}^{\theta}w\text{-}\roman{ord}_k,\\
\Gamma_k(\xi_k) & = {}^{\theta}\Gamma_k({}^{\theta}\xi_k) \hskip.36in \& \hskip.1in
\theta(\underline{\roman{Max}}\ \Gamma_k) = \underline{\roman{Max}}\ {}^{\theta}\Gamma_k\\
&\text{if\ }\max\ w\text{-}\roman{ord}_k = 0\\ 
t_k(\xi_k) &= {}^{\theta}t_k({}^{\theta}\xi_k) \hskip.38in \& \hskip.1in
\theta(\underline{\roman{Max}}\ t_k) = \underline{\roman{Max}}\ {}^{\theta}t_k\\ 
&\text{if\ the\
sequence\ satisfies\ condition\ }(\heartsuit),\\
\endaligned\right.$$
where ${}^{\theta}ord_k, {}^{\theta}w\text{-}ord_k, {}^{\theta}\Gamma_k, {}^{\theta}t_k$ are
the
$ord, w\text{-}ord, \Gamma, t$-invariants on $({}^{\theta}{\Cal
F}_k,({}^{\theta}W_k,{}^{\theta}E_k))$ (defined with respect to the
$\theta$-counterpart of the sequence).

If the sequence is $\theta$-equivariant, then $\theta$ lifts to an action on $({\Cal
F}_k,(W_k,E_k))$ with a $\dim W_k - \dim W_0 + d$-dimensional structure and we have
$$\left\{\aligned
\roman{ord}_k(\xi_k) & = \roman{ord}_k({}^{\theta}\xi_k), \\
w\text{-}\roman{ord}_k(\xi_k) & = w\text{-}\roman{ord}_k({}^{\theta}\xi_k) \hskip.1in \&
\hskip.1in \theta(\underline{\roman{Max}}\ w\text{-}ord_k) = \underline{\roman{Max}}\
w\text{-}ord_k,\\
\Gamma_k(\xi_k) & = \Gamma_k({}^{\theta}\xi_k) \hskip.36in \& \hskip.1in
\theta(\underline{\roman{Max}}\ \Gamma_k) = \underline{\roman{Max}}\ \Gamma_k\\
&\text{if\ }\max\ w\text{-}\roman{ord}_k = 0\\ 
t_k(\xi_k) &= t_k({}^{\theta}\xi_k) \hskip.38in \& \hskip.1in
\theta(\underline{\roman{Max}}\ t_k) = \underline{\roman{Max}}\ t_k\\ 
&\text{if\ the\
sequence\ satisfies\ condition\ }(\heartsuit).\\
\endaligned\right.$$
\endproclaim

\demo{Proof}\enddemo If $\xi_k \in F_k \cap W_k^{\lambda} = \roman{Sing}({\goth
a}_k^{\lambda},b^{\lambda})$ is a point of a chart $(\widetilde{W_k^{\lambda}},({\goth
a}_k^{\lambda},b^{\lambda}),\widetilde{E_k^{\lambda}})$, then ${}^{\theta}\xi_k 
\in {}^{\theta}F_k \cap
{}^{\theta}(W_k^{\lambda}) =
\roman{Sing}({}^{\theta}{\goth a}_k^{\lambda},{}^{\theta}b^{\lambda})$ is a point of the
chart 
$({}^{\theta}\widetilde{W_k^{\lambda}},({}^{\theta}{\goth
a}_k^{\lambda},{}^{\theta}b^{\lambda}),{}^{\theta}\widetilde{E_k^{\lambda}})$.  Therefore,
the above assertions are easy consequences of Definition-Proposition 4-5 and Lemma 8-9. 
Note that stability (equivariance) of the key invariants under isomorphism of general basic
objects can be stated and proved in a similar manner.

\vskip.1in

\proclaim{Proposition 8-11 (Equivariance of the inductive algorithm)} Let
$$({\Cal F}_0,(W_0,E_0)) \leftarrow \cdot\cdot\cdot \leftarrow ({\Cal F}_k,(W_k,E_k))$$
be the sequence of transformations of general basic objects, obtained via the inductive
algorithm of Theorem 5-1, representing resolution of singularities of a general basic object
$({\Cal F}_0,(W_0,E_0))$ with a $d$-dimensional structure.

Let $\theta$ be an action on the general basic object $({\Cal F}_0,(W_0,E_0))$ with a
$d$-dimensional structure.

Then the sequence is $\theta$-equivariant.

In particular, the sequence is equivaraint.
\endproclaim

\demo{Proof}\enddemo We have only to prove that, having a $\theta$-equivariant sequence,
$$({\Cal F}_0,(W_0,E_0)) \leftarrow \cdot\cdot\cdot \leftarrow ({\Cal F}_k,W_k,E_k)),$$
the extension given by the process of the inductive algorithm of Theorem 5-1 is also
$\theta$-equivariant.

\vskip.1in

$\underline{\bold{P1}}$: The sequence already represents
resolution of singularities and there is no need for the extension.  The original sequence
is $\theta$-equivariant by assumption. 

$\underline{\bold{P2}}$: The extension of the sequence of
transformations we create via Corollary 4-9 is $\theta$-equivariant, since the centers
for the transformations 
$$({\Cal F}_i,(W_i,E_i)) \leftarrow ({\Cal
F}_{i+1},(W_{i+1},E_{i+1})) \text{\ for\ } i \geq k$$ 
are
$\theta$-invariant., i.e.,
$\theta(\underline{\roman{Max}}\
\Gamma_i) =
\underline{\roman{Max}}\ \Gamma_i$ by Lemma 8-9.

\vskip.1in

$\underline{\bold{P3}}$: We deal with the case of possibility $\underline{\bold{P3}}$ in
the following.

\vskip.1in

$\bold{Case\ A}$: We have $\theta(R(1)(\underline{\roman{Max}}\ t_k)) =
R(1)(\underline{\roman{Max}}\ t_k)$, since $\theta(\underline{\roman{Max}}\ t_k) =
\underline{\roman{Max}}\ t_k$ by Lemma 8-9.  

Thus the extended sequence by adding $({\Cal
F}_k,(W_k,E_k)) \leftarrow ({\Cal F}_{k+1},(W_{k+1},E_{k+1}))$ is $\theta$-equivariant.

$\bold{Case\ B}$: Note first that $\theta$ induces an action on the general basic object
$({\Cal F}_k,(W_k,E_k))$ (since the sequence is $\theta$-equivariant) and that by the
property
$\theta(\underline{\roman{Max}}\ t_k) =
\underline{\roman{Max}}\ {}^{\theta}t_k)$ of Lemma 8-10 (not only referring to the case of
the index
$k$ but also to the further extension) $\theta$ induces an action on the general basic
object
$({\Cal G}_k,(W_k,E_k''))$ via Lemma 5-4.  (Recall that the properties $(\alpha)$ and
$(\beta)$ provide a characterization of the general basic object $({\Cal
G}_k,(W_k,E_k''))$.  See the conclusion of the proof of Theorem 5-1 for the assertions in
$\bold{Case\ B}$ under possibility $\underline{\bold{P3}}$.)  The sequence of
transformations, constructed via the inductive algorithm and representing resolution of
singularities of
$({\Cal G}_k,(W_k,E_k''))$, is
$\theta$-equivariant by induction on the dimension of the structure $d$.  (We
leave the proof of $\theta$-equivariance in the case $d = 1$ to the reader as an
exercise.)  Therefore, the extended sequence adding
$$({\Cal F}_k,(W_k,E_k)) \leftarrow \cdot\cdot\cdot \leftarrow ({\Cal F}_{k+N},(W_{k+N},E_{k+N}))$$
is also $\theta$-equivariant.

\vskip.1in

This completes the proof of Proposition 8-11.

\vskip.1in

\proclaim{Corollary 8-12 (Equivaraint resolution of singularities of a general basic
object)} Let
$({\Cal F}_0,(W_0,E_0))$ be a general basic object with a $d$-dimensional structure.

Then there exists equivariant resolution of singularities of $({\Cal F}_0,(W_0,E_0))$
$$({\Cal F}_0,(W_0,E_0)) \leftarrow \cdot\cdot\cdot \leftarrow ({\Cal F}_k,(W_k,E_k))$$
satisfying condition $(\heartsuit')$
$$(\heartsuit')\ Y_{i-1} \subset \underline{\roman{Max}}\ t_{i-1} \subset \underline{\roman{Max}}\
w\text{-}\roman{ord}_{i-1} \text{\ if\ }\max\ w\text{-}\roman{ord}_{i-1} > 0 \text{\ for\ }i
= 1, ... , k.$$
\endproclaim

\demo{Proof}\enddemo We only need to check that the sequence given by the inductive
algorithm to represent resolution of singularities, which satisfies condition
$(\heartsuit')$ by construction, is equivariant.  This is exactly the content of
Proposition 8-11.

\proclaim{Corollary 8-13 (Equivaraint resolution of singularities of a basic object)} Let
$(W_0,(J_0,b),E_0)$ be a basic object.  Then there exists equivariant resolution of singularities
of $(W_0,(J_0,b),E_0)$
$$(W_0,(J_0,b),E_0) \leftarrow \cdot\cdot\cdot \leftarrow (W_k,(J_k,b),E_k)$$
satisfying condition $(\heartsuit')$
$$(\heartsuit')\ Y_{i-1} \subset \underline{\roman{Max}}\ t_{i-1} \subset \underline{\roman{Max}}\
w\text{-}\roman{ord}_{i-1} \text{\ if\ }\max\ w\text{-}\roman{ord}_{i-1} > 0 \text{\ for\ }i
= 1,
... , k.$$
\endproclaim

\demo{Proof}\enddemo This is a direct consequence of Corollary 8-12 and Remark 4-2 (ii).

\proclaim{Corollary 8-14 (Resolution of singularities of a general basic
object (resp. basic object) over any field
$\bold{k}$ of characteristic zero)} Let $({\Cal F}_0,(W_0,E_0))$ (resp.
$(W_0,(J_0,b),E_0)$) be a general basic object with a $d$-dimensional structure (resp. basic object)
defined over a field
$k$ which is of characteristic zero but may not be algebraically closed.  Then there
exists resolution of singularities, satisfying condition
$(\heartsuit')$,
$$\align
&({\Cal F}_0,(W_0,E_0)) \leftarrow \cdot\cdot\cdot \leftarrow ({\Cal F}_k,(W_k,E_k))\\
(\text{resp.\ }& (W_0,(J_0,b),E_0) \leftarrow \cdot\cdot\cdot \leftarrow (W_k,(J_k,b),E_k)) \\
\endalign$$
which is defined over $k$.
\endproclaim

\demo{Proof}\enddemo Firstly remark that a general basic object $({\Cal F}_0,(W_0,E_0))$ is
defined over $k$ with a $d$-dimensional structure if, by definition, it has an open
covering $\{W_0^{\lambda}\}_{\lambda \in \Lambda}$ with charts
$\{\widetilde{W_0^{\lambda}},({\goth
a}_0^{\lambda},b^{\lambda}),\widetilde{E_0^{\lambda}})\}_{\lambda \in \Lambda}$ which are
defined over $k$ and that by $({\Cal F}_0,(W_0,E_0)) \times \roman{Spec}\ \overline{k}$ we
mean the general basic object having the open covering $\{W_0^{\lambda} \times
\roman{Spec}\ \overline{k}\}_{\lambda \in \Lambda}$ with charts
$\{\widetilde{W_0^{\lambda}},({\goth
a}_0^{\lambda},b^{\lambda}),\widetilde{E_0^{\lambda}}) \times \roman{Spec}\
\overline{k}\}_{\lambda \in \Lambda}$.  It should be warned, however, that the collection
(of sequences of transformations and smooth morphisms of pairs with specified closed
subsets) represented by $({\Cal F}_0,(W_0,E_0)) \times \roman{Spec}\ \overline{k}$ contains
more than those obtained by taking the Cartesian product of the sequences defined over $k$
in the collection represented by the original $({\Cal F}_0,(W_0,E_0))$ with $\roman{Spec}\
\overline{k}$.

Secondly note that, since we have been
assuming so far that the ambient space
$W_0$ to be irreducible, we need to generalize the theory to the case where $W_0$ may be
reducible and of pure dimension, as $W_0
\times \roman{Spec}\ \overline{k}$ may be.  This generalization can be made without any change in
the arguments.  

Thirdly note that the Galois group $\roman{Gal}(\overline{k}/k)$ acts on
the general basic object $({\Cal F}_0,(W_0,E_0)) \times
\roman{Spec}\ \overline{k}$ (resp. $(W_0,(J_0,b),E_0) \times \roman{Spec}\
\overline{k}$) in the sense of Definition 8-5 \linebreak
(resp. Definition 8-1).

\vskip.1in

We construct a sequence of
transformations representing resolution of singularities of $({\Cal F}_0,(W_0,E_0)) \times
\roman{Spec}\ \overline{k}$ (resp. $(W_0,(J_0,b),E_0) \times \roman{Spec}\
\overline{k}$) via the inductive algorithm of Theorem 5-1.  Since the
sequence is equivariant under
the action of the Galois group $\roman{Gal}(\overline{k}/k)$ by Proposition 8-11, the
centers
$\overline{Y_{i-1}}$ are defined over $k$, i.e., $\overline{Y_{i-1}} = Y_{i-1} \times
\roman{Spec}\ \overline{k}$ for some closed subscheme $Y_{i-1} \subset W_{i-1}$ defined
over $k$ and the sequence of transformations with centers
$Y_{i-1}$ provides resolution of singularities of $({\Cal F}_0,(W_0,E_0))$ (resp.
$(W_0,(J_0,b),E_0)$ over $k$.

\proclaim{Corollary 8-15 (Embedded resolution of singularities over any field
$\bold{k}$ of characteristic zero)} Let $X \subset W$ be a variety, embedded as a closed subscheme
(defined over $k$) of another variety $W$ smooth over a field $k$ which is of
characteristic zero but may not be algebraically closed.  Then there exists a sequence
defined over $k$ of blowups representing embedded resolution of singularities of $X \subset
W$, which is equivariant in the sense that for any automorphism $\theta:W
\overset{\sim}\to{\rightarrow} W$ with $\theta(X) = X$, we have all the centers of the
blowups being $\theta$-invariant.
\endproclaim

\demo{Proof}\enddemo This is an easy consequnce of Corollary 8-14 combined with our
construction of embedded resolution of singularities presented in Chapter 7.

\vskip.1in

\proclaim{Remark 8-16 (Equivariance of the sequence representing resolution of
singularities in Corollary 8-14 or Corollary 8-15)}\endproclaim

Let $\theta$ be an action on the general basic object $({\Cal F}_0,(W_0,E_0))$
with a $d$-dimensional structure (defined over $k$).  We claim that the sequence
constructed via the inductive algorithm of Theorem 5-1 is $\theta$-equivarinat, no
matter whether
$\theta$ is over
$k$ or not.  (Recall that
$\theta$ is an isomorphism as abstract varieties, satisfying certain conditions as described
in Definition 8-5, and that it is not necessarily over the base field $k$.) 

If $\theta:({\Cal F}_0,(W_0,E_0)) \overset{\sim}\to{\rightarrow} ({\Cal F}_0,(W_0,E_0))$ is
over $k$, then it extends to an action $\theta \times Id:({\Cal F}_0,(W_0,E_0)) \times
\roman{Spec}\ \overline{k} \overset{\sim}\to{\rightarrow} ({\Cal F}_0,(W_0,E_0)) \times
\roman{Spec}\ \overline{k}$.  Therefore, equivariance of $\theta$, i.e., $\theta(Y_{i-1}) =
Y_{i-1}$ follows that of $\theta \times Id$.

However, if $\theta:({\Cal F}_0,(W_0,E_0)) \overset{\sim}\to{\rightarrow} ({\Cal
F}_0,(W_0,E_0))$ is not over $k$, there is no obvious way that $\theta$ extends to an
action on $({\Cal F}_0,(W_0,E_0)) \times
\roman{Spec}$.  In other words, there is no obvious way to reduce the equivariance over
$k$ to that over $\overline{k}$ taking the Cartesion product $\times_{\roman{Spec}\
k}\roman{Spec}\ \overline{k}$. 

We just remark that the entire theory up to Chapter 7 and Chapter 8 can be developed over
any field of characteristic zero, without assuming that $k$ is algebraically closed and
that $\theta$-equivariance of the inductive algorithm, given any action whether it is
defined over $k$ or not, goes verbatim as in the proof of Proposition 8-11.  The essential
point is that the invariants $\roman{ord}$, $w\text{-}\roman{ord}$, $t$, and $\Gamma$
(of the original general basic object and of the auxiliary general basic objects
appearing in the inductive process (See the characterization via properties $(\alpha)$ and
$(\beta)$.)), which determine the inductive algorithm, depend only on the structure as
abstract varieties and not on the structure over $k$ and hence are preserved under any
isomorphism as abstract varieties whether it is over $k$ or not.  (We note that the only
place where we use the assumption of the base field
$k$ being algebraically closed is the definition of the extension $\Delta$ by making
explicit the partial derivatives
$\frac{\partial}{\partial x_i}$ via isomorphism
$\widehat{{\Cal O}_{W,p}}
\cong k[[x_1, ... , x_d]]$ for a choice of the system of regular parameters $(x_1, ... ,
x_d)$ of
$m_p$.  One can do this without looking at the isomorphism, and define the partial
derivatives and extension over any field $k$ of characteristic zero.  Then the rest of the
argument goes without any change.  The details are left to the reader as an exercise.) 

\vskip.1in

We finish this section stating the stability of our inductive algorithm (and hence that of
resolution process constructed via the inductive algorithm) under smooth morphisms.

\proclaim{Theorem 8-17 (Stability of the inductive algorithm under smooth morphism)}
Let
$\theta:({\Cal F}_0^{\theta},(W_0^{\theta},E_0^{\theta}))
\rightarrow ({\Cal F}_0,(W_0,E_0))$ be a smooth morphism of general basic objects
with a $(d^{\theta} = \dim W_0^{\theta} - \dim W_0 + d)$-structure and a
$d$-dimensional structure, respectively.  Let
$$({\Cal F}_0,(W_0,E_0)) \leftarrow \cdot\cdot\cdot \leftarrow ({\Cal F}_k,(W_k,E_k))$$
be a sequence of transformations of general basic objects satisfying condition
$(\heartsuit)$ and
$$({\Cal F}_0^{\theta},(W_0^{\theta},E_0^{\theta})) \leftarrow \cdot\cdot\cdot \leftarrow
({\Cal F}_k^{\theta},(W_k^{\theta},E_k^{\theta}))$$ 
be the sequence obtained by taking the Cartesian
product of the first with the smooth morphism
$\theta$.

Then the second is a sequence of transformations of general basic objects satisfing
condition $(\heartsuit)$, and the extension of the second sequence described by the
inductive algorithm is exactly the sequence
$$({\Cal F}_k^{\theta},(W_k^{\theta},E_k^{\theta})) \leftarrow \cdot\cdot\cdot \leftarrow
({\Cal F}_{k+N}^{\theta},(W_{k+N}^{\theta},E_{k+N}^{\theta}))$$
obtained by taking the Cartesian product of the smooth
morphism with the extension of the first
sequence described by the inductive algorithm
$$({\Cal F}_k,(W_k,E_k)) \leftarrow \cdot\cdot\cdot \leftarrow ({\Cal
F}_{k+N},(W_{k+N},E_{k+N}))$$
 (and by ignoring the trivial transformations
whenever the pull-backs of the centers are empty).

We say that the inductive algorithm is stable under smooth morphisms of general basic
objects.

\vskip.1in

In particular, the sequence representing resolution of singularities of a (general) basic
object constructed via the inductive algorithm is also stable under smooth morphisms.
\endproclaim 

\demo{Proof}\enddemo One can prove the invariance (stability) of the key invariants
under smooth morphisms in an identical manner to the one for proving the invariance
of the key invariants under actions.  Then the rest of the proof goes almost
verbatim as that of Proposition 8-11.  The details are left to the reader as an
exercise.

We remark that the inductive algorithm is stable under any field
extensions (of characteristic zero but may not be of finite type), which can be proved in an
identical manner.  

\newpage

$$\bold{CHAPTER\ 9.\ INVARIANTS\ REVISITED}$$

\vskip.1in

In this chapter, we construct the invariant $f^d$, based upon the key
invariants ($w\text{-}\roman{ord}$,
$\Gamma$ and $t$), associated to a (sequence of transformations of) general basic
object(s) with a
$d$-dimensional structure so that the centers of blowups in our inductive algorithm for
resolution of singularities are exactly the loci where the invariant $f^d$ attains its maximum. 
Since the invariant $f^d$ is easily seen to be stable under any action (or more generally any
smooth morphism), this will provide another easy proof for the equivariance (stability under
smooth morphisms) of the inductive algorithm.

Our invariant $f^d$ is slightly different from the one given in the paper ``A
course on constructive desingularization and equivariance" by Encinas and Villamayor, where their
invariant uses such global information as the global maximum of the $t$-invariant and hence it
is not stable under open immersions, much less so under general smooth morphisms. 

\vskip.1in

\proclaim{Definition-Construction 9-1 (Invariant $\bold{f^d}$)}\endproclaim
 
We define and construct the invariant $f^d$ by
induction on the dimension $d$ of the structure of a general basic object.

Let
$$({\Cal F}_0,(W_0,E_0)) \leftarrow \cdot\cdot\cdot \leftarrow ({\Cal F}_k,(W_k,E_k))$$
be a sequence of transformations of general basic objects satisfying condition
$(\heartsuit)$
$$(\heartsuit)\ Y_{i-1} \subset \underline{\roman{Max}}\ w\text{-}\roman{ord}_{i-1} \subset
F_{i-1}
\text{\ for\ }i = 1, ... , k.$$

\vskip.1in

$\boxed{\roman{Case}:\ d = 1.}$

\vskip.1in

When the dimension $d$ of the structure of the general basic objects is
equal to $1$, we define
$$f_k^1:F_k \rightarrow \{\{0\} \times \bold{\Gamma}^1\} \sqcup \{\bold{W}_{> 0} \times
\bold{T} \times \{\infty\}\}$$ 
in the following way: for $p \in F_k$
$$f_k^1(p) = \left\{\aligned
(w\text{-}\roman{ord}_k(p), \Gamma_k(p)) \hskip.3in \in \{0\} \times
\bold{\Gamma}^1  \hskip.52in &\text{\ if\ }\hskip.1in w\text{-}\roman{ord}_k(p) = 0 \\
(w\text{-}\roman{ord}_k(p), t_k(p), \infty) \hskip.15in \in \bold{W}_{> 0} \times \bold{T} \times
\{\infty\}&\text{\ if\ }\hskip.1in w\text{-}\roman{ord}_k(p) > 0\\
\endaligned\right.$$
with
$$\left\{\aligned
&w\text{-}\roman{ord}_k(p) \in \{0\} \sqcup \bold{W}_{> 0} = \bold{W} = \frac{1}{c!}{\Bbb
Z}_{\geq 0}
\\ &t_k(p) \in \bold{T} = \frac{1}{c!}{\Bbb Z}_{\geq 0} \times {\Bbb Z}_{\geq 0}\\
&\Gamma_k(p) \in \bold{\Gamma}^1 = ({\Bbb
Z}_{\geq -1} \times \frac{1}{c!}{\Bbb Z}_{\geq 0} \times {\Bbb Z}_{\geq 0}^1).\\
\endaligned\right.$$
Note that when $w\text{-}\roman{ord}_k(p) = 0$, by upper semi-continuity,
$w\text{-}\roman{ord}_k$ is zero in a neighborhood of $p$ and hence that $\Gamma_k(p) $ is
well-defined (cf. Definitin-Proposition 4-5 (iv)).  

We refer the reader to (GB-3) of Definition
4-1 of general basic objects for the meaning of the number $c$.

Note that we give the obvious lexicographical order (which is a total order) to the set $I_1 =
\{\{0\} \times
\bold{\Gamma}^1\}
\sqcup \{\bold{W}_{> 0} \times
\bold{T} \times \{\infty\}\}$, induced from the lexicographical orders on $\bold{W}$,
$\bold{\Gamma}^1$, and
$\bold{T}$.

(The superfluous-looking $\{\infty\}$ in the case $w\text{-}\roman{ord}_k(p) > 0$ is added to make
the invariant $f^d$ stable under smooth morphisms.)

\vskip.2in 

$\boxed{\roman{Case}:\ d = d\ \text{based\ upon\ }\roman{Case}:d = d-1 \text{\ by\ induction}.}$

\vskip.1in

Suppose we have already defined the invariant $f^{d-1}$ (with values in a totally
ordered set $I_{d-1}$) associated to a (sequence of) general basic object(s) with a \linebreak
$(d-1)$-dimensional structure. 

We define
$$f_k^d:F_k \rightarrow I_d = \{\{0\} \times \bold{\Gamma}^d\} \sqcup \{\bold{W}_{> 0} \times
\bold{T} \times \{\infty\}\} \sqcup \{\bold{W}_{> 0} \times \bold{T} \times \{0\} \times
I_{d-1}\}$$
in the following way: for $p \in F_k$
$$f_k^d(p) = \left\{\aligned
(w\text{-}\roman{ord}_k(p),&\Gamma_k(p)) \hskip.7in \in \{0\} \times \bold{\Gamma}^d \\
&\text{if\
}w\text{-}\roman{ord}_k(p) = 0 
\\ 
(w\text{-}\roman{ord}_k(p),&t_k(p),\infty) \hskip.53in \in \bold{W}_{> 0} \times \bold{T} \times
\{\infty\}
\\ &\text{if\ }w\text{-}\roman{ord}_k > 0
\text{\ and\ }R(1)(\{q
\in F_k;t_k(q) = t_k(p)\})_p
\neq
\emptyset
\\
(w\text{-}\roman{ord}_k(p),&t_k(p),0,{f''}_k^{d-1}(p)) \in \{\bold{W}_{> 0} \times \bold{T}
\times
\{0\}
\times I_{d-1}\}\\
&\text{if\ }w\text{-}\roman{ord}_k > 0
\text{\ and\ } R(1)(\{q
\in F_k;t_k(q) = t_k(p)\})_p = \emptyset\\ 
\endaligned\right.$$
with
$$\left\{\aligned
&w\text{-}\roman{ord}_k(p) \in \{0\} \sqcup \bold{W}_{> 0} = \bold{W} = \frac{1}{c!}{\Bbb
Z}_{\geq 0}
\\ &t_k(p) \in \bold{T} = \frac{1}{c!}{\Bbb Z}_{\geq 0} \times {\Bbb Z}_{\geq 0}\\
&\Gamma_k(p) \in \bold{\Gamma}^d = ({\Bbb
Z}_{\geq -d} \times \frac{1}{c!}{\Bbb Z}_{\geq 0} \times {\Bbb Z}_{\geq 0}^d),\\
\endaligned\right.$$
where $R(1)(\{q \in F_k;t_k(q) = t_k(p)\})_p$ is the codimension one part, i.e.,\linebreak
the $(d-1)$-dimensional part of the locus $\{q \in F_k;t_k(q) = t_k(p)\}$ passing through the
point
$p$.

The invariant ${f''}_k^{d-1}:G_{k,V} \rightarrow I_{d-1}$ is given by the
inductional assumption, defined on the general basic object $({\Cal
G}_{k,V},(V,{E''}_{k,V}))$ over $(G_{k,V},(V,{E''}_{k,V}))$ with a
$(d-1)$-dimensional structure, constructed in the following way: 

\vskip.1in

$\underline{\text{Construction\ and\ charactreization\ of\ the\ general\ basic\ object\
}({\Cal G}_{k,V},(V,{E''}_{k,V}))}$

\vskip.1in

Take an open neighborhood $V$ of $p$ such that
$$\align
w\text{-}\roman{ord}_k(p) &= \max \{w\text{-}\roman{ord}_k(q);q \in V \cap F_k\} \\
t_k(p) &= \max \{t_k(q);q \in V \cap F_k\}. \\
\endalign$$

In short, we follow the construction of $({\Cal G}_k,(W_k,E''_k))$ carried out in Lemma 5-3 and
Lemma 5-4 (for the conclusion of the proof for the assertions in $\bold{Case\ B}$ under
possibility $\underline{\bold{P\ 3}}$), locally over
$V$.  The general basic object we construct represents the collection of sequences of
transformations and smooth morphisms with the specifeid closed subsets being the loci
$\underline{\roman{Max}}\ t$ locally over
$V$.

\vskip.1in

We describe the construction more precisely in what follows.

\vskip.1in

Take the extension of the original sequence (only for the purpose of constructing the
general basic object $({\Cal G}_k,(W_k,E''_k))$ and defining the invariant ${f''}_k^{d-1}$)
$$({\Cal F}_0,(W_0,E_0)) \leftarrow \cdot\cdot\cdot \leftarrow ({\Cal F}_k,(W_k,E_k))
\leftarrow ({\Cal F}_{k+1},(W_{k+1},E_{k+1}))$$
where
$$({\Cal F}_k,(W_k,E_k))
\leftarrow ({\Cal F}_{k+1} = {\Cal F}_k|_V,(W_{k+1} = V,E_{k+1} = E_k|_V))$$
is induced by the open immersion $V \hookrightarrow W_k$. 

\vskip.1in

We construct a general basic object over $(G_{k,V},(V,{E''}_{k,V})) = (\underline{\roman{Max}}\
t_{k+1},(W_{k+1},E_{k+1}^+))$, with a
$d$-dimensional structure first, by specifying its charts of basic objects \linebreak
$\{(\widetilde{{W_{k+1}''}^{\lambda}}, ({{\goth a}_{k+1}''}^{\lambda},{b''}^{\lambda}),
\widetilde{{E_{k+1}''}^{\lambda}})\}$ of dimension $d$ in the following way:

Let $\{(\widetilde{W_{k+1}^{\lambda}}, ({\goth a}_{k+1}^{\lambda},b^{\lambda}),
\widetilde{E_{k+1}^{\lambda}})\}$ be the charts for the general basic objects $({\Cal F}_{k+1},
(W_{k+1},E_{k+1}))$ arising from the sequence (cf. Note 4-3)
$$(F_0,(W_0,E_0)) \leftarrow \cdot\cdot\cdot \leftarrow (F_k,(W_k,E_k)) \leftarrow ({\Cal
F}_{k+1},(W_{k+1},E_{k+1})).$$ 

We take $\widetilde{{W_{k+1}''}^{\lambda}} = \widetilde{W_{k+1}^{\lambda}}$.

If $\widetilde{{W_{k+1}''}^{\lambda}} \cap \underline{\roman{Max}}\ t_{k+1} = \emptyset$, then we
take the basic object $(\widetilde{{W_{k+1}''}^{\lambda}}, ({{\goth
a}_{k+1}''}^{\lambda},{b''}^{\lambda}),
\widetilde{{E_{k+1}''}^{\lambda}})$ to be

$$\left\{\aligned
\widetilde{{W_{k+1}''}^{\lambda}} &= \widetilde{W_{k+1}^{\lambda}} \\
{{\goth a}_{k+1}''}^{\lambda} &= {\Cal O}_{\widetilde{{W_{k+1}''}^{\lambda}}} \\
{b''}^{\lambda} &= 1 \\
\widetilde{{E_{k+1}''}^{\lambda}} &= \widetilde{E_{k+1}^{\lambda}}^+ = E_{k+1}^+ \cap
\widetilde{{W_{k+1}''}^{\lambda}}.\\
\endaligned\right.$$

If $\widetilde{{W_{k+1}''}^{\lambda}} \cap \underline{\roman{Max}}\ t_{k+1} \neq \emptyset$, then
we take the basic object $(\widetilde{{W_{k+1}''}^{\lambda}}, ({{\goth
a}_{k+1}''}^{\lambda},{b''}^{\lambda}),
\widetilde{{E_{k+1}''}^{\lambda}})$ to be

$$\left\{\aligned
\widetilde{{W_{k+1}''}^{\lambda}} &= \widetilde{W_{k+1}^{\lambda}} \\
{{\goth a}_{k+1}''}^{\lambda} &= ({\goth a}_{k+1}^{\lambda})'' \text{\ as\ constructed\ in\
Lemma\ 5-4}\\ 
{b''}^{\lambda} &= (b^{\lambda})'' \text{\ as\ constructed\ in\ Lemma\
5-4}\\
\widetilde{{E_{k+1}''}^{\lambda}} &= \widetilde{E_{k+1}^{\lambda}}^+ = E_{k+1}^+ \cap
\widetilde{{W_{k+1}''}^{\lambda}}.\\
\endaligned\right.$$

Let ${\goth C}_{G,V}$ be the collection of all the sequences of transformations and smooth
morphisms of pairs with specified closed subsets, starting with $(G_{k,V},(V,{E''}_{k,V}))$,
which satisfy condition (GB-1) with respect to the charts
$\{(\widetilde{{W_{k+1}''}^{\lambda}}, ({{\goth a}_{k+1}''}^{\lambda},{b''}^{\lambda}),
\widetilde{{E_{k+1}''}^{\lambda}})\}$.  Condition (GB-3) is trivially satisfied by the
construction, whereas condition (GB-0) is a consequence of the statement of Lemma 5-4
for $N = 0$ and condition (GB-2) a consequence of the statement of Lemma 5-4 for $N$
general.  (Note that we shift the starting point for the lemmas to the stage $k+1$.)

Therefore, the collection ${\goth C}_{G,V}$ is represented by a general basic object $({\Cal
G}_{k,V},(V,{E''}_{k,V})$ over
$(G_{k,V},(V,{E''}_{k,V}))$ with a $d$-dimensional structure. 

\vskip.1in

Now the ``Moreover" part of Lemma 5-4 and the key inducive lemma (Lemma 3-1) imply that the
general basic object $({\Cal G}_{k,V},(V,{E''}_{k,V}))$, which represents the collection ${\goth
C}_{G,V}$, has a
$(d-1)$-dimensional structure.

\vskip.1in

It also follows from Lemma 5-4 that the general basic object $({\Cal G}_{k,V},(V,{E''}_{k,V}))$
has properties $(\alpha)$ and $(\beta)$, which provide its characterization:

\vskip.1in

$(\alpha)$ With each sequence in ${\goth C}_{G,V}$
$$(G_{k,V},(V = V_k,{E''}_{k,V})) \overset{\pi_{k+1,V}''}\to{\leftarrow} \cdot\cdot\cdot
\overset{\pi_{k+N,V}''}\to{\leftarrow} (G_{k+N,V},(V_{k+N},{E''}_{k+N,V}))$$
satisfying the condition
$$G_{k+j,V} \neq \emptyset \text{\ for\ }j = 0, ... , N - 1,$$
there corresponds an extension of the original sequence of transformations and smooth
morphisms (where $\pi_{k+1}$ is the open immersion $V \hookrightarrow W_k$)
$$\align
&(F_0,(W_0,E_0)) \leftarrow \cdot\cdot\cdot \leftarrow (F_k,(W_k,E_k))
\leftarrow \\  
&(F_{k+1},(W_{k+1},E_{k+1})) \overset{\pi_{k+2}}\to{\leftarrow} \cdot\cdot\cdot
\overset{\pi_{k+N+1}}\to{\leftarrow} (F_{k+N+1},(W_{k+N+1},E_{k+N+1}))\\
\endalign$$
with condition
$$(\heartsuit') \hskip.1in \left\{\aligned &Y_{i-1} \subset \underline{\roman{Max}}\
t_{i-1} \subset
\underline{\roman{Max}}\ w\text{-}\roman{ord}_{i-1} (\subset F_{i-1}) \\
&\text{whenever\ } \pi_i \text{\ is\ a\ transformation\ with\ center\
}Y_{i-1}\\
\endaligned\right\} \text{\ for\ }i = k+2 , ... , k+N+1$$
satisfying the following conditions:

\vskip.1in

(i) $\pi_{k+j+1}''$ and $\pi_{k+j+2}$ are the transformations with the same
centers or the same smooth morphisms (as abstract varieties) for $j = 1, ... ,
N - 1$ with $V_{k+j+1} = W_{k+j+2}$ (which means, in particular, if
$\pi_{k+j+1}''$ is the transformation with center $Y_{k+j}'' \subset V_{k+j}$ which
is permissible for
$(G_{k+j,V}, (V_{k+j}, {E''}_{k+j,V}))$, then $Y_{k+j}''$ is also permissible for
$(F_{k+j+1}, (W_{k+j+1}, E_{k+j+1})$), 

(ii) we have 

\vskip.1in

\noindent either
$$\left\{\aligned
&\max\ t_{k+1} = \cdot\cdot\cdot = \max\ t_{k+N+1},
\text{\ and\ }\\
&G_{k+j,V} = \underline{\roman{Max}}\ t_{k+j+1} \text{\ for\ }j = 1,
\cdot\cdot\cdot, N\\
\endaligned\right.$$
or
$$\left\{\aligned
&\max\ t_{k+1} = \cdot\cdot\cdot = \max\ t_{k+N} >
\max\ t_{k+N+1}\\
&(\text{or\ }\max\ t_{k+1} = \cdot\cdot\cdot = \max\ t_{k+N}\ \&\ F_{k+N+1} =
\emptyset), \text{\ and}\\
&G_{k+j,V} = \underline{\roman{Max}}\ t_{k+j+1} \text{\ for\ }j = 0,
... , N-1 \hskip.1in \& \hskip.1in G_{k+N,V} = \emptyset.\\
\endaligned\right.$$

\vskip.1in

$(\beta)$ Conversely, with each extension of the original sequence of transformations and smooth
morphisms (where $\pi_{k+1}$ is the open immersion $V \hookrightarrow W_k$)
$$\align
&(F_0,(W_0,E_0)) \leftarrow \cdot\cdot\cdot \leftarrow (F_k,(W_k,E_k))
\leftarrow \\  
&(F_{k+1},(W_{k+1},E_{k+1})) \overset{\pi_{k+2}}\to{\leftarrow} \cdot\cdot\cdot
\overset{\pi_{k+N+1}}\to{\leftarrow} (F_{k+N+1},(W_{k+N+1},E_{k+N+1}))\\
\endalign$$
with condition
$$(\heartsuit') \hskip.1in \left\{\aligned &Y_{i-1} \subset \underline{\roman{Max}}\
t_{i-1} \subset
\underline{\roman{Max}}\ w\text{-}\roman{ord}_{i-1} (\subset F_{i-1}) \\
&\text{whenever\ } \pi_i \text{\ is\ a\ transformation\ with\ center\
}Y_{i-1}\\
\endaligned\right\} \text{\ for\ }i = k+2 , ... , k+N+1$$
and the condition
$$\max\ t_{k+1} = \cdot\cdot\cdot = \max\ t_{k+N},$$
there corresponds a sequence of transformations and smooth morphisms of general basic
objects starting from $(G_{k,V},(V = V_k,{E''}_{k,V}))$
$$(G_{k,V},(V = V_k,{E''}_{k,V})) \overset{\pi_{k+1,V}''}\to{\leftarrow} \cdot\cdot\cdot
\overset{\pi_{k+N,V}''}\to{\leftarrow} (G_{k+N,V},(V_{k+N},{E''}_{k+N,V}))$$
satisfying the condition 
$$G_{k+j,V} \neq \emptyset \text{\ for\ }j = 0, ... , N - 1$$
and conditions (i) and (ii) as in $(\alpha)$.

\vskip.1in

Once we construct the general basic object $({\Cal G}_{k,V},(V,{E''}_{k,V}))$ with a
$(d-1)$-dimensional structure, we have by induction the invariant ${f''}_k^{d-1}:G_{k,V}
\rightarrow I_{d-1}$ attached to this general basic object (considered as a trivial
sequence of general basic objects consisting only of itself).

Remark that the value ${f''}_k^{d-1}(p)$ is independent of the choice of the neighborhood $V$,
since if we choose a different open neighborhood $V'$, the general basic objects $({\Cal
G}_{k,V},(V,{E''}_{k,V}))$ and $({\Cal G}_{k,V'},(V',{E''}_{k,V'}))$ restrict to the same
general basic object $({\Cal G}_{k,V \cap V'},(V \cap V',{E''}_{k,V \cap V'}))$ (cf. Remark 4-2
(iv)(v)).
 
Note that we give the obvious lexicographical order (which is a total order) to the set
\linebreak
$I_d = \{\{0\} \times \bold{\Gamma}^d\}
\sqcup \{\bold{W}_{> 0} \times
\bold{T} \times \{\infty\}\} \sqcup \{\bold{W}_{> 0} \times \bold{T} \times \{0\} \times
I_{d-1}\}$ induced from the lexicographical orders on $\bold{W}, \bold{\Gamma}^d, \bold{T}$ and
$I_{d-1}$.

\vskip.1in

This completes the definition and construction of the invariant $f^d$.

\vskip.1in

\proclaim{Remark 9-2}\endproclaim

\vskip.1in

(i) In the definition above, we construct $f^d$ from BOTTOM UP based upon the construction of
$f^{d-1}$ inductively.  However, in reality, we can start writing down the invariant $f^d$ from TOP
DOWN without knowing what $f^{d-1}$ would be: First we compute
$w\text{-}\roman{ord}_k(p)$.  If
$w\text{-}\roman{ord}_k(p) = 0$, then go on to compute $\Gamma_k(p)$.  If
$w\text{-}\roman{ord}_k(p) > 0$, then go on to compute $t_k(p)$.  If $R(1)(q \in
F_k;\{t_k(q) = t_k(p)\})_p \neq
\emptyset$, then set the next factor to be $\infty$.  If $R(1)(q \in F_k;\{t_k(q) =
t_k(p)\})_p = \emptyset$, then set the next factor to be $0$ and construct $({\Cal
G}_{k,V},(V,{E''}_{k,V}))$.  Now with the general basic object $({\Cal
G}_{k,V},(V,{E''}_{k,V}))$ with a $(d-1)$-dimensional structue, we start writing down
$w\text{-}\roman{ord}$ and repeat the same procedure as above, and so on.

(ii) Though in the definition of $f^d$ above we used the charts
$\{(\widetilde{W_k^{\lambda}},({\goth
a}_k^{\lambda},b^{\lambda}),\widetilde{E_k^{\lambda}})\}$, the invariant $f^d$ is completely
determined only by the collection ${\goth C}_i$ of the sequences of transformations and smooth
morphisms with the specified closed subsets represented by the general basic object $({\Cal
F}_i,(W_i,E_i))$ for $i = 0, \cdot\cdot\cdot, k$, the number $d$ which refers to the dimension of
the structure of the general basic objects and by the original sequence   
$$({\Cal F}_0,(W_0,E_0)) \leftarrow \cdot\cdot\cdot \leftarrow ({\Cal F}_k,(W_k,E_k)).$$
In fact, by Hironaka's trick, the invariants $w\text{-}\roman{ord}_k$ and $\roman{ord}_k$
are determined by
${\goth C}_k$ and by the number $d$ referring to the dimension of the structure (cf.
Definition-Proposition 4-5 and Remark 4-7).  Therefore,
$\Gamma_k$ can also be determined by
${\goth C}_k$, the number $d$, and the original sequence, as it can be computed purely
by looking at
$\roman{ord}_k$, $F_k$ and $E_k$.  The invariant $t_k$ is determined by ${\goth C}_k$,
the number $d$, and the original sequence also, as it can be computed purely by looking
at
$w\text{-}\roman{ord}_k$ and
$E_k^-$.  Thus whether or not the codimension one part of the (local) maximum locus of
the invariant $t$ passes though a given point is also determined by ${\goth C}_k$, the
number $d$, and the original sequence.  Now notice that the general basic object $({\Cal
G}_{k,V},(V,{E''}_{k,V}))$ is also determined by ${\goth C}_k$, the number $d$, and the
original sequence, as the collection ${\goth C}_{G,V}$ is characterized by the loci
$\underline{\roman{Max}}\ t$ of the corresponding sequences in ${\goth C}_k$ (sequel to the
open immersion $({\Cal F}_k,(W_k,E_k)) \leftarrow ({\Cal F}_{k+1},(W_{k+1},E_{k+1}))$). 
Therefore, ${f''}_k^{d-1}$ is also determined by ${\goth C}_k$, the number $d$, and the
original sequence.  

\vskip.1in

(iii)  The invariant $f^d$ DOES depend on the number $d$ specifying the dimension of the
structure of your choice (of the general basic object) and is NOT purely determined by
the collection
${\goth C}_i$ of the sequences of transformations and smooth morphisms with the
specified closed subsets, represented by the general basic object $({\Cal
F}_i,(W_i,E_i))$ for $i = 0,
\cdot\cdot\cdot, k$ and by the original sequence   
$$({\Cal F}_0,(W_0,E_0)) \leftarrow \cdot\cdot\cdot \leftarrow ({\Cal F}_k,(W_k,E_k)),$$
since so does the invariant $w\text{-}\roman{ord}$ (cf. Remark 4-7) and hence also the
invariant $t$.

\vskip.1in

(iv) (Stability under smooth morphism) Let
$\theta:({\Cal F}_{0}^{\theta},(W_{0}^{\theta},E_{0}^{\theta}))
\rightarrow ({\Cal F}_0,(W_0,E_0))$ be a smooth morphism of general basic objects of relative
dimension $r$, so that the dimension of the structure of $({\Cal
F}_{0}^{\theta},(W_{0}^{\theta},E_{0}^{\theta}))$, induced by that of $({\Cal
F}_0,(W_0,E_0))$ with a
$d$-dimensional structure, is equal to $d + r$.  
 
Let
$$({\Cal F}_0,(W_0,E_0)) \leftarrow \cdot\cdot\cdot \leftarrow ({\Cal F}_k,(W_k,E_k))$$
be a sequence of transformations of general basic objects with a $d$-dimensional structure
satisfying condition
$(\heartsuit)$ and
$$({\Cal F}_{0}^{\theta},(W_{0}^{\theta},E_{0}^{\theta})) \leftarrow \cdot\cdot\cdot
\leftarrow ({\Cal F}_{k}^{\theta},(W_{k}^{\theta},E_{k}^{\theta}))$$ 
be the sequence obtained by taking the Cartesian
product of the first with the smooth morphism
$\theta$.

Then the second sequence is a sequence of transformations of general basic objects
with a $(d+r)$-dimensional sructure satisfying condition $(\heartsuit)$, and we have
$$f_k^{d+r}(p^{\theta}) = f_k^d(p) \text{\ for\ any\ point\ }p \in F_k \text{\ and\
}p^{\theta}
\in F_{k}^{\theta} \cap \theta^{-1}(p).$$
(The verification is straightforward identifing the factors from TOP DOWN, and left to
the reader as an exercise.  (cf. Theorem 8-17)).

(v) (Stability under (analytic) localization) The invariant $f^d$ is stable under
(analytic) localization in the following sense: Let
$$({\Cal F}_0,(W_0,E_0)) \leftarrow \cdot\cdot\cdot \leftarrow ({\Cal F}_k,(W_k,E_k))$$
be a sequence of transformations of general basic objects with $d$-dimensional structures
satisfying condition $(\heartsuit)$.  Let $p = p_k \in F_k$ be a point and $p_i \in
F_i$ its image in $W_i$.

Suppose we have another sequence of transformations of general basic objects with
$d$-dimensional structures 
$$({\Cal F}_0',(W_0',E_0')) \leftarrow \cdot\cdot\cdot \leftarrow ({\Cal F}_k',(W_k',E_k'))$$
with a point $p' = p_k' \in F_k'$ and its image $p_i' \in F_i'$ in $W_i'$.

Suppose we can find open neighborhoods $V_i$ of $p_i$ (resp. $V_i'$ of $p_i'$) such that we have a
commutative diagram of sequences of transformations of general basic objects, restricted
to the open subsets, with vertical arrows being isomorphisms of general basic objects

$$\CD
({\Cal F}_0|_{V_0},(V_0,E_0|_{V_0})) @. \leftarrow @. \cdot\cdot\cdot @. \leftarrow @. ({\Cal
F}_k|_{V_k},(V_k,E_k|_{V_k}))
\\ @VVV @.@.@. @VVV \\
({\Cal F}_0'|_{V_0'},(V_0',E_0'|_{V_0'})) @. \leftarrow @. \cdot\cdot\cdot @. \leftarrow @. ({\Cal
F}_k'|_{V_k'},(V_k',E_k'|_{V_k'})).
\\
\endCD$$

Or more generally suppose we have a
commutative diagram of sequences of transformations of basic objects, restricted to the
analytic neighborhoods, with vertical arrows being isomorphisms of (analytic) general
basic objects

$$\CD
({\Cal F}_0,(W_0,E_0)) \times \roman{Spec}\ \widehat{{\Cal O}_{W_0,p_0}}@. \leftarrow @.
\cdot\cdot\cdot @.
\leftarrow @. ({\Cal F}_k,(W_k,E_k)) \times \roman{Spec}\ \widehat{{\Cal O}_{W_k,p_k}}
\\
@VVV @.@.@. @VVV \\
({\Cal F}_0',(W_0',E_0')) \times \roman{Spec}\ \widehat{{\Cal O}_{W_0',p_0'}}@. \leftarrow @.
\cdot\cdot\cdot @. \leftarrow @. ({\Cal F}_k,(W_k,E_k)) \times \roman{Spec}\ \widehat{{\Cal
O}_{W_k',p_k'}}.\\
\endCD$$

Then we have
$$f_k^d(p) = f_k^d(p').$$

The verification is straightforward and left to the reader as an exercise.

\vskip.1in

(vi) (Invariance under action) Let
$$({\Cal F}_0,(W_0,E_0)) \leftarrow \cdot\cdot\cdot \leftarrow ({\Cal F}_k,(W_k,E_k))$$
be a sequence of transformations of general basic objects with $d$-dimensional structures
satisfying condition
$(\heartsuit)$.  Let $\theta$ be an action on the general basic object $({\Cal F}_0,(W_0,E_0))$. 
Suppose the sequence is $\theta$-equivariant.  Then $f_k^d(p) = f_k^d(\theta(p))$ for
any $p \in F_k$.

This is an immediate consequence of Lemma 8-10.

\vskip.1in

\proclaim{Theorem 9-3} Let
$$({\Cal F}_0,(W_0,E_0)) \leftarrow \cdot\cdot\cdot \leftarrow ({\Cal F}_k,(W_k,E_k))$$
be the sequence of transformations of general basic objects with $d$-dimensional
structures, obtained via the inductive algorithm (Theorem 5-1), representing resolution
of singularities of a general basic object
$({\Cal F}_0,(W_0,E_0))$.

Then the centers of the transformations are exactly the loci where the invariants $f^d$ take their
maxima, i.e.,
$$Y_{i-1} = \underline{\roman{Max}}\ f_{i-1}^d \text{\ for\ } i = 1, ... , k.$$

\endproclaim

\demo{Proof}\enddemo This is a straightforward consequence of the description of the process of the
inductive algorithm on how to choose the centers for the sequence
representing resolution of singularities in Theorem 5-1 and the way we define the
invariant $f^d$.
 
\vskip.1in

Note that the invariance of the invariant $f^d$ and its stablity under smooth morphism
(cf. Remark 9-2 (vi) and (iv)) give an alternative proof (or rather to say, a different
presentation of the same proof), for the corresponding statements (cf. Proposition 8-11
Theorem 8-17) for the sequence of transformations representing resolution of
singularities obtained via the inductive algorithm.

\newpage

${}\hskip1in$$\bold{CHAPTER\ 10.\ NON\text{-}EMBEDDED}$ \linebreak
${}\hskip1in$$\bold{RESOLUTION\ OF\ SINGULARITIES}$

\vskip.1in

In this chapter, we prove ((non-embedded) resolution of singularities), achieving Main Theme
0-1.

A variety $X$ can be covered by a finite number of open subsets $\{U_s\}_{s \in S}$ which are
embedded, as closed subschemes, into smooth varieties $W_{U_s}$, i.e., $U_s \subset W_{U_s}$. 
By choosing a number
$d$ sufficiently large and replacing $W_{U_s}$ with $W_{U_s} \times {\Bbb A}^{d -
\dim W_{U_s}}$, we may assume that all the ambient smooth varieties $W_{U_s}$ are of the same
dimension
$d$. 

We observe then, on the intersections $U_s \cap U_{s'}$ of the open subsets of the covering
$\{U_s\}_{s \in S}$, not only that the invariants
$f^d_{U_s}$ and
$f^d_{U_{s'}}$, defined on the singular loci \linebreak
$U_s \subset W_{U_s}$ and $U_{s'}
\subset W_{U_{s'}}$ of the basic objects
$(W_{U_s},({\Cal I}_{U_s},1),\emptyset)$ and
$(W_{U_{s'}},({\Cal I}_{U_{s'}},1), \emptyset)$ (and on their transformations) as in
Chapter 9, coincide and hence give rise to a global invariant $f^d_X$ on $X$ (and
its transformations), but also that the ideals defining the loci
$\{f^d_{U_s} = \max f^d_X\} \subset W_{U_s}$ and $\{f^d_{U_{s'}} = \max f^d_X\} \subset
W_{U_{s'}}$, restricted to $U_s \subset W_{U_s}$ and $U_{s'} \subset W_{U_{s'}}$, coincide and
hence give rise to a global ideal of the center of blowup on $X$ (and on its transformations).

Choosing the center(s) of blowup(s) this way, though based upon the method of embedded resolution
of singularities, we obtain a sequence representing non-embedded resolution of singularities of
$X$.

The
sequence thus obtained is independent of the choice of the number $d$ (though the invariant
$f^d_X$ is dependent of the number $d$) or the choice of the covering $\{U_s\}_{s
\in S}$.  

The sequence is
equivariant with respect to any automorphism $\theta:X \overset{\sim}\to{\rightarrow} X$ in
the sense that it lifts to an automorphism of the sequence.

\vskip.1in

\proclaim{Theorem 10-1 ((Non-embedded) Resolution of singularities)} Let $X$ be a variety over
a field $k$ of characteristic zero.  Take an open covering $\{U\}$ of $X$ so that the open
subsets $U$ are embedded, as closed subschemes, into varieties $W_U$ of dimension $d$ smooth over
$k$.  (The number $d$ is common to all the varieties $W_U$.)  

Then the invariants $f^d_U$, defined as in Chapter 9, on the
singular loci $U \subset W_U$ of the basic objects $(W_U,({\Cal I}_U,1),\emptyset)$ (as well
as the invariants on their transformations) patch together to give rise to a global invariant
$f^d_X$ on $X$ (as well as to the global invariants on the transformations of $X$). 
Moreover, the ideals, restricted to $U$, defining the loci $\{f^d_U = \max f^d_X\}
\subset W_U$ (as well as the ideals, restricted to the transformations of $U$,
defining the maximum loci of the invariants on the transformations of the original
basic objects $(W_U,({\Cal I}_U,1),\emptyset)$) patch together to determine the
global ideal of the center of blowup on
$X$ (as well as the global ideals of the centers of blowups on the
transformations of $X$).

This provides an algorithm to choose the centers of blowups, which lie over the singular locus
$\roman{Sing}(X)$ of $X$, for constructing a sequence of (non-embedded) resolution of
singularities of
$X$.

The sequence thus obtained is
independent of the choice of the number $d$ (though the invariant
$f^d_X$ is dependent of the number $d$) or the choice of the covering $\{U\}$. 

The sequence thus obtained is
equivariant with respect to any automorphism \linebreak
$\theta:X
\overset{\sim}\to{\rightarrow} X$ in the sense that it lifts to an automorphism of
the sequence.  

\vskip.1in 

We give a more precise description of our algorithm in the following:

\vskip.1in

Inductively, we construct (a part of, i.e., up to the $k$-th stage of) the sequence
representing (non-embedded) resolution of singularities
$$X = X_0 \overset{\pi_1}\to{\leftarrow} X_1 \overset{\pi_2}\to{\leftarrow} \cdot\cdot\cdot
\overset{\pi_{k-1}}\to{\leftarrow} X_{k-1} \overset{\pi_k}\to{\leftarrow} X_k$$
with centers $Y_{i-1} \subset X_{i-1}$ for $i = 1, ... , k$ and the invariant $f_{X_i}^d$ on
$X_i$ for $i = 0, 1, ... , k$ with the
following properties:

(i) For each closed point $p \in X$ and an open subset $U$ containing $p$ and taken
from the open covering, there exists an open neighborhood $p \in U_p \subset U$ with
an induced embedding
$U_p \subset W_{U_p}$ such that
$$\psi_{i-1}(Y_{i-1}) \cap U_p \neq \emptyset \text{\ if\ and\ only\ if\
}\psi_{i-1}^{-1}(p) \cap Y_{i-1} \neq \emptyset \text{\ for\ }i = 1, ... , k$$ 
where $\psi_{i-1} = \pi_1 \circ \pi_2 \circ \cdot\cdot\cdot \circ \pi_{i-1}$.

(ii) The first $l_{p,k}$-stages of the sequence of transformations representing
resolution of singularities of the basic object $(W_{U_p},({\Cal I}_{U_p},1),
\emptyset)$ 
$$(W_{U_p}, ({\Cal I}_{U_p},1), \emptyset) = ((W_{U_p})_0,((J_{U_p})_0,1), (E_{U_p})_0) \leftarrow
\cdot\cdot\cdot \leftarrow ((W_{U_p})_{l_{p,k}},((J_{U_p})_{l_{p,k}},1), (E_{U_p})_{l_{p,k}}),$$
where
$$l_{p,k} = \#\{i; \psi_{i-1}^{-1}(p) \cap Y_{i-1} \neq \emptyset, i = 1, ... ,
k\},$$ 
gives rise to a sequence representing (the first $l_{p,k}$-stages of) resolution of
singularities of
$U_p$
$$U_p = (U_p)_0 \leftarrow \cdot\cdot\cdot \leftarrow (U_p)_{l_{p,k}}$$ 
where $(U_p)_j$ is the strict transform of $U_p = (U_p)_0$ on $(W_{U_p})_j$ for $j = 1, ...
, l_{p,k}$.  

This sequence coincides with 
$$X = X_0 \overset{\pi_1}\to{\leftarrow} X_1 \overset{\pi_2}\to{\leftarrow} \cdot\cdot\cdot
\overset{\pi_{k-1}}\to{\leftarrow} X_{k-1} \overset{\pi_k}\to{\leftarrow} X_k$$
restricted over $U_p$ where we ignore the trivial transformations $\pi_i$ whenever \linebreak
$\psi_{i-1}^{-1}(U_p)
\cap Y_{i-1} = \emptyset$.

(We remark that, in the above sequence of basic objects starting from $(W_{U_p},({\Cal
I}_{U_p},1),
\emptyset)$, we only consider the neighborhoods of the strict transforms $(U_p)_j$ of $(U_p)_0 =
U_p$.  That is to say, in the sequence representing resolution of singularities of $(W_{U_p},({\Cal I}_{U_p},1),
\emptyset)$, we ignore the transformations whose centers of blowups are away from the strict
transforms.)

(iii) The invariant $f_{X_k}^d(p_k)$ for $p_k \in \psi_k^{-1}(U_p) \subset X_k$, where
$\psi_k^{-1}(U_p) = (U_p)_{l_{p,k}}$ with
$\psi_k = \pi_1 \circ \cdot\cdot\cdot \circ \pi_k$, is defined to be equal to the value
$f_{U_p,l_{p,k}}^d(p_k)$, where the invariant $f_{U_p,l_{p,k}}^d$ is attached to the 
$l_{p,k}$-th stage
of the sequence of basic objects given in (ii) and where the invariant $f_{U_p,l_{p,k}}^d$ is
defined on the singular locus of the basic object
$((W_{U_p})_{l_{p,k}},((J_{U_p})_{l_{p,k}},1), (E_{U_p})_{l_{p,k}})$, which contains
$(U_p)_{l_{p,k}}$.

The invariant $f_{U_p,l_{p,k}}^d(p_k)$ is independent of the choice of $U$
or
$U_p$ (justifying the omission of reference to $U_p$ or $U$ in the notation
$f_{X_k}^d$).

(Note that the invariant $f^d$ is stable (invariant)
under open immersion.  Thus the invariant
$f_{U_p}^d$ is just the restriction of $f_U^d$ over $U_p$.) 

(iv) The defining ideal ${\Cal I}_{Y_k}$ of the center $Y_k \subset X_k$, which may
not be smooth, reduced or irreducible in general, is taken so that ${\Cal
I}_{Y_k}|_{\psi_k^{-1}(U_p)}$ coincides with the ideal, restricted to $(U_p)_{l_{p,k}} =
\psi_k^{-1}(U_p)$, defining the locus $\{f_{U_p,l_{p,k}}^d = \max f^d_{X_k}\}$ inside of
$(W_{U_p})_{l_{p,k}}$.

(v) The center $Y_k$ lies over the singular locus $\roman{Sing}(X)$ of $X$.

\endproclaim

\proclaim{Remark 10-2}\endproclaim

(i) As we stated in the footnote to Main Theme 0-1, we do not require that the
centers $Y_i \subset X_i \hskip.1in (i = 0, ... , k)$ to be smooth or reduced (and we also allow
the centers to be reducible).  In fact, our algorithm produces centers which may not be
smooth or reduced.  Therefore, though it is true that set-theoretically we have
$$\roman{Supp}(Y_i) = \underline{\roman{Max}}\ f^d_{X_i} = \{p \in
X_i;f^d_{X_i}(p) = \max f^d_{X_i}\},$$ this description of the center is
not enough to determine its scheme-theoretic structure.  This feature is in
clear contrast to the situiation where, if we required the centers to be (smooth
and) reduced, the set-theoretic description of them as the maximum loci of the
invariants would suffice.

(ii) We remark that it is NOT sufficient merely to prove, in order to construct a global sequence
of non-embedded resolution of singularities, that
$X$ has an open covering $\{U\}$ with embeddings $U \subset W_U$ into varieties smooth over $k$ and
that the process of embedded resolution of
$U
\subset W_U$ restricted to $U \cap V$ ``coincides in the naive sense" with the
process of embedded resolution of
$V
\subset W_V$ restricted to $U \cap V$ for any two open subsets $U$ and $V$ of the open
covering.  The reason, which involves the interpretation of the words ``coincides in
the naive sense", is the following:  

When we restrict the process to a smaller open subset $U \cap V$, we ignore
the trivial transformations blowing up the centers outside of $U
\cap V$.  So even if we prove the processes, obtained by restricting those over $U$ and $V$,
coincide after ignoring those trivial transformations, we would be at loss, without introducing an
invariant, about how to patch the processes for $U$ and $V$ together including the transformations
with centers outside of $U
\cap V$, in order to obtain a global order of choosing centers.   

\vskip.1in

\demo{Proof}\enddemo Let $\overline{k}$ be the algebraic closure of the field $k$. 
If we prove the assertion over $\overline{k}$, i.e., for $X \times \roman{Spec}\
\overline{k}$ with its open covering
$\{U
\times
\roman{Spec}\ \overline{k}\}$ and embeddings $U \times
\roman{Spec}\ \overline{k} \hookrightarrow W_U \times
\roman{Spec}\ \overline{k}$ (We may lose the irreducibility assumption on ``varieties", but
since the theory remains valid without any change in the argument, we ignore this point.),
together with the equivariance assertion (which obviously implies the equivariance under the
action of the Galois group
$\roman{Gal}(\overline{k}/k)$), then the assertion over $k$ follows.  Therefore, we may assume
that $k$ is algebraically closed in what follows.

\vskip.1in

We check the inductive construction of our algorithm (stated in ``a more precise
description" in the statement of Theorem 10-1), condition by condition, starting with
the $(k = 0)$-th stage. 

\newpage

$\boxed{\roman{Case}\ k = 0}$ 

\vskip.1in

We look at the case when $k = 0$.

\vskip.1in

Condition (i) is obvious at the stage $k = 0$, choosing $p \in U_p \subset U$ where $U$ is
an open subset containg $p$ and taken from the covering and where $U_p$ is any open neighborhood
$p
\in U_p
\subset U$ 
\vskip.1in

Condition (ii) is clear at the stage $k = 0$.

\vskip.1in

We verify condition (iii) in the following. 

\vskip.1in

Let $e = \dim m_{X,p}/m_{X,p}^2$ be the embedding dimension of $X$ at $p$.

Take a small open neighborhood $U_p$ so that we
may assume that it is embedded into a smooth affine variety $W_{U_p} = \roman{Spec}\ A(W_{U_p})$ of
dimension
$d$ with affine coordinate ring $A(W_{U_p})$ and that there exists a system of regular parameters
$(x_1, ... , x_{d - e}, y_1, ... , y_e)$ satisfying the following properties:

\vskip.1in

$(\alpha)\ \langle x_1, ... , x_{d - e}\rangle \subset \Gamma(W_{U_p},{\Cal I}_{U_p}) =
I_{U_p}$, and 

$(\beta)$ we have a commutative diagram
$$\CD
0 @. \hskip.1in \rightarrow @. {\Cal K} @. \rightarrow @. k[[Y_1, ... , Y_e]] @. \rightarrow
\hskip.1in @.
\widehat{{\Cal O}_{X,p}} @. \rightarrow @. 0 \\
@| @. @AAA @. @AAA @.@| \\
0 @. \rightarrow @. I_{U_p}/\langle x_1, ... , x_{d-e}\rangle \otimes \widehat{{\Cal
O}_{W_{U_p},p}} @.
\rightarrow @. A(W_{U_p})/\langle x_1, ... , x_{d-e}\rangle \otimes \widehat{{\Cal O}_{W_{U_p},p}}
@.
\rightarrow @. \widehat{{\Cal O}_{U_p,p}} @. \rightarrow @. 0 \\
\endCD$$
where the second middle vertical arrow is an isomorphism sending $y_1, ... , y_e$ to $Y_1,
... , Y_e$, regarded as independent variables. 

This induces an isomorphism of (analytic) basic objects
$$\align
(W_{U_p}, ({\Cal I}_{U_p},1), \emptyset)/&(x_1, ... , x_{d-e}) \times \roman{Spec}\
\widehat{{\Cal O}_{W,p}} \\
&\overset{\sim}\to{\rightarrow} ((\widehat{{\Bbb A}^e})_0 = \roman{Spec}\
k[[Y_1, ... , Y_e]], ({\Cal K}_0 = {\Cal K},1),({\Cal E}_p)_0 = \emptyset),\\
\endalign$$
where ``$/\langle x_1, ... , x_{d-e}\rangle$" denotes the restriction to the nonsingular subvariety
defined by the ideal $\langle x_1, ... , x_{d-e}\rangle$, by abuse of notation.

\vskip.1in

We now compute the invariant $f_{U_p,0}^d(p)$.

\vskip.1in

$\underline{\roman{Subcase}\ e - \dim X > 0}$: We consider the subcase $e - \dim X > 0$.

\vskip.1in

Suppose $d - e > 0$.

First, since $x_1 \in (J_{U_p})_0 = \overline{(J_{U_p})_0} = {\Cal I}_{U_p}$, we have
$$w\text{-}\roman{ord}_{U_p,0}^d(p) = 1 \text{\ and\ } t_{U_p,0}^d(p) = (1,0).$$
Since $d - e + e - \dim X > 1$, the next factor is $0$ (cf. Chapter 9).  Thus we have the pattern
$$(1, (1,0), 0).$$
Now we look at the general basic object with a $(d-1)$-dimensional structure as in
Definition-Construction 9-1, whose chart at
$p$ is given by $(W_{U_p}, ({\Cal I}_{U_p},1), \emptyset)/\langle x_1\rangle$ by construction (cf.
Lemma 5-3 and Lemma 5-4).  Note that, since $b = b_0 = b' = b'' = 1$, the coefficient ideal is
nothing but the restriction of the original ideal $(J_{U_p})_0 = \overline{(J_{U_p})_0} = {\Cal
I}_{U_p}$ to the smooth hypersurface $\{x_1 = 0\}$. 

\vskip.1in

Suppose $d - e - 1 > 0$.

First, since $x_2 \in (J_{U_p}/\langle x_1\rangle)_0 = \overline{(J_{U_p}/\langle x_1\rangle)_0} =
{\Cal I}_{U_p}/\langle x_1\rangle$, we have 
$$w\text{-}\roman{ord}_{U_p/\langle x_1\rangle,0}^{d-1}(p) = 1 \text{\
and\ }t_{U_p/\langle x_1\rangle,0}^{d-1}(p) = (1,0).$$ 
Since $d - e - 1 + e - \dim X > 1$, the next
factor is
$0$.  Thus we have the pattern
$$(1, (1,0), 0, 1, (1,0), 0).$$
Now we look at the general basic object with a $(d-2)$-dimensional structure as in
Definition-Construction 9-1,
whose chart at
$p$ is given by $(W_{U_p}, ({\Cal I}_{U_p},1), \emptyset)/\langle x_1,x_2\rangle$ by construction
(cf. Lemma 5-3 and Lemma 5-4).  Note that, since $b = b_0 = b' = b'' = 1$, the coefficient ideal is
nothing but the restriction of the ideal $(J_{U_p}/\langle x_1\rangle)_0 =
\overline{(J_{U_p}/\langle x_1\rangle)_0} = {\Cal I}_{U_p}/\langle x_1\rangle$ to the smooth
subvariety
$\{x_1 = x_2 = 0\}$.

\vskip.1in

Inductively, carrying out the same argument with $x_3, ... , x_{d-e}$, we conclude that
after repeating the pattern
$(1,(1,0),0)$ for
$(d-e)$-times
$$(1, (1,0), 0, 1, (1,0), 0, \cdot\cdot\cdot, 1, (1,0),0)$$
we reach the general basic object with a $(e = d - (d-e))$-dimensional structure, whose chart at
$p$ is given by $(W_{U_p}, ({\Cal I}_{U_p},1), \emptyset)/\langle x_1, \cdot\cdot\cdot,
x_{d-e}\rangle$, which is analytically isomorphic to $((\widehat{{\Bbb A}^e})_0,({\Cal
K}_0,1),({\Cal E}_p)_0)$.  

It is straightforward to see that the analytic basic object $((\widehat{{\Bbb A}^e})_0,({\Cal
K}_0,1),({\Cal E}_p)_0)$ is purely determined by $\widehat{{\Cal O}_{X,p}}$, and independent of the
choice of $U$, $U_p \subset U$, or a system of regular parameters $(x_1, ... , x_{d-e}, y_1, ... ,
y_e)$.  In fact, let $p \in U_p' \subset U'$ be another choice of open subsets with a system
regular parameters $(x_1', ... , x_{d-e}', y_1', ... , y_e')$ as above, which leads to a
commutative diagram 
$$\CD
0 @. \hskip.1in \rightarrow @. {\Cal K}' @. \rightarrow @. k[[Y_1', ... , Y_e']] @. \rightarrow
\hskip.1in @.
\widehat{{\Cal O}_{X,p}} @. \rightarrow @. 0 \\
@| @. @AAA @. @AAA @.@| \\
0 @. \rightarrow @. I_{U_p}/\langle x_1', ... , x_{d-e}'\rangle \otimes \widehat{{\Cal
O}_{W_{U_p},p}} @.
\rightarrow @. A(W_{U_p})/\langle x_1', ... , x_{d-e}'\rangle \otimes \widehat{{\Cal
O}_{W_{U_p},p}} @.
\rightarrow @. \widehat{{\Cal O}_{U_p,p}} @. \rightarrow @. 0. \\
\endCD$$
Then there exists an isomorphism (though non-canonical) 
$$\phi:k[[Y_1, ... , Y_e]] \rightarrow
k[[Y_1', ... , Y_e']]$$ 
which makes the following diagram commute
$$\CD
0 @. \hskip.1in \rightarrow \hskip.1in @. {\Cal K} @. \hskip.1in \rightarrow \hskip.1in @. k[[Y_1,
... , Y_e]] @.
\hskip.1in \rightarrow
\hskip.1in @.
\widehat{{\Cal O}_{X,p}} @. \hskip.1in \rightarrow \hskip.1in @. 0 \\
@| @. @VVV @. @VV{\phi}V @.@| \\
0 @. \hskip.1in \rightarrow \hskip.1in @. {\Cal K}' @. \hskip.1in \rightarrow \hskip.1in @.
k[[Y_1', ... , Y_e']]
\hskip.1in @.
\hskip.1in \rightarrow
\hskip.1in @.
\widehat{{\Cal O}_{X,p}} @. \hskip.1in \rightarrow \hskip.1in @. 0. \\
\endCD$$
Thus
denoting the invariant attached to the analytic basic object $((\widehat{{\Bbb A}^e})_0,({\Cal
K}_0,1),({\Cal E}_p)_0)$ by $\widehat{f_{p,0}^e}$, we conclude

$$f_{X_0}^d(p) = f_{U_p,0}^d(p) = (1, (1,0), 0, 1, (1,0), 0, \cdot\cdot\cdot, 1, (1,0),0,
\widehat{f_{p,0}^e}(p)),$$

where the pattern $(1, (1,0), 0)$ is repeated for $(d-e)$-times.  This is independent of the choice
of $U$, $U_p \subset U$, or a system of regular parameters $(x_1, ... , x_{d-e}, y_1, ... , y_e)$.

\vskip.1in

Suppose $d - e = 0$.  Then we conclude
$$f_{X_0}^d(p) = f_{U_p,0}^d(p) = \widehat{f_{p,0}^e}(p).$$

\vskip.1in

$\underline{\roman{Subcase}\ e - \dim X = 0}$: We consider the remaining subcase $e - \dim X = 0$.

\vskip.1in

In the case $e - \dim X =
0$, we conclude by a similar consideration 

$$f_{X_0}^d(p) = f_{U_p,0}^d(p) = (1, (1,0), 0, 1, (1,0), 0, \cdot\cdot\cdot, 1, (1,0), 0, 1,
(1,0),\infty),$$

where the pattern $(1, (1,0), 0)$ is repeated for $(d-e-1)$-times.  This is also independent of the
choice of $U$, $U_p \subset U$, or a system of regular parameters $(x_1, ... , x_{d-e}, y_1, ... ,
y_e)$.  Observe that in this case, i.e., when $p$ is a nonsingular point of $X$, the value
$f_{X_0}^d(p)$ is minimum, i.e.,
$$f_{X_0}^d(p) \leq f_{X_0}^d(q) \hskip.1in \forall q \in X.$$   

\vskip.2in

Now we look at condition (iv), using the same notation as for the verification of condition (iii).

\vskip.1in

The locus $\{q \in X_0;f^d_{U_p,0}(q) = \max\ f_X^d\} \subset W_{U_p} = (W_{U_p})_0$ can be
identified, after taking the product
$\times \roman{Spec}\ \widehat{{\Cal O}_{W,p}}$, with the locus \linebreak
$M_{0,p} = \{q \in (\widehat{{\Bbb
A}^e})_0;(1,(1,0),0, \cdot\cdot\cdot ,1,(1,0),0,\widehat{f_{p,0}^e}(q)) = \max f_X^d\}$ where the
pattern
$(1,(1,0),0)$ is repeated for $(d-e)$-times, via the inclusion
$(\widehat{{\Bbb A}^e})_0 \subset \roman{Spec}\ \widehat{{\Cal O}_{W,p}}$, i.e.,
$$\CD
(\widehat{{\Bbb
A}^e})_0 @. \hskip.1in \hookrightarrow \hskip.1in @. \roman{Spec}\ \widehat{{\Cal O}_{W,p}} \\
\cup @.  @. \cup \\
M_{0,p} @. \hskip.1in = \hskip.1in @.
\{q \in X_0;f_{(U_p)_0}(q) = \max f_X^d\}
\times
\roman{Spec}\
\widehat{{\Cal O}_{W,p}}. \\
\endCD$$
Therefore, the defining ideal ${\Cal I}_{Y_0}$ of the center $Y_0 \subset X_0$, is characterized
analytically locally as the ideal defining the locus \linebreak
$M_{0,p} = \{q \in (\widehat{{\Bbb
A}^e})_0;(1,(1,0),0, \cdot\cdot\cdot ,1,(1,0),0,\widehat{f_{p,0}^e}(q)) = \max f_X^d\}$,
restricted to $\roman{Spec}\ \widehat{{\Cal O}_{X,p}} \subset (\widehat{{\Bbb A}^e})_0$.  As
before, it is straightforward to see that this characterization is independent of the choice of
$U_p$, $U_p \subset U$, or a system of regular parameters $(x_1, ... , x_{d-e}, y_1, ... , y_e)$,
in the following sense: For any other choice of open subsets or a system of regular parameters,
there exists an isomorphism $\phi$, as in the verification of condition (iii), which identifies
the ideals defining $\roman{Spec}\ \widehat{{\Cal O}_{X,p}}$ in $\roman{Spec}\ \widehat{{\Cal
O}_{W,p}}$ (i.e., ${\Cal K}$ and ${\Cal K}'$) as closed subschemes, as well as the ideals
defining the centers.  Therefore, via $\phi$, we identify the two a priori different defining
ideals of the centers, restricted to
$\roman{Spec}\ \widehat{{\Cal O}_{X,p}}$, as one.

\vskip.1in

This finishes the verification of condition (iv) at the $(k = 0)$-th stage.

\vskip.1in

Condition (v) follows from the observation in Subcase $e - \dim X = 0$ that $f_{X_0}^d(p)$ is
minimum when
$p$ is a nonsingular point and hence that when the center contains $p$ the entire $X$ has to be an
irreducible component of the center, which implies $X$ is already nonsingular.

\vskip.1in

This completes checking of the conditions at the $(k = 0)$-th stage.

\vskip.1in

$\boxed{\roman{Inductional\ Assumption\ of\ Case}\ k = k}$

\vskip.1in

Before proving the assertions for the case $k = k + 1$, based upon the assertions for the case
$k = k$, we make some extra inductional assumption explicit for the case $k = k$ aside from
condition (i) through (iv).

\vskip.1in

$\underline{\roman{Subcase}\ e - \dim X > 0}$: In fact, we assume inductively, starting with
$$\align
(W_{U_p}, ({\Cal I}_{U_p},1), &\emptyset)/(x_1, ... , x_{d-e}) \times \roman{Spec}\
\widehat{{\Cal O}_{W,p}}\\
&= ((W_{U_p})_0, ((J_{U_p})_0,1), (E_{U_p})_0)/(x_{1,0}, ... ,
x_{d-e,0}) \times
\roman{Spec}\
\widehat{{\Cal O}_{W,p}} \\
(\widehat{{\Bbb A}^e},({\Cal K},1), \emptyset) &= ((\widehat{{\Bbb A}^e})_0,({\Cal K}_0,1), ({\Cal
E}_p)_0), \\
\endalign$$
we have a commutative diagram between the sequences of
transformations of (analytic) basic objects 
$$\CD
((W_{U_p})_0, ((J_{U_p})_0,1), (E_{U_p})_0)/(x_{1,0}, ... , x_{d-e,0}) \times
\roman{Spec}\
\widehat{{\Cal O}_{W,p}} @. \overset{\sim}\to{\rightarrow} @. ((\widehat{{\Bbb A}^e})_0,({\Cal
K}_0,1), ({\Cal E}_p)_0)  \\
@AAA @.@AAA \\
\cdot @.@. \cdot \\
\cdot @.@. \cdot \\
\cdot @.@. \cdot \\
@AAA @. @AAA \\
((W_{U_p})_{l_{p,k}}, ((J_{U_p})_{l_{p,k}},1),
(E_{U_p})_{l_{p,k}})/(x_{1,l_{p,k}}, ... , x_{d-e,l_{p,k}})
\times
\roman{Spec}\
\widehat{{\Cal O}_{W,p}} @. \overset{\sim}\to{\rightarrow} @. ((\widehat{{\Bbb
A}^e})_{l_{p.k}},({\Cal K}_{l_{p,k}},1), ({\Cal E}_p)_{l_{p,k}})\\
\endCD$$
where the $x_{j,i}$ denote the strict transform of $x_{j,0} = x_j$, i.e., the equations defining the
strict transforms of the smooth hypersurfaces $\{x_j = 0\}$.  

(Remark that only at this place of
the notes we write down the sequences of transformations vertically, due to the limitation of
space.)  

Note that this commutative diagram leads us immediately to the computation of the
invariant $f_{X_k}^d(p_k)$ for $p_k \in \psi_k^{-1}(p)$ as was done in the case $k = 0$
$$f_{X_k}^d(p_k) = f_{U_p,l_{p,k}}^d(p_k) = (1, (1,0), 0, \cdot\cdot\cdot, 1, (1,0), 0,
\widehat{f_{p,l_{p,k}}^d}(p_k)),$$

where the pattern $(1, (1,0), 0)$ is repeated for $(d-e)$-times and where
$\widehat{f_{p,l_{p,k}}^d}$ denotes the invariant attached to the (analytic) basic object $((\widehat{{\Bbb
A}^e})_{l_{p.k}},({\Cal K}_{l_{p,k}},1), ({\Cal E}_p)_{l_{p,k}})$ obtained from the sequence
above.

\vskip.1in

$\underline{\roman{Subcase}\ e - \dim X = 0}$: Note that in this case $X$ is smooth in a neighborhood of
$p$ and we have not touched this neighborhood in the process of
(non-embedded) resolution of singularities.  Just as at the $0$-th stage, we have 

$$f_{X_k}^d(p_k) = f_{U_p,l_{p,k}}^d(p_k) = (1, (1,0), 0, \cdot\cdot\cdot, 1, (1,0), 0, 1, (1,0),
\infty),$$

where the pattern $(1, (1,0), 0)$ is repeated for $(d-e-1)$-times.  Observe that in this case the
value
$f_{X_k}^d(p)$ is minimum, i.e.,
$$f_{X_k}^d(p) \leq f_{X_k}^d(q) \hskip.1in \forall q \in X_k.$$

\vskip.1in

$\boxed{\roman{Case}\ k = k + 1}$

\vskip.1in

We look at the assertions for the case $k = k+1$.

\vskip.1in

If $\psi_k^{-1}(p_k) \cap Y_k \neq \emptyset$, then we choose $p \in U_p \subset U$ just as we did
at the $k$-th stage.  If $\psi_k^{-1}(p_k) \cap Y_k = \emptyset$, then we shrink $U_p$ so that
$U_p \cap \psi_k(Y_k) = \emptyset$.  This shrinking can be done, since $\psi_k(Y_k)$ is closed by
properness of $\psi_k$ and since $p_k \not\in \psi_k(Y_k)$.  For this choice of $p \in U_p
\subset U$, condition (i) is satisfied at the stage $k = k  + 1$. 

Condition (ii) immediately follows from the construction, since in case $\psi_k^{-1}(p_k) \cap
Y_k \neq \emptyset$, the center $Y_k$ over $(U_p)_{l_{p,k}}$ is given as the restriction to
$(U_p)_{l_{p,k}}$ of the center of the transformation of basic objects
$$((W_{U_p})_{l_{p,k}},((J_{U_p})_{l_{p,k}},1),(E_{U_p})_{l_{p,k}}) \leftarrow
((W_{U_p})_{l_{p,{k+1}}},((J_{U_p})_{l_{p,{k+1}}},1),(E_{U_p})_{l_{p,{k+1}}}).$$
Note that in case $\psi_k^{-1}(p_k) \cap
Y_k = \emptyset$ we have $l_{p,k} = l_{p,k+1}$.

In condition (iii), the only thing to check is that the invariant $f_{U_p,l_{p,k+1}}^d(p_{k+1})$
is independent of the choice of $U$ or $U_p$.  This can be seen easily, once one realizes that
the vertical commutative diagram in the inductional assumption of the case $k = k$ can be extended
to the case $k = k + 1$, giving an anlytic characterization of the invariant, and that for
another choice $p \in U_p' \subset U'$ the isomorphism $\phi$ (discussed in the verification of
the assertions at the stage $k = 0$) can be extended to an isomorphism between the two sequences
of analytic basic objects all the way to the stage $k = k + 1$.  

Condition (iv) indicates how to choose the next center $Y_{k+1} \subset X_{k+1}$ of blowup, which
follows from Theorem 9-3 and from the analytic characterization of the ideal defining the
center.   This analytic characterization can be verified via an argument identical to the one for
the verification of condition (iii).

\vskip.1in
 
Condition (v) follows from the observation in Subcase $e - \dim X = 0$ that $f_{X_k}^d(p)$ is
minimum when
$p$ is a nonsingular point of $X$ and hence that when the center contains $p$ the entire strict
transform has to an irreducible component of the center, which implies $X_k$ is already
nonsingular.

\vskip.1in

This finishes checking the inductive construction of our algorithm.

\vskip.2in

The sequence obtained through our algorithm is independent of the covering $\{U\}$, as it is
purely determined by the analytic basic objects which are independent of the covering.
 
If we choose a number $d'$ (say $d' > d$) for the common dimension of the smooth ambient
varieties $W_U$, then the invariant $f_{X_k}^{d'}$ only differs from $f_{X_k}^d$ by repeating the
pattern
$(1, (1,0), 0)$ for $(d' - d)$-times.  However, the analytic basic objects $(\widehat{\Bbb
A},({\Cal K},1),({\Cal E}_p)_0)$ (and their transformations), being of the same dimension as
the embedding dimension, remain unchanged.  These analytic objects are the only ingredients to
determine the centers of the transformations of the sequence.  Therefore, the sequence obtained
through our algorithm is independent of the number $d$.

Finally, we can see that any automorphism $\theta:X
\overset{\sim}\to{\rightarrow} X$ can be lifted to an automorphism of the sequence, once we
realize that $\{\theta(U)\}$ gives another open covering of $X$ and then that the argument as
above showing the independence of the sequence of the choice of the covering also shows the
lifting of the automorphism. 

\vskip.1in

This completes the proof for Theorem 10-1.
  
\newpage

$$\bold{CHAPTER\ 11.\ EXAMPLES}$$

\vskip.1in

In this chapter, we present examples, some of which demonstrate a couple of essential points of
the inductive algorithm and some of which simply demonstrate how it works.  Many of them are
communicated to the author by Profs. Encinas and Villamayor, and/or taken directly from the
lectures delivered by the latter at Purdue University.  (However, any inaccuracy in the
presentation is solely the responsibility of the author.)

\proclaim{Example 11-1 (Why do we need to keep the history ?)}\endproclaim

In our inductive algorithm of resolution of singularities of a (general) basic object, the
$t$-invariant plays a key role.  The second factor $n_k$ of the $t$-invariant \linebreak
$t_k =
(w\text{-}\roman{ord}_k,n_k)$ depends upon the ``history" of the process.  Namely we have to
look at the sequence from the $0$-th stage up to the
$k$-th stage, finding when the maximum of the invariant
$w\text{-}\roman{ord}$ changed in the past (cf. Definition 1-10).

Do we really have to keep track of the history of the process of resolution of singularities or
principalization ?

\vskip.1in

The following simple example shows that the answer is \it yes \rm, not only in our inductive
algorithm but also in \it any \rm algorithm (which looks only at the weak transforms of the
ideal), in the following sense:

\vskip.1in

Suppose we look for an algorithm of principalization, which assigns a uniquely determined
sequence to a gievn ideal ${\Cal I} \subset {\Cal O}_W$ on a smooth variety, satisfying
conditions (i), (ii), and (iii) in Main-Theme 0-3, and the following extra requirements:

$(\alpha)$ the algorithm is equivariant with respect to an action, and

$(\beta)$ the algorithm is stable with respect to truncation and localization.

\vskip.1in

The example below shows that there is NO such algorithm.  

\vskip.1in

In other words, condition
$(\beta)$ inevitably leads to an infinite loop in the algorithm.  Hence, in order to
guarantee that an algorithm of principalization come to an end after finitely many steps, we
have to give up condition $(\beta)$.

\vskip.1in

Note that the precise meaning of condition $(\beta)$ is:

\vskip.1in

Let
$$W_0 \overset{\pi_1}\to{\leftarrow} W_1 \overset{\pi_2}\to{\leftarrow} \cdot\cdot\cdot
\overset{\pi_{l-1}}\to{\leftarrow} W_{l-1} \overset{\pi_l}\to{\leftarrow} W_l$$
and
$$W_0' \overset{\pi_1}\to{\leftarrow} W_1' \overset{\pi_2}\to{\leftarrow} \cdot\cdot\cdot
\overset{\pi_{l'-1}}\to{\leftarrow} W_{l'-1}' \overset{\pi_{l'}}\to{\leftarrow} W_{l'}'$$
be sequences constructed according to the algorithm for principalization of ideals ${\Cal I} =
{\Cal I}_0 = \overline{{\Cal I}_0}
\subset {\Cal O}_W$ on a smooth variety $W = W_0$ and ${\Cal I}' = {\Cal I}_0' =
\overline{{\Cal I}_0'} \subset {\Cal O}_{W'}$ on
$W' = W_0'$.  Let
$\overline{{\Cal I}_l}$ and
$\overline{{\Cal I}'_{l'}}$ be the weak transforms of ${\Cal I}$ and ${\Cal I}'$, respectively
(cf. Remark 1-11 (ii)).

Suppose there exist open subsets $U \subset W_l$ and $V \subset W'_{l'}$ such that there is an
isomorphism $U \overset{\sim}\to{\rightarrow} V$ which induces an isomorphism of ideals
$\overline{{\Cal I}_l}|_{U} \overset{\sim}\to{\rightarrow} \overline{{\Cal I}'_{l'}}|_{V}$.

Then the extensions of the sequences of principalization of ideal, constructed according to the
algorithm, coincide over
$U$ and over $V$ (after ignoring the trivial transformations whose centers lie outside of the
loci over $U$ or $V$).

\vskip.1in

We translate the truncation property as looking only at
the present (situation of the weak transform).  Therefore, we interpret the necessity to give up
condition
$(\beta)$ as the need to look into the history.

\vskip.1in

Let
$${\Cal I} = \langle x_1, x_2x_3\rangle \subset {\Cal O}_W \text{\ where\ }W = \roman{Spec}\
k[x_1, x_2, x_3].$$
Since $J$ is not principal, we have to choose a center $Y \subset W$.  By condition (i) of
Main Theme 0-3, we have $Y \subset \roman{Supp}\ {\Cal O}_W/{\Cal I}$.  By condition (ii) of
Main Theme 0-3, the center $Y$ has to be smooth.  There is an obvious
action of ${\Bbb Z}_2$ on $W$, switiching $x_2$ and $x_3$, under which ${\Cal I}$ is
invariant.  Therefore, $Y$ has to be invariant under the ${\Bbb Z}_2$-action.

It is easy to see that the only center which satisfy all the above requirements and which
includes the origin is the origin itself.

Therefore, the first blowup (in a neighborhood of the origin) of principalization \it must \rm
be along the ideal
$\langle x_1, x_2, x_3\rangle$. 
  
However, over the open subset $U$ of $W_1$ with the system of regular parameters $(t_1, x_2,
t_3)$ with
$$x_1 = t_1x_2, x_2 = x_2, x_3 = t_3x_2,$$
the weak transform $\overline{{\Cal I}_1}|_U$ is in the identical form to the ideal ${\Cal I}_0
= \overline{{\Cal I}_0}$, i.e.,
$$\overline{{\Cal I}_1}|_U = \langle t_1, x_2t_3 \rangle
\overset{\sim}\to{\rightarrow} \overline{{\Cal I}_0}|_V$$ where we set $V = W = W_0$.

Now it is clear that condition $(\beta)$  leads to an infinite loop of the process induced by
the algorithm.

\vskip.1in

\proclaim{Example 11-2 (Fundamental obstruction to carry out our algorithm in positive
characteristic)}\endproclaim

The success of our inductive algorithm of resolution of singularities depends in an essential
way on finding a hypersurafce of maximal contact (cf. Remark 1-5, Lemma 3-1 (key inductive
lemma)).

The following example shows that a hypersurface of maximal contact does \it not \rm always exist
in positive characteristic, and hence that there is a fundamental obstruction to carry out our
algorithm in positive characteristic.  (It is brought to the attention of the author by Prof.
J. W{\l}odarczyk via communication with Prof. P. Milman.  The reader is also encouraged to look
at Moh [1].)

\vskip.1in

Consider a hypersurafec singularity
$$0 \in \{f = 0\} \subset {\Bbb A}^4 = \roman{Spec}\ k[x_1, x_2, x_3, x_4]$$
where
$$f = x_4^2 + x_1^3x_2 + x_2^3x_3 + x_3^7x_1.$$
Suppose that the characteristic of the base field $k$ (which is assumed to be algebraically
closed for simplicity) is equal to 2, i.e.,
$$\roman{char}(k) = 2.$$
We look for the locus where the multiplicity of $f$ is equal to 2 (or more), where 2 is the
multiplicity of $f$ at the origin.

By substituting the following into $f$
$$\left\{\aligned
x_1 &= y_1 + a \\
x_2 &= y_2 + b \\
x_3 &= y_3 + c \\
x_4 &= y_4 + d, \\
\endaligned\right.$$
we see
$$\align
f &= \{d^2 + a^3b + b^3c + c^7a\} + \{(a^2b + c^7)y_1 + (b^2c + a^3)y_2 + (ac^6 + b^3)y_3\} \\
&+\ \text{higher\ terms}.\\
\endalign$$
It is straightforward to see from this that the locus of multiplicity 2 (or more) has the
parametrization
$$\left\{\aligned
a &= t^{15} \\
b &= t^{19} \\
c &= t^7 \\
d &= t^{32}. \\
\endaligned\right.$$
The embedding dimension of the curve parametrized as above at the origin is 4, and hence it can
never be contained in a smooth hypersurface in a neighborhood of the origin.  Therefore, there
is no hypersurface of maximal contact at the origin for this example.

\vskip.1in

\proclaim{Example 11-3 (Our algorithm of resolution of singularities of a general basic object
DOES depend on the specification of the dimension $d$ of its structure.)}\endproclaim

Our inductive algorithm of resolution of singularities of a general basic object $({\Cal F},
(W,E))$ is determined by the invariant $f^d$ (cf. Chapter 9) where $d$ is the number specifying
the dimension of its structure.  Sometimes the general basic object
can have a $d$-dimensional structure as well as a $d'$-dimensional structure for two different
numbers $d \neq d'$.  That is to say, we can have two different sets of charts
$\{(\widetilde{W^{\lambda}},({\goth
a}^{\lambda},b^{\lambda}),\widetilde{E^{\lambda}})\}_{\lambda \in \Lambda}$ and
$\{(\widetilde{W^{\mu}},({\goth b},c^{\mu}),\widetilde{E^{\mu}})\}_{\mu \in M}$ being of
different dimensions $d$ and $d'$, i.e., $\dim \widetilde{W^{\lambda}} = d \neq d' =
\widetilde{W^{\mu}}$, but giving rise to the same collection ${\goth C}$ of sequences of smooth
morphisms and transformations of pairs with specified closed subsets, represented by the
general basic object $({\Cal F},
(W,E))$.

Since the invariants $f^d$ and $f^{d'}$ could be different (cf. Remark 4-7 and Remark 9-2
(iii)), it is natural to suspect that our algorithm depends on the specification of
the dimension of the structure of one's choice of a general basic object.

The following example, taken directly from Encinas [1], shows that this is indeed the
case, demonstrating a general basic object having two different sequences of resolution of
singularities, though both are prescribed by our algorithm, depending on two different
specifications $d$ and $d'$. 

\vskip.1in

Consider the following basic object $(W,(J,b),E)$ of dimension 4 where
$$\left\{\aligned
W &= {\Bbb A}^4 = \roman{Spec}\ k[x_1,x_2,x_3,x_4], \\
J &= \langle f \rangle \hskip.1in \text{with} \hskip.1in f = x_4^2 + x_3^3 + x_2x_3^2 + x_1^3,\\
b &= 2,\\
E &= \{H\} \hskip.1in \text{with} \hskip.1in H = \{x_3 = 0\}.\\
\endaligned\right.$$
As in Remark 4-2 (ii), the basic object $(W,(J,b),E)$ defines a general basic object $({\Cal
F},(W,E))$ with a $4$-dimensional structure.
 
Also consider the following basic object $(X,({\goth a},c),F)$ of dimension 3 where
$$\left\{\aligned
X &= {\Bbb A}^3 = \roman{Spec}\ k[x_1,x_2,x_3], \\
{\goth a} &= \langle g \rangle \hskip.1in \text{with} \hskip.1in g = x_3^3 + x_2x_3^2 + x_1^3,\\
c &= 2,\\
F &= \{H_V\} \hskip.1in \text{with} \hskip.1in H_V = \{x_3 = 0\}.\\
\endaligned\right.$$
From Giraud's Lemma (cf, Claim 3-4) and from a view point of looking at $f$ as a polynomial
in $x_4$, it follows immediately that, via the closed immersion of pairs
$(X =
\{x_4 = 0\},F)
\hookrightarrow (W,E)$, the basic object provides a (global) 3-dimensional chart to the general
basic object $({\Cal F},(W,E))$.

Therefore, the general basic object has two different dimensions, namely 3 and 4, for its
structure.

\vskip.2in

(i) Resolution of singularities of $({\Cal F},(W,E))$ with the $4$-dimensional structure:

\vskip.1in

We apply our inductive algorithm to the basic object $(W,(J,b),E)$ of dimension 4.  

\vskip.1in

Via direct computation, we see that
$$\align
\roman{Sing}(J,b) &= V(x_1, x_3, x_4), \\
\max\ t &= \max\ (w\text{-}\roman{ord},n) = (1,1),\\
\underline{\roman{Max}}\ t &= V(x_1,x_2,x_3).\\
\endalign$$
Since $\roman{codim}_W\underline{\roman{Max}}\ t > 1$, accoprding to the algorithm, we proceed
to construct a basic object $(W'',(J'',b''),E'')$ where
$$\left\{\aligned
W'' &= W, \\
J'' &= \overline{J} + \langle x_3^2 \rangle = J + \langle x_3^2 \rangle = \langle x_4^2 + x_1^3,
x_3^2
\rangle,\\ b'' &= b = 2,\\
E'' &= \emptyset.\\
\endaligned\right.$$
Now out of the basic object $(W'',(J'',b''),E'')$, knowing $x_4 \in \Delta(J'')$, we construct a
basic object $(\widetilde{W''},(C(J''),b''!),\widetilde{E''})$ of dimension 3 where
$$\left\{\aligned
\widetilde{W''} &= \{x_4 = 0\}, \\
C(J'') &= \langle x_1^3, x_3^2 \rangle + \langle x_1^2, x_3 \rangle^2 = \langle
x_1^3,x_3^2,x_1^2x_3 \rangle,\\ 
b''! &= 2! = 2,\\
\widetilde{E''} &= E''|_{\widetilde{W''}} = \emptyset.\\
\endaligned\right.$$
Denote $(\widetilde{W''},(C(J''),b''!),\widetilde{E''})$ by
$(W^{(3)},(J^{(3)},b^{(3)}),E^{(3)})$ and the associated $t$-invariant by $t^{(3)}$ and the
others by putting the superscript ${}^{(3)}$.  Via direct computation, we see that 
$$\align
\roman{Sing}(J^{(3)},b^{(3)}) &= V(x_1, x_3), \\
\max\ t^{(3)} &= \max\ (w\text{-}\roman{ord}^{(3)},n^{(3)}) = (1,0),\\
\underline{\roman{Max}}\ t^{(3)} &= V(x_1,x_3).\\
\endalign$$
Since $\roman{codim}_{W^{(3)}}\underline{\roman{Max}}\ t^{(3)} > 1$, according to the
algorithm, we proceed to construct a basic object
$({W^{(3)}}'',({J^{(3)}}'',{b^{(3)}}''),{E^{(3)}}'')$, which is nothing but
$(W^{(3)},(J^{(3)},b^{(3)}),E^{(3)})$ itself in this case.
Now out of the basic object $({W^{(3)}}'',({J^{(3)}}'',{b^{(3)}}''),{E^{(3)}}'')$, knowing 
$x_3 \in \Delta({J^{(3)}}'')$, we construct a basic object
$(\widetilde{{W^{(3)}}''},(C({J^{(3)}}''),{b^{(3)}}''!),\widetilde{{E^{(3)}}''})$ of dimension
2 where
$$\left\{\aligned
\widetilde{{W^{(3)}}''} &= \{x_3 = 0\}, \\
C({J^{(3)}}'') &= \langle x_1^3 \rangle + \langle x_1^2 \rangle^2 = \langle
x_1^3 \rangle,\\ 
{b^{(3)}}''! &= 2! = 2,\\
\widetilde{{E^{(3)}}''} &= {E^{(3)}}''|_{\widetilde{{W^{(3)}}''}} = \emptyset.\\
\endaligned\right.$$
Denote $(\widetilde{{W^{(3)}}''},(C({J^{(3)}}''),{b^{(3)}}''!),\widetilde{{E^{(3)}}''})$ by
$(W^{(2)},(J^{(2)},b^{(2)}),E^{(2)})$ and the associated $t$-invariant by $t^{(2)}$ and the
others by putting the superscript ${}^{(2)}$.  Via direct computation, we see that 
$$\align
\roman{Sing}(J^{(2)},b^{(3)}) &= V(x_1), \\
\max\ t^{(2)} &= \max\ (w\text{-}\roman{ord}^{(2)},n^{(2)}) = (\frac{3}{2},0),\\
\underline{\roman{Max}}\ t^{(2)} &= V(x_1).\\
\endalign$$

Since this time
$$\roman{codim}_{W^{(2)}}\underline{\roman{Max}}\ t^{(2)} = 1,$$
according to the algorithm, we finally decide that the center $Y_0$ of the first transformation
must be $Y_0 = R(1)(\underline{\roman{Max}}\ t^{(2)})$, i.e., via the inclusion $W^{(2)} =
V(x_3,x_4) \subset W^{(4)} = W$ we have the description of the center
$$Y_0 = V(x_1,x_3,x_4) \subset W.$$

\vskip.1in

(ii) Resolution of singularities of $({\Cal F},(W,E))$ with the $3$-dimensional structure: 

\vskip.1in

We
apply our inductive algorithm to the basic object $(X,({\goth a},c),F)$ of dimension 3.

\vskip.1in

Via direct computation, we see that
$$\align
\roman{Sing}({\goth a},c) &= V(x_1, x_3), \\
\max\ t &= \max\ (w\text{-}\roman{ord},n) = (\frac{3}{2},1),\\
\underline{\roman{Max}}\ t &= V(x_1,x_2,x_3).\\
\endalign$$
Since $\roman{codim}_X\underline{\roman{Max}}\ t > 1$, we have to proceed constructing the
auxiliary basic objects.  However, since we know that the center $Y_0'$ that we take for the
first transformation must satisfy the condition $Y_0' \subset \underline{\roman{Max}}\ t$, we
conclude that via the inclusion $X = V(x_4) \subset W$ we have the description of the center
$$Y_0' = V(x_1,x_2,x_3,x_4) \subset W.$$

\vskip.1in

Comparing (i) and (ii), we see that the two sequences of resolution of singularities of the
general basic object $({\Cal F},(W,E))$, one with a $4$-dimensional structure and the other
with a $3$-dimensional structure, have two different centers for the first transformations
$$Y_0 = V(x_1,x_3,x_4) \neq Y_0' = V(x_1,x_2,x_3,x_4) \subset W.$$
Therefore, the two sequences are obviously different.

\vskip.1in

\proclaim{Remark 11-4}\endproclaim

(i) Example 11-3 should not be confused with the fact that our algorithm for non-embedded
resolution of singularities of a variety $X$ does NOT depend on the choice of
the number $d$ which represents the common dimension of the ambient smooth varieties $W_U$, into
which the open subsets $U$ (in an open covering $\{U\}$ of $X$) are embedded (cf. Theorem
10-1).  

It should be noted that the algorithm for non-embedded resolution of singularities of a variety
$X$ described in Encinas [1] is different from our algorithm in Chapter 10.  Therefore, there
is no contradiction between our Theorem 10-1 and the claim in Encinas [1] that the algorithm for
non-embedded resolution of singularities of a variety $X$ DOES depend on the choice of
the number $d$ which represents the common dimension of the ambient smooth varieties $W_U$, into
which the open subsets $U$ (in an open covering $\{U\}$ of $X$) are embedded.

\vskip.1in

(ii) Let $({\Cal F},(W,E))$ be a general basic object with a $d$-dimensional structure,
having charts $\{(\widetilde{W^{\lambda}},({\goth
a},b^{\lambda}),\widetilde{E^{\lambda}})\}_{\lambda \in \Lambda}$ of
$\dim \widetilde{W^{\lambda}} = d \hskip.1in \forall \lambda \in \Lambda$.  Assume that $({\Cal
F},(W,E))$ has a $d'$-dimensional structure with $d' < d$.

Suppose that $({\Cal F},(W,E))$ is simple in the sense that $(\widetilde{W^{\lambda}},({\goth
a},b^{\lambda}),\widetilde{E^{\lambda}})$ is a simple basic object for all $\lambda \in
\Lambda$.

Suppose further that $E = \emptyset$.

Then the sequence representing resolution of singularities of $({\Cal
F},(W,E))$ constructed according to our algorithm with the specified dimension of the structure
being
$d$, coincides with the sequence constructed according to our algorithm with the specified dimension of the structure
being
$d'$.  

This can be seen as follows: Firstly observe that for a simple basic
object (and their transformations) the invariants $\roman{ord}$ and $w\text{-}\roman{ord}$
coincide and, if further the boundary divisor is empty, then $w\text{-}\roman{ord}$ and the
$t$-invariant coincide (until the maximum drops).  Secondly, based upon the first observation,
observe that the auxiliary general basic object of dimension $d'$ that we construct in the
process prescribed by our algorithm must coincide with the original general basic object with
a $d'$-dimensional structure.

In Example 11-3, the general basic object $({\Cal F},(W,E))$ is simple, arising from the
simple basic object $(W,(J,b),E)$.  Therefore, setting $E = \{H\} \neq \emptyset$ is
essential in order to get two different sequences of resolution of singularities, depending
on $d = 4$ and $d' = 3$.

The author does not know an example of a general basic object $({\Cal F},(W,E))$ with $E =
\emptyset$, having structures of two different dimensions, for which our algorithm gives rise
to two different sequences of resolution of singularities depending on the specified dimensions
of the structures.  Such a general basic object cannot be simple.

To reveal more ignorance, the author does not know an example of a basic object of dimension $d$
which is not simple and which has a $(d-1)$-dimensional structure as a general basic object (or
even if such a basic object exists).

\vskip.1in

\proclaim{Example 11-5 (The centers for non-embedded resolution may not be smooth.)}\endproclaim

\vskip.1in

The following example of a sequence of embedded resolution of singularities of a variety
$X
\subset W$ (embedded as a closed subscheme in a smooth variety)
$$X = X_0 \subset W = W_0 \overset{\pi_1}\to{\leftarrow} X_1 \subset
W_1 \overset{\pi_2}\to{\leftarrow} \cdot\cdot\cdot \overset{\pi_{l-1}}\to{\leftarrow} X_{l-1}
\subset W_{l-1} \overset{\pi_l}\to{\leftarrow} X_l
\subset W_l,$$
shows that the center $Y_{i-1}
\subset W_{i-1}$, chosen by our inductive algorithm, is always smooth inside of the ambient
variety
$W_{i-1}$ (cf. Theorem 7-1) by construction, but that 

($\alpha$) the center $Y_{i-1}$ may not be contained in the strict transform $X_{i-1}$, and/or

($\beta$) the intersection of the center with the strict transform $Y_{i-1} \cap X_{i-1}$ may

not be smooth.

(Remark that in such an example $X$ must not be a hypersurface (cf. Remark 7-2 (ii)) and hence
that the defining ideal ${\Cal I}_X$ of $X$ in $W$ has to have two or more generators even
locally.)

\vskip.1in

Since the sequence representing non-embedded resolution of singularities of $X$, constructed
according to our inductive algorithm, is based upon the one representing embedded resolution of
singularities (cf. Chapter 10), this example also shows that the centers for the sequence
representing non-embedded resolution of singularities may not be smooth.

\vskip.1in

Let $X \subset W = \roman{Spec}\ k[x,y,z,w]$ be a subvariety defined by the ideal
$${\Cal I}_X = \langle x^2 + y^2 + z^2 + w^2, x^6 + y^6 + z^6 + w^6 \rangle.$$

In order to obtain the sequence representing embedded resolution of singularities, we consider
(cf. Chapter 10) a basic object $(W,(J,b),E)$ where
$$\left\{\aligned
W &= {\Bbb A}^4 = \roman{Spec}\ k[x,y,z,w] \\
J &= {\Cal I}_X \\
b &= 1 \\
E &= \emptyset.\\
\endaligned\right.$$
It is straightforward to see that
$$w\text{-}\roman{ord}(p) = \left\{\aligned
2 \hskip.2in &\text{\ if\ } p = 0 \\
1 \hskip.2in &\text{\ if\ } p \neq 0.\\
\endaligned\right.$$
Since the center $Y_0 \subset W_0$ of the first transformation of basic objects 
$$(W,(J,b),E) =
(W_0,(J_0,b),E_0)
\overset{\pi_1}\to{\leftarrow} (W_1,(J_1,b),E_1)$$
must be contained in the maximum locus of $w\text{-}\roman{ord} = w\text{-}\roman{ord}_0$,
i.e., \linebreak
$Y_0 \subset \underline{\roman{Max}}\ w\text{-}\roman{ord}_0 = \{0\}$, we conclude that
$$Y_0 = \{0\}.$$
Consider the affine open subset of $W_1$ obtained by inverting $x$ (We use the same letter
$W_1$ for the affine open subset by abuse of notation.), with a system of regular
parameters $(x_1,y_1,z_1,w_1)$ with
$x = x_1, y = xy_1, z = xz_1, w = xw_1$.

By direct computation we see that
$$\align
J_1 &= \langle x_1(1 + y_1^2 + z_1^2 + w_1^2), x_1^5(1 + y_1^6 + z_1^6 + w_1^6) \rangle \\
\overline{J_1} &= \langle 1 + y_1^2 + z_1^2 + w_1^2, x_1^4(1 + y_1^6 + z_1^6 + w_1^6) \rangle \\
{\Cal I}_{X_1} &= \langle 1 + y_1^2 + z_1^2 + w_1^2, 1 + y_1^6 + z_1^6 + w_1^6 \rangle \\
\endalign$$
and hence that
$$\align
\roman{Sing}(J_1,b) &= V(x_1) \cup V(\langle 1 + y_1^2 + z_1^2 + w_1^2, 1 + y_1^6 + z_1^6 +
w_1^6 \rangle) \\
w\text{-}\roman{ord}_1(p_1) &= \left\{\aligned
0 \hskip.25in &\text{\ if\ }p_1 \in \roman{Sing}(J_1,b) \setminus V(1 + y_1^2 + z_1^2 + w_1^2)
\\ 
1 \hskip.25in &\text{\ if\ }p_1 \in \roman{Sing}(J_1,b) \cap V(1 + y_1^2 + z_1^2 + w_1^2)\\
\endaligned\right. \\
t_1(p_1) &= \left\{\aligned
(0,1) &\text{\ if\ }p_1 \in \roman{Sing}(J_1,b) \setminus V(1 + y_1^2 + z_1^2 + w_1^2) \\
(1,0) &\text{\ if\ }p_1 \in (\roman{Sing}(J_1,b) \cap V(1 + y_1^2 + z_1^2 + w_1^2)) \setminus
V(x_1) \\ 
(1,1) &\text{\ if\ }p_1 \in \roman{Sing}(J_1,b) \cap V(1 + y_1^2 + z_1^2 + w_1^2) \cap V(x_1).
\\
\endaligned\right.\\ 
\endalign$$
Therefore, we have
$$\align
\max\ t_1 &= (1,1) \\
\underline{\roman{Max}}\ t_1 &= V(x_1) \cap V(1 + y_1^2 + z_1^2 + w_1^2). \\
\endalign$$
Now since
$$\roman{codim}_{W_1}\underline{\roman{Max}}\ t_1 > 1,$$
our inductive algorithm tells us (cf. Theorem 5-1) to construct a basic object
$(W_1'',(J_1'',b''),E'')$ (cf. Lemma 5-3 and lemma 5-4) where
$$\left\{\aligned
W_1'' &= W_1,\\
J_1'' &= \overline{J_1} + \langle x_1 \rangle = \langle 1 + y_1^2 + z_1^2 + w_1^2, x_1^4(1 +
y_1^6 + z_1^6 + w_1^6), x_1 \rangle, \\ 
b'' &= b = 1, \\
E_1'' &= E_1^+ = \emptyset.\\
\endaligned\right.$$
The basic object $(W_1'',(J_1'',b''),E'')$ has a $3 (= 4-1)$-dimensional structure, whose chart
is given by
$$((W_1'')_h,(C(J_1''),b''!),(E_1'')_h) = (\{x_1 = 0\}, (\langle 1 + y_1^2 + z_1^2 + w_1^2
\rangle, 1),
\emptyset).$$
We see that
$$\roman{Sing}(C(J_1''),b''!) = V(1 + y_1^2 + z_1^2 + w_1^2)$$
and that
$$w\text{-}\roman{ord}(q) = 1, t(q) = (1,0) \hskip.1in \forall q \in
\roman{Sing}(C(J_1''),b''!).$$
Therefore, we have
$$R(1)(\underline{\roman{Max}}\ {t_1''}^{(3)}) = \underline{\roman{Max}}\ {t_1''}^{(3)} =
\roman{Sing}(C(J_1''),b''!) = V(1 + y_1^2 + z_1^2 + w_1^2).$$
Therefore, according to our inductive algorithm, the center $Y_1 \subset W_1$ for the second
transformation of basic objects
$$(W_1,(J_1,b),E_1) \overset{\pi_2}\to{\leftarrow} (W_2,(J_2,b),E_2)$$
has to be taken so that
$$Y_1 = V(x_1,1 + y_1^2 + z_1^2 + w_1^2) \subset W_1.$$
Although $Y_1$ itself is smooth, the intersection with the strict transform is described by
$$Y_1 \cap X_1 = V(x_1,1 + y_1^2 + z_1^2 + w_1^2,1 + y_1^6 + z_1^6 + w_1^6),$$
whose singular locus is characterized by the condition
$$\roman{rank}\left[\matrix
2y_1 & 2z_1 & 2w_1 \\
6y_1^5 & 6z_1^5 & 6w_1^5 \\
\endmatrix\right] < 2$$
and hence contains a point, e.g.,
$$(x_1,y_1,z_1,w_1) = (0,0,0,i).$$
Therefore, $Y_1 \cap X_1$ is SINGULAR.

\proclaim{Remark 11-6}\endproclaim

In general, we observe (cf. the footnote to Main Theme 0-1, Remark 10-2 (ii)) that

($\alpha$) the center $Y_{i-1}$ may not be contained in the strict transform $X_{i-1}$, and/or

($\beta$) the intersection of the center with the strict transform $Y_{i-1} \cap X_{i-1}$ may

not be smooth or even reduced.

We leave it as an exercise to the reader to produce an example where $Y_{i-1} \cap X_{i-1}$ is
actually non-reduced.

\vskip.1in

\proclaim{Example 11-7 (Resolution of singularities of a monomial basic object.)}\endproclaim

Finally, we give an example demonstrating how to construct a sequence representing resolution
of singularities of a monomial basic object, explaining how the tie breaker works  and how the
convention in Definition 1-8 (iii) works.  It could be of help to the reader trying to fill in
the proof to Proposition 2-5.

\vskip.1in

We start with a monomial basic object $B = (W,(J,b),E)$ where
$$\left\{\aligned
W &= {\Bbb A}^2 = \roman{Spec}\ k[x,y] \\
J &= I(H_1)^3I(H_2)^3 \\
b &= 2 \\
E &= \{H_1,H_2\},\\
\endaligned\right.$$
where $H_1 = V(x)$ and $H_2 = V(y)$.

It is straightforward to see that
$$\align
\roman{Sing}(J,b) &= H_1 \cup H_2, \\
\Gamma_{B_0}(p) &= \left\{\aligned
(-1,\frac{3}{2},(1,0)) &\text{\ if\ }p \in H_1 \setminus H_2 \\
(-1,\frac{3}{2},(2,0) &\text{\ if\ }p \in H_2.\\
\endaligned\right.\\
\endalign$$
Therefore, we have
$$\align
\max\ \Gamma_{B_0} &= (-1,\frac{3}{2},(2,0)) \\
\underline{\roman{Max}}\ \Gamma_{B_0} &= H_2.\\
\endalign$$
Observe that the indices of $H_1$ and $H_2$ work as a tiebreaker.

According to our algorithm, we choose $Y_0 = \underline{\roman{Max}}\ \Gamma_{B_0} = H_2$ to be
the center for the first transformation
$$(W,(J,b),E) = B = B_0 = (W_0,(J_0,b),E_0) \hskip.1in \overset{\pi_1}\to{\leftarrow} \hskip.1in B_1 =
(W_1,(J_1,b),E_1).$$
Since $Y_0$ is a divisor, $W_0 \overset{\pi_1}\to{\leftarrow} W_1$ is an isomorphism, whereas,
according to the convention in Definition 1-8 (iii), $H_2$ is now called $H_3$ and hence we have
$$J_1 = I(H_1)^3I(H_3) \hskip.1in \text{and} \hskip.1in E_1 = \{H_1,H_3\}.$$
We see that
$$\align
\roman{Sing}(J_1,b) &= H_1, \\
\Gamma_{B_1}(p) &= (-1,\frac{3}{2},(1,0)) \hskip.1in \forall p \in H_1.\\
\endalign$$
Therefore, we have
$$\align
\max\ \Gamma_{B_1} &= (-1,\frac{3}{2},(1,0)) \\
\underline{\roman{Max}}\ \Gamma_{B_1} &= H_1.\\
\endalign$$
According to our algorithm, we choose $Y_1 = \underline{\roman{Max}}\ \Gamma_{B_1} = H_1$ to be
the center for the second transformation
$$B_1 = (W_1,(J_1,b),E_1) \hskip.1in \overset{\pi_2}\to{\leftarrow} \hskip.1in B_2 =
(W_2,(J_2,b),E_2).$$
Since $Y_1$ is a divisor, $W_1 \overset{\pi_2}\to{\leftarrow} W_2$ is an isomorphism, whereas,
according to the convention in Definition 1-8 (iii) again, $H_1$ is now called $H_4$ and hence
we have
$$J_2 = I(H_3)I(H_4) \hskip.1in \text{and} \hskip.1in E_2 = \{H_3,H_4\}.$$
We see that
$$\align
\roman{Sing}(J_2,b) &= H_3 \cap H_4, \\
\Gamma_{B_2}(p) &= (-2,1,(4,3)) \hskip.1in \forall p \in H_3 \cap H_4. \\
&\hskip.1in \text{(Note
that\ }H_1
\cap H_4 \text{\ consists\ of\ a\ point.)}\\
\endalign$$
Therefore, we have
$$\align
\max\ \Gamma_{B_2} &= (-2,1,(4,3)) \\
\underline{\roman{Max}}\ \Gamma_{B_2} &= H_3 \cap H_4.\\
\endalign$$
According to our algorithm, we choose $Y_2 = \underline{\roman{Max}}\ \Gamma_{B_2} = H_3 \cap
H_4$ to be the center for the third transformation
$$B_2 = (W_2,(J_2,b),E_2) \hskip.1in \overset{\pi_3}\to{\leftarrow} \hskip.1in B_3 =
(W_3,(J_3,b),E_3).$$
We observe then that
$$\roman{Sing}(J_3,b) = \emptyset$$
and the sequence stops here achieving resolution of singularities.  

Note that in the process we
have the invariant $\Gamma$ strictly decreasing
$$\max\ \Gamma_{B_0} = (-1,\frac{3}{2},(2,0)) > \max\ \Gamma_{B_1} = (-1,\frac{3}{2},(1,0)) > 
\max\ \Gamma_{B_2} = (-2,1,(4,3)).$$

\newpage

$$\bold{REFERENCES}$$

\vskip.1in

We should emphasize that the purpose of these notes is to give a self-contained
exposition of the results of our seminar on the specified inductive algorithm by
Encinas and Villamayor, and not to write a treatise on the subject of resolution of
singularities in general.  Accordingly, our references are very restricted and limited
to those which are directly related to the inductive algorithm we concentrate our focus
on and which happen to fall upon the eyes of the author, and we do not even pretend
to try to list a part of the vast literature which may probably be connected to the
algorithm in an explicit or implicit way.  (For example, we list no references to the
recent development of ``weak resolution" by Abramovich-de Jong, Bogomolov-Pantev and
others, based upon the work of de Jong on alterations.)  The readers are encouraged to
look at the references of the papers listed below for those references that we miss
here and beyond.

\vskip.1in

Bierstone, E. and Milman, P.

\hskip.3in [1] \it Canonical desingularization in characteristic zero by blowing-up 

\hskip.3in the
maximal strata of a local invariant, \rm Invent. Math. $\bold{128}$ (2) (1997), 

\hskip.3in 207-302

\hskip.3in [2] \it Desingularization algorithms I; Role of exceptional divisors, 

\hskip.3in \rm
mathAG/0207098

Encinas, S.

\hskip.3in [1] \it On properties of constructive desingularization, \rm preprint

\vskip.1in

Encinas, S. and Villamayor, O.

\hskip.3in [1] \it A course on constructive desingularization and equivariance, 

\hskip.3in \rm in Resolution of singularities (Obergurgl, 1997), Progress in Math.

\hskip.3in $\bold{181}$, Birkh\"auser (2000), 147-227

\hskip.3in [2] \it A new theorem of desingularization over fields of characteristic
zero,

\hskip.3in \rm preprint (1999)

Hauser, H.

\hskip.3in [1] \it The Hironaka theorem on resolution of singularities (Or: a proof 

\hskip.3in we
always wanted to understand), \rm preprint (2002)

Hironaka, H.

\hskip.3in [1] \it Resolution of singularities of an algebraic variety over a field 

\hskip.3in of
characteristic zero, \rm Ann. of Math. $\bold{79}$ (1964), 109-326

\hskip.3in [2] \it Idealistic exponents of singularity, \rm Algebraic Geometry 

\hskip.3in The
Johns Hopkins Cent. Lect., Johns Hopkins Univ. Press (1977), 52-125

Lipman, J.

\hskip.3in [1] \it Introduction to resolution of singularities, \rm ``Algebraic
Geometry"

\hskip.3in (Proc. Sympos. Pure Math. $\bold{29}$), Amer. Math. Soc., Providence R.I.

\hskip.3in (1975), 187-230
 
Moh, T.T.

\hskip.3in [1] \it On a Newton polygon approach to the uniformization of singularities
of

\hskip.3in characteristic $p$, \rm Progress in Mathematics $\bold{134}$ (1996), 49-93

Villamayor, O.

\hskip.3in [1] \it Constructiveness of Hironaka's resolution, \rm Ann. Scient. Ecole
Norm.

\hskip.3in Sup. Paris 4 $\bold{22}$ (1989), 1-32

\enddocument